\documentclass{article}[11pt]

\usepackage{hyperref}  % hyperlinks
\usepackage{geometry}  % page geometry
\usepackage{amssymb, amsbsy}  % maths symbols
\usepackage{mathtools}  % misc maths enhancements and tools
\usepackage{xfrac}  % compact inline fractions
\usepackage{xspace}  % smart space character
\usepackage{float}  % floating environments
\usepackage{graphicx}  % enhanced graphics support
\usepackage{subfigure}
\usepackage{natbib}  % extended citation commands + bib styling
\usepackage[capitalise]{cleveref}  % smart cross-references [always capitalise]
\usepackage{appendix}  % appendix
\usepackage{enumitem}  % customisable enumerate environment
\usepackage{multicol}  % multiple column environment
\usepackage[noend]{algpseudocode}  % pseudo-code commands [hide end statements]
\usepackage{algorithm}  % algorithm float environment
\usepackage{acro}  % automatic acronyms
\usepackage{xr}  % cross referencing between documents
\usepackage{bm} % bold greek symbols in maths
\usepackage{tabularx}  % more flexible tabular environment
\usepackage{etoolbox}  % tools for macros
\usepackage{booktabs}
\usepackage{multirow}

\hypersetup{colorlinks=true, allcolors=blue}

\makeatletter
 
  \@addtoreset{equation}{section}
 \makeatother
\usepackage{bm}
\usepackage{amsfonts}
\usepackage{amssymb}
\usepackage{mathrsfs}
\usepackage{amsmath,amsthm}
\usepackage{authblk}

\usepackage{comment}
\usepackage{graphicx}
\usepackage{color}
\usepackage{indentfirst}
\usepackage{here}
\usepackage{natbib} 

\newtheorem{axiom}{Axiom}
\newtheorem{claim}[axiom]{Claim}
\newtheorem{theorem}{Theorem}[section]
\newtheorem{proposition}[theorem]{Proposition}
\newtheorem{lemma}[theorem]{Lemma}
\newtheorem{remark}[theorem]{Remark}
\newtheorem{assump}[theorem]{Assumption}
\newtheorem*{ex}{Example}
\usepackage{mathtools}
\usepackage{booktabs} 

\DeclarePairedDelimiter\floor{\lfloor}{\rfloor} 

\makeatletter
\newcommand*{\addFileDependency}[1]{% argument=file name and extension
\typeout{(#1)}
\@addtofilelist{#1}
\IfFileExists{#1}{}{\typeout{No file #1.}}
}\makeatother

\usepackage{multirow}
\usepackage{docmute}
\usepackage{xr} 
\providecommand{\keywords}[1]{%
    \textbf{Keywords:} #1
}

\title{A Closed-Form Transition Density Expansion for Elliptic and Hypo-Elliptic SDEs}

\author[1]{Yuga Iguchi}
\author[2]{Alexandros Beskos}

\affil[1]{Lancaster University, Lancaster, UK}
\affil[2]{University College London, London, UK}
 
\begin{document}

\maketitle 

\begin{abstract}
We introduce a closed-form expansion for the transition density of elliptic and hypo-elliptic multivariate Stochastic Differential Equations (SDEs), over a period $\Delta\in (0,1)$, in terms of powers of $\Delta^{j/2}$, $j\ge 0$. Our methodology provides approximations of the transition density, easily evaluated via any software that performs symbolic calculations. 
A major part of the paper is devoted to an analytical control of the remainder in our expansion for fixed $\Delta\in(0,1)$. 
The obtained error bounds validate theoretically the methodology, by 
characterising the size of the distance from the true value. It is the first time that such a closed-form expansion becomes available for the important class of hypo-elliptic SDEs, to the best of our knowledge. For elliptic SDEs, closed-form expansions are available, with some works identifying the size of the error for fixed $\Delta$, as per our contribution. Our methodology allows for a uniform treatment of elliptic and hypo-elliptic SDEs, when earlier works are intrinsically restricted to an elliptic setting. We show numerical applications highlighting the effectiveness of our method, by carrying out parameter inference for hypo-elliptic SDEs that do not satisfy stated conditions. The latter are sufficient for controlling the remainder terms, but the closed-form expansion itself is applicable in general settings.  
\end{abstract}


\section{Introduction} \label{sec:intro}  
Stochastic Differential Equations (SDEs) constitute an effective tool for modelling %continuous-time
non-linear dynamics that arise in numerous application fields, including, e.g., finance, physics and neuroscience \citep{kloe:92}. Over the past few decades, a large amount of research has contributed to methodological and theoretical advances
on the theme of parameter inference for SDEs. An overarching challenge is that the transition density of a non-linear SDE is in general intractable, thus appropriate proxies must be formulated to conduct likelihood-based inference. 
%The objective of this work is to 
We propose a new closed-form (CF) transition density expansion for SDEs, which can approximate the true density with high precision.  
%The development of such an approximation can lead to an effective parameter estimation for SDEs incorporated within various likelihood-based inference methods. 
In contrast to previous approaches, one of the novelties of our methodology is that it covers a broad class of diffusion processes, including \emph{hypo-elliptic} SDEs, 
i.e.~processes with a degenerate diffusion matrix and a transition law that still admits a density with respect to (w.r.t.) the Lebesgue measure. 
Hypo-elliptic SDEs appear in broad areas of applications (including physics, neuroscience) and parameter inference for these models has been a very active area of research in the last years.  
%\noindent 
%\textit{Model classes.}  

Let $B_t = (B_{1,t}, \ldots, B_{d, t})$, $t \ge 0$, be the standard $d$-dimensional Brownian motion, $d\ge 1$, defined upon the filtered probability space $(\Omega, \mathcal{F}, \{\mathcal{F}_t \}_{t \ge 0}, \mathbb{P})$. We consider $N$-dimensional SDEs, $N \ge 1$, of the following general form:
\begin{align} 
\label{eq:sde}
d X_t = V_0 (X_t, \theta)\,dt 
+ \sum_{1 \le j \le d} V_j (X_t, \theta) \, d B_{j,t}, 
\quad X_0 = x_0 \in \mathbb{R}^N,
\end{align}
for parameter vector $\theta \in \Theta \subseteq \mathbb{R}^{N_\theta}$, $N_\theta \ge 1$, and functions  $ V_j : \mathbb{R}^N \times \Theta \to \mathbb{R}^N, \,  0 \le j \le d$.  We set  
$\sigma = [V_{1}, \ldots, V_{d}]$ and $a = \sigma  \sigma^\top$.
% 
% In our contribution, the matrix $a=\sigma \sigma^{\top}\in\mathbb{R}^{N\times N}$, with $\sigma= [V_1, \ldots, V_d]\in\mathbb{R}^{N\times d}$,
% is not required to be of full rank. 
Our work focuses on two model classes, covering a large set of SDEs used in applications. The first class is the \emph{elliptic} one, where we consider SDEs of the following form: 
\begin{align} \label{eq:ellip}
\begin{aligned}
d X_{t} = V_{R, 0} (X_t, \theta) dt + \sum_{1 \le j \le d} V_{R, j} (X_t, \theta) d B_{j, t}, \quad X_0 = x_0 \in \mathbb{R}^N,  
\end{aligned} \tag{\textrm{E}}
\end{align}
so that $V_{j} = V_{R, j}$, $0 \le j \le d$. We set  
$\sigma_R = [V_{R,1}, \ldots, V_{R, d}]$, $a_R = \sigma_R  \sigma_R^\top$, and assume that $a_R=a_R (x, \theta)$ is positive definite for all $(x, \theta) \in \mathbb{R}^N \times \Theta$. Thus, w.l.o.g.~here $d=N$.
Class (\ref{eq:ellip}) 
% refers to the class of elliptic SDEs, which 
includes a multitude of models used in applications, see e.g.~\cite{kloe:92}. 
%Under uniform ellipticity and further regularity assumptions, it can be shown that the law of the solution to (\ref{eq:ellip}) at a finite time $t>0$ admits a density w.r.t.~the Lebesgue measure. 
%the Heston stochastic volatility model \citep{heston:93} and the stochastic Lorenz system in physics. %a%re The solution process characterised via this SDE is referred to as \emph{elliptic} diffusion. For instance, Heston stochastic volatility model \citep{heston:93} in finance and the stochastic Lorenz system in physics are classified within this class.  
% 
% \item[2.] \textbf{Hypo-elliptic Class.}  
The second model class we work with is the \emph{hypo-elliptic} one, where the SDE in (\ref{eq:sde}) now splits into smooth and rough components as $X_t = (X_{S, t}, X_{R, t}) \in \mathbb{R}^{N_S+ N_R}$, 
so that $N=  N_S + N_R$, $N_S \ge 1$, $N_R \ge 1$, and we re-express (\ref{eq:sde}) as:  
\begin{align} 
\begin{aligned} 
d X_{S,t} & = V_{S, 0} (X_t, \theta) \, dt; \qquad d X_{R,t}  = V_{R, 0} (X_t, \theta) \, dt + \sum_{1 \le j \le d} V_{R, j} (X_t, \theta) \, dB_{j, t},\\ 
X_0  &= x_0 = (x_{S,0}, x_{R,0}) \in \mathbb{R}^{N_S+ N_R}. 
\end{aligned}
%\qquad  X_0  = x_0 = (x_{0, S}, x_{0, R}) \in \mathbb{R}^{N_S+ N_R}. 
\label{eq:hypo} 
\tag{\textrm{H}} 
\end{align}
% 
% with an initial value $X_0 = (x_S, x_R) \in \mathbb{R}^{N_S} \times \mathbb{R}^{N_R}$. 
In (\ref{eq:hypo}), the involved functions are defined as  % 
$
V_{S, 0} : \mathbb{R}^N \times \Theta \to \mathbb{R}^{N_S}, \, 
V_{R, j} : \mathbb{R}^N \times \Theta \to \mathbb{R}^{N_R}, \,  0 \le j \le d. 
$ 
Model class (\ref{eq:hypo}) stems from the generic form (\ref{eq:sde}), where we now have that, for $(x, \theta) \in \mathbb{R}^N \times \Theta$:
\begin{align*}
V_0 (x, \theta) 
=  
\bigl[ 
V_{S, 0} (x, \theta)^\top,   
V_{R, 0} (x, \theta)^\top  
\bigr]^\top , \qquad
V_j (x, \theta) 
=
\bigl[
\mathbf{0}_{N_S}^\top, 
V_{R, j} (x, \theta)^\top 
\bigr]^\top 
, \qquad 1 \le j \le d. 
\end{align*}
Notice that component $X_{S, t}$ is not driven by the Brownian motion, and consequently class (\ref{eq:hypo}) requires a separate treatment from (\ref{eq:ellip}). Later on, we introduce sufficient requirements associated with the \emph{weak H\"ormander's condition}, so that $V_{S, 0} (X_t, \theta)$ depends on $X_{R, t}$, thus Brownian noise propagates 
into the smooth component, and the law of $X_{t}$, $t>0$, admits a density w.r.t.~the Lebesgue measure. 
% Such a process is referred to as a \emph{hypo-elliptic} SDE. 
Hypo-elliptic models are used in several application fields, including, e.g.: the FitzHugh-Nagumo SDE \citep{devil:05} and the Jansen-Rit neural mass SDE \citep{abl:17} in neuroscience;  the underdamped or generalised Langevin equation \citep{pavl:14} in physics.  
%For this model class, we later introduce sufficient assumptions, associated with the \emph{weak H\"ormander's condition}, so that the law of the solution of SDE (\ref{eq:hypo}) similarly admits a Lebesgue density. 

%\noindent \textit{Overview.} 
We consider parameter inference for SDEs given a collection of  discrere-time data $\{X_t\}_{t \in \mathbb{T}_n}$, for the set of time instances $\mathbb{T}_n = \{t_0, t_1, \ldots,  t_n \}$, $n\ge 1$. 
For simplicity, we consider equidistant step-sizes, with $\Delta: = t_i - t_{i-1}$. %is assumed to be fixed throughout, for simplicity. 
The likelihood function  is given as: 
$$ L_n (\{X_t\}_{t \in \mathbb{T}_n}; \theta) = p (X_{t_0}; \theta) \prod_{1 \le i \le n} p_{\Delta}^X (X_{t_{i-1}}, X_{t_i}; \theta),$$
for some initial law $p(\cdot;\theta)$, 
where 
$x'\mapsto p_{\Delta}^X (x, x'; \theta)$ is the transition density of SDE (\ref{eq:sde}), with the latter being  in general unavailable in closed form. A practical standard approach to circumvent this intractability is by introducing a {time-discretisation} scheme and using the induced CF approximate transition density as a proxy for the true one. For instance, a common scheme is the \emph{Euler-Maruyama} one, which yields a conditionally Gaussian approximate density upon application to elliptic SDEs. However, it is well-understood that such an approximation 
cannot correctly capture the true non-linear dynamics unless the step-size $\Delta$ is  
close to $0$. Thus,  parameter estimation relying on such a simple Gaussian approximation requires a \emph{high-frequency} observation regime, with  $\Delta \ll 1$. In practice, the step-size of available data is usually fixed, and its value may not be too small. 

In the context of fixed~$\Delta$, the prominent work of \cite{ait:02} proposes an elaborate approximation of the transition density for time-homogeneous univariate elliptic SDEs via a CF Hermite-series expansion. Roughly, the expansion for the transition density is of the structure:
\begin{align} 
\label{eq:expansion_intro}
{p_\Delta^X (x, y; \theta)}
\approx
%\quad 
{q_\Delta (x, y; \theta)  \times 
\bigl\{ 1 + (\mathrm{correction \, term}) \bigr\}}.
\end{align} 
Here, $y \mapsto q_\Delta (x, y; \theta)$ is a `baseline' tractable density. The `correction term' is given in closed-form, and includes Hermite polynomials up to a degree $J \ge 1$, obtained via working with $q_\Delta (x, y; \theta)$. 
The correction term plays a key role in capturing non-linear/non-Gaussian effects in the true transitions. 
%thus the CF expansion delivers a more accurate approximation of $p_\Delta^X (x, y; \theta)$, compared with using only the {baseline} Gaussian density. 
In detail, \cite{ait:02} constructs the 
CF-expansion by first applying an 1--1 `Lamperti transform' \citep{robe:01}, thus replacing the original scalar  $X_t$ with a process $Y_t$ of unit diffusion coefficient, and then obtaining the Hermite series expansion for the transition density of~$Y_t$. 
We refer to this line of research as the \emph{Hermite approach}.   
\cite{ait:02} proves convergence of the 
CF-expansion to the true density for fixed $\Delta\in (0,1)$ as the degree of Hermite polynomials, $J$, grows to infinity. The result is a qualitative one, as no order of convergence is provided. 
% After the groundbreaking work of \cite{ait:02}, numerous extensive CF expansions in the form of (\ref{eq:expansion_intro}) have been proposed.
The Hermite approach works only for the sub-class of `reducible' elliptic SDEs for which the  Lamperti transform is applicable. Also, as stated in \cite{ait:08}, convergence of the Hermite series expansion is not guaranteed when back-transforming onto the original density of~$X_t$. 
%\cite{ait:08} proposes an alternative approach to treat a wider class of non-reducible 
%multivariate elliptic SDEs.
%where the Lamperti transform is unavailable, i.e., \emph{irreducible diffusions}.  
To treat a wider class of non-reducible 
multivariate elliptic SDEs, \cite{ait:08} utilises the Kolmogorov backward/forward equations (PDEs) to construct a series expansion in $\Delta$ and $y - x$. No analytical results are provided for fixed~$\Delta$.
%contains a sanity-check type of result, by proving convergence of the approximate likelihood built from the CF expansion in probability as $\Delta \to 0$, for fixed datasize $n$. 
We refer to this contribution as the \emph{PDE approach}. 
%{However, we note that the above result is an asymptotic analysis, and thus, the effect of expansion on the precision of approximation is not clarified.} 
\cite{li:13} develops a \emph{probabilistic approach}, by making use of Malliavin calculus
and carrying out an asymptotic analysis of Wiener functionals \citep{wata:87, yoshi:92-2} to obtain a CF-expansion, accompanied by an analytic bound for the approximation error, for fixed $\Delta$. The expansion is given in terms of powers $\Delta^{j/2}$, $j\ge 0$, for $\Delta \in (0,1)$. More precisely: 
\begin{align} 
\label{eq:delta_expansion}
p_\Delta^X (x, y; \theta)
=  
q_\Delta (x, y; \theta)  \times 
\Bigl\{ 1 + \sum_{1 \le j \le J} 
\Delta^{j/2} \cdot e^{(j)}_{\Delta} (x, y; \theta) 
\Bigr\} 
+ R (\Delta, x, y; \theta), \qquad J \ge 1, 
\end{align}
for tractable coefficients $e_{\Delta}^{(j)}(\cdot)$, $j\ge 1$, and a remainder term $R(\cdot)$. \cite{li:13} proves under conditions that the remainder is of size $\mathcal{O} (\Delta^{(J+1-N)/2})$. % uniformly in $x, y,\theta$. 
%over an appropriate state and parameter space, where recall that $N$ is the dimension of $X_t$. The expansion coefficient $e^{(k)}$ involves conditional expectations of iterated It\^o-integrals, and these have an explicit form as a linear combination of Hermite polynomials, as pointed out in \cite{yang:19}. 
The probabilistic approach is extended to elliptic SDEs with jumps in \cite{li:16}. For time-inhomogeneous elliptic SDEs, \cite{choi:15} develops a CF-expansion via the PDE approach, similarly to \cite{ait:08}. \cite{yang:19}  use  It\^o-Taylor expansions and obtain a series of the form (\ref{eq:delta_expansion}) that involves Hermite polynomials, with explicit bounds provided on residuals as in \cite{li:13}.
%{where they also gave a non-asymptotic error bound as \cite{li:13} did}. 
Even if alternative approaches have been followed in the literature, the produced expansions are closely related to each other. 
E.g., one can obtain a series expansion as in (\ref{eq:delta_expansion}) involving Hermite polynomials via the two different approaches in \cite{li:13, yang:19}. 
Furthermore, \cite{lee:14} show that the Hermite expansion of \cite{ait:02} can be expressed in the form  (\ref{eq:delta_expansion}) by rearranging terms in the expansion w.r.t.~powers $\Delta^{j/2}$, $j\ge 1$.

Importantly, for developed CF-expansions to be \emph{theoretically validated}, the remainder terms should be controlled and vanish. This property guarantees convergence of the expansion, with a rate in $\Delta$ that grows when more terms are used in the expansion.
As mentioned, such an elaborate analysis has been carried out in  \cite{li:13, yang:19} in the context of elliptic SDEs. 
% The detailed connection between existing CF expansions in the literature is explained in \cite{yang:19}. 

The aforementioned works also demonstrate the effective use of a CF-expansion within parameter inference procedures. In particular, the approaches provide an approximate Maximum Likelihood Estimator (MLE). Obtained numerical results showcase that: the proxy MLEs stays close to the true ones even when the step-size $\Delta$ is not too small; the CF-expansions outperform proxy methods based on Gaussian-type quasi-likelihoods. 
\cite{chang:11} provide analytical consistency and convergence rate results for the proxy MLE, and demonstrate good performance of their CF-expansion by clarifying the effect of the length of  the expansion and of the fixed step-size $\Delta \in (0,1)$. In the context of Bayesian inference for SDEs, \cite{stra:10} utilise the expansion-based likelihood and show advantages over the Gaussian-type 
(Euler-Maruyama-based) likelihood.  

Critically, the development of CF-expansions in the literature is so far restricted to elliptic SDEs {and a limited class of hypo-elliptic ones, e.g., with linear drift and constant diffusion coefficient \citep{bar:17, Hab:19}, thus general hypo-elliptic SDEs specified as (\ref{eq:hypo}) have yet to be covered} 
% and 
% does not cover hypo-elliptic ones, 
even though the latter are widely used in applications.
%e.g.~in molecular dynamics, physics, neuroscience and ecology. 
% This limitation is intrinsic, in the sense that 
Available methods for elliptic SDE build upon steps that cannot be readily extended to the hypo-elliptic setting.
In brief, one limitation derives from the definition of the reference Gaussian density $q_\Delta (x, y; \theta)$ in the expansion relying on positive definiteness of its covariance matrix, when such a property is violated within the hypo-elliptic class  (\ref{eq:hypo}). Our work develops a novel CF-expansion that covers both elliptic and hypo-elliptic SDEs in a unified framework. To this end, we consider a \emph{non-degenerate} baseline Gaussian density $q_\Delta (x, y; \theta)$ that is well-defined for both SDE classes, (\ref{eq:ellip}) and (\ref{eq:hypo}). We then construct a CF-expansion in the form of (\ref{eq:delta_expansion}) based on such a well-posed $q_\Delta (x, y; \theta)$. We emphasise that the error analysis is much more challenging in the hypo-elliptic setting than in the elliptic one, due to varying scales across the SDE co-ordinates. 
We manage to provide analytical error estimates for the proposed expansion by utilising a recent result on estimates of the transition density for degenerate SDEs \citep{piga:22}, thus theoretically validating our CF-expansion both within the elliptic and the hypo-elliptic classes of SDEs.  

Beyond the above-mentioned literature on CF-expansions for elliptic SDEs with fixed $\Delta\in (0,1)$, our work is also motivated by several recent developments in the area of parametric inference for hypo-elliptic SDEs, 
albeit in a \emph{high-frequency observation regime}, i.e.~$n \to \infty$, $\Delta=\Delta_n \to 0$, $n \Delta_n \to \infty$, together with an extra `design' condition on $\Delta_n$. 
Indicatively, \cite{dit:19, mel:20, glot:21, pil:24-II} propose contrast estimators, under the design condition $\Delta_n = o (n^{-1/2})$. The latter is weakened to $\Delta_n = o (n^{-1/3})$ and $\Delta_n = o (n^{-1/p})$, for a general integer $p \ge 2$, by \cite{igu:BJ} and \cite{igu:23}, respectively. 
\cite{igu:SPA} also treat a class of `highly degenerate' hypo-elliptic SDEs.

Our main contributions are briefly summarised as follows:     
\begin{enumerate}
\item[a.] We propose a CF-expansion for the transition density of both elliptic and hypo-elliptic SDEs, in (\ref{eq:ellip}) and (\ref{eq:hypo}), respectively. 
Within the elliptic class, a starting point for developing the CF-expansion is motivated by the work of one of the co-authors in \cite{igu:21-2}. This latter work lies in the area of numerical methods for SDEs and looks at the development of approximation schemes for elliptic SDEs of improved weak order of convergence.  
%by extending a methodology proposed in \citep{igu:21-2}. 
%
%Our expansion is applied to both elliptic and hypo-elliptic diffusions specified as (\ref{eq:ellip}) and (\ref{eq:hypo}). 
To the best of our knowledge, this is the first time in the literature that a CF-expansion is obtained for hypo-elliptic SDEs with a general form of coefficients.   

%the construction of density expansion covering the hypo-elliptic case  is novel.
% 
\item[b.] Our proposed CF-expansion involves a linear combination of differential operators acting on an appropriately chosen baseline Gaussian density, thus is easily computable via available software with symbolic calculations. 
Though we initially obtain an expression of different structure from~(\ref{eq:delta_expansion}), we later show that our CF-expansion indeed  takes up the form of (\ref{eq:delta_expansion}), i.e.~a series expansion in powers of $\sqrt{\Delta}$. Thus, our 
CF-expressions align with existing works for elliptic SDEs. 
\item[c.] We theoretically validate our CF expansions by proving analytically,  under appropriate conditions, that the residuals are of size $\mathcal{O} (\Delta^{K/2})$ for a step-size  $\Delta \in (0,1)$, where $K \ge 1$ is an integer differing between the elliptic and hypo-elliptic classes, and which depends on the model dimension $N$. In particular, the effect of the dimensionality varies amongst the two SDE classes. %the hypo-elliptic SDEs differs from that of elliptic SDEs. 
\item[d.] We present numerical results showcasing that the use of the proposed CF-expansion leads to effective parameter estimation for SDEs, with an emphasis on hypo-elliptic models. In particular, we conduct Bayesian inference for a real dataset and show that the posterior distribution is accurately estimated by the proposed density expansion. 
% compared to the case of using a Gaussian-type approximation with no extra correction terms. 
\end{enumerate}

%As a prelude, Fig.~\ref{fig:ae_I} shows a rapid decrease in absolute error for the CF-expansion
%when moving from the baseline Gaussian density (corresponding here to $J=2$) to $J=5$, for 3 choices of $\Delta$. The model used is the hypo-elliptic \emph{FitzHugh-Nagumo} (FHN) SDE, with full details given in Section \ref{sec:num}.

%\noindent {\it Other related works.} 
%\noindent 
%{\it Organisation of the paper.}
The structure of the paper is as follows.
In Section \ref{sec:CF} we outline our strategy for constructing a CF-expansion which covers both (\ref{eq:ellip}), (\ref{eq:hypo}), and then proceed with the development of the expansion. Section \ref{sec:error} provides a rigorous error analysis for the proposed CF-expansion, separately for classes (\ref{eq:ellip}) and (\ref{eq:hypo}). Section \ref{sec:num} shows numerical applications, and the codes that reproduce the results are available at \url{https://github.com/YugaIgu/CF-density-expansion}. Section~\ref{sec:conclusion} provides a summary and conclusions. Most proofs are collected in a Supplementary Material.

\noindent \textit{Notation.} 
We set $\mathbb{Z}_{\ge 0} = \mathbb{N} \cup \{0\}$. 
For a multi-index $\alpha \in \mathbb{Z}_{\ge 0}^k$,  $k \in \mathbb{N}$, we write 
$\| \alpha \|= k$,  
$| \alpha | = \textstyle{\sum_{1 \le j \le k} \alpha_j}$,   
$\textstyle \| \alpha \|_\infty=  \max_{1 \le i \le k} \alpha_i$, \  $\alpha! = \textstyle{\prod_{j = 1}^k} \alpha_j.$ 
For $\alpha \in \mathbb{Z}_{\ge 0}^m$ and a sufficiently smooth $\varphi : \mathbb{R}^m \to \mathbb{R}, \, m \in \mathbb{N}$, 
we write $\partial^\alpha \varphi (y) 
= \partial^{\alpha_1}_{y_1} \cdots \partial^{\alpha_m}_{y_m} \varphi (y)$, where 
$\partial^{\alpha_i}_{y_i} \equiv \partial^{\alpha_i} / \partial y^{\alpha_i}_{i}$. 
We often write 
$\partial^\alpha_y \varphi (y) \equiv \partial^\alpha \varphi (y)$ to emphasise the argument upon which the derivative acts. The generator associated with SDE (\ref{eq:sde}) writes as: 
\begin{align}
\mathscr{L}_\theta \varphi (x) = \sum_{1 \le i \le N} V_0^i (x, \theta)\,\partial_i \varphi (x) 
+ \tfrac{1}{2} \sum_{1 \le i_1, i_2 \le N} 
\sum_{1 \le j \le d} V_j^{i_1} (x, \theta) V_j^{i_2} (x, \theta) \,
\partial_{i_1 i_2} \varphi (x), \label{eq:generator}
\end{align}  
$(x, \theta) \in \mathbb{R}^N \times \Theta$, 
for $\varphi : \mathbb{R}^N \to \mathbb{R}$, where we use integer superscripts to indicate co-ordinates in vectors. For $x \in \mathbb{R}$, we write $\floor*{x} = \max \{m \in \mathbb{Z} \, | \, m \le x \}$. 
For differential operators $D_1, D_2$, we define the \emph{commutator} as 
$
\mathrm{ad}_{D_1} (D_2) = [D_1, D_2] 
\equiv D_1 D_2 - D_2 D_1. 
$ 
The $k$-times iteration of the commutator writes as
${\mathrm{ad}_{D_1}^k (D_2) = [D_1, \mathrm{ad}_{D_1}^{k-1} (D_2)]}$, $k\ge 1$,  
%$   \mathrm{ad}_{D_1}^k (D_2) = \underbrace{[D_1, [D_1, \ldots [D_1}_{k-\mathrm{times}}, D_2]]]$.    
% \begin{gather*}
%     \mathrm{ad}_{D_1}^k (D_2) = \underbrace{[D_1, [D_1, \ldots [D_1}_{k-\mathrm{times}}, D_2]]]  
% \end{gather*}  
% 
with $\mathrm{ad}_{D_1}^0 (D_2)= D_2$. 
% 
%\begin{figure}[!h] 
%\centering
%\includegraphics[width=12cm]{Contour-density-abs-err-DE-I-set2.pdf}
%%\begin{}
%%  \centering
%%  \includegraphics[width=1\linewidth]{Contour-density-abs-err-DE-I-delta-0.05-set2.png}
%%\end{subfigure}
%%\begin{subfigure}
%%  \centering
%%  \includegraphics[width=1\linewidth]{Contour-density-abs-err-DE-I-delta-0.02-set2.png}
%%\end{subfigure}
%\caption{\small Heatmap of the absolute error for the CF-expansion in the case of the hypo-elliptic FHN model (see Section \ref{sec:num}). Rows correspond to 3 choices $\Delta=(0.1,0.05,0.02)$ and columns to 4 choices $J=2,3,4,5$. 
%% The row corresponds to the step-size $\Delta= 0.1, 0.05, 0.02$ from top to bottom, and the column does to the order of density expansion $J = 2, 3, 4, 5$ from left to right.
%}
%\label{fig:ae_I}
%\end{figure}

% 
% 
\section{Closed-Form Transition Density Expansion}
\label{sec:CF}
We will present a new CF transition density expansion for a wide class of It\^o processes in (\ref{eq:sde}), including the family of hypo-elliptic SDEs specified in (\ref{eq:hypo}). We write the transition density of $X_{t + \Delta}$ given $X_t = x \in \mathbb{R}^N$ as
$ y \mapsto p^X_\Delta (x, y; \theta) :=  
\mathbb{P} ( X_{t + \Delta} \in dy \,| \, X_t = x )/dy 
$, with $t \ge 0$, $\Delta > 0$.
\subsection{Conditions for Closed-Form Expansion} 
\begin{assump} 
\label{ass:diff} 
For both model classes (\ref{eq:ellip}) and (\ref{eq:hypo}), the maps $x \mapsto V_j (x, \theta)$, $0 \le j \le d$, are infinitely differentiable for any $\theta \in \Theta$. 
\end{assump}
\noindent For a vector-valued $V:\mathbb{R}^{N}\to \mathbb{R}^{N}$, we make use of the standard correspondence $V \leftrightarrow \textstyle\sum_{i=1}^{N} V^{i}\partial_{i}$.
\begin{assump} 
\label{ass:hor}
We distinguish between model classes (\ref{eq:ellip}) and (\ref{eq:hypo}).
\begin{itemize}
\item[I.] For class (\ref{eq:ellip}), it holds that $a_R (x, \theta) = (\sigma_R \sigma_R)^{\top} (x, \theta) $ is positive definite for all $(x, \theta) \in \mathbb{R}^N \times \Theta$. This is equivalent to 
$
\mathrm{Span}\{ V_{R, j} (x, \theta), \, 1 \le j \le d \} = \mathbb{R}^N, 
$ 
for all $(x, \theta) \in \mathbb{R}^N \times \Theta$. 
\item[II.] For class (\ref{eq:hypo}), it holds that:
\begin{gather*}
\mathrm{Span} \bigl\{ V_{R, j} (x, \theta), \,  1 \le j \le d  \bigr\} = \mathbb{R}^{N_R}, \qquad 
\mathrm{Span} \Bigl\{ \bigl\{ V_{j} (x, \theta), \, [\widetilde{V}_{0}, V_j] (x, \theta) \bigr\}, \,  1 \le j \le d  \Bigr\} = \mathbb{R}^{N},
\end{gather*}
for all $(x, \theta) \in \mathbb{R}^N \times \Theta$, where $\widetilde{V}_0: \mathbb{R}^N \times \Theta \to \mathbb{R}^N$ is the drift function when the It\^o SDE (\ref{eq:hypo}) is written in a Stratonovich form, namely
$\textstyle \widetilde{V}_0 (x, \theta) 
= V_0 (x, \theta) - \tfrac{1}{2} \sum_{j = 1}^d \sum_{i = 1}^N V_j^i (x, \theta) \partial_{x_i} V_j (x, \theta).
$  
% = 
% \begin{bmatrix}
% V_{S, 0} (x, \theta) \\[0.1cm]  
% V_{R, 0} (x, \theta) - \tfrac{1}{2} 
% \sum_{j = 1}^d \sum_{i = 1}^N V_j^i (x, \theta) \partial_{x_i} V_{R, j} (x, \theta)    
% \end{bmatrix}. 
% \label{eq:tilde_V0} 
% \end{align}  
%
\end{itemize}
\end{assump}
Assumption \ref{ass:hor} is related to \emph{H\"ormander's condition} (it suffices for H\"ormander's condition to hold) and implies that the law of $X_t$, $t>0$, admits a Lebesgue density. For the hypo-elliptic class (\ref{eq:hypo}), Assumption \ref{ass:hor}-II guarantees that the %Brownian %(Brownian motion $B_t$) 
noise in the rough component $X_{R,t}$ (of size $\sqrt{t}$ for a period of length $t$) propagates into the smooth component $X_{S, t}$.
%with the time-scale $\sqrt{t}$ in the rough component $X_{R,t}$ propagates into the smooth component $X_{S, t}$. 
Inclusion of the vector fields $[\widetilde{V}_0, V_j] (x, \theta)$, $1 \le j \le d$, relates to the appearance of terms $\textstyle \int_0^t B_s^j ds$ (of a different scale $\sqrt{t^3}$) in the smooth component $X_{S, t}$ after an
It\^o-Taylor expansion of $\textstyle\int_0^t V_{S, 0} (X_u, \theta) du$. Thus, the transition law of SDE~(\ref{eq:hypo}) is non-degenerate and admits a Lebesgue density even if not all coordinates are directly driven by the Brownian motion. 
A precise  definition of H\"ormander's condition can be found, e.g., in~\cite{nual:06}.      

\subsection{Background Idea} 
\label{sec:background}
Before presenting the CF-expansion we explain an idea that underpins its development -- more precisely the starting point of the latter. 
Consider the elliptic class (\ref{eq:ellip})
and the Euler-Maruyama (EM) scheme which  approximates the transition dynamics of $X_{t + \Delta} | X_t = x$, with $x\in \mathbb{R}^N$, $t > 0$, $\Delta > 0$, so that:
\begin{align}
\label{eq:EM}
\bar{X}_{t + \Delta}^{\mathrm{EM}, \theta} 
:= x + V_{R,0} (x, \theta) \Delta + \sigma_{R} (x, \theta) \bigl( B_{t + \Delta} - B_t  \bigr).  
\end{align}
Under regularity conditions on $x \mapsto V_{R,j} (x, \theta)$, $0 \le j \le d$, and the requirement that the matrix 
$a_R  (x, \theta) = (\sigma_R \sigma_R^\top)(x, \theta)$ 
is positive definite for all $(x, \theta) \in \mathbb{R}^N \times \Theta$, 
the EM scheme gives rise to a well-defined baseline Gaussian transition density, $y \mapsto p_\Delta^{\bar{X}^{\mathrm{EM}}} (x, y; \theta)$.
In the present elliptic setting \cite{igu:21-2} constructed a CF transition density approximation of the following form: 
\begin{align}
\label{eq:IY-CF}
p^X_\Delta (x, y; \theta)
\approx 
p^{\bar{X}^{\mathrm{EM}}}_\Delta (x, y; \theta) \times \bigl(1 + (\mathrm{correction \, term}) \bigr).
\end{align}   
The tools utilised in \cite{igu:21-2} to derive the approximation include Taylor expansion, Kolmogorov backward/forward equations, use of the infinitesimal generators for the target SDE and its EM approximation.  
In the above expression, the `correction term' involves $\Delta$, partial derivatives of the SDE coefficients and Hermite polynomials obtained via differentiating the transition density of the EM scheme. As mentioned in the Introduction, other approaches are also available, including the ones developed in \cite{ait:02, ait:08, li:13, yang:19}, and all such works also assume invertibility of the matrix $a_R$, thus are not relevant for the hypo-elliptic 
class (\ref{eq:hypo}). 

To construct a CF transition density expansion for a broader family of SDEs that includes hypo-elliptic SDEs, it is critical to choose an appropriate  reference  Gaussian density which is \emph{non-degenerate} for the target class of models.
To achieve this, 
%we consider a baseline Gaussian transition density 
%(rather than EM-type density) 
%that is well-defined for wider classes of SDEs, including~(\ref{eq:hypo}),  
% Motivated by \cite{mel:20}, we 
we consider the \emph{local drift linearisation (LDL) scheme}, which, upon application on the general SDE model in (\ref{eq:sde}), is defined via the following expression, for each given $t\ge 0$, $\Delta>0$, and for  $\bar{X}_{t}^{\theta}=x\in\mathbb{R}^{N}$: 
% 
%\begin{align} 
%\label{eq:LDL}
%d\bar{X}_{t}^{\theta} 
%    = 
%    \bigl( 
%        A_{x, \theta} \bar{X}_{t}^\theta 
%        + b_{x, \theta} 
%    \bigr) ds 
%    + \sigma (x, \theta)\, 
%    dB_{t}, \qquad \bar{X}_{0}^{\theta} =  x, 
%\end{align} 
%
\begin{align} 
\label{eq:LDL}
\bar{X}_{t + \Delta}^{\theta} 
= x + \int_t^{t + \Delta}
\bigl( 
A_{x, \theta} \bar{X}_{s}^\theta 
+ b_{x, \theta} 
\bigr) ds 
+ \sigma (x, \theta) 
\bigl(B_{t + \Delta} - B_{t} \bigr), 
\end{align} 
where $A_{x, \theta}\in \mathbb{R}^{N\times N}$
and $b_{x, \theta}\in\mathbb{R}^{N}$ are specified as follows: 
\begin{align*}
A_{x, \theta} = \bigl[ \partial_{x_j}  V_0^i (x, \theta) \bigr]_{1 \le i, j \le N}, \qquad 
b_{x, \theta} = V_0 (x, \theta) 
- A_{x, \theta} \, x.  
\end{align*} 
That is, (\ref{eq:LDL}) is obtained from a 1st-order Taylor expansion of the drift $V_0$ about the initial position $x$ and $\sigma(\cdot)$ fixed at its initial value. Expression (\ref{eq:LDL}) corresponds to a linear SDE, with a solution for $\bar{X}_{t + \Delta}^\theta | \bar{X}_{t}^\theta = x$ that has the following explicit form:  
\begin{align*}
\bar{X}_{t + \Delta}^\theta 
= e^{\Delta A_{x, \theta}} x 
+ \int_t^{t + \Delta} e^{(t + \Delta - s ) A_{x, \theta}} b_{x, \theta}\,ds
+ 
\int_t^{t + \Delta} e^{(t + \Delta - s ) A_{x, \theta}} \sigma (x, \theta)\,dB_s. 
\end{align*} 
Thus, $\bar{X}_{t + \Delta}^\theta | \bar{X}_{t}^\theta = x$ follows a Gaussian law, with mean $\mu(\Delta, x, \theta)$ and covariance $\Sigma (\Delta, x, \theta)$ given as: 
\begin{align}
& \mu (\Delta, x, \theta)
= e^{\Delta A_{x, \theta}} \hat{\mu}^x (\Delta, x, \theta), \qquad 
\hat{\mu}^z (\Delta, x, \theta) 
:=  x + \int_0^{\Delta} e^{- s A_{z, \theta}} b_{z, \theta} ds, \quad z \in \mathbb{R}^N; \label{eq:mu} \\
& \Sigma (\Delta, x, \theta) 
=
e^{\Delta A_{x, \theta}} \hat{\Sigma} (\Delta, x, \theta) e^{\Delta A_{x, \theta}^\top}, \qquad   
\hat{\Sigma} (\Delta, x, \theta)  
:= \int_0^\Delta e^{- s A_{x, \theta}} a (x, \theta) e^{- s A_{x, \theta}^\top} ds,   
\label{eq:Sigma} 
\end{align} 
{where we recall $a = \sigma \sigma^\top$.} The introduction of the extra argument $z$ in $\hat{\mu}^z (\Delta, x, \theta)$ will be of use in later developments.

%  

% for every $(x, \theta) \in \mathbb{R}^N \times \Theta$ as:
% % 
% \begin{align*}
% p_\Delta^{\bar{X}^{\mathrm{EM}}} (x, y; \theta) 
% = \varphi_{\mathcal{N}} 
% \bigl(y  \, ; \, x + \Delta V_0 (x, \theta),
%    \Delta a (x, \theta)
% \bigr),
% \end{align*} 
% where $\varphi_{\mathcal{N}} (\cdot; \mu, \Sigma)$ is the density of Normal distribution with mean $\mu$ and covariance $\Sigma$.   
% 
%   
%
% 

% 
%The differential form of the LDL scheme (\ref{eq:LDL}), i.e., 
% 
%\begin{align} 
%\label{eq:LDL_SDE}
%d \bar{X}_{t}^{\theta} = \bigl( A_{x, \theta} \bar{X}_{t}^{\theta} + b_{x, \theta}  \bigr) dt 
%+ \sum_{1 \le j \le d} \sigma_j (x, \theta) dB_{j, t}, \qquad \bar{X}_0^\theta = x \in \mathbb{R}^N, 
%\end{align}
% 

% 
% \begin{align} 
% \label{eq:LDL_cov}
% \hat{\mu}^z (t, x, \theta) 
%  :=  x + \int_0^{t} e^{- s A_{z, \theta}} b_{z, \theta} ds, \qquad 
%  \hat{\Sigma} (t, x, \theta)  
%  := \int_0^t e^{- s A_{x, \theta}} a (x, \theta) e^{- s A_{x, \theta}^\top} ds. 
% \end{align} 
% 
%

%{Exponential matrix can explode with a general param space. }
In the sequel we show that $\hat{\Sigma}$ (thus also $\Sigma$) is positive definite for both model classes (\ref{eq:ellip}), (\ref{eq:hypo}) under Assumption \ref{ass:hor}. 
Then, under regularity conditions on the SDE coefficients together with the invertibility $\Sigma$, we appropriately expand upon the direction followed by \cite{igu:21-2} to  construct a CF transition density approximation that covers the model class (\ref{eq:hypo}) and writes as:  
\begin{align} \label{eq:generic}
p^X_\Delta ( x, y; \theta) 
\approx 
p^{\bar{X}}_\Delta (x, y; \theta) 
\times \bigl(1 + (\mathrm{correction \, term}) \bigr), 
\end{align}  
where $y \mapsto p^{\bar{X}}_\Delta (x, y; \theta)$ is the transition density of the LDL scheme (\ref{eq:LDL}).
% given as: 
% % 
% \begin{align} \label{eq:density_LDL}
% p_\Delta^{\bar{X}} (x, y; \theta) 
% = \tfrac{1}{\sqrt{(2\pi)^N \det {\Sigma} (\Delta, x, y)}} \exp 
% \left( - \tfrac{1}{2} 
% \bigl(
% e^{- \Delta A_{x, \theta}} y - \hat{\mu} (\Delta, x, \theta)  \bigr)^\top 
% \hat{\Sigma}^{-1} (\Delta, x, \theta) 
% \bigl( 
% e^{- \Delta A_{x, \theta}} y - \hat{\mu} (\Delta, x, \theta) 
% \bigr) 
% \right)
% \end{align}
% % 
% for $x, y \in \mathbb{R}^N$ and $\theta \in \Theta$. 
% % 
%    
Similarly to the case of the expansion for elliptic diffusions, the correction term appearing in (\ref{eq:generic}) involves partial derivatives of SDE coefficients w.r.t.~the state argument and Hermite polynomials now defined via partial derivatives of the non-degenerate Gaussian density $p_\Delta^{\bar{X}} (x, y; \theta)$. 

\begin{remark}
\cite{igu:21-2} work in an elliptic setting to develop Monte-Carlo estimators of improved weak order of convergence for $\mathbb{E}[\varphi(X_{T})]$, $\varphi:\mathbb{R}^{N}\to \mathbb{R}$, $T>0$, and their expansion in the form of (\ref{eq:IY-CF}) is used for such a purpose. In brief, they use samples from the baseline $p_\Delta^{\bar{X}^{\mathrm{EM}}} (x, y; \theta)$, weighted by $(1 + (\mathrm{correction \, term}))$, in an iterative procedure over $\floor*{T/\Delta}$ steps.
Even if the initial derivations in the 
CF-expansion we develop here resemble steps followed in \cite{igu:21-2}, our objectives 
and, consequently, the structure of the 
CF-expansion and its theoretical analysis 
(and, in general, the overall contribution) 
fully deviate from \cite{igu:21-2}. 
\end{remark}
\begin{remark} \label{rem:vsLL} 
LDL scheme (\ref{eq:LDL}) differs from the so-called Local Linearisation (LL) scheme (its definition can be found, e.g., in \cite{jim:17}) in the sense that the latter applies a first order Taylor expansion for both drift and diffusion coefficients. As shown in the next subsection, in particular in Lemma \ref{lemma:LDL_hor}, the LDL scheme follows conditionally a non-degenerate Gaussian distribution that admits a transition density for (\ref{eq:ellip}) and (\ref{eq:hypo}) under Assumptions \ref{ass:diff}-\ref{ass:hor}, leading to the development of (\ref{eq:generic}). The key idea here is that the noise in the rough component $X_R$ propagates to the smooth component $X_S$ via the locally linearised drift. Similarly, the LL scheme can be shown to admit a well-defined tractable Lebesgue density, thus it could also form the basis of a transition density expansion like (\ref{eq:generic}). In our analysis below, we employ the LDL scheme since its density admits a simpler expression and suffices to build a tractable density expansion. Also, such a non-degenerate Gaussian approximation can be constructed via a drift linearisation for partial coordinates (not full) so that the matrix $A_{x, \theta}$ is upper-triangular, thus reducing the computation cost of calculating $\exp\{A_{x, \theta}\}$. We briefly discuss a practical choice of $A$ in the numerical experiment Section \ref{sec:num_DE} and in Section \ref{sec:discussion} as well. 
\end{remark} 
 
\subsection{Non-Degeneracy of the LDL Scheme}
As  mentioned in Section \ref{sec:background}, existence of a Lebesgue density for the transition dynamics  of the LDL scheme (\ref{eq:LDL}) is essential for the  construction of our CF-expansion for both model classes (\ref{eq:ellip}) and (\ref{eq:hypo}). 
In this section we show that such an existence is implied by Assumption \ref{ass:hor}. 
That is, 
%We use the above discussion to study the LDL scheme (\ref{eq:LDL}). 
we show under 
Assumption~\ref{ass:hor}, i.e.~H\"ormander's condition for the target model classes (\ref{eq:ellip}) and (\ref{eq:hypo}), that the vector fields defined via the LDL scheme (\ref{eq:LDL}) also satisfy H\"ormander's condition. Specifically, the vector fields defined from the coefficients of (\ref{eq:LDL}) coincide with those of the original SDE, upon fixing the argument $x$ to the initial condition of (\ref{eq:LDL}). %Thus, matrix $\hat{\Sigma}$ defined in (\ref{eq:Sigma}) is positive definite and  the transition density of the LDL scheme is well-defined.  
We thus have the following result whose proof is given in Appendix~\ref{app:pf_LDL_hor}. 

\begin{lemma} \label{lemma:LDL_hor}
Let $(\Delta, x, \theta) \in (0,\infty) \times \mathbb{R}^N \times \Theta$. Under Assumptions \ref{ass:diff}-\ref{ass:hor}, the law of $\bar{X}_{t + \Delta}^\theta | \bar{X}_t^\theta  = x$ defined in (\ref{eq:LDL}) admits a smooth Lebesgue density for model classes (\ref{eq:ellip}) and (\ref{eq:hypo}). 
\end{lemma}
% 
% The proof is postponed in Appendix \ref{app:pf_LDL}, where we show that the Hormander's condition holds for vector fields defined from the coefficients of the LDL scheme. Lemma \ref{lemma:LDL_hor} implies that under Assumption \ref{ass:hor} the matrix $\Sigma (\Delta, x, \theta)$ defined in (\ref{eq:Sigma}) is positive definite for any $(t, x, \theta) \in (0, \infty) \times \mathbb{R}^N \times \Theta$, and then the Gaussian density given in (\ref{eq:density_LDL}) is well-defined. % 
%
% In the next subsection, we will show that the covariance $\Sigma$ is indeed positive definite and then the density of LDL scheme is well-defined under Assumption \ref{ass:hor} for both model classes (\ref{eq:ellip}) and (\ref{eq:hypo}). 
We describe a sub-class from (\ref{eq:hypo}) where the LDL scheme delivers well-posed Lebesgue densities for SDE transition dynamics while the Euler-Maruyama scheme provides degenerate distributions. 
%\vspace{-0.8cm}
% 
\begin{ex}[Underdamped Langevin Equation]
We consider the following bivariate SDE: 
\begin{align}
\begin{aligned} \label{eq:ULE}
d X_t^1 & = X_t^2\,dt; \qquad d X_t^2 & = \bigl( - V' (X_t^1) - \alpha X_t^2 \bigr) dt + \sigma\,d B_{1, t},  
\end{aligned}
\end{align} 
for parameters $\alpha, \sigma > 0$ and potential $V : \mathbb{R} \to \mathbb{R}$. Such dynamics are used  to describe the motion of a particle on the real line $\mathbb{R}$, with $X_t^1$ and $X_t^2$ representing position and momentum, respectively.  The coefficients in SDE (\ref{eq:ULE}) correspond to the following vector fields, for $x = (x_1, x_2) \in \mathbb{R}^2$: 
% for 
\begin{align*}
V_0  \equiv \widetilde{V}_0 = x_2 \partial_{x_1} 
+ \bigl( - V' (x_1) - \alpha x_2 \bigr) \partial_{x_2}, \qquad V_1 = \sigma \partial_{x_2}. 
\end{align*}
The diffusion matrix is degenerate, so (\ref{eq:ULE}) belongs to class (\ref{eq:hypo}). Also, Assumption \ref{ass:hor}-II is satisfied as: 
\begin{align}
V_1 (x) = 
[0,  \sigma]^\top, 
\qquad 
[\widetilde{V}_0, V_1] (x) = \widetilde{V}_0 V_1 (x) - V_1 \widetilde{V}_0 (x) 
=  [- \sigma, \sigma\alpha]^\top,  
\label{eq:vf_true}
\end{align}
and, given $\sigma, \alpha > 0$, we get that 
$ \mathrm{Span} \{ V_1 (x), [V_0, V_1] (x) \} = \mathbb{R}^2,
$
for all $x \in \mathbb{R}^2$. The Euler-Maruyama scheme for $\bar{X}_{t+\Delta}^{\mathrm{EM}}|\bar{X}_t^{\mathrm{EM}} = x$ writes as: 
\begin{align}
%\begin{aligned}
\bar{X}_{t + \Delta}^{\mathrm{EM}, 1} 
= x_1 
+ x_2 \Delta; \qquad  
\bar{X}_{t + \Delta}^{\mathrm{EM}, 2} 
= x_2  + 
\bigl( - V' (x_1) - \alpha x_2  \bigr) \Delta 
+ \sigma (B_{1, t+ \Delta} - B_{1, t}).  
%\end{aligned} 
\label{eq:em} 
\end{align}
% \begin{align} 
% \begin{aligned}
% \begin{bmatrix}
% \bar{X}_{t + \Delta}^{\mathrm{EM}, 1}  \\
% \bar{X}_{t + \Delta}^{\mathrm{EM}, 2} 
% \end{bmatrix} 
% =
% \begin{bmatrix}
% x_1 
% + x_2 \Delta   \\
% x_2  + 
% \bigl( - V' (x_1) - \alpha x_2  \bigr) \Delta 
% \end{bmatrix}  
% + 
% \begin{bmatrix}
% 0 \\
% \sigma (B_{1, t+ \Delta} - B_{1, t}) 
% \end{bmatrix}. 
% \end{aligned}
% \end{align} 
% 
% 
So, the law of 
$\bar{X}_{t+\Delta}^{\mathrm{EM}} | \bar{X}_t^{\mathrm{EM}}$ involves a degenerate covariance matrix.
%and then its transition density can not be expressed as a well-defined Lebesgue density. 
In contrast, in this setting the LDL scheme~(\ref{eq:LDL}) contains the $2 \times 2$ matrix $A_{x}$ and the vector $b_{x}$ specified as follows: %for $x= (x_1, x_2) \in \mathbb{R}^2$,   
\begin{align*}
[A_x]_{11} = 0, \ \ 
[A_x]_{12} = 1,  \ \
[A_x]_{21} = - V'' (x_1), \ \
[A_x]_{22} = - \alpha,
\quad 
b^1_x = 0, \ \ 
b^2_x = - V' (x_1) + V'' (x_1) x_1.  
% \begin{bmatrix}
% 0 & 1 \\
% - V'' (x_1) & - \alpha 
% \end{bmatrix},
% \qquad 
% b_x = 
% \begin{bmatrix}
% 0 \\ 
% - V' (x_1) + V'' (x_1) x_1
% \end{bmatrix}, 
% \quad x= (x_1, x_2) \in \mathbb{R}^2. 
\end{align*}  
% 
% \begin{align} \label{eq:ex_LDL}
% d \bar{X}_t = \bigl( A_x \bar{X}_t + b_{x} \bigr) dt + 
% %\begin{bmatrix}
%   [0, \sigma]^{\top}
% %\end{bmatrix} 
% d B_{1 , t},  \qquad  \bar{X}_0 = x \in \mathbb{R}^2, 
% \end{align} 
% 
% where 
% 
% 
The vector fields associated with the coefficients of the LDL scheme are given as follows, 
for $x, z \in \mathbb{R}^2$:  
\begin{align*}
\bar{V}_0^z = \sum_{i = 1, 2} 
[A_z x  +  b_z]_i \partial_{x_i}, 
\qquad 
\bar{V}_1 = \sigma \partial_{x_2}, 
\end{align*} 
where (with some abuse of notation) we introduce $z \in \mathbb{R}^2$ to represent the initial condition for (\ref{eq:LDL}), thus distinguish the latter from argument $x \in \mathbb{R}^2$ upon which the linear drift in (\ref{eq:LDL}) applies.
%differential operators for argument `x', to stress that such operators will not act on $A_z$ and $b_z$. 
For the above vector fields, H\"ormander's condition holds via:  
\begin{align}
\bar{V}_1 (x) = 
\bigl[ 
0, \sigma]^\top, 
\qquad 
[\bar{V}_0^z, \bar{V}_1] (x) 
= - \bar{V}_1 \bar{V}_0^z (x)   
= \bigl[- \sigma, \sigma\alpha \bigr]^\top.  
\label{eq:vf_LDL}
\end{align}
Thus, the law of $\bar{X}_{t + \Delta} | \bar{X}_t$ admits a well-defined (Gaussian) transition density. Note that the vector fields in (\ref{eq:vf_LDL}) coincide with the ones in (\ref{eq:vf_true})  defined for the original SDE (\ref{eq:ULE}), so the hypo-ellipticity of the target SDE is inherited by its induced LDL scheme. 
% \begin{align*}
% \bar{X}_{t + \Delta} = e^{A_{\bar{X}_t}} \bar{X}_{t} 
% + \int_t^{t + \Delta} e^{(t + \Delta - s) A_{\bar{X}_t}} 
% \cdot 
% \begin{bmatrix}
% 0 \\
% - V' (\bar{X}_t^1) + V''(\bar{X}_t^1) \bar{X}_t^1 
% \end{bmatrix}
% ds 
% + \int_t^{t + \Delta} 
% e^{(t + \Delta - s) A_{\bar{X}_t}} \cdot 
% \begin{bmatrix}
% 0 \\
% \sigma
% \end{bmatrix}dB_{1, s},  
% \end{align*}
% 
\end{ex}

\subsection{Transition Density CF Expansion}
\label{sec:CF_expansion}
% 
% 
% \subsubsection{Nondegeneracy of LDL scheme}
% 
% 
\subsubsection{Preliminaries}
% We introduce some notation and auxiliary results to construct a CF expansion of transition density for both (\ref{eq:ellip}) and (\ref{eq:hypo}).  
We prepare some ingredients for the construction of our CF-expansion. We define the LDL scheme starting from a point $x \in \mathbb{R}^N$ with its coefficients being frozen at a point $z \in \mathbb{R}^N$ as: 
\begin{align} \label{eq:LDL_frozen}
d\bar{X}_{t}^{\theta, z} 
= 
\bigl( 
A_{z, \theta} \bar{X}_{t}^{\theta, z}  
+ b_{z, \theta} 
\bigr)dt 
+ \sigma (z, \theta) 
dB_t, \qquad \bar{X}_{0}^{\theta,z}= x_0\in \mathbb{R}^{N}.
\end{align}
%\begin{align} \label{eq:LDL_frozen}
%   \bar{X}_{t + \Delta}^{\theta, x, z} 
%  = x + \int_t^{t + \Delta}
% \bigl( 
%    A_{z, \theta} \bar{X}_{s}^{\theta, x, z}  
%    + b_{z, \theta} 
%\bigr) ds 
%+ \sigma (z, \theta) 
%\bigl(B_{t + \Delta} - B_{t} \bigr).  
%\end{align}
% 
Notice that $\big[\,\bar{X}_{t}^{\theta,z} |_{z=x}\,\big|\, \bar{X}_{0}^{\theta,z}=x\,\big]\equiv [\bar{X}_{t}^{\theta}|\bar{X}_{0}^{\theta}=x]$. The generator corresponding to (\ref{eq:LDL_frozen}) is given as:  
\begin{align*}
\mathscr{L}_\theta^{0, z} \varphi (x)
&
=  \sum_{1 \le i \le N} 
\bigl[A_{z, \theta} x + b_{z, \theta} \bigr]_i 
\partial_i \varphi (x) + \tfrac{1}{2} \sum_{1 \le i_1, i_2 \le N} \sum_{1 \le j \le d} V_j^{i_1} (z, \theta) V_j^{i_2} (z, \theta) 
\partial_{i_1 i_2} \varphi (x),  
%\\
% & = \sum_{i = 1}^N \bigl\{ V_0^i (z, \theta) +  \partial_z^\top V_0^i (z, \theta) \cdot (x - z) \bigr\} 
% \partial_i \varphi (x) + \tfrac{1}{2} \sum_{i_1, i_2 = 1}^N \sum_{j = 1}^d V_j^{i_1} (z, \theta) V_j^{i_2} (z, \theta) 
% \partial_{i_1 i_2} \varphi (x),  
\end{align*} 
for $\varphi : \mathbb{R}^N \to \mathbb{R}$, $x, z \in \mathbb{R}^N$, $\theta \in \Theta$.   
Notice that $\bar{X}^{\theta, z}_{t + \Delta} | \bar{X}_{t}^{\theta, z} = x \sim \mathscr{N} \big(e^{\Delta A_{z, \theta}} \hat{\mu}^z (\Delta, x, \theta), \Sigma (\Delta, z, \theta) \big)$, where $\hat{\mu}^z$ and $\Sigma$ are defined in (\ref{eq:mu}) and (\ref{eq:Sigma}), respectively. We write the density of $\bar{X}^{\theta, z}_{t + \Delta} | \bar{X}_{t}^{\theta, z} = x$ as
$
y \mapsto p^{\bar{X}^z}_\Delta (x, y; \theta) 
$ %\label{eq:ds_LDL_z}
% where $\hat{\Sigma} : [0, \infty) \times \mathbb{R}^N \times \Theta \to \mathbb{R}^{N \times N}$ is defined in (\ref{eq:LDL_cov}) and 
% %  
% \begin{gather}
%     \hat{\mu}^z (\Delta, x, \theta) 
%     = x  + \int_0^{\Delta} e^{- s A_{z, \theta}} b_{z, \theta} ds. 
% \end{gather} 
% 
and note that $p^{\bar{X}^z}_\Delta (x, y; \theta) |_{z=x} \equiv p^{\bar{X}}_\Delta (x, y; \theta)$, where the right-hand-side is the transition density of the LDL scheme (\ref{eq:LDL}).  

We introduce semi-groups $\{P_t^\theta\}_{t \ge 0}$ and $\{\bar{P}_t^{\theta, z} \}_{t \ge 0}$ associated with the Markov processes $\{X_t\}_{t \ge 0}$ and $\{\bar{X}_t^{\theta, z}\}_{t \ge 0}$, respectively as follows:  
\begin{gather}
P_{t}^{\theta} \varphi (x) 
= \int_{\mathbb{R}^N} \varphi (y) p^X_t (x, y; \theta) dy, \qquad 
\bar{P}_{t}^{\theta, z} \varphi (x) 
= \int_{\mathbb{R}^N} \varphi (y) 
p^{\bar{X}^z}_t (x, y; \theta) dy, \quad z \in \mathbb{R}^N, 
\end{gather}   
for $\varphi : \mathbb{R}^N \to \mathbb{R}$ and $(t, x, \theta) \in (0, \infty) \times \mathbb{R}^N \times \Theta$. For notational simplicity, we introduce:  
\begin{gather}
\widetilde{\mathscr{L}}^{z}_\theta 
:= \mathscr{L}_\theta - \mathscr{L}_\theta^{0, z}, 
% x \mapsto \Psi_{s}^{\theta, z, y} (x) 
% \equiv \widetilde{\mathscr{L}}_\theta^{z} \bigl\{ p^{\bar{X}^z}_s (\cdot, y; \theta) \bigr\} (x), \qquad z, y \in \mathbb{R}^N, \, s \in (0, \Delta].    
\end{gather} 
where we recall that $\mathscr{L}_\theta$ is the generator associated with the target SDE, given in (\ref{eq:generator}). 
%Our work follows the approach of \cite{igu:21-2}, 
The first steps in the derivation of our CF-expansion are provided in the following two results whose proofs are provided in Appendices \ref{app:pf_aux_1} and \ref{app:pf_aux_2}:% 
% \YI{Explain the motivation to introduce the following result and give its proof in this main text.}

\begin{lemma} \label{lemma:aux}
Let $t > 0$, $x,y \in \mathbb{R}^N$ and $\theta \in \Theta$. Also, let $\varphi \in C^\infty (\mathbb{R}^N; \mathbb{R})$. It holds that: 
\begin{align} 
p^X_t (x, y; \theta) 
& = p^{\bar{X}^z}_t (x, y; \theta) |_{z=x}
+ \int_0^t P_s^\theta
\widetilde{\mathscr{L}}_\theta^z
p^{\bar{X}^z}_{t - s} (\cdot, y; \theta) (x) 
|_{z =x} ds; \label{eq:density_app} \\  
P_t^{\theta} \varphi (x) 
& = \bar{P}_t^{\theta, z} \varphi (x) |_{z =x} 
+ 
\int_0^t P_{s}^{\theta} 
\widetilde{\mathscr{L}}^z_\theta  
\bar{P}_{t-s}^{\theta, z} \varphi (x) |_{z=x} ds. \label{eq:semi_group}    
\end{align}
% 
% In particular, for a sufficiently smooth test function $\varphi : \mathbb{R}^N \to \mathbb{R}$, we have that:  
% % 
% \begin{align} \label{eq:semi_group}
% P_t^{\theta} \varphi (x) 
% = \bar{P}_t^{\theta, z} \varphi (x) |_{z =x} 
% + 
% \int_0^t P_{s}^{\theta} 
% \widetilde{\mathscr{L}}^z_\theta  
% \bar{P}_{t-s}^{\theta, z} \varphi (x) |_{z=x} ds.  
% \end{align} 
\end{lemma}
\begin{lemma} \label{lemma:aux2}
Let $0 < s < t$, $\theta \in \Theta$ and $z, x, y \in 
\mathbb{R}^N$. Then it holds that, for any $K \in \mathbb{N}$, 
\begin{align} 
\bar{P}_s^{\theta, z} 
\widetilde{\mathscr{L}}_\theta^z 
\bar{P}_{t-s}^{\theta, z} \varphi (x)
& = \sum_{0 \le k \le K} \frac{s^k}{k!} 
\bigl\{ \mathrm{ad}^k_{\mathscr{L}_\theta^{0,z}} (\widetilde{\mathscr{L}}_\theta^z) 
\bigr\} 
\bar{P}_t^{\theta, z} \varphi (x) 
+ \mathscr{R}^{K+1, z} (s, t, x; \theta); \label{eq:expansion_step2_1}  
\\
\bar{P}_s^{\theta, z} \widetilde{\mathscr{L}}_\theta^z 
p^{\bar{X}^z}_{t-s}  (\cdot, y; \theta) (x) 
& = \sum_{0 \le k \le K} \frac{s^k}{k!} 
\bigl\{ 
\mathrm{ad}^k_{\mathscr{L}_\theta^{0,z}} (\widetilde{\mathscr{L}}_\theta^z) 
\bigr\} 
\{ p^{\bar{X}^z}_t (\cdot , y; \theta) \} (x)
+ \widetilde{\mathscr{R}}^{K + 1, z} (s, t, x, y; \theta)  
\label{eq:expansion_step2_2}   
\end{align}     
where the remainder terms ${\mathscr{R}}^{K + 1, z}$ and  $\widetilde{\mathscr{R}}^{K + 1, z}$ are specified in the proof of Lemma \ref{lemma:aux2} in Appendix \ref{app:pf_aux_2}. 
% 
% \begin{align} 
% \mathscr{R}_2^{J_1+1, z} (s, t, x; \theta)
% & = {s^{J_1+1}} \int_0^1 \bar{P}_{su}^{\theta, z} 
% \bigl\{ 
% \mathrm{ad}^{J_1 +1 }_{\mathscr{L}_\theta^{0,z}} (\widetilde{\mathscr{L}}_\theta^z) 
% \bigr\} 
% \bar{P}_{t-su}^{\theta, z} \varphi (x) 
% \tfrac{(1-u)^{J_1}}{J_1!} du; \label{eq:R2} \\ 
% \widetilde{\mathscr{R}}_2^{J_2 + 1, z} (s, t, x; \theta) 
% & = {s^{J_2+1}} \int_0^1 \bar{P}_{su}^{\theta, z}
% \mathrm{ad}^{J_2 + 1}_{\mathscr{L}_\theta^{0,z}} (\widetilde{\mathscr{L}}_\theta^z) 
% \bigl( 
% p^{\bar{X}^z}_{t-su} (\cdot, y; \theta) 
% \bigr) (x)
% \tfrac{(1-u)^{J_2}}{J_2!} du.  
% \end{align}   
% 
% 
%  
\end{lemma} 
% % % 
% 
\noindent Results similar to the above, for the elliptic case and for the Euler-Maruyama scheme used as a baseline transition density, are obtained in \cite{igu:21-2}. %but we provide their derivation here (within proofs) for completeness.
% Furthermore, we introduce: for $0 \le s_j \le s_{j-1} \le \cdots \le s_1 \le s_0, \, j \in \mathbb{N}$, 
% % 
% \begin{align}
% \mathbb{R}^N \ni x \mapsto \bar{P}_{(s_j, s_{j-1}, \ldots, s_0)}^{\theta, y, z} (x) 
% \equiv \bar{P}_
% \end{align}
% 
\subsubsection{Construction of 
CF-Expansion}
\label{sec:construction}
Based on the auxiliary results (Lemma \ref{lemma:LDL_hor}, \ref{lemma:aux} and \ref{lemma:aux2}) in the previous subsections, we construct a CF transition density expansion in the following three steps:  

\noindent
{\textit{Step 1.}} We recursively apply formula (\ref{eq:semi_group}) within (\ref{eq:density_app}), from Lemma \ref{lemma:aux}, to obtain for any $M \in \mathbb{N}$: 
\begin{align}
p^X_\Delta (x, y; \theta)  
&= p^{\bar{X}^z}_\Delta (x, y; \theta) |_{z=x} 
+ \int_0^\Delta 
{P}_{s_1}^{\theta} 
\widetilde{\mathscr{L}}_\theta^z p_{\Delta - s_1}^{\bar{X}^z} (\cdot, y; \theta) (x) |_{z=x}ds_1       \nonumber 
%\\[0.2cm]  \nonumber 
%& = p^{\bar{X}^z}_\Delta (x, y; \theta) |_{z=x} 
%+ \int_0^\Delta \bar{P}_{s_1}^{\theta, z} 
%\widetilde{\mathscr{L}}_\theta^z p_{\Delta - s_1}^{\bar{X}^z} (\cdot, y; \theta) (x) |_{z=x}ds_1  \\ & \qquad \qquad\qquad\qquad\qquad\qquad + \int_0^\Delta \int_0^{s_1} {P}_{s_2}^{\theta} \widetilde{\mathscr{L}}_\theta^{z} 
%\bar{P}_{s_1 - s_2}^{\theta, z} 
%\widetilde{\mathscr{L}}_\theta^z p_{\Delta - %s_1}^{\bar{X}^z} (\cdot, y; \theta) (x) |_{z=x} ds_2 ds_1
%\nonumber  
\\%[0.2cm] 
& =  
p^{\bar{X}^x}_\Delta (x, y; \theta) 
+  
\sum_{1 \le j \le M-1} 
\mathscr{T}^{j} (\Delta, x, y; \theta)
\ \ 
+  
\ \  
\mathscr{R}_1^{M} (\Delta, x, y; \theta),  
\label{eq:expansion_step1}
\end{align} 
where we have set: 
\begin{align}
& \mathscr{T}^{\, j} (\Delta, x, y; \theta) 
=
\int_{I(s_{1:j})}
% \int_0^\Delta \int_0^{s_1} \cdots \int_0^{s_{j-1}}
\bar{P}_{s_j}^{\theta, z} 
\widetilde{\mathscr{L}}_\theta^z 
\bar{P}_{s_{j-1} - s_j}^{\theta, z} \cdots 
\widetilde{\mathscr{L}}_\theta^z \bar{P}_{s_1 - s_2}^{\theta, z} 
\widetilde{\mathscr{L}}_\theta^z p_{\Delta - s_1}^{\bar{X}^z} (\cdot, y; \theta) (x) \big|_{z=x} 
\,  ds_j \cdots ds_1; \nonumber \\
& \mathscr{R}_1^{M} (\Delta, x, y; \theta) 
= %\int_0^\Delta \int_0^{s_1} \cdots \int_0^{s_{M-1}}
\int_{I(s_{1:M})}   
{P}_{s_{M}}^{\theta} \widetilde{\mathscr{L}}_\theta^z \bar{P}_{s_{M-1} - s_{M}}^{\theta, z} \cdots 
\widetilde{\mathscr{L}}_\theta^z \bar{P}_{s_1 - s_2}^{\theta, z} \widetilde{\mathscr{L}}_\theta^z p_{\Delta - s_1}^{\bar{X}^z} (\cdot, y; \theta) (x) 
\big|_{z=x} \, ds_{M} \cdots ds_1, 
\label{eq:first_err}  \\[0.2cm]  \nonumber 
&\qquad\qquad \qquad \quad I(s_{1:k}) := \{s_{1:k}=(s_1,\ldots, s_k) : 0 \le s_k \le \cdots \le s_1\le \Delta\}, \quad k\ge 0, 
\end{align}
%
%for 
%
%$$
%I(s_{1:k}) = \{s_{1:k}=(s_1,\ldots, s_k) : 0 \le s_k \le \cdots \le s_1\le \Delta\}, \quad k\ge 0,  
%$$
% 
with the convention $s_0 \equiv \Delta$.  

\noindent 
{\textit{Step 2.}} Since $\mathscr{T}^j (\Delta, x,y; \theta)$, $1 \le j \le M$, 
is not tractable, we obtain a computable quantity for it via use of Lemma \ref{lemma:aux2}. 
Let $\mathscr{T}_{s_{1:j}} (\Delta, x, y; \theta)$ be the integrand of  $\mathscr{T}^j (\Delta, x, y; \theta)$, so that $ \mathscr{T}^j (\Delta, x, y; \theta) = \textstyle \int_{I(s_{1:j})} \mathscr{T}_{s_{1:j}} (\Delta, x, y; \theta) ds_{j} \cdots ds_1$. 
%with convention ${s_{1:j}} \equiv (s_1, \ldots, s_j)$.  
% % 
Recursive application of Lemma \ref{lemma:aux2} to  $\mathscr{T}_{s_{1:j}}$, $1 \le j \le M-1$, gives: 
%for any $J \in \mathbb{N}$, 
% 
\begin{align} \label{eq:tau_j}
\mathscr{T}_{s_{1:j}} (\Delta, x, y; \theta) 
= 
\sum_{\substack{\alpha \le \beta^{[j]}}}
\frac{\prod_{k = 1}^j (s_k)^{\alpha_k}}{\alpha!} 
\cdot 
\mathscr{D}_\alpha^{z, \theta} 
\bigl\{ p^{\bar{X}^z}_{\Delta} (\cdot, y; \theta) \bigr\} (x) \bigr|_{z=x}
+ \mathscr{E}^{\beta^{[j]}}_{s_{1:j}}
(\Delta, x, y; \theta), 
\end{align} 
for a multi-index $\beta^{[j]} = (\beta_1^{[j]}, \ldots, \beta_j^{[j]}) \in \mathbb{Z}_{\ge 0}^{j}$, where we have defined:   
\begin{align} 
\mathscr{D}_{\alpha}^{z, \theta} 
& \equiv 
\Bigl(
\mathrm{ad}_{\mathscr{L}_\theta^{0, z}}^{\alpha_j} (\widetilde{\mathscr{L}}_{\theta}^z) 
\Bigr)
\Bigl(\mathrm{ad}_{\mathscr{L}_\theta^{0, z}}^{\alpha_{j-1}} 
(\widetilde{\mathscr{L}}_{\theta}^z)  
\Bigr)
\cdots 
\Bigl(\mathrm{ad}_{\mathscr{L}_\theta^{0, z}}^{\alpha_{1}} (\widetilde{\mathscr{L}}_{\theta}^z) \Bigr), \qquad z \in \mathbb{R}^N, \, \theta \in \Theta,  \label{eq:D} 
\end{align}
and the residual $\mathscr{E}^{\beta^{[j]}}_{s_{1:j}}
(\Delta, x, y; \theta)$ is given in (\ref{eq:E_def}) of Supplementary Material.  
% and appropriately bounded under assumptions. 
% 
% In (\ref{eq:tau_j}), $\mathcal{E}^{[J]}_{s_{1:j}} (\Delta, x, y; \theta)$ is an error term such that 
% % 
% \begin{align} 
% \mathscr{R}_2^{j, [J]} (\Delta, x, y, \theta) \equiv \int_0^\Delta \int_0^{s_1} \cdots \int_0^{s_{j-1}} \mathcal{E}^{[J]}_{s_{1:j}}
% (\Delta, x, y; \theta) ds_j \cdots ds_1
% \end{align}
% has an appropriate upper bound under assumptions (this is shown later).  % 
We now obtain, for $1 \le j \le M-1$: 
\begin{align}
&\mathscr{T}^{j} (\Delta, x, y; \theta)
=
\sum_{\substack{\alpha \le \beta^{[j]}}}
\Delta^{|\alpha| + j} 
\cdot
\tfrac{K (\alpha)}{\alpha!}
\cdot 
\mathscr{D}_\alpha^{z, \theta}
\{ p^{\bar{X}^z}_\Delta (\cdot, y; \theta) \} (x) |_{z = x} 
+ \mathscr{R}_2^{\, j, \beta^{[j]}} (\Delta, x, y, \theta),  
\label{eq:expansion_step2} \\[0.2cm]
&\qquad\qquad\qquad K(\alpha) := 
\int_{0 \le s_{j} \le \cdots \le s_1 \le 1}
\prod_{1 \le k \le j} (s_k)^{\alpha_k} ds_j \cdots ds_1,
\nonumber  \\   
&\qquad\qquad\qquad \mathscr{R}_2^{\, j,  \beta^{[j]}} (\Delta, x, y, \theta) 
:= 
% \int_{0 \le s_{j} \le \cdots \le s_1 \le \Delta} 
\int_{I(s_{1:j})} \mathscr{E}^{\beta^{[j]}}_{s_{1:j}}
(\Delta, x, y; \theta) ds_j \cdots ds_1. 
\label{eq:second_err}
\end{align} 
{\textit{Step 3.}} From (\ref{eq:expansion_step1}) in {\textit{Step 1.}} and 
(\ref{eq:expansion_step2}) in {\textit{Step 2.}}, we obtain the following CF-expansion for the true (intractable) transition density. For any $M\ge 1$ and multi-indices $\beta^{[j]} \in \mathbb{Z}_{\ge 0}^j$, $1 \le j \le M-1$:
% /
\begin{align}
p^X_\Delta (x, y; \theta)
& = 
p^{\bar{X}}_\Delta (x, y; \theta)  
+ 
\sum_{1 \le j \le M-1}
\sum_{\substack{ \alpha \le \beta^{[j]}}}
\Delta^{|\alpha| + j} \cdot 
\tfrac{K (\alpha)}{\alpha!}
\cdot \mathscr{D}_\alpha^{z, \theta}
\{ p^{\bar{X}^z}_\Delta (\cdot, y; \theta) \} (x) \Bigr|_{z = x} %\label{eq:app} 
\nonumber \\[0.2cm]
& \qquad\qquad\qquad \qquad \qquad  
+ \mathscr{R}_1^{M} (\Delta, x, y; \theta) 
+ \sum_{1 \le j \le M-1} \mathscr{R}_2^{\, j, \beta^{[j]}} (\Delta, x, y, \theta), 
\label{eq:err} 
\end{align} 
where $\mathscr{R}_1^{M}$, $\mathscr{R}_2^{\, j, \beta^{[j]}}$ are defined in (\ref{eq:first_err}), (\ref{eq:second_err}), respectively. Under assumptions, we show in Section \ref{sec:main}  that the remainder terms are of size $\mathcal{O} 
(\Delta^{p})$ for an arbitrary $p > 0$ by choosing a large enough $M$ and appropriate $\beta^{[j]}$, $1 \le j \le M-1$. The double sum in (\ref{eq:err}) involves tractable terms and can be utilised as a proxy for the true transition density. In particular, the expansion is well-defined for both model classes (\ref{eq:ellip}) and (\ref{eq:hypo}) since the Gaussian density $p^{\bar{X}^z} (x, y; \theta)$ and its partial derivatives (involved in $\mathscr{D}_\alpha^{z, \theta} \bigl\{ p_\Delta^{\bar{X}^z} (\cdot, y; \theta) (x)  \bigr\} |_{z=x}$) are well-defined from Lemma \ref{lemma:LDL_hor}. We note that the tractable double sum in  (\ref{eq:err}) is regarded as a 
CF-expansion,  but the current form of the expansion does not yet correspond  to the 
`$\Delta$'-expansion (\ref{eq:delta_expansion}). For instance, the exponents of the step-size $\Delta$ are integers in (\ref{eq:err}), while they are given as $k / 2$, $k \in \mathbb{N}$, in (\ref{eq:delta_expansion}). However, we emphasise that (\ref{eq:err}) will be ultimately expressed as a $\Delta$-expansion of the form in  (\ref{eq:delta_expansion}) after carefully working with the terms $\mathscr{D}_\alpha^{z, \theta}
\{ p^{\bar{X}^z}_\Delta (\cdot, y; \theta) \} (x)$. 
Indeed, taking partial derivatives of $x \mapsto p^{\bar{X}^z} (x, y; \theta)$, will give Hermite polynomials and powers $\Delta^{- k/2}$, where the integer $k$ depends on the number of derivatives. We explain this in detail 
% of (\ref{eq:app}) into the form of delta expansion (\ref{eq:delta_expansion}) 
later on in Section \ref{sec:series}. 
\section{Error Analysis for the 
CF-Expansion} \label{sec:error}
In Section \ref{sec:CF_expansion} we have constructed a CF-expansion (\ref{eq:err}) for the true transition density. Our objective now is to provide rigorous error estimates for this expansion, thus theoretically justifying its derivation, similarly to results obtained by a few  earlier works in the case of the elliptic class (\ref{eq:ellip}). We also describe that the obtained expansion (\ref{eq:cf}) can be given in the form (\ref{eq:delta_expansion}), namely a series in  powers of $\Delta$. 
As the error estimates vary for classes (\ref{eq:ellip}), (\ref{eq:hypo}), we make use of the notation $w \in \{\ref{eq:ellip}, \ref{eq:hypo} \}$ and 
write $p^{X, (w)}\equiv p^{X}$, $p^{\bar{X}, (w)} \equiv p^{\bar{X}}$, $\mathscr{R}_1^{M, (w)} \equiv \mathscr{R}_1^{M}$ and $\mathscr{R}_2^{\, j, \beta^{[j]}, (w)} \equiv \mathscr{R}_2^{\, j, \beta^{[j]}}$ to indicate the class under consideration.  
\subsection{Main Result}
\label{sec:main}
We derive upper bounds for the residuals of the CF expansion $\mathscr{R}_1^M$ and  $\mathscr{R}_2^{j, \beta^{[j]}}$ specified in (\ref{eq:first_err}) and (\ref{eq:second_err}), respectively.  We will need the  following additional assumptions.
%with which we work to obtain non-asymptotic bounds for the residuals.
% 
\begin{assump} \label{ass:param}
The parameter space $\Theta$ is compact. Also, for each $x \in \mathbb{R}^N$, the function $\theta \mapsto V_j (x, \theta)$ is continuous, $0\le j\le d$. 
\end{assump} 
\begin{assump} 
\label{ass:coeff}
Let $x \in \mathbb{R}^N$ be the initial state of the transition dynamics. The SDE coefficients satisfy the following properties: 
\begin{itemize}
\item[1](Boundedness of drift at initial state): There exists a constant $\kappa > 0$ such that  $|x| + |V_0(x, \theta)| \le \kappa$ for all $\theta \in \Theta$; 
\item[2](Uniform boundedness of diffusion coefficients): There exists a constant $C > 0$ such that 
\mbox{$ 
| V_j (y, \theta)| \le C$}, $1 \le j \le d 
$,
for all $(y, \theta) \in \mathbb{R}^N \times \Theta$;   
\item[3] (Uniform boundedness of derivatives): There is a constant $C > 0$ such that 
$
\textstyle\sum_{j = 0}^d 
|\partial^\alpha_y V_j (y, \theta)|
\le C 
$ 
for all $\alpha \in \mathbb{Z}^N_{\ge 0}$ with $|\alpha| \ge 1$ and all $(y, \theta) \in \mathbb{R}^N \times \Theta$.  
\end{itemize} 
\end{assump}

\begin{assump} \label{ass:dim}
For the hypo-elliptic model class (\ref{eq:hypo}), $N_S = N_R = d$. 
\end{assump}

Assumptions \ref{ass:param}--\ref{ass:dim} suffice for obtaining appropriate bounds for the residuals of the CF-expansion. The uniform boundedness for the derivatives of the SDE coefficients is a standard assumption for the existence of a smooth transition density, when combined with  H\"ormander's condition. Such uniform boundedness is also assumed in \cite{li:13} to control the residual of the expansion developed therein for elliptic SDEs. Assumptions \ref{ass:coeff}--\ref{ass:dim}  are used mainly in the proof of Theorem \ref{thm:bd_r1}, where we need {an} upper bound for the true density $y \mapsto p^{X}_\Delta (x, y; \theta)$. \cite{piga:22} shows that under  Assumptions \ref{ass:coeff}--\ref{ass:dim} the true transition density has a Gaussian-type bound as given later at (\ref{eq:G_bd_hypo}). Based on this result, we show that the errors are appropriately bounded, analogously to  \cite{yang:19} who also used a Gaussian-type bound for the true density to control the residuals of an expansion for inhomogeneous elliptic SDEs. We stress that Assumptions \ref{ass:param}--\ref{ass:dim} are not necessary for the construction of the CF-expansion, in the sense that our formulae can still be evaluated for SDEs with coefficients whose partial derivatives exhibit, e.g., polynomial growth as assumed in earlier works \citep{ait:08, yang:19} for elliptic SDEs. {Assumption \ref{ass:dim} can be weakened as there is a possibility to obtain an upper bound for the true transition density by carefully developing the arguments in \cite{piga:22}. However, this is not straightforward and is beyond the scope of the present work.} The relaxation of our conditions is left for future research. 
%
% \begin{remark}
% \end{remark} 
%

To provide a statement of our main result, we introduce some notation. We set:   
\begin{align} \label{eq:m}
m(\mathrm{\ref{eq:ellip}}) := N, 
\qquad  
m(\mathrm{\ref{eq:hypo}}) := 4 d. 
\end{align}
% \begin{align}
%     m (w) \equiv 
%     \begin{cases}
%        N  & w = \ref{eq:ellip}; \\
%        4d & w = \ref{eq:hypo}. 
%     \end{cases}
% \end{align}
% 
Also, 
$\mathscr{G}^{(w)} : (0, \infty) \times \mathbb{R}^N \times \mathbb{R}^N \times \Theta \to \mathbb{R}$, $w \in \{\ref{eq:ellip}, \ref{eq:hypo} \}$, is a mapping  characterised as follows. There exist constants $C, c > 0$ such that: 
\begin{align}
\bigl| 
\mathscr{G}^{(\ref{eq:ellip})} (t, x, y, \theta) 
\bigr| 
& \le 
C t^{-\tfrac{m(\ref{eq:ellip})}{2}} \times 
\exp \Bigl( - c \tfrac{|y-x|^2}{t} \Bigr); 
\label{eq:G_bd_ell} \\[0.2cm]  
\bigl| 
\mathscr{G}^{(\ref{eq:hypo})} (t, x, y, \theta) 
\bigr| 
& \le 
C t^{- \tfrac{m(\ref{eq:hypo})}{2}} \times  
\exp 
\left( 
- c \Bigl(\tfrac{|y_S-x_S - V_{S, 0} (x, \theta) t |^2}{t^3} + \tfrac{|y_R-x_R|^2}{t} \Bigr) 
\right),  \label{eq:G_bd_hypo}  
\end{align}
% 
% \begin{align} 
% \bigl| 
% \mathscr{G}^{(w)} (t, x, y, \theta) 
% \bigr| 
% \le 
% \begin{cases}
% C t^{-\tfrac{m(\ref{eq:ellip})}{2}} \times 
% \exp \Bigl( - c \tfrac{|y-x|^2}{t} \Bigr),  
% & w = \ref{eq:ellip}; \\[0.3cm]  
% C t^{- \tfrac{m(\ref{eq:hypo})}{2}} \times  
% \exp 
% \left( 
% - c \Bigl(\tfrac{|y_S-x_S - V_{S, 0} (x, \theta) t |^2}{t^3} + \tfrac{|y_R-x_R|^2}{t} \Bigr) 
% \right), 
% & w = \ref{eq:hypo}, 
% % , \quad N_S = N_R = d, 
% \end{cases}
% \label{eq:G_bd_hypo} 
% \end{align}
for all $(t, x, y, \theta) \in (0, \infty) \times \mathbb{R}^N \times \mathbb{R}^N\times \Theta$.  
Notice that for some constant $C > 0$, for $\Delta \in (0, 1)$:
\begin{align}
\begin{aligned} \label{eq:G_bd}
& \sup_{(x, y, \theta) \in \mathbb{R}^N \times \mathbb{R}^N \times \Theta} 
\,  \bigl| \mathscr{G}^{(w)} (\Delta, x, y, \theta) \bigr| 
\le C \Delta^{-\tfrac{m(w)}{2}}, \\ 
& \textcolor{black}{\sup_{(x, \theta) \in \mathbb{R}^N \times \Theta} 
\int_{\mathbb{R}^N} \bigl| \mathscr{G}^{(w)} (\Delta, x, y, \theta) \bigr| dy 
\le C,}
\end{aligned} 
\qquad w \in \{\ref{eq:ellip}, \ref{eq:hypo} \}, 
\end{align}
{which implies that the size of $\mathscr{G}^{(w)}$ in $L_1$-norm is $\mathcal{O}(1)$ irrespective of the model classes \eqref{eq:ellip} or \eqref{eq:hypo}.} 
% 
%The two remainder terms in the CF expansion (\ref{eq:err}) are bounded as follows.
% 
\begin{theorem}[Bound for $\mathscr{R}_1^{M, (w)}$] \label{thm:bd_r1} 
Let $x \in \mathbb{R}^N$ be the initial state of the transition dynamics and $M \ge 1$.  Under Assumptions \ref{ass:diff}--\ref{ass:dim}, there exists a constant $C> 0$ such that for all $(\Delta, y, \theta) \in (0,1) \times \mathbb{R}^N \times \Theta$:   
\begin{align*}
\bigl| 
\mathscr{R}_1^{M, (w)} (\Delta, x, y; \theta) \bigr| 
\leq 
C \Delta^{\tfrac{M}{2}} \times 
{\bigl| \mathscr{G}^{(w)} (\Delta, x, y, \theta) \bigr|}, 
% \le C_2 \Delta^{\tfrac{M -  m(w)}{2}}, 
\qquad w \in \{\ref{eq:ellip}, \ref{eq:hypo} \}.   
\end{align*} 
\end{theorem}

\begin{theorem}[Bound for $\mathscr{R}_2^{\, j, \beta^{[j]}, (w)}$]
\label{thm:bd_r2}
Let $x \in \mathbb{R}^N$ be the initial state of the transition dynamics, and let $1 \le j \le M-1, \, M \in \mathbb{N}$, $\beta^{[j]} = (\beta^{[j]}_1, \ldots, \beta^{[j]}_j) \in \mathbb{Z}_{\ge 0}^j$. Under Assumptions \ref{ass:diff}--\ref{ass:dim}, there exists  constant $C > 0$ such that for all $(\Delta, y, \theta) \in (0,1) \times \mathbb{R}^N \times \Theta$:
\begin{gather*} 
\bigl| \mathscr{R}_2^{\, j, \beta^{[j]}, (w)} (\Delta, x, y; \theta) \bigr| 
\le C \Delta^{K^{[j], w} (\beta^{[j]})} \times \bigl| {\mathscr{G}^{(w)} (\Delta, x, y, \theta) \bigr|},
% \le C_2 \Delta^{K^{[j], w} (\beta^{[j]}) - \tfrac{m(w)}{2}},  
\quad  w \in \{\ref{eq:ellip}, \ref{eq:hypo}\}, \\[0.2cm] 
%\end{align*}
% 
% 
%\begin{align*}
\qquad  K^{[j], \ref{eq:ellip}} (\beta^{[j]})
:= \min_{1 \le  i \le j} \tfrac{\beta_i^{[j]}}{2} + \tfrac{j}{2}, \qquad 
K^{[j], \mathrm{\ref{eq:hypo}}} (\beta^{[j]}) 
:= \min_{1 \le i \le j} \tfrac{1}{2} 
\bigl( 
\floor{\tfrac{\beta_i^{[j]}}{2}} - \mathbf{1}_{\beta_i^{[j]} \ge 2} 
\bigr)  + \tfrac{j}{2}. 
\end{gather*}
% \begin{align}
% K^{[j], (w)} (\beta^{[j]})
% \equiv 
% \begin{cases}
% \min_{1 \le  i \le j} \tfrac{\beta_i^{[j]}}{2} + \tfrac{j}{2}, & w = \ref{eq:ellip}; \\[0.2cm] 
% \min_{1 \le i \le j} \tfrac{1}{2} 
% \bigl( 
% \floor{\tfrac{\beta_i^{[j]}}{2}} - \mathbf{1}_{\beta_i^{[j]} \ge 2} 
% \bigr)  + \tfrac{j}{2}, & w = \ref{eq:hypo}. 
% \end{cases}
% \end{align}
\end{theorem}
\noindent  The proofs of Theorems \ref{thm:bd_r1}--\ref{thm:bd_r2} are given in {Section}  \ref{app:pfs_main} of Supplementary Material. From Theorem \ref{thm:bd_r2}, by selecting the multi-index $\beta^{[j], (w)} \in \mathbb{Z}_{\ge 0}^j$ so that $K^{[j], w} (\beta^{[j], (w)}) 
\ge \tfrac{M}{2}$, we have that: 
{
\begin{align} 
\begin{aligned} 
& p^{X, (w)}_\Delta (x, y; \theta) 
=  
p^{\bar{X}, (w)}_\Delta (x, y; \theta)  \\ 
&\quad + 
\sum_{1 \le j \le M-1}
\sum_{\substack{ \alpha \le \beta^{[j], (w)}}}
\Delta^{|\alpha| + j} \cdot 
\tfrac{K (\alpha)}{\alpha!}
\cdot \mathscr{D}_\alpha^{z, \theta}
\{ p^{\bar{X}^z, (w)}_\Delta (\cdot, y; \theta) \} (x) \Bigr|_{z = x}
+ \mathscr{E}^{(w)} (\Delta, x, y; \theta), 
% + \mathcal{O} \bigl( \Delta^{\tfrac{M - m(w)}{2}} \bigr).
\end{aligned} 
 \label{eq:cf}  
\end{align}
% . 
for a residual $\mathscr{E}^{(w)}$ so that for any initial state $x \in \mathbb{R}^N$ there exists a constant $C> 0$ such that for all $(\Delta, y, \theta) \in (0,1) \times \mathbb{R}^N \times \Theta$:   
\begin{align*}
\bigl| 
\mathscr{E}^{(w)} (\Delta, x, y; \theta) \bigr| 
\leq 
C \Delta^{\tfrac{M}{2}} \times 
\bigl| \mathscr{G}^{(w)} (\Delta, x, y, \theta) \bigr|, 
% \le C_2 \Delta^{\tfrac{M -  m(w)}{2}}, 
\qquad w \in \{\ref{eq:ellip}, \ref{eq:hypo} \}.   
\end{align*}
} 
\subsection{Series Expansion in $\Delta$} 
\label{sec:series}
We study the CF-expansion given in (\ref{eq:cf}) in detail. Expression (\ref{eq:cf}) involves $\Delta^{|\alpha|+j}$ in front of each summand, but an additional $\Delta^{-K_\alpha/2}$, for some  $K_\alpha \in \mathbb{N}$, is produced from $\mathscr{D}_\alpha^{z, \theta}
\{ p^{\bar{X}^z}_\Delta (\cdot, y; \theta) \} (x)$, where we recall that $\mathscr{D}_\alpha^{z, \theta}$ is the differential operator defined in (\ref{eq:D}). We show that upon rearrangement of terms in powers of $\Delta$, the right-hand-side of (\ref{eq:cf}) attains the form of the $\Delta$-expansion in (\ref{eq:delta_expansion}), i.e.~a series expansion in (positive) powers of $\sqrt{\Delta}$.
%, involving Hermite polynomials defined via the auxiliary Gaussian density $p^{\bar{X}^z, (w)} (x, y; \theta)$. 
In particular, we clarify below that differentiating the Gaussian density $p^{\bar{X}^z} (x, y; \theta)$ w.r.t.~the initial state $x \in \mathbb{R}^N$ produces additional powers $\Delta^{-{K}/{2}}$, for $K\ge 1$ depending on the number of derivatives. For the model class (\ref{eq:hypo}), the value of $K$ varies depending on whether the differentiation acts on smooth or rough components.
% 
% \begin{enumerate}
% \item[II.] The operator 
% $\mathscr{D}_\alpha^{z, \theta}$ involves appropriate number of differentiation, and then $\mathscr{D}_\alpha^{z, \theta} 
% \{ p^{\bar{X}^z}_\Delta (\cdot, y; \theta) \} (x) |_{z = x}$ produces additional $\Delta^{-\tfrac{K_\alpha}{2}}$, 
% where $K_\alpha$ is a positive integer satisfying $1 \le K_\alpha \le 2 |\alpha| + j$. 
% \end{enumerate}
% 
We define, for $\alpha \in \mathbb{Z}_{\ge 0}^N$: 
\begin{align}
\| \alpha \|_{\ref{eq:ellip}}: = \tfrac{1}{2} | \alpha |, \qquad 
\| \alpha \|_{\ref{eq:hypo}} = \tfrac{3}{2} | \alpha_S|  + \tfrac{1}{2} | \alpha_R|,  \label{eq:mi_norm}
\end{align}
% 
% \begin{align}
% \| \alpha \|_{w} =
% \begin{cases}
% \tfrac{1}{2} | \alpha |   &   w = \ref{eq:ellip};  \\[0.1cm] 
% \tfrac{3}{2} | \alpha_S|  + \tfrac{1}{2} | \alpha_R|,  &  w = \ref{eq:hypo},  
% \end{cases} \label{eq:mi_norm}
% \end{align} 
% 
where, for class (\ref{eq:hypo}), we interpret 
$\alpha = (\alpha_{S,1},\ldots,\alpha_{S,N_S},\alpha_{R,1},\ldots,\alpha_{R,N_R})$ for given $\alpha_S \in \mathbb{Z}_{\ge 0}^{N_S}$,  $\alpha_R  \in \mathbb{Z}_{\ge 0}^{N_R}$.
%$\alpha  = \alpha_S\ast\alpha_R$
%with $\alpha_S \in \mathbb{Z}_{\ge 0}^{N_S}$,  $\alpha_R  \in \mathbb{Z}_{\ge 0}^{N_R}$ for the hypo-elliptic model (\ref{eq:hypo}). 
We then have the following key result whose proof is provided in Section  \ref{sec:pf_ldl_deriv} of Supplementary Material: 
\begin{lemma} \label{lemma:deriv_LDL}
Let $x, y \in \mathbb{R}^N$, $\theta \in \Theta$, $\Delta > 0$. Also, let $\alpha \in \mathbb{Z}_{\ge 0}^N$. Under Assumptions \ref{ass:diff}--\ref{ass:dim}, we have that: 
\begin{align*}
\partial^{\alpha}_x \, p_\Delta^{\bar{X}^z, (w)} (x, y; \theta) |_{z=x} 
& = \Delta^{- \| \alpha \|_w} \times 
\mathscr{H}^{(w)}_{\alpha} (\Delta, x ,y ;\theta) \, p_\Delta^{\bar{X}, (w)} (x, y; \theta), \qquad w \in \{ \ref{eq:ellip}, \ref{eq:hypo} \}, 
\end{align*} 
where $y \mapsto \mathscr{H}^{(w)}_{\alpha} (\Delta, x, y ;\theta)$ is obtained explicitly defined via differentiation of $x \mapsto p_\Delta^{\bar{X}^z, (w)} (x, y; \theta)$ and is characterised as follows. There exists a constant $C > 0$ such that for all $\Delta > 0$, $x, y \in \mathbb{R}^N$, $\theta \in \Theta$:
\begin{gather*}
\bigl| \mathscr{H}^{(w)}_{\alpha} (\Delta, x ,y ;\theta) \, p_\Delta^{\bar{X}, (w)} (x, y; \theta) 
\bigr|   \le 
C \, 
\bigl| 
\mathscr{G}^{(w)} (\Delta, x, y, \theta)
\bigr|. 
% ; \\ 
% \int_{\mathbb{R}^N} \Bigl| \mathscr{H}^{(w)}_{\alpha} (\Delta, x ,y ;\theta) \, p_\Delta^{\bar{X}, (w)} (x, y; \theta) \Bigr| dy 
% \le C_2. 
% \\[0.2cm] 
% \left\| 
% \mathscr{H}^{(w)}_{\alpha} (\Delta, x , \cdot ; \theta) 
% \, p_\Delta^{\bar{X}, (w)} (x, \cdot ; \theta) 
% \right\|_{L^1} 
% & < C_2. 
\end{gather*} 
% 
%for all $\Delta > 0$, $x, y \in \mathbb{R}^N$ and $\theta \in \Theta$. In particular, there exists a constant $C > 0$ such that 
% 
%\begin{align}
%\int_{\mathbb{R}^N} \Bigl| \mathscr{H}^{(w)}_{\alpha} (\Delta, x ,y ;\theta) \, p_\Delta^{\bar{X}, (w)} (x, y; \theta) \Bigr| dy 
%\le C, \qquad w \in \{\ref{eq:ellip}, \ref{eq:hypo} \}.  
%\end{align} 
%% 
%for all $\Delta > 0$, $x \in \mathbb{R}^N$ and $\theta \in \Theta$. 
% 
\end{lemma} 
\noindent In brief, Lemma \ref{lemma:deriv_LDL} states the following. 
% 
%\begin{itemize}
%\item 
For model class \eqref{eq:ellip}, taking $k \in \mathbb{N}$ partial derivatives of $x \mapsto p_\Delta^{\bar{X}^z, (\ref{eq:ellip})} (x, y; \theta)$ yields a term $\Delta^{-\tfrac{k}{2}}$.
%\item 
For class \eqref{eq:hypo}, taking $k \in \mathbb{N}$ partial derivatives of $p_\Delta^{\bar{X}^z, (\ref{eq:hypo})} (x, y; \theta)$ w.r.t.~the smooth components (resp.~the rough components) produces the term $\Delta^{-\tfrac{3}{2}k}$ (resp.~the term $\Delta^{-\tfrac{k}{2}}$).
%\end{itemize}
% 
Based upon Lemma \ref{lemma:deriv_LDL}, the form of the CF-expansion is determined from the expression of the differential operator $\mathscr{D}_\alpha^{z, \theta}$ and the number of derivatives involved therein. A detailed characterisation for the differential operator $\mathscr{D}_\alpha^{z, \theta}$ is provided in Supplementary Material. In particular, Lemma \ref{lemma:diff_step_2} in Supplementary Material states that the operator admits the following expression. For $\varphi \in C^\infty (\mathbb{R}^N; \mathbb{R})$, $\alpha \in \mathbb{Z}^j_{\ge 0}, \, j \in \mathbb{N}$ and $(x, \theta) \in \mathbb{R}^N \times \Theta$:  
\begin{align} \label{eq:do_main}
\mathscr{D}_\alpha^{z, \theta, (w)} \varphi (x) |_{z=x} 
= 
\sum_{\gamma \in \mathscr{J}_{w}(\alpha)}   
\mathscr{W}_\gamma^{[\alpha]} (x, \theta) \partial^\gamma \varphi (x), \qquad  
w \in \{\ref{eq:ellip}, \ref{eq:hypo}\}{,}  
\end{align}
where $\mathscr{J}_w (\alpha)$ is a set of multi-indices $\mathbb{Z}_{\ge 0}^N$ defined in (\ref{eq:set_mi}) in Supplementary Material and $\mathscr{W}_\gamma^{[\alpha]}: \mathbb{R}^N \times \Theta \to \mathbb{R}$ is explicitly determined from products of partial derivatives of the SDE coefficients and can be evaluated in applications using software performing symbolic calculations. 
%from symbolic calculation and explicitly determined from products of partial derivatives of the SDE's coefficients. 
% 
%\begin{align*}
%\mathscr{J}_w (\alpha) 
%\equiv  
%\begin{cases}
%\bigl\{  \gamma \in \mathbb{Z}_{\ge 0}^N \ \Bigl| \  1 \le |\gamma| \le | \alpha |  + 2 j - \| \alpha \|_{0}  \bigr\}, & w = \ref{eq:ellip}; \\[0.2cm] 
%\bigl\{  \gamma \in \mathbb{Z}_{\ge 0}^N \ \Bigl| \  1 \le |\gamma_R | + 2 |\gamma_S|  \le | \alpha |  + 2 j  - \| \alpha \|_{0}, 
%\  \| \gamma \|_{\ref{eq:hypo}} \le \tfrac{1}{2} \sum_{1 \le i \le j} \floor{\tfrac{3 \alpha_i}{2}} + \tfrac{3}{2}j - \tfrac{1}{2} \| \alpha \|_{0, 1, 2} \bigr\}, & w = \ref{eq:hypo},  
%\end{cases}
%\end{align*} 
%% 
%with 
%% 
%\begin{align*}
%\| \alpha \|_{0} \equiv \sum_{i = 1}^j \mathbf{1}_{\alpha_i = 0}, \qquad 
%\| \alpha \|_{0, 1, 2} \equiv \sum_{i = 1}^j \bigl\{ 2 \times \mathbf{1}_{\alpha_i = 0} +  \mathbf{1}_{\alpha_i \in \{ 1, 2\} } \bigr\}. 
%\end{align*} 
% 
Due to (\ref{eq:do_main}) and Lemma \ref{lemma:deriv_LDL}, we have that for $w \in \{\ref{eq:ellip}, \ref{eq:hypo}\}$:
\begin{align}
&
\sum_{1 \le j \le M-1}
\sum_{\substack{ \alpha \le \beta^{[j]}}}
\Delta^{|\alpha| + j} \cdot 
\tfrac{K (\alpha)}{\alpha!}
\cdot \mathscr{D}_\alpha^{z, \theta, (w)}
\{ p^{\bar{X}^z, (w)}_\Delta (\cdot, y; \theta) \} (x) \Bigr|_{z = x}
\nonumber \\[0.2cm]
& 
\quad = 
\sum_{1 \le j \le M-1}
\sum_{\substack{\alpha \le \beta^{[j]}}} 
\sum_{\gamma \in \mathscr{J}_w (\alpha)}
\Delta^{|\alpha| + j - \| \gamma \|_w } \,
\tfrac{K (\alpha)}{\alpha!} \, 
\mathscr{W}_\gamma^{[\alpha]} (x, \theta) \,  
\mathscr{H}_\gamma^{(w)} (\Delta, x, y; \theta) \, 
p^{\bar{X}, (w)}_\Delta (x, y; \theta) 
\nonumber \\ 
& \quad \equiv \sum_{1 \le k \le M-1} 
\Delta^{\tfrac{k}{2}} \cdot e_k^{(w)} (\Delta, x, y; \theta)\cdot  p^{\bar{X}, (w)}_\Delta (x, y; \theta) 
+ \mathscr{R}_3^{(w)} (\Delta, x, y; \theta) \, 
p_\Delta^{\bar{X}, (w)} (x, y; \theta),  
\label{eq:arg}
\end{align}
where in the last line, we rearrange the sum in ascending order in powers of $\sqrt{\Delta}$ and have defined:  
\begin{align}
& e_k^{(w)} (\Delta, x, y; \theta)  \nonumber  \\
& \quad  := 
\sum_{1 \le j \le M-1}
\sum_{\substack{\alpha \le \beta^{[j]}}} 
\sum_{\gamma \in \mathscr{J}_w (\alpha)} 
\tfrac{K (\alpha)}{\alpha!} \, 
\mathscr{W}_\gamma^{[\alpha]} (x, \theta)   
\mathscr{H}_\gamma^{(w)} (\Delta, x, y; \theta)
\cdot \mathbf{1}_{|\alpha| + j - \|\gamma\|_w = \tfrac{k}{2}};  
\label{eq:e} 
\\ 
& \mathscr{R}_3^{(w)} (\Delta, x, y; \theta)   \nonumber \\ 
&\quad  :=
\sum_{1 \le j \le M-1}
\sum_{\substack{\alpha \le \beta^{[j]}}} 
\sum_{\gamma \in \mathscr{J}_w (\alpha)} 
\Delta^{|\alpha| + j - \| \gamma \|_w } \,
\tfrac{K (\alpha)}{\alpha!} \, 
\mathscr{W}_\gamma^{[\alpha]} (x, \theta) \,  
\mathscr{H}_\gamma^{(w)} (\Delta, x, y; \theta)
\cdot \mathbf{1}_{|\alpha| + j - \|\gamma\|_w \ge \tfrac{M}{2}}.   \nonumber 
\end{align}
Under Assumptions \ref{ass:diff}--\ref{ass:dim}, from Lemma \ref{lemma:deriv_LDL}, there exists a constant $C > 0$ such that:  
\begin{align} \label{eq:bd_r3}
\bigl| 
\mathscr{R}_3^{(w)} (\Delta, x, y; \theta)
p_\Delta^{\bar{X}, (w)} (x, y; \theta) 
\bigr| \le C \Delta^{\tfrac{M}{2}} 
\bigl| \mathscr{G}^{(w)} (\Delta, x, y, \theta) \bigr|,
\end{align}
for all $(\Delta, y, \theta) \in (0,1) \times \mathbb{R}^N \times \Theta$. 
Working with Theorems \ref{thm:bd_r1}--\ref{thm:bd_r2}, (\ref{eq:cf}), (\ref{eq:arg}) and (\ref{eq:bd_r3}), we finally obtain the following `minimal' representation of density expansion for diffusion models (\ref{eq:ellip}) and (\ref{eq:hypo}). 
\begin{theorem}[$\Delta$-Expansion] \label{thm:main}
Let $(\Delta, x, y, \theta) \in (0,1) \times \mathbb{R}^N \times\mathbb{R}^N \times \Theta$. 
Under Assumptions \ref{ass:diff}--\ref{ass:dim}, the transition density admits the following expansion. For every $J \in \mathbb{N}$ and $w \in \{\ref{eq:ellip}, \ref{eq:hypo}\}$:  
\begin{align} \label{eq:series}
\begin{aligned}
p^{X, (w)}_\Delta (x, y; \theta) 
= p_\Delta^{\bar{X}, (w)} (x, y ; \theta)   
\cdot \Bigl\{
1 + \sum_{1 \le j \le J} \Delta^{\tfrac{j}{2}} \cdot e_{j}^{(w)} (\Delta, x, y ; \theta) 
\Bigr\} 
+  
\mathscr{R}^{J, (w)} (\Delta, x, y; \theta).   
\end{aligned}
\end{align}
The coefficients $e_j^{(w)}$ are explicitly determined in (\ref{eq:e}). Also, for $1 \le j \le J$, there exists a constant $C > 0$ such that for all $(\Delta, y, \theta) \in (0,1) \times \mathbb{R}^N \times \Theta$, 
\begin{gather*}
\bigl| 
e_j^{(w)} (\Delta, x, y ; \theta) \times p_\Delta^{\bar{X}, (w)} (x, y; \theta) 
\bigr| 
\le C \, 
\bigl|
\mathscr{G}^{(w)} (\Delta, x, y, \theta)
\bigr|.  
% \int_{\mathbb{R}^N} \bigl|  
% e_k^{(w)} (\Delta, x, y ; \theta) \times p_\Delta^{\bar{X}, (w)} (x, y; \theta)    
% \bigr| dy \le C_2.   
\end{gather*}
% 
% and 
% \begin{align*}
% \int_{\mathbb{R}^N} \bigl|  
% e_k^{(w)} (\Delta, x, y ; \theta) \times p_\Delta^{\bar{X}, (w)} (x, y; \theta)    
% \bigr| dy \le C_2.   
% \end{align*}  
% 
For the residual $\mathscr{R}^{J, (w)}$, there exist constants $C_1, C_2, C_3 > 0$ such that for all $(\Delta, y, \theta) \in (0,1) \times\mathbb{R}^N \times \Theta$: 
\begin{gather}
\label{eq:resi}
\bigl| 
\mathscr{R}^{J, (w)} (\Delta, x, y; \theta) 
\bigr| 
\leq 
C_1 \Delta^{\tfrac{J+1}{2}}
\bigl| 
\mathscr{G}^{(w)} (\Delta, x, y, \theta) 
\bigr|
\leq 
C_2 \Delta^{\tfrac{J+1}{2} - \tfrac{m(w)}{2}}
\end{gather}
{and 
\begin{align}
 \int_{\mathbb{R}^N}
 \bigl| 
 \mathscr{R}^{J, (w)} (\Delta, x, y; \theta) \bigr| dy 
 \leq C_3 \Delta^{\tfrac{J+1}{2}}. \label{eq:R_L1}   
\end{align} }
% % 
% 
% 
% 
\end{theorem}  
{We note that the pointwise error bound (\ref{eq:resi}) differs across model classes \eqref{eq:ellip} and \eqref{eq:hypo} due to $m(w)$ taking a larger value in the latter case. In brief,  this is due to $X_S$ in (\ref{eq:hypo}) being a smooth component, driven by a Gaussian noise $\int_0^\Delta B_s ds$ of size $\mathcal{O}(\Delta^{3/2})$ rather than by $B_\Delta$ of size $\mathcal{O} (\Delta^{1/2})$ in the case of $X_{R}$, thus $X_S$ has a smaller variance for a fixed $\Delta \in (0, 1)$. I.e., the existence of the smooth component $X_S$ in (\ref{eq:hypo}) leads to a sharper density and/or concentration around the mode in the $X_S$ coordinate. However, in terms of $L_1$-error, its order only depends on the choice of $J$ and not on the model class because $y \mapsto \mathscr{G}^{(w)} (x, y, \theta)$ can be treated as a Gaussian density for any $(x, \theta, w) \in \mathbb{R}^d \times \Theta \times \{\ref{eq:ellip}, \ref{eq:hypo} \}$; recall also (\ref{eq:G_bd}).} 
\begin{remark}
\label{rem:sub}
Let $\textstyle \pi^{[J],(w)} ( \Delta, x, y; \theta) := \sum_{1 \le j \le J} \Delta^{\tfrac{j}{2}} \cdot e_{j}^{(w)} (\Delta, x, y ; \theta)$.
To avoid negative values for $\pi^{[J],(w)}$, we use a standard technique (see e.g.~\cite{stra:10} for a related approach)
%
%\begin{align*}
where 
$1+\xi = \textstyle \exp\big\{\log (1+\xi)\big\} = 
\exp\{ T_{J'}(\xi)\}\, \exp\{ R_{J'}(\xi) \}$,
%\end{align*}
%
for $J'\ge 1$, with $T_{J'}(\xi):=\textstyle\sum_{j=1}^{J'} (-1)^{j+1} \tfrac{\xi^{j}}{j}$ the $J'$-order Taylor expansion of $\xi\mapsto \log(1+\xi)$ and $R_{J'}(\xi)$ its residual. Via simple arguments, for $|\xi|<\delta<1$ one has $\big|(1+\xi)-\exp\{T_{J'}(\xi)\}\big|\le C \delta^{{J'}+1}$, for $C>0$. The above suggests the use of the following proxy: 
%
%\begin{align*}
%1+\pi^{[J],(w)} ( \Delta, x, y; \theta) \longleftrightarrow \exp\Big\{  
%T_{J'}\big(\pi^{[J],(w)} ( \Delta, x, y; \theta) \big)
%\Big\}.
%\end{align*}
%
\begin{align}
\label{eq:practice}
\widetilde{p}_\Delta^{(w)} (x, y; \theta) :=  p_\Delta^{\bar{X}, (w)} (x, y ; \theta)   
\cdot \exp\Big\{  
T_{J'}\big(\pi^{[J],(w)} ( \Delta, x, y; \theta) \big)
\Big\}.
\end{align}
Thus, $\pi^{[J],(w)}$ includes  powers $\Delta^{1/2},\ldots,\Delta^{J/2}$ (assuming non-zero $e_j$'s), so is of size $\delta=\mathcal{O}(\Delta^{1/2})$ and the residual in (\ref{eq:series}) is $\mathcal{O}(\Delta^{(J+1)/2})$ -- in the sense of the first bound in (\ref{eq:resi}). For the replacement by the Taylor approximation to only affect terms of size $\mathcal{O}(\Delta^{(J+1)/2})$, one should select $J'$ as the smallest even integer so that $J'\ge J$. An even $J'$ guarantees integrability of the  density proxy.
\end{remark} 

\section{Numerical Experiments} \label{sec:num}
%We present numerical results to illustrate the effectiveness of our CF expansion in an applied setting. 
We focus on the
bivariate \emph{FitzHugh-Nagumo (FHN)} SDE used in neuroscience. This model writes as:  
\begin{align} 
\begin{aligned}  \label{eq:FHN} 
d V_t & = \tfrac{1}{\varepsilon} \bigl( V_t - V_t^3 - U_t - s \bigr)  dt; \qquad  
d U_t & = \bigl( \gamma V_t  - U_t + \beta \bigr) dt
+ \sigma \, d B_{1, t}, 
\end{aligned}  
\end{align} 
with $V$ describing the membrane potential of a single neuron and the recovery variable $U$ expressing the ion channel kinetics.
Also, $s$ is the magnitude of the stimulus current and is often controlled 
%in practice, 
and $\theta = (\epsilon, \gamma, \beta, \sigma)$ is the parameter.
This SDE does not satisfy the boundedness conditions in Assumption~\ref{ass:coeff} as there is a non-Lipschitz term in the drift.
Statistical inference for the FHN SDE is an important topic from a theoretical and a practical viewpoint, see \cite{dit:19, mel:20, sam:24}. SDE (\ref{eq:FHN}) belongs in class (\ref{eq:hypo}) as the weak H\"ormander's condition (specifically, Assumption \ref{ass:hor}-II)  holds in this case. The transition density is intractable and we approximate it with the CF-expansion given in Section \ref{sec:series}. We  investigate the accuracy of the 
CF-expansion in Section~\ref{sec:num_DE} and 
% and observe that the approximation error is gradually removed as we include more expanded terms, i.e., as the number `$M$' in the proposed density expansion (\ref{eq:series}) increases. 
use  the expansion to carry out Bayesian inference with real data in Section \ref{sec:num_HMC}.  
\subsection{Accuracy of the CF-Expansion} \label{sec:num_DE} 
%We here present numerical results to validate the proposed density expansion given in the form (\ref{eq:series}). 
We produce two expansions via use of different baseline Gaussian densities. In particular, for a given initial value $x = (x_1, x_2) \in \mathbb{R}^2$, $\Delta > 0$ and $J \in \mathbb{N}$, we work with the CF-expansions: 
\begin{align}
y \mapsto p_\Delta^{(\mathrm{\iota}), [J]} (x, y; \theta) := \bar{p}_\Delta^{\, (\mathrm{\iota})} (x, y; \theta) \times \Bigl\{1 + \sum_{1 \le k \le J}
\Delta^{\tfrac{k}{2}}  \cdot  e^{(\iota)}_k (\Delta, x , y)  \Bigr\}, 
\qquad \iota \in \{ \mathrm{I}, \mathrm{II} \},   
\end{align}
with the reference density $\bar{p}^{(\iota)}_{\Delta}(\cdot)$, $\iota\in \{ \mathrm{I}, \mathrm{II} \}$, corresponding to the following `full' (for $\iota=\mathrm{I}$) or `partial' (for $\iota=\mathrm{II}$) LDL scheme: 
\begin{align}
d \bar{X}_t^{(\iota)} = 
\bigl( 
A^{(\mathrm{\iota})}_{x, \theta} \bar{X}_t^{(\iota)} + b_{x, \theta}^{(\iota)}
\bigr) dt  
+ 
\begin{bmatrix}
0 \\ 
\sigma 
\end{bmatrix} 
% [0,\sigma]^{\top}
d B_{1, t}, 
\label{eq:ldl_fhn}
\end{align} 
where we consider the following two choices:
\begin{gather*}
A^{(\mathrm{I})}_{x, \theta} 
:=  
\begin{bmatrix} 
\bigl( 1 - 3 (x_1)^2 \bigr) / \varepsilon
& - 1/\varepsilon \\[0.1cm]
\gamma & -1 
\end{bmatrix}, \qquad 
A^{(\mathrm{II})}_{x, \theta} 
:= 
\begin{bmatrix} 
\bigl( 1 - 3 (x_1)^2 \bigr) / \varepsilon
& - 1 / {\varepsilon} \\[0.1cm]
0 & -1 
\end{bmatrix}; \\
b_{x, \theta}^{(\iota)} = 
\begin{bmatrix}
\bigl(x_1 - (x_1)^3 - x_2 + s \bigr) / \varepsilon \\[0.1cm]
\gamma x_1 - x_2 + \beta 
\end{bmatrix}
- A_{x, \theta}^{(\iota)} \, x.  
\end{gather*} 
The $e_k^{(\iota)}$'s are found starting from expansion (\ref{eq:err}), with $\mathscr{L}_\theta^{0, z}$ (used by the differential operator~$\mathscr{D}_{\alpha}^{z,\theta}$) corresponding to the generator associated with (\ref{eq:ldl_fhn}) for $\iota \in \{ \mathrm{I}, \mathrm{II} \}$, and then re-arranging terms in powers of $\sqrt{\Delta}$ as described in Section \ref{sec:series}. Matrix $A^{(\mathrm{II})}$ is upper-triangular, so the baseline  $\bar{p}^{(\mathrm{II})}_{\Delta}$ takes a simpler form compared to when using $A^{({\mathrm{I}})}$. In both cases, the reference Gaussian laws are non-degenerate. 
%We recover the $w_k^{(\iota)}$'s after sorting the produced CF-expansion terms 
%w.r.t.~$\Delta^{k/2}$, $k = 1, \ldots, 5$, as explained in Section \ref{sec:series}. 
To calculate the $e_k^{(\iota)}$'s we use \texttt{Mathematica} with full expressions given in Section~\ref{app:supp_num} in Supplementary Material. Due to the SDE noise being additive, we have $w_{k}^{(\iota)} = 0$, for $k = 1,2$, $\iota \in \{\mathrm{I}, \mathrm{II}\}$. 
Thus, the CF-expansions with $J\in\{1,2\}$ coincide with the baseline.
% in this model example, which implies that the above two drift linearisation already approximates the target well. However, as we shall observe later, adding more expanded terms leads to more accurate density estimation. 
The reference density for $\iota=\mathrm{I}$ involves full  linearisation, so the $e_k^{(\mathrm{I})}$'s have simpler expressions than the $e_k^{(\mathrm{II})}$'s, see Section~\ref{app:cf_fhn} in Supplementary Material for details. 

We choose $s = 0.01$, initial value $x = (V_0, U_0) = (-0.1, 0.2)$ and  $\theta = (\varepsilon, \gamma, \beta, \sigma) = (0.1, 1.2, 0.3, 0.8)$. We consider $\Delta \in\{ 0.1, 0.05, 0.02\}$ and compute 
CF-expansions using the transform $\widetilde{p}_\Delta$ described in Remark \ref{rem:sub}, 
which we denote here
$\widetilde{p}_\Delta^{\, (\iota),[J]}$, 
$\iota \in \{\mathrm{I}, \mathrm{II}\}$. We try $J=2,3,4,5$, and for $\widetilde{p}_\Delta$ we set $J'=2$, as 
%That is, the CF expansion ultimately employed writes as 
the correction term includes powers $\Delta^{3/2},\ldots,\Delta^{J/2}$, $J\le 5$, and the transform can only affect terms of size $\mathcal{O}((\Delta^{3/2})^{(J'+1)})=\mathcal{O}(\Delta^{9/2})$. %i.e.~smaller than $\mathcal{O}(\Delta^{3})$. 
% 
%\begin{align} \label{eq:non_neg_de}
%\widetilde{p}_\Delta^{\, (\iota), [J]} (x, y; \theta) := 
%\bar{p}_\Delta^{\, (\mathrm{\iota})} (x, y; \theta) \times \exp \Bigl( G \bigl( \pi^{(\iota), [J]} ( \Delta, x, y; \theta) \bigr) \Bigr), \quad \iota \in \{ \mathrm{I}, \mathrm{II} \},  
%\end{align} 
% 
%with $\xi \mapsto  G(\xi) := \xi - \xi^2/2$, $\xi\in\mathbb{R}$, and $\pi^{(\iota), [J]} ( \Delta, x, y; \theta) := \textstyle\sum_{k = 1}^J 
%\Delta^{\tfrac{k}{2}}  \cdot  w^{(\iota)}_k (\Delta, x, y)$. 
%
We find the benchmark `true' density via a simulation that: (i) uses $2\times 10^7 $ samples from the FHN SDE at $\Delta$ via an  EM scheme with discretisation step $\Delta/800$; (ii) applies a standard Kernel Density Estimator (KDE) approach to reconstruct the density.  
We write the benchmark density as $p_\Delta^{(\mathrm{B})}$. 
Densities are evaluated on a regular $51\times 51$ grid  $D = \bigl\{ (s_i, r_j) \, | \, 0 \le i,j \le 50 \bigr\} \subset \mathbb{R}^2$ for reals $r_0 < \cdots < r_{50}$ and $s_0 < \cdots < s_{50}$ defined in an apparent way. 

\begin{figure}
\centering 
\includegraphics[width=12cm]{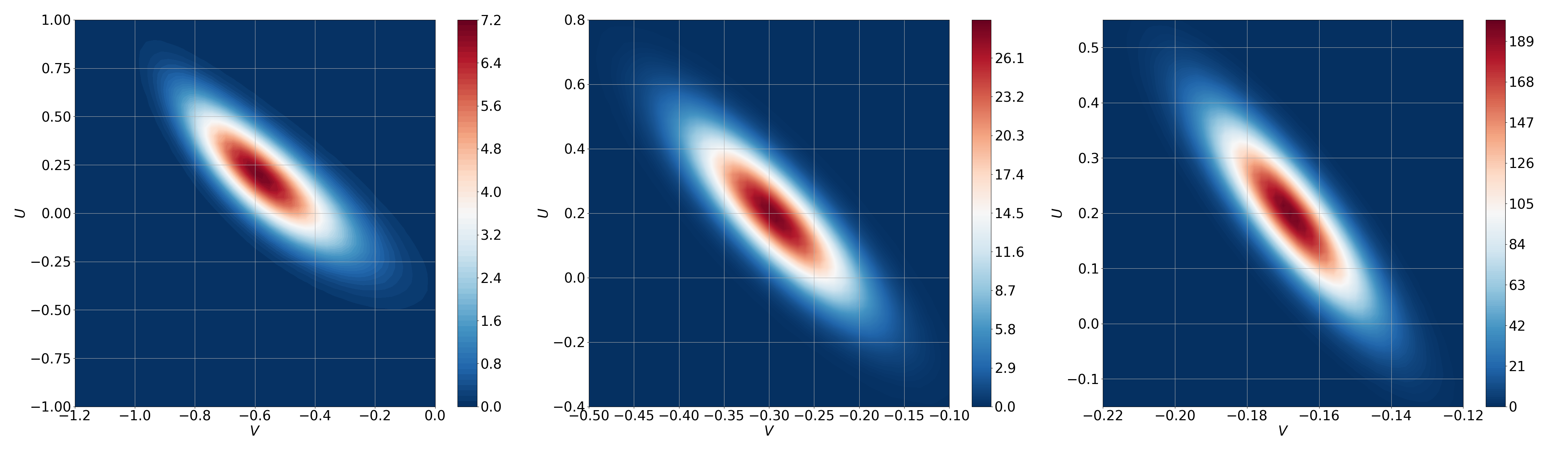}   
\caption{Contours of the benchmark densities $p_\Delta^{(\mathrm{B})}$. Left: $\Delta = 0.1$. Middle: $\Delta = 0.05$. Right: $\Delta = 0.02$.} 
\label{fig:benchmark} 
\end{figure}

\begin{figure} 
\centering
\subfigure[DE-I.]{\includegraphics[width=11.5cm]{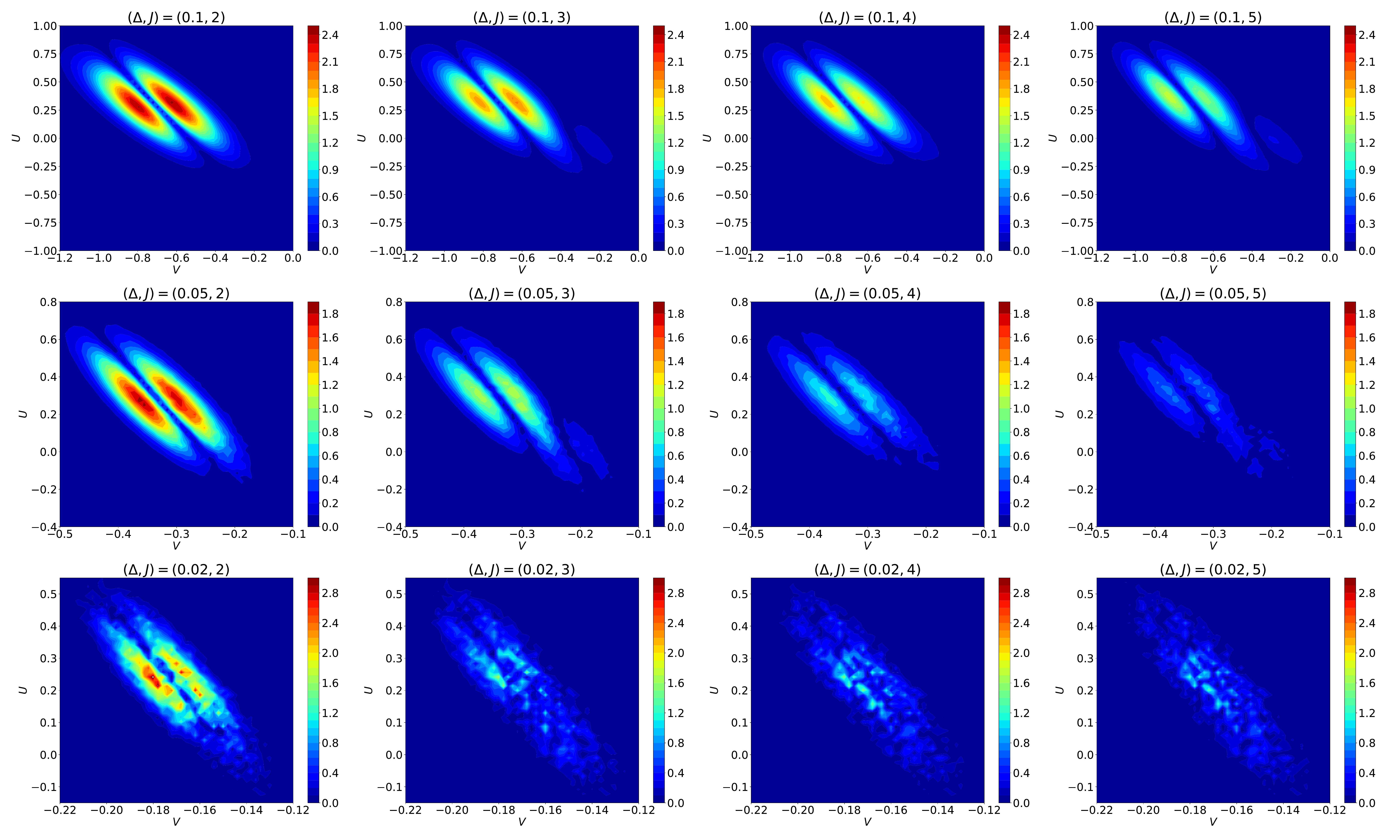}} 
\subfigure[DE-II.]{\includegraphics[width=11.5cm]{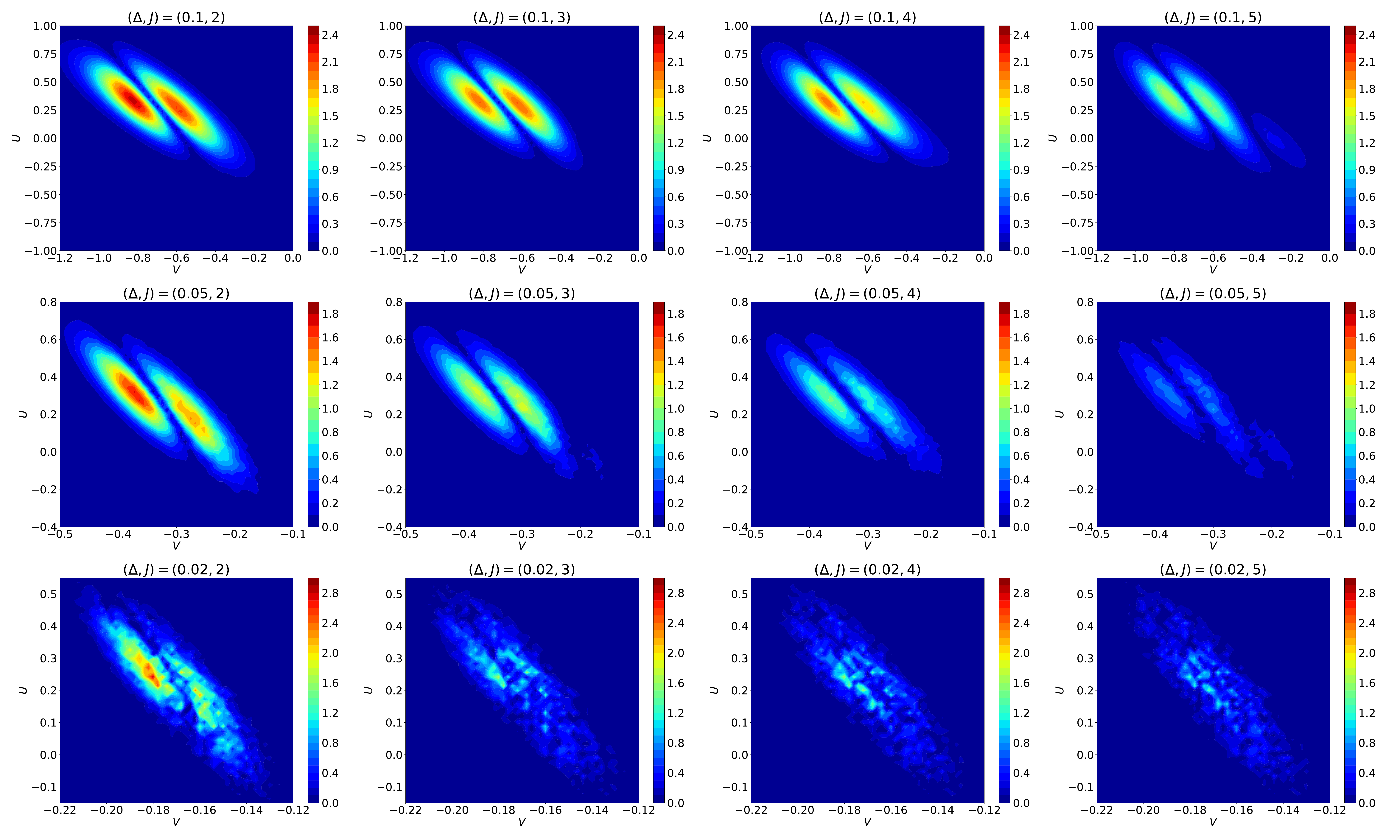}}  
%\begin{subfigure}
%  \centering  \includegraphics[width=1\linewidth]{Contour-density-abs-err-DE-II-delta-0.1-set2.png}
%\end{subfigure}
%\begin{subfigure}
%  \centering
%  \includegraphics[width=1\linewidth]{Contour-density-abs-err-DE-II-delta-0.05-set2.png}
%\end{subfigure}
%\begin{subfigure}
%  \centering
%  \includegraphics[width=1\linewidth]{Contour-density-abs-err-DE-II-delta-0.02-set2.png}
%\end{subfigure}
\caption{\small Heatmap of the absolute error for the CF-expansion in the case of the hypo-elliptic FHN model (see Section \ref{sec:num}). Rows correspond to 3 choices $\Delta=(0.1,0.05,0.02)$ and columns to 4 choices $J=2,3,4,5$.} 
\label{fig:ae}
\end{figure}

\begin{figure}
\centering
\subfigure[$L_1$-error]{\includegraphics[width=0.35 \textwidth]{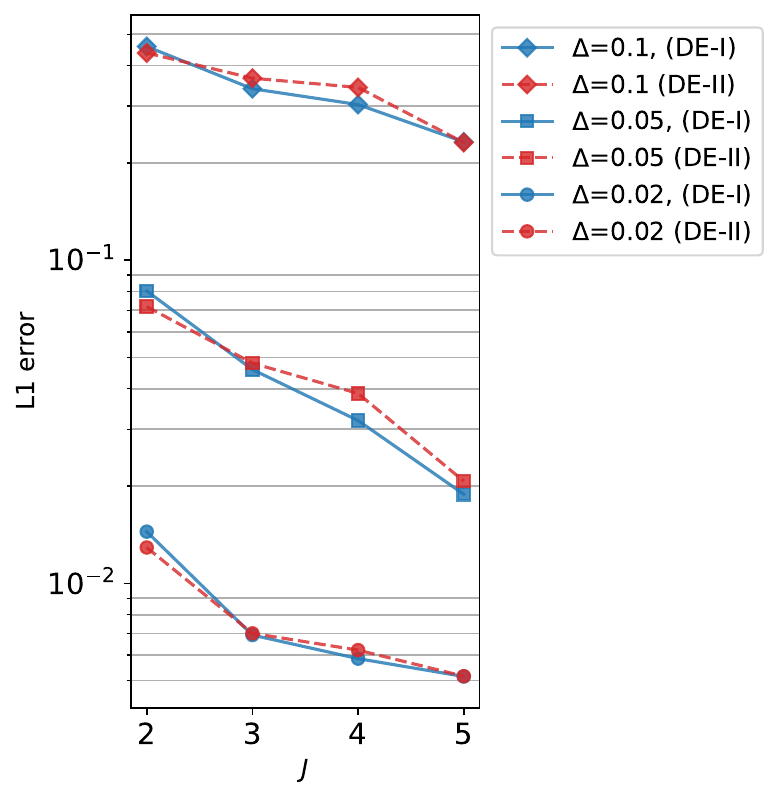} \label{fig:L1_err}}
\hspace{1cm}
\subfigure[Computational time per grid]{\includegraphics[width=0.4\textwidth]{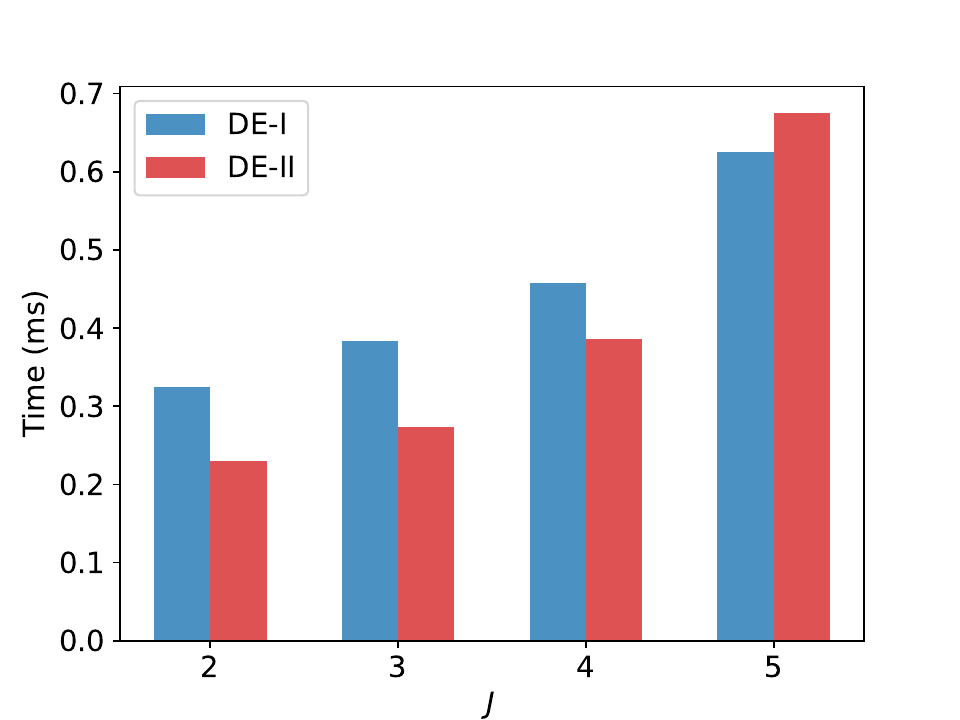}  \label{fig:ct}}
\caption{Summary for performance of density expansions.}
\label{fig:summary} 
\end{figure}

Fig.~\ref{fig:benchmark} shows the contours of the benchmark  $p_\Delta^{(\mathrm{B})}$, {which indicate the unimodality of the target transition densities.} Fig.~\ref{fig:ae} plots the absolute errors, 
$$
\mathscr{E}_\Delta^{(\mathrm{\iota}), [J]} (x, y; \theta) :=  \bigl| p_\Delta^{(B)} (x, y; \theta)  
-  \widetilde{p}_\Delta^{(\iota), [J]} (x, y; \theta) \bigr|, \quad y \in D, 
$$ 
between $p_\Delta^{(\mathrm{B})}$ and the 
CF-expansions of order $J= 2, 3, 4, 5$.  
% 
%$
%\mathscr{E}_\Delta^{(\mathrm{\iota}), [J]} (x, y; \theta) :=  \bigl| p_\Delta^B (x, y; \theta)  
% -  \widetilde{p}_\Delta^{(\iota), [J]} (x, y; \theta) \bigr|$, $y \in D$, $\iota \in 
% \{\mathrm{I}, \mathrm{II} \}$. 
Fig.~\ref{fig:summary} summarises the overall performance of the 
CF-expansions. Fig. \ref{fig:L1_err} gives the $L_1$-error of the CF-expansions, defined as   
$
\mathrm{L}^{(\iota), [J]}_1 (\Delta, x ;  \theta) := \textstyle \sum_{y\in D} \mathscr{E}_\Delta^{(\iota), [J]} (x, y ; \theta) \times \delta_V \times \delta_U 
$, where $\delta_V = (s_{50} - s_0) / 50, \, \delta_V = (r_{50} - r_0) / 50$. Fig.~\ref{fig:ct} shows the average running time of DE-I and DE-II (denoting the two density expansions for $\iota\in\{\mathrm{I},\mathrm{II}\}$), with the average taken from the 3  choices of $\Delta$. Fig.~\ref{fig:ae} demonstrates that absolute errors diminish as $J$ increases. %(plots moving from left to right) 
%for both DE-I and DE-II.
{In particular, the error of single mode and variance (or higher order moments) between the benchmark and approximate transition densities gradually diminishes as $J$ increases}. We observe a similar decrease in $L_1$-error in Fig.~\ref{fig:L1_err}. Note that the errors by the two CF-expansions with $J=5$ are less than half of those with $J=2$. Also, errors decrease for smaller $\Delta$. 
In terms of computing cost, for the DE with $J=2$, i.e.~the baseline Gaussian density without correction, DE-II %(CF-expansion based on the modified LDL scheme) 
is computationally cheaper due to the simpler expression in the matrix exponential.  Costs are similar between DE-I with $J=2$ and DE-II with $J=4$. 
Costs grow as $J$ increases, but the growth rate seems faster in DE-II since the latter makes use of the simpler but slightly less accurate baseline density, thus involves more correction terms as $J$ grows.      

\subsection{Application to Bayesian Inference} \label{sec:num_HMC} 

\subsubsection{MCMC via CF-expansion -- Design of Posterior}
\label{sec:mcmc_CF}
We use our CF-expansion to carry out Bayesian inference for SDEs. 
In this subsection we consider general SDEs rather than just the FHN SDE as the approach is relevant in a wide setting. 
Via the CF-expansion we obtain a posterior law that can be integrated within well-established MCMC methodologies, including  centred/non-centred parameterisations \citep{papa:07}, Particle MCMC and Particle Gibbs algorithms \citep{andrieu:10}.
Note that early literature \citep{stra:10} investigated the use of CF-expansions (for elliptic models) within a standard Metropolis-Hastings method under centred parametrisation, thus the options provided were limited.
%We obtain a posterior use  involving CF-expansion is well-combined with various computational MCMC methodologies, including both so-called \emph{centred}/\emph{non-centred parameterisation} approaches, as we shall see later. 
%While early literature \citep{stra:10} studied the use of CF expansion in the context of standard Metropolis-Hastings algorithm under centred-parametrisation and so the computational options of MCMC are extremely limited, the discussion below opens up wider options, including the particle-based such as Particle MCMC or Particle Gibbs \citep{andrieu:10}, which requires sampling from the tractable (approximate) transition density.
%
%the developed methodology can be applied to a broad class of SDEs and describe how to integrate the CF-expansion into the construction of tractable posterior density. 
%In particular, we clarify that the
%The obtained posterior involving CF-expansion is well-combined with various computational MCMC methodologies, including both so-called \emph{centred}/\emph{non-centred parameterisation} approaches, as we shall see later. While early literature \citep{stra:10} studied the use of CF expansion in the context of standard Metropolis-Hastings algorithm under centred-parametrisation and so the computational options of MCMC are extremely limited, the discussion below opens up wider options, including the particle-based such as Particle MCMC or Particle Gibbs \citep{andrieu:10}, which requires sampling from the tractable (approximate) transition density.
Particle-based MCMC methods require sampling from the SDE transition density, i.e.~in our case from the CF-expansion used as its proxy.
It is typically difficult to simulate from the CF-expansion. However, notice  that the expansion writes as `Gaussian density' $\times$ `correction term'. Thus, particle-based MCMC and general Sequential Monte Carlo (SMC) methodology can be implemented using the baseline density (which we can sample from) with the correction term being attached in the `weights' within the algorithm.
%and the weights of particles are computed from some test function (e.g. a tractable density of assumed measurement noise) with the correction term being attached. 
Furthermore, the CF-expansion structure of `Gaussian density' $\times$ `correction term' permits a non-centred approach -- such an algorithm turns out to be the most effective one in our numerics in the next section. We provide more details on the mentioned algorithms directly below.

Consider the data $\mathscr{Y}_n = \{ Y_{t_k} \}_{0 \le k \le n}$ at instances $t_k$, $0 \le k \le n$, for which we assume an equidistant step-size $\Delta$. We consider the setting of noisy observations, so that there is a density $p(Y_{t_k}|X_{t_k})$, assumed known.
Under a \emph{data augmentation} approach, 
we set $\mathbf{q} := (\theta, \mathscr{X}_n) \in \mathbb{R}^{d_\theta} \times \mathbb{R}^{N \times (n+1)}$, where $\mathscr{X}_n := \{X_{t_k}\}_{1 \le k \le n}$.
%
%We write the available data set as $\mathscr{Y}_n \equiv \{ Y_{t_k}^{\mathrm{obs}} \}_{0 \le k \le n}$, observed at $(n+1)$ time instances $t_k, \, 0 \le k \le n$ with an equidistant step-size $\Delta$. To treat a general inferential problem in literature, we assume that the data $\mathscr{Y}_n$ is obtained from some partial coordinates of the SDE $X_{t_k} \equiv (Y_{t_k}, Z_{t_k}) \in \mathbb{R}^N$, depending on a parameter vector $\theta \in \Theta \subseteq  \mathbb{R}^{d_\theta} $,  with measurement noise as :  
% 
%$ 
%Y_{t_k}^{\mathrm{obs}} = Y_{t_k} + \epsilon_{t_k}, \,  0 \le k \le n, 
%$
% 
%where $\{ \epsilon_{t_k} \}_{0 \le k \le n}$ is the independent sequence of measurement noise following a tractable distribution, and then the Lebesgue density of the conditional law of $Y_{t_k}^{\mathrm{obs}} | Y_{t_k} = y$, expressed as $\xi \mapsto f (\xi| y)$, is analytically available. 
%Because the marginal likelihood of $\mathscr{Y}_n$ given the parameter $\theta$ is highly complex, we adopt \emph{Data Augmentation} (DA), i.e., data is augmented with the hidden components $\mathscr{X}_n \equiv \{X_{t_k}\}_{1 \le k \le n}$, and then design the (approximate) tractable posterior distribution defined on the space of augmented variable $\mathbf{q} \equiv (\theta, \mathscr{X}_n) \in \mathbb{R}^{d_\theta} \times \mathbb{R}^{N \times (n+1)}$. If the true density of transition dynamics $X_{t_{k+1}}| X_{t_k} = x$, denoted as $y \mapsto p_\Delta^X (x, y; \theta)$, is analytically available, 
The posterior density on the augmented state $\mathbf{q} = (\theta, \mathscr{X}_n)$ writes as: 
\begin{align}
\label{eq:posterior}
{P} ( \mathbf{q} \, | \, \mathscr{Y}_n ) 
\propto \Bigl\{  \prod_{0 \le k \le n} p (Y_{t_k} | X_{t_k}) \Bigr\} 
\times \Bigl\{ \prod_{1 \le k  \le n} p_\Delta^X (X_{t_{k-1}}, X_{t_k}; \theta) \Bigr\} \times p_0 (X_0) \times p_\theta (\theta),  
\end{align}
where $p_0$, $p_\theta$ denote priors on the initial value $X_0$ and the parameter $\theta$, respectively. 
%We remark that the above posterior is generally unavailable since the transition density is. 
% \begin{align} 
% \ell \bigl( \mathcal{Z}_n \,  | \, \theta, X_0  \bigr) 
% \equiv \sum_{k = 0}^n \log  f (Y_{t_k} | V_{t_k}) 
% + \sum_{k = 1}^n \log p_\Delta^X ( X_{t_{k}}, X_{t_{k-1}}; \theta). 
% \end{align}
% 
We replace the true transition density with the CF-expansion as given in (\ref{eq:practice}), that is: 
\begin{align}
\label{eq:a-posterior}
\prod_{1 \le k \le  n} p_\Delta^X (X_{t_{k-1}}, X_{t_{k}}; \theta) 
\approx
\prod_{1 \le k \le n} 
p_\Delta^{\bar{X}} (X_{t_{k-1}}, X_{t_k}; \theta)
\times 
\prod_{1 \le k \le n} \exp \Bigl( T_{J'} \bigl( \pi^{[J]} ( \Delta, X_{t_{k-1}}, X_{t_k}; \theta) \bigr) \Bigr).
\end{align} 
The approximate posterior obtained via (\ref{eq:posterior})-(\ref{eq:a-posterior}) can now be used within standard or particle-based MCMC methods:
(i) For standard MCMC, the `correction terms' can be treated as a part of the likelihood function, so that a-priori the dynamics of the $X$-process are determined by the baseline density. This allows for application of centred/non-centred algorithms, as in the latter case one can use as latent components the standard Gaussian noise that generates samples from the baseline density; (ii) For particle-based methods, the `correction terms' can become part of the weights and one can apply, e.g., particle filters by sampling from the tractable baseline density.
{\begin{remark} \label{rem:davie}
In the above, we have discussed the use of baseline Gaussian density as a `proposal' within the standard/particle-based MCMC computational framework with `correction terms' becoming part of the weights, rather than directly sampling from the approximate density of the form (baseline Gaussian density) $\times$ (1 + (correction)), as the latter approach is in general unavailable.  However, one may employ a methodology proposed by \cite{davie:22} who constructed a tractable sampling scheme via a corresponding density expansion in an elliptic setting, in a way so that the used expansion preserves a high order proximity in Wasserstein distance. Extension to hypo-elliptic SDEs is not straightforward and could be an interesting future research direction. 
\end{remark}}

\subsubsection{Experimental Design and MCMC Results}
\label{sec:results}
We apply our CF-expansion  %obtained transition density expansion to the inferential problem and, in particular, 
to carry out Bayesian inference for the FHN SDE (\ref{eq:FHN}) with the real dataset used in \cite{sam:24}. The data are available at \url{ https://data.mendeley.com/ datasets/ybhwtngzmm/1} which provides 20 neural recordings of the 5th lumbar dorsal rootlet from a single adult female rat with time length $250$ms and equidistant step-size $0.02$ms. 
In our study we choose a particular dataset, specifically the file  1554.mat from the above URL, which was obtained while the 5th lumber dermatome was stimulated. We  subsample the first $40$ms of data with a step-size $0.08$ms, i.e.~we have $(T, \Delta) = (40, 0.08)$ and a  number of datapoints $n = 501$, so  $\Delta$ is relatively large.
As in  \cite{sam:24}, we set $s = 0$ and focus on the parameter $\theta = (\epsilon, \gamma, \beta, \sigma)$. We assume that the data $\mathscr{Y}_n = \{ Y_{t_k} \}_{0 \le k \le n}$ are observed with a small measurement noise as $Y_{t_k} = V_{t_k} + \epsilon_{t_k}$, with $V$ the smooth coordinate in the FHN SDE and $\epsilon_{t_0}, \ldots, \epsilon_{t_n} \overset{i.i.d}{\sim} \mathscr{N} (0, 0.01^2)$. 
We adopt a non-centred parametrisation, assign log-normal priors on $\theta$, i.e., $\log \varepsilon, \log \beta, \log \gamma, \log \sigma \overset{\mathrm{i.i.d.}}{\sim} \mathscr{N} (0,1)$ and set $V_0 \sim \mathscr{N} (0, 0.1^2), \, U_0 \sim \mathscr{N} (0, 0.2^2)$ for the initial state. We employ Hybrid Monte Carlo (HMC) to sample from the posterior, using the Python package Mici (\url{https://pypi.org/project/mici/}) which offers a variety of MCMC methods based on Hamiltonian dynamics. We use a dynamic integration-time HMC implementation \citep{beta:17} with a dual-averaging algorithm \citep{hoff:14} to adapt the step-size of the leapfrog integrator. The mass matrix is set to identity. 
%We use the leapfrog integrator with St\"ormer-Verlet splitting for the time-discretisation of the Hamiltonian dynamics and set the mass matrix to identity. 

We consider 3 designs of tractable posteriors: [\textbf{P0}] Benchmark `true' posterior. This is constructed via a local Gaussian (LG) transition density scheme \citep{glot:21}, which provides an approximation of the transition density of the hypo-elliptic SDE (\ref{eq:hypo}) for a sufficiently small step-size. A data augmentation step is applied, whereby $d_M = 100$ signal points are added in-between observation pairs to eliminate the bias. 
% The dimension of augmented variable $\mathbf{q}$ is $d_q  = d_\theta +  (n-1)  \times d_M  \times d_w^{\mathrm{LG}} + N  = 100,006$, where $d_w^{\mathrm{LG}} \equiv 2$ is the number of Gaussian variables to produce one-step integration of LG scheme. 
The obtained posterior values are treated as the benchmark true values; 
%[\textbf{P1}] Posterior based on the LG scheme without any time-disretisation between data points, i.e.~P0 with $d_M = 1$;   
% $d_q = d_\theta + (n-1) \times d_w^{\mathrm{LG}} + N = 1,006$.  
%
[\textbf{P1}] Posterior based on the partial LDL scheme given in (\ref{eq:ldl_fhn}) with $\iota = \mathrm{II}$;  
% This is defined as (\ref{eq:nc_posterior}) with the baseline Gaussian density being the density of the modified LDL scheme and  $\widetilde{\pi}_n \equiv 1$.  $d_q = 1,006$.  
[\textbf{P2}] Posterior produced via implementation of the non-centred parameterisation of the initial target given by (\ref{eq:posterior})-(\ref{eq:a-posterior}), based on the CF-expansion around the partial LDL scheme with $J=3$. %For the correction term, the non-negative transformation is applied as indicated in (\ref{eq:non_neg_de}); %$d_q = 1,006$.      
%[\textbf{P4}] Posterior as in P3, with $J=4$. 
%(\ref{eq:nc_posterior}) based on the CF-expansion around the modified LDL scheme with $J=4$.
For each posterior, we run two HMC chains of 8,000 iterations with the first 4,000 iterations used as an adaptive warm-up phase. 

\begin{figure}
\centering
%\subfigure[\footnotesize P1 (Local Gaussian scheme)]{\includegraphics[width=7cm]{FHN-Posterior-local-gaussianvsbenchmark_overlaid.pdf} \label{fig:P1}}
%\hspace{1cm}
\subfigure[\footnotesize P1 (Partial LDL scheme)]{\includegraphics[width=7cm]{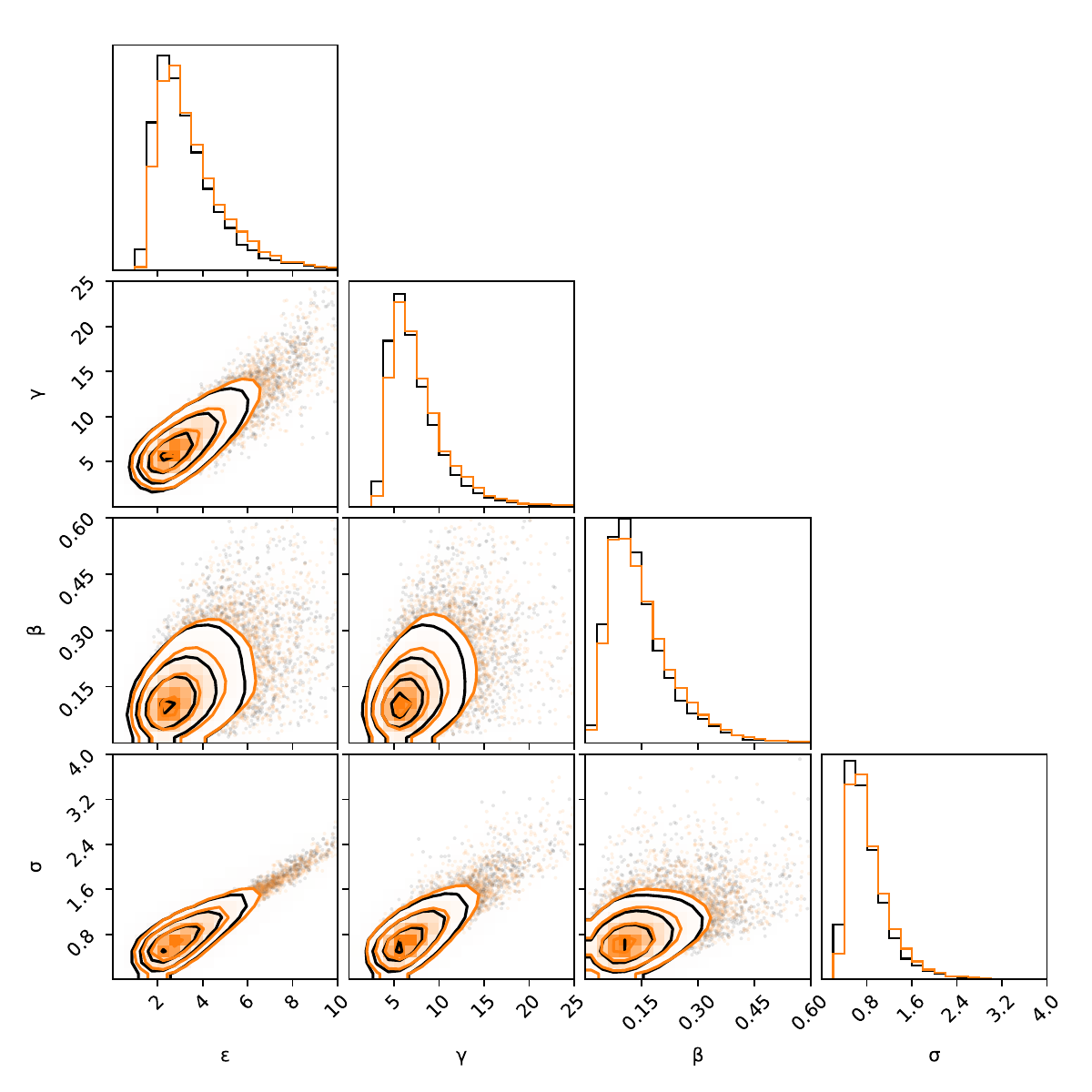}  \label{fig:P1}}
\subfigure[\footnotesize P2 (CF-expansion with $J = 3$)]{\includegraphics[width=7cm]{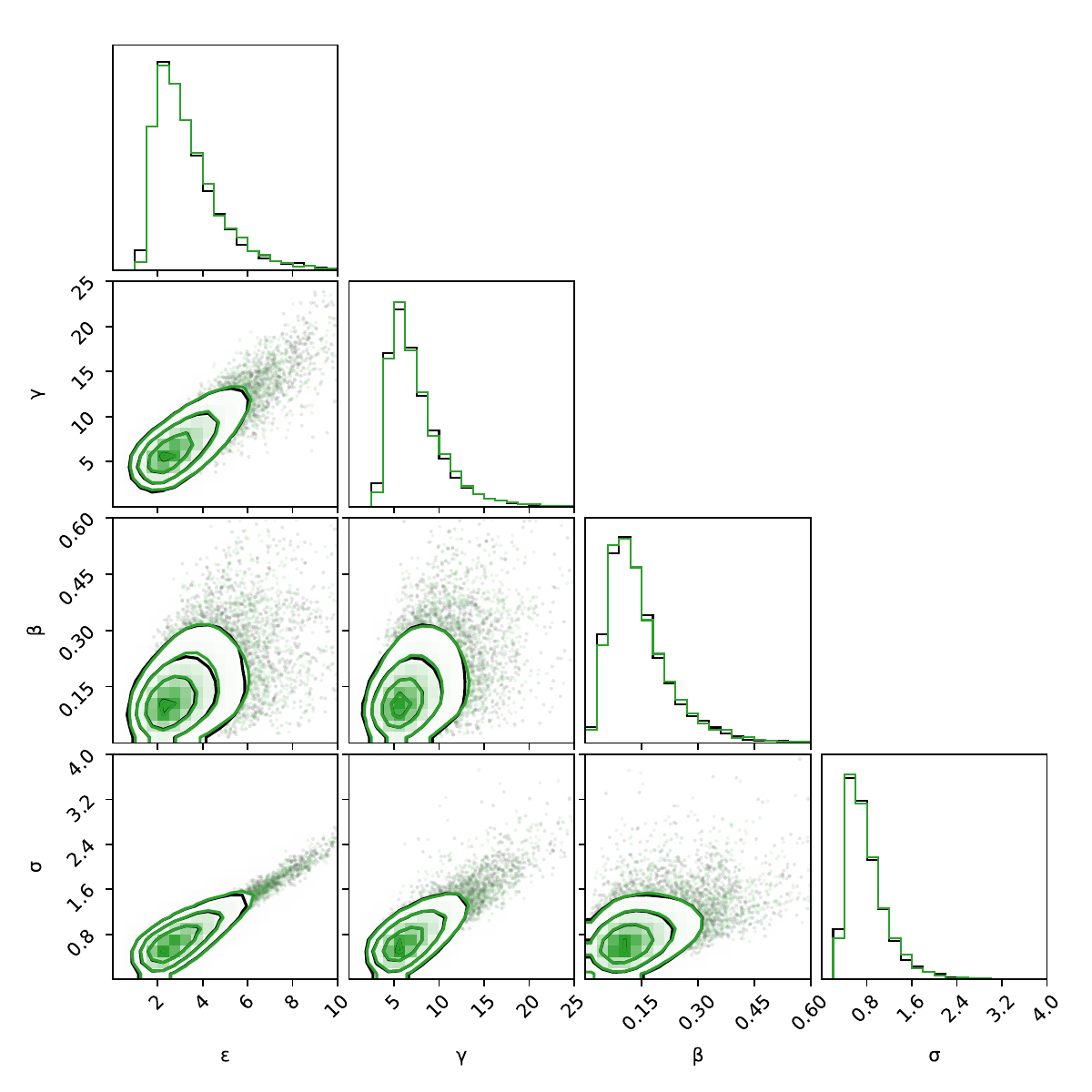} \label{fig:P2}}
%\hspace{1cm} 
%\subfigure[\footnotesize P4 (CF-expansion with $J = 4$)]{\includegraphics[width=7cm]{FHN-Posterior-DE-J-4vsbenchmark_overlaid.pdf}  \label{fig:P4}}
\caption{Posterior estimates. P0 (benchmark posterior) is overlaid in each figure in black.}
\label{fig:posterior} 
\end{figure}

Fig.~\ref{fig:posterior}  shows  results for targets P0-P2.  %presented in black, blue, orange, green and red colour, respectively. 
Results for the true posterior P0 are given in black and are overlaid in the sub-figures to observe the accuracy of posteriors P1-P2.
% and the gap between P0 and the others is interpreted as bias. 
Table \ref{table:hmc_fhn} shows average running times per iteration from two chains. 
Additional convergence diagnostics provided in Table \ref{table:hmc_suppl} of Section~\ref{sec:results_suppl} in Supplementary Material show similarly good convergence performance for all $3$ cases, we can thus conclude that the posteriors shown in Fig.~\ref{fig:posterior} are reliable. 
%From Fig.~\ref{fig:P1} we see that the contours of P1, i.e.~scheme without correction terms, deviate from P0. %though the gap is reduced for P2. 
In Fig.~\ref{fig:posterior} it is clear that P2 (i.e.~scheme based upon the CF-expansion) captures P0 more accurately than P1 (i.e.~scheme without correction terms) does. Thus, the inclusion of the  correction term  eliminates the bias even for $J=3$, with the algorithm targeting P2 having a computing cost approximately 10 times smaller than that of the benchmark P0 (see Table \ref{table:hmc_fhn}). 
%We remark that P1  achieves the smallest computational cost but fails to capture the benchmark P0. 
Our experiment implies that, in this case, the CF-expansion is effective both from the perspectives of computational cost and estimation accuracy. We remark that a centred-parametrisation led to MCMC chains with very poor convergence performance.% and great challenges at obtaining accurate posterior estimates.  
\begin{table}[h]
\caption{Computational cost of MCMC chains for the FHN model. Schemes: P0$\rightarrow$ benchmark; %P1$\rightarrow$ Local Gaussian; 
{P1}$\rightarrow$ Modified LDL; {P2} $\rightarrow$ CF-Expansion, $J=3$.} %P4$\rightarrow$ CF-Expansion, $J=4$.}
\label{table:hmc_fhn}
\centering 
\begin{tabular}{l|ccc}
\toprule 
scheme &  P0  & 
P1  & P2  %& 
%P3  & P4 
\\
time(sec)/iter & 4.745 & %0.118 & 
0.237 & 0.460 %& 0.548 
\\ 
\midrule
\end{tabular}
\end{table}  

{
\section{Discussion on Practical Perspectives}  \label{sec:discussion}
\begin{enumerate}
\item \textit{Data augmentation between data points}.  We have developed a density expansion and discussed its use in statistical inference for the setting where the step-size $\Delta$ between observations is less than $1$. In a general setting where $\Delta \ge 1$, one can still employ the density expansion in a \emph{Data Augmentation} framework (\cite{papa:13}), i.e., by imputing the latent variables between data points via time-discretisation. This can indeed be realised by generating a Markov chain of the baseline Gaussian scheme with the products of correction terms attached to the test function as a weight. As illustrated at the experiments in Section \ref{sec:num_HMC}, the use of the correction terms can lead to efficient inference with a smaller number of discretisations in-between data points, compared to the case without the corrections. Similar efficiency gains were studied and illustrated in \cite{igu:BJ} where a weak second-order  sampling scheme is compared with the Euler-Maruyama (weak first order) one in a Bayesian data augmentation framework. 
\\ 

\item \emph{Numerical properties of the LDL scheme.} 
Motivated by the use of the LDL scheme in a data augmentation framework, one may be interested in its numerical properties such as stability or (geometric) ergodicity. 
Though a full investigation is beyond the scope of this paper, we will make some comments below on the preservation of ergodicity by the LDL scheme. Let us consider the following Langevin-type equation: 
\begin{align} \label{eq:Langevin}
d X_t  = b (X_t) dt + \Sigma dB_t, 
\end{align}
where $X_t \in \mathbb{R}^d$, $\Sigma \in \mathbb{R}^{d \times d}$ and $\{ B_t \}_{t \ge 0}$ is a $d$-dimensional Wiener process. We assume    
standard sufficient conditions for (\ref{eq:Langevin}) to be ergodic, specifically, (i) minorisation and (ii) Lyapunov condition, see e.g. Lemma 2.3.  and Assumption 2.2., respectively, in \cite{matt:02}. We will check if such conditions are inherited by the LDL scheme applied to (\ref{eq:Langevin}). We first notice that this is not generally true for the Euler-Maruyama scheme unless the drift function is globally Lipschitz. In particular, the second condition can break down when the drift is only locally Lipschitz, while the minorisation condition can still hold, see, e.g., the proof of \cite[Corollary 7.4.]{matt:02}, regardless of the growth of the drift.  In a similar manner, the minorisation condition should hold for the LDL scheme as well, thus we focus on the Lyapunov condition. Here, we will show that such a preservation can occur for the LDL scheme applied to the following 1-dimensional SDE with non-globally Lipschitz drift: 
\begin{align} 
\label{eq:toy_ex}
d X_t = - X_t^3 dt + dB_{1, t}, \quad X_0  = x \in \mathbb{R}.  
\end{align}
We note that \cite[Lemma 6.3.]{matt:02} proved that for any step-size and initial state, the Euler-Maruyama scheme applied to (\ref{eq:toy_ex}) can be unstable with a positive probability, thus does not preserve the ergodic property. Then, the LDL scheme $\{Y_k \}_{k \ge 0}$ is defined as follows, for $0 < t_{n-1} < t_{n}$ and $\Delta = t_{n} - t_{n-1}$:
\begin{align} \label{eq:ldl_toy}
Y_{n} = \exp (- 3 Y_{n-1}^2 \Delta) Y_{n-1} + 2 Y_{n-1}^3 \int_0^\Delta \exp (- 3 Y_{n-1}^2 s) ds + \int_{t_{n-1}}^{t_n} \exp (- 3 Y_{n-1}^2 (t_{n} - s)) dB_{1,s}.   
\end{align}
%    
%Interestingly, \eqref{eq:ldl_toy} satisfies the Lyapunov condition as we show below. 
We use the Lyapunov function $V(x) = x^2$ and write $\mathcal{F}_{n-1}$ as the $\sigma$-algebra generated by the Markov chain $\{Y_{t_k}\}_{k \le n-1}$. We have that: 
\begin{align*}
\mathbb{E} [ V (Y_{n}) | \mathcal{F}_{t_{n-1}}] 
= m (\Delta, Y_{n-1})^2 + \int_0^\Delta \exp (- 6 Y_{n-1}^2 s) ds \le m (\Delta, Y_{n-1})^2 + \Delta, 
\end{align*}
where $m (\Delta, Y_{n-1}) = \exp (- 3 Y_{n-1}^2 \Delta) Y_{n-1} + 2 Y_{n-1}^3 \int_0^\Delta \exp (- 3 Y_{n-1}^2 s) ds$. We derive an  upper bound for $m (\Delta, Y_{n-1})^2$ and for a fixed $\varepsilon > 0$ consider the following three cases separately: (a) $Y_{n-1} = 0$; (b) $|Y_{n-1}| > \varepsilon$; (c) $|Y_{n-1}| \le \varepsilon, Y_{n-1} \neq 0$. For (a), we immediately see that $m (\Delta, Y_{n-1}) = 0$. We note that, when $Y_{n-1} \neq 0$: 
\begin{align*}
m (\Delta, Y_{n-1}) = \exp (- 3 Y_{n-1}^2 \Delta) Y_{n-1} + 2 Y_{n-1}^3  \times \tfrac{1 - \exp (- 3 Y_{n-1}^2 \Delta)}{3 Y_{n-1}^2} 
= \Bigl( \tfrac{2}{3}  + \tfrac{1}{3} \exp (- 3 Y_{n-1}^2 \Delta) \Bigr) Y_{n-1}. 
\end{align*}
Therefore, for case (b), we have that:  
\begin{align*} 
m (\Delta, Y_{n-1})^2 \le \Bigl( \tfrac{2}{3}  + \tfrac{1}{3} \exp (- 3 \varepsilon^2 \Delta) \Bigr)^2 \times Y_{n-1}^2  \equiv \rho \times  V (Y_{n-1}), 
\end{align*} 
with $\rho \in (0,1)$. In case (c), it also follows that:  
\begin{align*}
m (\Delta, Y_{n-1})^2 
& = \rho \times Y_{n-1}^2  
+ 2 \sqrt{\rho} \Bigl( \exp (- 3 Y_{n-1}^2 \Delta) - \exp (- 3 \varepsilon^2 \Delta )\Bigr) Y_{n-1}^2 \\
& \quad +  \Bigl( \exp (- 3 Y_{n-1}^2 \Delta) - \exp (- 3 \varepsilon^2 \Delta )\Bigr)^2  Y_{n-1}^2 \\ 
& \le \rho \times V (Y_{n-1}) + (2 \sqrt{\rho} + 1) \varepsilon^2.  
\end{align*} 
We thus conclude  that the discrete-time Lyapunov condition holds for (\ref{eq:ldl_toy}), i.e., there exists $\alpha \in (0,1)$ and $\beta \ge 0$ s.t.: 
\begin{align*} 
%\label{eq:DL}
\mathbb{E} [V (Y_{n}) | \mathcal{F}_{n-1} ] \le \alpha V (Y_{n-1}) + \beta, \qquad \forall n \in \mathbb{N}. 
\end{align*} 
In summary, due to minorisation and discrete-time Lyapunov conditions, the LDL scheme (\ref{eq:ldl_toy}) preserves (geometric) ergodicity for the SDE (\ref{eq:toy_ex}) with locally Lipschitz drift. This example is used as an indication that the LDL scheme can preserve ergodicity for general SDEs with non-globally Lipschitz drift. Detailed analysis is left as future work. 
\\ 
\item \emph{Design of local drift linearisation -- choice of matrix $A$}. In the development of the density expansion, we considered full-drift linearisation, i.e., first-order Taylor expansion of the drift for all coordinates to define the matrix $A$. However, as mentioned in Section \ref{sec:num_DE}, one can also consider a partially linearised drift approximation, e.g., with the matrix $A$ being upper-triangular, in order to reduce the computational cost of calculating $\exp (A)$, as long as the vector field defined in the baseline scheme satisfies H\"ormander's condition as in Lemma \ref{lemma:LDL_hor}. For preservation of hypoellipticity, at least linearisation of $V_{S, 0}$ (drift of the smooth component $X_S$) w.r.t. $X_R$ is required so that the noise in the rough component $X_R$ is lifted to $X_S$. The optimal way of linearisation (choice of $A$) would depend on the model at hand, but if a user considers a lower level of density expansion `$J$', say, $J = 2, 3, 4$, which is indeed sufficient to see improvements in estimation accuracy, then the partial drift linearisation will be a better option in terms of computational cost; recall e.g., DE-II in Figure \ref{fig:summary}-(b).
\\  
\item \emph{Computational cost w.r.t. the state dimension.} We stress that our CF-expansion converges exponentially fast with $J\ge 1$, so small values of $J$ will typically provide accurate proxies. Such a consideration counterbalances the computing cost for increasingx state dimension $N$. Following the analytical expressions of the $e_k$'e for the FHN model in Section~\ref{app:supp_num} of Supplementary Material, in the case of additive noise, one has $e_1=e_2=0$, while the calculation of $e_5$ requires all 3rd order derivatives of the baseline Gaussian density, at a cost of $\mathcal{O}(N^3)$. An extra derivative is added in the calculation when increasing $k$ in $e_k$ by one. Note that calculations involving just the baseline Gaussian transition density will typically already involve costs of $\mathcal{O}(N^3)$ due to matrix inversions, so in the additive noise setting using $J=5$ will not increase computing costs vs $J=0$ as an order of $N$.      
\end{enumerate} 
}

\section{Conclusion} 
\label{sec:conclusion}
We propose a new CF-expansion for the transition density of multivariate SDEs over a time interval with fixed length $\Delta \in (0,1)$, of the form `baseline Gaussian density' $\times$ `correction term', where the `correction term' involves %partial derivatives of SDE coefficients and the reference Gaussian density with 
quantities of size 
$\Delta^{j/2}$, $j = 0, \ldots, J$, for $J\ge 1$. Analytical 
expressions can be obtained via any software that carries out symbolic calculations. 
We have shown analytically that the error has a size of  $\mathcal{O} (\Delta^{(J+1)/2})$.
The proposed CF-expansion covers hypo-elliptic classes of SDEs, whereas most of the developments in earlier works are dedicated to elliptic SDEs. In the numerical studies in 
Section~\ref{sec:num} the errors from our CF-expansion are fastly eliminated as $J$ increases for a fixed $\Delta \in (0,1)$.

%This work contains many novelties and contributions to the field of statistical inference of diffusion processes. 
%First, the proposed CF-expansion covers hypo-elliptic classes of SDEs. 

%In terms of further connections and remarks
We also mention the following.
First, we take the direction described in the paper to produce our CF-expansion because potential alternative approaches used in the literature for the elliptic class (involving, e.g., Malliavin calculus) are arguably much more challenging in terms of producing a practical and theoretically validated methodolology.
Second, several recent works on the theme of parameter inference for hypo-elliptic SDEs produce methodology and analytical results in the high-frequency scenario $\Delta\rightarrow 0$, see 
e.g.~\cite{dit:19,mel:20,glot:21, pil:24-II,igu:23,igu:SPA,igu:BJ}. Then, numerical experiments are used to check the precision for a fixed $\Delta>0$ given in practice. In contrast, our contribution assumes a fixed $\Delta$, thus is expected to be more robust in deviations of $\Delta$ from $0$ than high-frequency approaches.

\subsection*{Acknowledgements}
{We are grateful to the Editor-in-Chief, Associate Editor, and two anonymous referees for their constructive suggestions that enhanced the quality of this paper. We also thank Dr Karen Habermann for suggesting related references on closed-form heat kernel expansions of hypo-elliptic diffusions.}

\subsection*{Funding}
YI acknowledges the supports from the Additional Funding Programme for Mathematical Sciences, delivered by EPSRC (EP/V521917/1) and the Heilbronn Institute for Mathematical Research. YI is also supported by the EPSRC grant Prob AI (EP/Y028783/1).

\appendix

\subsection*{Appendix}

\section{Proof of Lemma \ref{lemma:LDL_hor}}
\label{app:pf_LDL_hor} 

We write the LDL scheme (\ref{eq:LDL_frozen}) with a frozen variable $z \in \mathbb{R}^N$ as the following differential form:
\begin{align*}
d \bar{X}_t^{\theta, z} = {V}_{0}^z (\bar{X}_t^{\theta, z}, \theta) dt + \sum_{j = 1}^d V_{j}^z (\bar{X}_t^{\theta, z}, \theta) dB_{j, t}, \quad X_0 = x, 
\end{align*} 
for $(x, \theta) \in \mathbb{R}^N \times \Theta$, where  we have set: 
\begin{align}
V_0^z (x, \theta) \equiv A_{z, \theta} x + b_{z, \theta}, 
\qquad 
V_j^z (x, \theta) \equiv V_j (z, \theta), \quad  1 \le j \le d. 
\end{align} 
Since the diffusion coefficients are independent of the state $\bar{X}_t^{\theta, z}$, the above It\^o-type SDE identifies with the Stratonovich one. We show under Assumption \ref{ass:hor} that the vector fields determined from the coefficients of the above SDE satisfy H\"ormander's condition for each model class \ref{eq:ellip} and \ref{eq:hypo}.  
\\ 

\noindent  
\textit{Elliptic model \ref{eq:ellip}.} 
We immediately have from Assumption \ref{ass:hor}-I that 
\begin{align}
 \mathrm{Span} \bigl\{ V_j^z (x, \theta) |_{z = x},  \ \ 
 1 \le j \le d  \bigr\} = \mathbb{R}^N, 
\end{align} 
for all $(x, \theta) \in \mathbb{R}^N \times \Theta$.
\\ 

\noindent  
\textit{Hypo-elliptic model \ref{eq:hypo}.}  
We firstly note that 
\begin{align*}
\widetilde{V}_{S, 0} (x, \theta) \equiv {V}_{S, 0} (x, \theta), \quad 
V_{j}^z (x, \theta) \equiv  \bigl[ \mathbf{0}_{N_S}^\top, V_{R, j} (z, \theta)^\top \bigr]^\top, \qquad  1 \le j \le d.
\end{align*} 
Then Assumption \ref{ass:hor}-II is equivalent to the following condition:
\begin{align} 
\begin{aligned}
\mathrm{Span} \bigl\{ V_{R, j} (x, \theta), \, 1 \le j \le d \bigr\} & = \mathbb{R}^{N_R}; \\ 
\mathrm{Span} \bigl\{ \mathrm{Proj}_{1, N_S} \{[ \widetilde{V}_0, V_j] (x, \theta)  \}, \, 1 \le j \le d \bigr\} 
& = \mathrm{Span} \bigl\{ 
\partial_{x_R}^\top V_{S, 0} (x, \theta) V_{R, j} (x, \theta), \, 1 \le j \le d \bigr\} = \mathbb{R}^{N_S}.  \label{eq:span_I} 
\end{aligned}
\end{align} 
The condition (\ref{eq:span_I}) leads to: 
\begin{align*}
\mathrm{Span} \bigl\{ V_{R, j} (z, \theta) |_{z = x }, \, 1 \le j \le d  \bigr\} = \mathbb{R}^{N_R}, 
\end{align*} 
and
\begin{align*} 
& \mathrm{Span} \bigl\{ \mathrm{Proj}_{1, N_S} 
\bigl\{ 
[V_0^z, V_j^z] (x, \theta) |_{z = x} \bigr\}, \, 1 \le j \le d  \bigr\}
= \mathrm{Span} \bigl\{ \mathrm{Proj}_{1, N_S} 
\{ A_{z, \theta} V_j^z (x, \theta) |_{z = x} \}, \, 1 \le j \le d  \bigr\} \\
& \qquad 
= \mathrm{Span} \bigl\{ 
 \partial_{x_R}^\top V_{S, 0} (z, \theta) V_{R,j} (z, \theta) |_{z = x}, \, 1 \le j \le d  \bigr\}  = \mathbb{R}^{N_S}. 
\end{align*} 
Thus, we obtain: 
\begin{align} 
\mathrm{Span} \Bigl\{ \bigl\{ V_j^z (x, \theta), [V_0^z, V_j^z] (x, \theta) \bigr\} |_{z = x}, \, 1 \le j \le d \Bigr\} 
= \mathbb{R}^{N}, 
\end{align}
for all $(x, \theta) \in \mathbb{R}^N \times \Theta$.  The proof is now complete.  

\section{Proof of Lemma \ref{lemma:aux}} \label{app:pf_aux_1}

We define 
\begin{align*}
G(s) \equiv \int_{\mathbb{R}^N} p^{\bar{X}^z}_{\Delta -s} (\xi, y; \theta) p^X_s (x, \xi; \theta) |_{z = x} \, 
d \xi, \qquad s \geq 0. 
\end{align*} 
Noticing that $p^X_0 (x, y; \theta) = \delta_y (x)$ and $p^{\bar{X}}_0 (x, y; \theta) = \delta_y (x)$, we have 
\begin{align*}
p^X_\Delta (x, y; \theta) - p^{\bar{X}^z}_\Delta (x, y; \theta) |_{z=x}
= G(\Delta) - G(0) = \int_0^\Delta G'(s) ds. 
\end{align*} 
Also, note that the transition densities $p^X_\Delta (x, y ; \theta)$ and $p^{\bar{X}^z}_\Delta (x, y ; \theta)$ satisfy the following backward/forward Kolmogorov equations: 
\begin{gather} 
 \partial_t \, p^X_t (x, y; \theta) 
= \mathscr{L}_\theta \{ p^X_t ( \cdot, y; \theta) \} (x), \qquad 
\partial_t \, p^X_t (x, y; \theta) 
= \mathscr{L}_\theta^* \{ p^X_t ( x, \cdot ; \theta) \} (y);  
\\  
 \partial_t \, p^{\bar{X}^z}_t ( x, y; \theta)
= \mathscr{L}_\theta^{0, z} 
\{ p^{\bar{X}^z}_t (\cdot, y; \theta) \} (x).  \label{eq:FPEs}
%\qquad 
%\partial_t \, p^{\bar{X}^z}_t ( x, y; \theta) 
%= \mathscr{L}_\theta^{0, z, *} 
%\{ p^{\bar{X}^z}_t (x, \cdot; \theta) \} (y).  
\end{gather}  
% 
%where $\mathscr{L}$ is the generator given by (\ref{eq:generator}), and $\mathscr{L}_\theta^\ast$, $\mathscr{L}_\theta^{0, z, \ast}$ are $L^2 (\mathbb{R}^N)$-adjoint operators of $\mathscr{L}_\theta$, $\mathscr{L}_\theta^{0, z}$, respectively. 
% 
% where we note that the operator $\mathscr{L}_\theta$ and $\mathscr{L}_{\theta}^{0, z}$ act on the current state $x \in \mathbb{R}^N$ of the transition.
It then follows that:   
\begin{align*}
G'(s)  
& = - \int_{\mathbb{R}^N} \mathscr{L}_\theta^{0, z} \{ p^{\bar{X}^z}_{\Delta -s} (\cdot, y; \theta)  \} (\xi) 
p^X_s (x, \xi; \theta) |_{z=x} d \xi  \\ 
& \qquad \qquad \qquad  +  \int_{\mathbb{R}^N} 
p^{\bar{X}^z}_{\Delta - s} (\xi, y; \theta) \mathscr{L}_\theta^\ast \{ p^X_s (x, \cdot; \theta) \} (\xi) |_{z=x} d \xi  \nonumber \\ 
& = - \int_{\mathbb{R}^N} \mathscr{L}_\theta^{0, z} \{ p^{\bar{X}^z}_{\Delta - s} (\cdot, y; \theta)  \} (\xi) 
p^X_s (x, \xi; \theta) |_{z=x} d \xi  \\  
& \qquad  \qquad \qquad +  \int_{\mathbb{R}^N}
\mathscr{L}_\theta \{ p^{\bar{X}^z}_{\Delta - s} 
(\cdot, y; \theta) \} (\xi)  p^X_s (x, \xi; \theta) |_{z=x} d \xi.
\end{align*}
% 
%where we used integration by parts on the second equation. 
The proof is now complete. 

\section{Proof of Lemma \ref{lemma:aux2}} \label{app:pf_aux_2} 

We focus on showing the formula (\ref{eq:expansion_step2_1}), and
 (\ref{eq:expansion_step2_2}) is obtained from a similar argument. We make use of the approach used in \citet[Proposition 2.1]{igu:21-2}. We define 
 \begin{align} \label{eq:f_s}
    g^{\widetilde{\mathscr{L}}^{z}_\theta } (s) 
    & \equiv \bar{P}_s^{\theta, z} 
    \widetilde{\mathscr{L}}^{z}_\theta  \bar{P}_{t-s}^{\theta, z}  \varphi (x) 
    = \int_{\mathbb{R}^N} \widetilde{\mathscr{L}}^{z}_\theta   
    \bar{P}_{t-s}^{\theta, z} \varphi ( {\xi_1})
    \, p^{\bar{X}^z}_s (x, \xi_1; \theta) 
    \, d \xi_1 
 \end{align}
 and consider the Taylor expansion of $g^{\widetilde{\mathscr{L}}^{z}_\theta } (s)$ at $s=0$: 
 \begin{align*}
    g^{\widetilde{\mathscr{L}}^{z}_\theta } (s) 
    &  = \sum_{k = 0}^{J} \partial^k 
    g^{\widetilde{\mathscr{L}}^{z}_\theta } (s)|_{s = 0} \times \frac{s^k}{k!} 
    + \mathscr{R}^{J+1, z} (s, t, x; \theta),  \qquad J \in \mathbb{N},
 \end{align*} 
where 
$\textstyle 
\mathscr{R}^{J+1, z} (s, t, x; \theta) 
\equiv {s^{J +1}} \int_0^1 
\partial^{J+1} g^{\widetilde{\mathscr{L}}^{z}_\theta } (su) 
\tfrac{(1-u)^{J}}{J!} du
$. 
We have from (\ref{eq:f_s}) that 
\begin{align*}
& \partial g^{\widetilde{\mathscr{L}}^{z}_\theta } (s) 
= \int_{\mathbb{R}^N} \widetilde{\mathscr{L}}^{z}_\theta
\bar{P}_{t-s}^{\theta, z} \varphi (\xi_1) \,  \partial_s p^{\bar{X}^z}_s (x, \xi_1; \theta) \, d \xi_1 
+ \int_{\mathbb{R}^N} 
\widetilde{\mathscr{L}}^{z}_\theta  
\partial_s 
\bar{P}_{t-s}^{\theta, z} \varphi (\xi_1) \, 
p^{\bar{X}^z}_s (x, \xi_1; \theta) \,  d \xi_1 \\
& = \int_{\mathbb{R}^N} 
\widetilde{\mathscr{L}}^{z}_\theta  
\bar{P}_{t-s}^{\theta, z} \varphi (\xi_1) \,
 (\mathscr{L}_\theta^{0,z})^\ast 
\{ p^{\bar{X}^z}_s (x, \cdot; \theta)  \} (\xi_1) \, d \xi_1  
- \int_{\mathbb{R}^N} \widetilde{\mathscr{L}}^{z}_\theta \mathscr{L}_\theta^{0,z} 
\bar{P}_{t-s}^{\theta, z} \varphi (\xi_1) \, 
p^{\bar{X}^z}_s (x, \xi_1; \theta) \,  d \xi_1  \\ 
& 
= \int_{\mathbb{R}^N} 
[ \mathscr{L}_\theta^{0, z}, \widetilde{\mathscr{L}}^{z}_\theta]  
\bar{P}_{t-s}^{\theta, z} \varphi (\xi_1)  \, 
p^{\bar{X}^z}_s (x, \xi_1; \theta) \,  d \xi_1  \\
& = g^{[ \mathscr{L}_\theta^{0, z}, \widetilde{\mathscr{L}}^{z}_\theta]} (s), 
\end{align*} 
where we made use of (\ref{eq:FPEs}) and integration by parts in the second and third lines, respectively. Thus, the higher-order derivatives of $g^{\widetilde{\mathscr{L}}^{z}_\theta } (s)$ are given as: 
\begin{align*}
    \partial^k g^{\widetilde{\mathscr{L}}^{z}_\theta} (s) 
= g^{(\mathrm{ad}_{\mathscr{L}_\theta^{0, z}})^k (\widetilde{\mathscr{L}}_\theta^{z})} (s) 
= \bar{P}_s^{\theta, z} 
\bigl\{ 
    \bigl(
        \mathrm{ad}_{\mathscr{L}_\theta^{0, z}} 
    \bigr)^k
    (\widetilde{\mathscr{L}}_\theta^{z}) 
 \bigr\}
 \bar{P}_{t-s}^{\theta, z} \varphi (x),  \qquad 0 \le k \le J, 
\end{align*}
and then, 
\begin{align*} 
    \partial^k g^{\widetilde{\mathscr{L}}^{z}_\theta} (s) |_{s = 0} = \bigl\{ 
        \bigl(
            \mathrm{ad}_{\mathscr{L}_\theta^{0, z}} 
        \bigr)^k
        (\widetilde{\mathscr{L}}_\theta^{z}) 
     \bigr\}
     \bar{P}_{t}^{\theta, z} \varphi (x). 
\end{align*}
The proof is now complete.

\subsection*{Supplementary Material}     
This material contains supporting information for the manuscript ``{A closed-form transition density expansion for elliptic and hypo-elliptic SDEs}''. The material is organised as follows: 
In Section \ref{app:pfs_main}, we prove Lemma \ref{lemma:deriv_LDL} and the central analytic results provided in the main text, i.e., Theorems \ref{thm:bd_r1} and \ref{thm:bd_r2}, error estimates for the proposed CF expansion. 
Section \ref{app:DO} studies the expression of the differential operator $\mathscr{D}_\alpha^{z, \theta}$ given in (\ref{eq:D}). Section \ref{app:supp_num} provides supporting information for the implementation of density expansion used in the numerical experiments and additional experiments of Bayesian inference for FHN model.   

%
% Substituting $z=x$ in (), we obtain the desired statement of Lemma \ref{lemma:deriv_LDL}. 
\section{Proof of Lemma \ref{lemma:deriv_LDL} and Theorems \ref{thm:bd_r1}-\ref{thm:bd_r2}}
\label{app:pfs_main}
This section is organised as follows: In Section \ref{app:pre}, we collect notations and basic results that are frequently used throughout Section \ref{app:pfs_main}. In Section \ref{sec:pf_ldl_deriv}, we show Lemma \ref{lemma:deriv_LDL}. Section \ref{app:aux} provides five auxiliary results (Propositions \ref{prop:density_bd}, \ref{prop:bd_E_1}, \ref{prop:bd_E_2} and Lemmas \ref{lemma:conv}, \ref{lemma:poly_bd}) to prove the main analytic results (Theorems \ref{thm:bd_r1} and \ref{thm:bd_r2}). Propositions \ref{prop:bd_E_1} and \ref{prop:bd_E_2} in particular play a central role to show Theorems \ref{thm:bd_r1} and \ref{thm:bd_r2}, respectively. Proposition \ref{prop:density_bd} and Lemmas \ref{lemma:conv}--\ref{lemma:poly_bd} are required to show Propositions \ref{prop:bd_E_1}--\ref{prop:bd_E_2}. Given Propositions \ref{prop:bd_E_1} and \ref{prop:bd_E_2}, we prove Theorems \ref{thm:bd_r1} and \ref{thm:bd_r2} in Section \ref{app:pf_bd_r1} and \ref{app:pf_bd_r2}, respectively. 
\subsection{Preliminaries} 
\label{app:pre}
We introduce and recall notations frequently used throughout Section \ref{app:pfs_main}. 
\begin{itemize}
	\item For $x, y, z \in \mathbb{R}^N, \, \theta \in \Theta$, $0 \le s \le t < \infty$ and $w \in 
	\{\ref{eq:ellip}, \ref{eq:hypo} \}$, 
	\begin{gather}
		% \Gamma_{t, (w)} 
		% & \equiv  
		% \begin{cases}
			% \mathrm{diag} \bigl( 
			% \underbrace{\sqrt{t}, \ldots, \sqrt{t}}_{N} \bigr), 
			% & w  = \ref{eq:ellip}; \\ 
			% \mathrm{diag} 
			% \bigl(
			% \underbrace{\sqrt{t^3}, \ldots, \sqrt{t^3}}_{N_S}, 
			% \underbrace{\sqrt{t}, \ldots, \sqrt{t}}_{N_R}
			% \bigr), & w  = \ref{eq:hypo}, 
			% \end{cases}
		% \qquad t \ge 0;  \\ 
		\Gamma_{t, (\ref{eq:ellip})} 
		\equiv \mathrm{diag} \bigl( 
		\underbrace{\sqrt{t}, \ldots, \sqrt{t}}_{N} \bigr), \quad 
		\Gamma_{t, (\ref{eq:hypo})} \equiv \mathrm{diag} 
		\bigl(
		\underbrace{\sqrt{t^3}, \ldots, \sqrt{t^3}}_{N_S}, 
		\underbrace{\sqrt{t}, \ldots, \sqrt{t}}_{N_R}
		\bigr);  \label{eq:Gamma} \\ 
		\nu^{(w), z} (t, y, x, \theta) 
		\equiv  
		\Gamma_{t, (w)}^{-1}
		\Bigl(
		e^{- t A_{z, \theta}} y 
		- x - 
		\int_0^t e^{-s A_{z, \theta}} b_{z, \theta} ds  
		\Bigr);   \label{eq:nu}   \\ 
		v (t, s, y, x, \theta)  
		\equiv x + V_0 (x, \theta) s - e^{-(t-s) A_{x,\theta}} y + \int_0^{t-s} e^{-u A_{x, \theta}} b_{x, \theta} du. 
		\label{eq:v_func} 
	\end{gather}  
	\item Recall: 
	$$
	m(\ref{eq:ellip}) \equiv N, \quad m (\ref{eq:hypo}) \equiv  4 d  
	$$
	and for $\alpha \in \mathbb{Z}_{\ge 0}^N$, 
	\begin{align}
		\| \alpha \|_{\ref{eq:ellip}} \equiv \tfrac{1}{2} | \alpha |, \qquad 
		\| \alpha \|_{\ref{eq:hypo}} \equiv \tfrac{3}{2} | \alpha_S|  + \tfrac{1}{2} | \alpha_R|  \label{eq:mi_norm_supp}
	\end{align} 
	with convention $\alpha_S = (\alpha_{1}, \ldots, \alpha_{N_S})$ and $\alpha_R = (\alpha_{N_S + 1}, \ldots, \alpha_{N})$  for hypo-elliptic model $\eqref{eq:hypo}$. 
	% 
%	\item For $\alpha \in \mathbb{Z}_{\ge 0}^j, \, j \in \mathbb{N}$, 
%	\begin{align}
%		\| \alpha \|_{0} \equiv  \sum_{1 \le i \le j} \mathbf{1}_{\alpha_i = 0}, \qquad 
%		\| \alpha \|_{0, 1, 2} \equiv  \sum_{1 \le i \le j} \bigl\{ 2 \times \mathbf{1}_{\alpha_i = 0} +  \mathbf{1}_{\alpha_i \in \{ 1, 2\} } \bigr\}. 
%	\end{align}
	\item Let  $\mathscr{F}^{(w)}, \, w \in \{\ref{eq:ellip}, \ref{eq:hypo} \}$ be a space of functions $\psi^{(w)} : (0, 1) \times \mathbb{R}^N \to \mathbb{R}, \, w \in \{\ref{eq:ellip}, \ref{eq:hypo}\}$ satisfying the following property: for any $\alpha \in \mathbb{Z}_{\ge 0}^N$, there exists a constant $C > 0$ such that for all $(t, \xi) \in (0,1) \times \mathbb{R}^N$, 
	\begin{align}
		\bigl| \partial_\xi^\alpha \psi^{(w)} (t, \xi) \bigr| 
		\le C t^{-\| \alpha \|_w} \times \bigl| \widetilde{\psi}^{(w)} (t, \xi) \bigr|,  
	\end{align}
	for some function $\widetilde{\psi}^{(w)}: (0,1) \times \mathbb{R}^N \to \mathbb{R}$. 
	\item For $\kappa > 0$ defined in Assumption \ref{ass:coeff}, 
	\begin{align}
		\mathscr{Z}_\kappa 
		\equiv  
		\bigl\{\xi \in \mathbb{R}^N \, \bigl| \, 
		|\xi| + \bigl| V_0 (\xi, \theta) \bigr| \le \kappa,  \quad  \forall \ \theta \in \Theta
		\bigr\}.  
		\label{eq:Z}
	\end{align}
	In particular, under Assumption \ref{ass:coeff}, the initial state $x$ is an element of $\mathscr{Z}_\kappa$. 
	\item Let $\mathscr{G}^{(w)}: (0, \infty) \times \mathbb{R}^N \times \mathbb{R}^N \times \Theta \to \mathbb{R}$ and $\overline{\mathscr{G}}^{(w), z} : (0,\infty) \times \mathbb{R}^N \times \mathbb{R}^N \times \Theta \to \mathbb{R}, \, w \in \{\ref{eq:ellip}, \ref{eq:hypo} \} $, $z \in \mathscr{Z}_\kappa$ be generic notations to represent functions having the following property: there exist constants $C, \lambda > 0$ such that for all $(t, x, y, z, \theta) \in (0,\infty) \times \mathbb{R}^N \times \mathbb{R}^N \times \mathscr{Z}_\kappa \times \Theta$,
	\begin{align}
		\bigl|  
		{\mathscr{G}}^{(w)} (t, x, y, \theta) 
		\bigr|   
		& \le  C t^{- \tfrac{m(w)}{2}} \exp \Bigl( - \lambda \bigl| \Gamma_{t, (w)}^{-1} 
		\bigl(   y - x - V_0 (x, \theta) t  \bigr) \bigr| ^2 \Bigr);  
		\nonumber \\  
		\bigl| 
		\overline{\mathscr{G}}^{(w), z} (t, x, y, \theta) 
		\bigr| 
		& \le C t^{-  \tfrac{ m(w)}{2}}  
		\exp \Bigl( - \lambda \bigl| \nu^{(w), z} (t, y, x, \theta) \bigr|^2 \Bigr),  
	\end{align}
	where $\Gamma_{t, (w)}$ and $\nu^{(w), z}$ are defined in (\ref{eq:Gamma}) and (\ref{eq:nu}), respectively.  
    \item For a differential operator $D$ acting on the functions over $\mathbb{R}^N$, we write $D^\ast$ as its adjoint defined via the formula: $\textstyle \int_{\mathbb{R}^N} D f (x) g (x) dx  = \int_{\mathbb{R}^N} f (x) D^\ast g (x) dx$, for sufficiently smooth functions $f,g \in  L^2 (\mathbb{R}^N)$ vanishing at infinity. 
\end{itemize}

The following basic results/estimates are frequently used throughout Section \ref{app:pfs_main}: 
\begin{enumerate}
\item[I.] For $\alpha \in \mathbb{Z}_{\ge 0}^N$, it holds that 
\begin{align}
\partial^\alpha_y p_t^{\bar{X}^z}  (x, y; \theta) 
= \partial^\alpha_x p_t^{\bar{X}^z}  (x, y; \theta) G_\alpha^z (t, \theta), 
\qquad (t, x, y, z, \theta) \in (0, \infty) \times \mathbb{R}^N \times \mathbb{R}^N \times \mathscr{Z}_\kappa \times \Theta,   
\label{eq:basic} 
\end{align}
where $G_\alpha^z: (0, \infty) \times \Theta \to \mathbb{R}$ satisfies the following property under Assumptions \ref{ass:param}--\ref{ass:coeff}:  
\begin{align*}
\sup_{(t, z, \theta) \in (0, 1) \times \mathscr{Z}_\kappa  \times \Theta } 
\bigl| G_\alpha^z (t, \theta) \bigr| < \infty. 
\end{align*}
\item[II.] For two multidimensional Gaussian densities $\mathbb{R}^N \ni \xi \mapsto \varphi (\xi; a_i, A_i), \, i = 1, 2$ with mean $a_i \in \mathbb{R}^N$ and covariance $A_i \in \mathbb{R}^{N \times N}$, it follows (see e.g. \cite{vin:04}) that: 
\begin{align} \label{eq:standard}
\begin{aligned}
\int_{\mathbb{R}^N} \varphi_1 (\xi; a_1, A_1) \varphi_2 (\xi; a_2, A_2) d \xi  
= \tfrac{1}{\sqrt{(2\pi)^N \det \bigl( A_1 + A_2 \bigr)}} \exp \Bigl( - \tfrac{1}{2} (a_1 - a_2)^\top 
(A_1 + A_2)^{-1} (a_1 - a_2) \Bigr).  
\end{aligned}
\end{align} 
\item[III.] Let $\lambda_1 > 0$. For any $\alpha \in \mathbb{Z}_{\ge 0}^N$, there exist constants $c, \lambda_2 > 0$ such that for all $\xi \in \mathbb{R}^N$,  
\begin{align}
|\xi^\alpha| \times 
\exp \bigl( - \lambda_1 |\xi|^2 \bigr)
\le C  \exp \bigl( - \lambda_2 |\xi|^2 \bigr).  
\label{eq:exp_bd}
\end{align}
We note that this type of inequality is also used to prove the bound of the remainder term of density expansion for elliptic diffusions in \cite{yang:19}.  
\end{enumerate} 

\subsection{Proof of Lemma \ref{lemma:deriv_LDL}} \label{sec:pf_ldl_deriv}
We will show the following statement: it holds under Assumptions \ref{ass:diff}--\ref{ass:dim} that for the initial state $x \in \mathbb{R}^N$,  $(\Delta, y, z, \theta) \in (0, \infty) \times \mathbb{R}^N \times \mathscr{Z}_\kappa \times \Theta$,  
\begin{align} \label{eq:derive_LDL_z}
\partial^{\alpha}_x \, p_\Delta^{\bar{X}^z, (w)} (x, y; \theta) 
& = \Delta^{- \| \alpha \|_w} \times 
\mathscr{H}^{(w), z}_{\alpha} (\Delta, x ,y, \theta) \, p_\Delta^{\bar{X}^z, (w)} (x, y; \theta), \quad \alpha \in \mathbb{Z}_{\ge 0}^N, \quad w \in \{ \ref{eq:ellip}, \ref{eq:hypo} \}, 
\end{align} 
for some $\mathscr{H}_\alpha^{(w), z}: (0, \infty) \times \mathbb{R}^N \times \mathbb{R}^N \times \Theta \to \mathbb{R}$ satisfying that there exists a constant $C > 0$ such that for all 
$(\Delta , x, y, z , \theta) \in (0,1) \times  \mathbb{R}^N \times  \mathbb{R}^N  \times  \mathscr{Z}_\kappa \times \Theta$, 
\begin{gather*}
\bigl| \mathscr{H}^{(w), z}_{\alpha} (\Delta, x ,y, \theta) \, p_\Delta^{\bar{X}^z, (w)} (x, y; \theta) 
\bigr|   \le 
C \Delta^{- \tfrac{m(w)}{2}} \exp \Bigl( - \lambda | \nu^{(m), z} (\Delta, y, x, \theta) |^2 \Bigr) .  
\end{gather*} 
Note that the statement of Lemma \ref{lemma:deriv_LDL} is immediately obtained by setting $z=x$ in (\ref{eq:derive_LDL_z}). 

We first consider the case $| \alpha |  = 1$, i.e. $\partial_x^\alpha \equiv \partial_{x_j}, \, 1 \le j \le N$. We have  that
\begin{align}
p_\Delta^{\bar{X}^z, (w)} (x,y; \theta) 
= \tfrac{1}{\sqrt{(2\pi)^N \det {\Sigma} (\Delta, z, y)}} \exp 
\Bigl(- \tfrac{1}{2} 
\nu^{(w), z} (\Delta, y, x, \theta)^\top   
\bigl( \mathscr{M}^{(w)} (\Delta, z, \theta) 
\bigr)^{-1}   
\nu^{(w), z} (\Delta, y, x, \theta)
\Bigr)    \label{eq:ds_LDL_z_2} 
\end{align} 
with 
$
\mathscr{M}^{(w)} (\Delta, z, \theta) 
\equiv 
\Gamma_{\Delta, (w)}^{-1}
\hat{\Sigma} (\Delta, z, \theta) 
\Gamma_{\Delta, (w)}^{-1} 
$%,
, where the matrices ${\Sigma}, \hat{\Sigma}$ are defined in (\ref{eq:Sigma}) in the main text and $\nu^{(w), z}$ in (\ref{eq:nu}). Then, 
\begin{align} \label{eq:deriv_first_z}
\partial_{x_j} p_\Delta^{\bar{X}^z, (w)} (x,y; \theta)  
= \Delta^{-\|(j)\|_w} \times \mathscr{H}_{(j)}^{(w), z} (\Delta, x, y, \theta) \times p_\Delta^{\bar{X}^z, (w)} (x,y; \theta)    
\end{align} 
with 
\begin{align} \label{eq:H_j}
\mathscr{H}_{(j)}^{(w), z} (\Delta, x, y, \theta) 
= \sum_{1 \le k \le N} 
\bigl[ \bigl( \mathscr{M}^{(w)} (\Delta, z, \theta) 
\bigr)^{-1} \bigr]_{jk} \bigl[ \nu^{(w), z} (\Delta, y, x, \theta) \bigr]_k.  
\end{align}
We will
% study the bound for $\bigl( \mathscr{M}^{(w)} (\Delta, z, \theta)  \bigr)^{-1}$ and in particular, 
show under Assumptions \ref{ass:hor} and \ref{ass:coeff} that there exists a constant $C > 0$ such that for all $(t, z, \theta) \in (0,1)\times \mathscr{Z}_\kappa \times \Theta$, 
\begin{align} \label{eq:inv_M_bd}
\Bigl| 
\bigl[
\bigl( \mathscr{M}^{(w)} (t, z, \theta) \bigr)^{-1}  
\bigr]_{i_1 i_2} 
\Bigr| 
\le C, 
\qquad 1 \le i_1, i_2 \le N, 
\quad w \in \{\ref{eq:ellip}, \ref{eq:hypo} \}.
\end{align}
To check (\ref{eq:inv_M_bd}), we study the bound for the elements of the matrix $\hat{\Sigma}$. Due to the expansion $e^{- s A_{z, \theta}} = I_N - s A_{z, \theta} + \mathscr{O} (s^2)$, we have under Assumption \ref{ass:coeff} that  
\begin{align*}
& \hat{\Sigma} (t, z, \theta) 
\equiv \int_0^t e^{- s A_{z, \theta}} \sigma (z, \theta) \sigma (z, \theta)^\top  e^{- s A_{z, \theta}^\top} ds  \\ 
& = t a (z, \theta)^\top 
- \tfrac{t^2}{2} \bigl( A_{z, \theta} a (z, \theta) + a (z, \theta) A_{z, \theta}^\top \bigr)
+ \tfrac{t^3}{3} A_{z, \theta} a (z, \theta) A_{z, \theta}^\top + \widetilde{\Sigma} (t, z, \theta), 
\end{align*}
where the remainder term $\widetilde{\Sigma} (t, z, \theta)$ satisfies: $\sup_{(z, \theta) \in \mathscr{Z}_\kappa \times \Theta} \bigl| \widetilde{\Sigma} (t, z, \theta)  \bigr| \le C t^4$ for some constant $C > 0$ independent of $t \in (0,1)$. Note that under model class (\ref{eq:hypo}), the matrix $a (z, \theta)$ is degenerate, and also the matrix 
$$t a (z, \theta)^\top 
- \tfrac{t^2}{2} \bigl( A_{z, \theta} a (z, \theta) + a (z, \theta) A_{z, \theta}^\top \bigr) + \tfrac{t^3}{3} A_{z, \theta} a (z, \theta) A_{z, \theta}^\top
$$ 
is positive definite under Assumption \ref{ass:hor}. Thus, there exists a constant $C > 0$ such that for all $(t, z, \theta) \in (0, 1) \times \mathscr{Z}_\kappa \times \Theta$, 
\begin{itemize}
\item under model class (\ref{eq:ellip}), 
\begin{align} \label{eq:v_scale_e}
[\hat{\Sigma} (t, z, \theta)]_{i_1 i_2} \le C t, \qquad  1 \le i_1, i_2 \le N,  
\end{align}
\item under model class (\ref{eq:hypo}), 
\begin{align}  \label{eq:v_scale_h} 
[\hat{\Sigma} (t, z, \theta)]_{i_1 i_2}  
\le 
\begin{cases}
C t^3, & 1 \le i_1, i_2 \le N_S; \\ 
C t, & N_S + 1 \le i_1, i_2 \le N; \\ 
C t^2, & (\mathrm{otherwise}). 
\end{cases}
\end{align}
\end{itemize}
Then, the inequality (\ref{eq:inv_M_bd}) is deduced from (\ref{eq:v_scale_e}) and (\ref{eq:v_scale_h}). We also have from the positive definiteness of $\hat{\Sigma}$,  (\ref{eq:v_scale_e}) and (\ref{eq:v_scale_h}) that: there exists a constant $c > 0$ such that for all $(t, z, \theta) \in (0,1) \times \mathscr{Z}_\kappa \times \Theta$, 
\begin{align} \label{eq:det_LDL}
%\det \Sigma^{(w)} (t, x, \theta) 
%\ge c \, t^{m(w)}, \qquad w \in \{\ref{eq:ellip}, \ref{eq:hypo} \}. 
\det \Sigma (t, z, \theta) \ge c \times 
\begin{cases}
	t^{N}  & \mathrm{under \, model \, class \, \eqref{eq:ellip}}; \\ 
	t^{4d}  & \mathrm{under \, model \,  class \, \eqref{eq:hypo}}.  
\end{cases}
\end{align} 
Thus, it follows from (\ref{eq:H_j}), (\ref{eq:inv_M_bd}) and (\ref{eq:det_LDL}) that: 
\begin{align*}
 \Bigl| \mathscr{H}_{(j)}^{(w), z} (\Delta, x , y, \theta) \times p_\Delta^{\bar{X}^z, (w)} (x, y; \theta) \Bigr| 
 & \le
 C_1  \bigl| \nu^{(w), z} (\Delta, y, x, \theta) \bigr| 
 \times \Delta^{- \tfrac{m(w)}{2}} \exp \Bigl( - \lambda_1 \bigl| \nu^{(w), z} (\Delta, y, x, \theta)   \bigr|^2 \Bigr)  \\ 
 & \le C_2 \Delta^{-\tfrac{m(w)}{2}} \exp \Bigl( - \lambda_2 \bigl| \nu^{(w), z} (\Delta, y, x, \theta)   \bigr|^2 \Bigr)  \quad \bigl( \because \eqref{eq:exp_bd} \bigr)
\end{align*}
for some constants $C_1, C_2, \lambda_1, \lambda_2  > 0$ that are independent of $\Delta, x, y, z$ and $\theta$. Hence, the formula (\ref{eq:derive_LDL_z}) holds with $|\alpha| = 1$. (\ref{eq:derive_LDL_z}) with $| \alpha | \ge 1$ can be shown via iterative application of (\ref{eq:deriv_first_z}) with a similar argument presented above, and we conclude. \qed  
\subsection{Auxiliary reesults for Theorems \ref{thm:bd_r1} and \ref{thm:bd_r2}} \label{app:aux}
\begin{proposition}
Let Assumptions \ref{ass:diff}--\ref{ass:dim} hold. For any $\alpha \in \mathbb{Z}_{\ge 0}^N$, there exist constants $C, \lambda > 0$ such that for all $(\Delta, y, \theta) \in (0,1) \times \mathbb{R}^N \times \Theta$, 
\begin{align}
\bigl| \partial^\alpha_y p_\Delta^{X, (w)} (x, y; \theta) \bigr| 
& \le  
C \Delta^{-\|\alpha\|_w - \tfrac{m(w)}{2}}
\times \exp \Bigl( - \lambda 
\bigl|
\Gamma_{\Delta, (w)}^{-1}
\bigl( 
y - x - V_0 (x, \theta) \Delta
\bigr) 
\bigr|^2 
\Bigr), \qquad w \in \{\ref{eq:ellip}, \ref{eq:hypo} \}, 
\label{eq:true_density_bd} 
\end{align}
where we recall $\| \cdot \|_w$ is given in (\ref{eq:mi_norm_supp}). 
% 
% \begin{align}
% \bigl| \partial^\alpha_y p_\Delta^{X, (w)} (x, y; \theta) \bigr| 
% & \le 
% \begin{cases}
% C \Delta^{- \| \alpha \|_{\ref{eq:ellip}} -\tfrac{N}{2}} \exp \Bigl( - \lambda \cdot \tfrac{|y - x |^2}{\Delta} \Bigr),  & w = \ref{eq:ellip};
% \\[0.2cm] 
% C \Delta^{- \| \alpha \|_{\ref{eq:hypo}} -\tfrac{4d}{2}}  \exp \Bigl( - 
% \lambda  
% \Bigl( 
% \tfrac{|y_S - x_S - V_{S,0} (x,\theta) \Delta |^2}{\Delta^3} 
% + 
% \tfrac{|y_R - x_R|^2}{\Delta} 
% \Bigr) 
% \Bigr),   & w = \ref{eq:hypo}. 
% \end{cases}   
% \label{eq:true_density_bd}
% \end{align} 
% 
% 
Also, for any $\alpha \in \mathbb{Z}_{\ge 0}^N$, there exists a constant $C, \lambda > 0$ such that for all $(\Delta, x, y, z, \theta) \in (0, 1) \times \mathbb{R}^N \times \mathbb{R}^N \times \mathscr{Z}_\kappa \times \Theta$, 
\begin{align} 
\Bigl| 
\partial^\alpha_y p_\Delta^{\bar{X}^z, (w)} (x, y; \theta) 
\Bigr| 
& \le C \Delta^{-\| \alpha \|_w - \tfrac{m (w)}{2}} 
\times 
\exp
\Bigl( 
- \lambda 
\bigl| \nu^{(w), z} (\Delta, y, x, \theta)  
\bigr|^2
\Bigr),   \qquad w \in \{ \ref{eq:ellip}, \ref{eq:hypo} \}.   
\label{eq:LDL_density_bd} 
\end{align}
% 
% where we have set:
%
% In particular, 
% % 
% \begin{align} 
% \Bigl| 
% \partial^\alpha_y p_\Delta^{\bar{X}^z, (w)} (x, y; \theta) 
% \bigr|_{z=x}
% \Bigr| 
% & \le C_2 \Delta^{-\| \alpha \|_w} 
% \times 
% \bigl|
% \mathscr{G}^{(w)} (\Delta, x, y; \theta)
% \bigr|.   
% \end{align} 
% We recall that the definition of $\|\alpha\|_w$ is given in (\ref{eq:mi_norm}). 
\label{prop:density_bd}
\end{proposition} 

\noindent
\textit{Proof of Proposition \ref{prop:density_bd}.}
We first show (\ref{eq:true_density_bd}). 
We have from \citep[Theorem 2.1.]{piga:22} that:  under Assumptions \ref{ass:diff}--\ref{ass:dim}, there exist constants $C, \lambda >0$ such that for all $(\Delta, y, \theta) \in (0, 1) \times \mathbb{R}^N \times \Theta$, 
\begin{align}
\bigl| 
\partial^\alpha_y p_\Delta^{X, (w)} (x, y; \theta) 
\bigr|  
\le 
\begin{cases}
C \Delta^{- \| \alpha \|_{\ref{eq:ellip}} -\tfrac{N}{2}} \exp \Bigl( - \tfrac{1}{C} \cdot \tfrac{|y - \xi_\Delta (x, \theta) |^2}{\Delta} \Bigr),  & w = \ref{eq:ellip};
\\[0.2cm] 
C \Delta^{- \| \alpha \|_{\ref{eq:hypo}} -\tfrac{4d}{2}}  \exp \Bigl( - 
\tfrac{1}{C} 
\Bigl( 
\tfrac{|y_S - \xi_{S, \Delta}(x, \theta)|^2}{\Delta^3} 
+ 
\tfrac{|y_R - \xi_{R, \Delta}(x, \theta)|^2}{\Delta} 
\Bigr) 
\Bigr),   & w = \ref{eq:hypo},  
\end{cases}
\label{eq:piga}
\end{align}
where $\xi_t(x, \theta) = 
\bigl(\xi_{S, t} (x, \theta), \xi_{R, t} (x, \theta) \bigr) \in \mathbb{R}^{N_S} \times \mathbb{R}^{N_R}, \, t \ge 0$ ($\xi_t (x, \theta) = \xi_{R, t} (x, \theta)$ when $w = \ref{eq:ellip}$) is defined via the following ODE: 
\begin{align} \label{eq:ode}
d \xi_t (x, \theta)  
= \widetilde{V}_0 
\bigl(\xi_t (x, \theta), \theta \bigr) dt, \qquad 
\xi_0 (x, \theta) = x,   
\end{align}
with $\widetilde{V}_0 : \mathbb{R}^N \times \Theta \to \mathbb{R}^N$ defined as:
\begin{align}
\widetilde{V}_0 (x, \theta) \equiv 
\Bigl[
V_{S, 0} (x, \theta) \Delta,
\, 
\Bigl( {V}_{R, 0} (x, \theta)  
- \tfrac{1}{2} \sum_{1 \le j \le d} 
\sum_{1 \le i \le N_R} V_{R, j}^i (x, \theta) \partial_{x_{R,i}} V_{R, j} (x, \theta)
\Bigr) \Delta 
\Bigr]^\top.  
\end{align}
Taylor expansion of the solution of ODE (\ref{eq:ode}) gives: for model class  (\ref{eq:ellip})
\begin{align} \label{eq:ode_E}
\xi_\Delta (x, \theta) 
= x + \widetilde{V}_0 (x, \theta) \Delta + r_{\ref{eq:ellip}} (\Delta, x, \theta), 
\end{align}
and for (\ref{eq:hypo}),  
\begin{align} \label{eq:ode_H}
\xi_\Delta (x, \theta) 
= \begin{bmatrix}
\xi_{S, \Delta} (x, \theta) \\[0.1cm] 
\xi_{R, \Delta} (x, \theta)  
\end{bmatrix}
= 
\begin{bmatrix}
x_S \\
x_R 
\end{bmatrix}
+ 
\widetilde{V}_0 (x, \theta) \Delta 
+ r_{\ref{eq:hypo}} (\Delta, x, \theta),   
\end{align} 
where the remainders $r_{\ref{eq:ellip}}$ and $r_{\ref{eq:hypo}}$ have the following properties under Assumptions \ref{ass:param}--\ref{ass:coeff}: 
\begin{align}
\sup_{(x, \theta) \in \mathbb{R}^N \times \Theta}
\bigl| r_w (\Delta, x, \theta) \bigr| \le C \Delta^2, \qquad  
w \in \{\ref{eq:ellip}, \ref{eq:hypo} \}, 
\end{align}
for some constant $C > 0$ independent of $\Delta$. We then deduce (\ref{eq:true_density_bd}) from (\ref{eq:piga}), (\ref{eq:ode_E}) and (\ref{eq:ode_H}).  

The inequality (\ref{eq:LDL_density_bd}) is immediately obtained from (\ref{eq:derive_LDL_z}) established in the proof of Lemma \ref{lemma:deriv_LDL}. The proof of Proposition \ref{prop:density_bd} is now complete. 
\qed  \\ 
\begin{lemma} \label{lemma:conv}
Let Assumptions \ref{ass:diff}--\ref{ass:dim} hold. Let  $w \in \{\ref{eq:ellip}, \ref{eq:hypo}\}$ and $x \in \mathbb{R}^N$. There exist constants $C, \lambda > 0$ such that for all $( y, \theta) \in \mathbb{R}^N \times \Theta$ and $0 \le s \le t < 1$, 
\begin{align}  
\Bigl|
\int_{\mathbb{R}^N}  
\overline{\mathscr{G}}^{(w), x} (t-s, \xi, y, \theta)
{\mathscr{G}}^{(w)} (s, x, \xi, \theta) 
d \xi 
\Bigr| 
\le C  t^{- \tfrac{m(w)}{2}} 
\exp \Bigl( - \lambda 
\bigl| 
\Gamma^{-1}_{t, (w)} 
\bigl( y - x - V_0 (x, \theta) t \bigr)
\bigr|^2 \Bigr).   
\label{eq:conv_1} 
\end{align}
Also, there exist constants $C, \lambda > 0$ such that for all $(x, y, z, \theta) \in \mathbb{R}^N \times \mathbb{R}^N  \times \mathscr{Z}_\kappa \times \Theta$ and $0 \le s \le t < 1$, 
\begin{align}
\Bigl|
\int_{\mathbb{R}^N}  
\overline{\mathscr{G}}^{(w), z} (t-s, \xi, y, \theta)
\overline{\mathscr{G}}^{(w), z} (s, x, \xi, \theta) 
d \xi 
\Bigr| 
\le C t^{- \tfrac{m(w)}{2}} 
\exp \Bigl( - \lambda 
\bigl| 
\nu^{(w), z} (t, y, x, \theta)
\bigr|^2 \Bigr),
\label{eq:conv_2}  
\end{align}
where $\nu^{(w), z}$ is defined in (\ref{eq:nu}). 
\end{lemma}

\noindent 
\textit{Proof of Lemma \ref{lemma:conv}.} 
%During the proof, we omit ``$(w)$" for notational simplicity because the argument applies to both cases $w = \ref{eq:ellip}$ and $w = \ref{eq:hypo}$. 
We first show the bound (\ref{eq:conv_1}). We set $\widetilde{\Gamma}_{t}^{(w)} \equiv 
({\Gamma}_{t, (w)})^\top {\Gamma}_{t, (w)},
\ t \in (0,1)$. Proposition \ref{prop:density_bd} with the multi-index $\alpha = \mathbf{0}$ gives: 
\begin{align}
 \Bigl| 
& \int_{\mathbb{R}^N}  
\overline{\mathscr{G}}^{(w), x} (t-s, \xi, y, \theta)
{\mathscr{G}}^{(w)} (s, x, \xi, \theta)  
d \xi \Bigr|
\nonumber \\ 
& \le C_1 \int_{\mathbb{R}^N} 
(t-s)^{- \tfrac{m(w)}{2}} \exp 
\Bigl( - \lambda_1
\bigl|
\nu^{(w), x} (t-s, y, \xi, \theta)
\bigr|^2 \Bigr) 
\times 
s^{- \tfrac{m(w)}{2}} 
\exp 
\Bigl( - \lambda_2 
\bigl| 
\Gamma_{s, (w)}^{-1} 
\bigl( \xi - x - V_0 (x, \theta) s \bigr)  
\bigr|^2 \Bigr)  
d\xi \nonumber \\
& \le \frac{C_2}{\sqrt{\det \bigl(c_1 \widetilde{\Gamma}_{t-s}^{(w)} + c_2 \widetilde{\Gamma}_{s}^{(w)} \bigr)}}
\exp \Bigl( - \lambda_3 v (t, s, y, x, \theta)^\top 
\bigl(c_1 \widetilde{\Gamma}_{t-s}^{(w)} + c_2 \widetilde{\Gamma}_{s}^{(w)} \bigr)^{-1} 
v (t, s, y, x, \theta) 
\Bigr) 
\qquad \bigl( \because (\mathrm{\ref{eq:standard}}) \bigr)
\nonumber \\
& \le \frac{C_3}{\sqrt{\det \widetilde{\Gamma}_{t}^{(w)}}}
\exp \Bigl( - \lambda_4 \, v (t, s, y, x, \theta)^\top 
\big(\widetilde{\Gamma}_{t}^{(w)} \bigr)^{-1}
v (t, s, y, x, \theta) 
\Bigr)  \nonumber \\
& = C_3 t^{- \tfrac{m(w)}{2}} 
\exp \Bigl( - \lambda_4 
\bigl| \Gamma_{t, (w)}^{-1} v (t, s, y, x, \theta) \bigr|^2 \Bigr),  
\label{eq:conv_3}
\end{align}
for some constants $C_1, C_2, C_3, c_1, c_2, \lambda_1, \lambda_2, \lambda_3, \lambda_4  > 0$ independent of $t, y, \theta$, where $v$ is defined in (\ref{eq:v_func}). 
% 
% Note that we have used (\ref{eq:standard}) in the derivation of the bound (\ref{eq:gauss_conv}). 
Under Assumptions \ref{ass:diff}-- \ref{ass:coeff}, we have that  
\begin{align} \label{eq:v}
v (t, s, y, x, \theta)   
& = - e^{-(t-s)A_{x, \theta}} 
\Bigl( y - e^{(t-s)A_{x, \theta}} 
\bigl( x + V_{0} (x, \theta) s \bigr) + \int_0^{t-s} e^{(t -s - u) A_{x, \theta}} b_{x, \theta} du
\Bigr) \nonumber \\ 
& = - e^{-(t-s)A_{x, \theta}} 
\Bigl( y - x - V_0 (x, \theta) t  
+ r (t, s, x, \theta)
\Bigr),
\end{align} 
where we considered the power series of the exponential matrices, and the term $r (t, s, x, \theta)$ satisfies: 
\begin{align}
| r (t, s, x, \theta) | \le C t^2, 
\end{align}
for some $C >0$ independent of $t, s, x, \theta$. We then deduce the bound (\ref{eq:conv_1}) from (\ref{eq:conv_3}) and (\ref{eq:v}). 

We next prove (\ref{eq:conv_2}). Using the same argument to derive (\ref{eq:conv_3}), typically the result (\ref{eq:standard}), we have that: 
\begin{align}
\Bigl|
\int_{\mathbb{R}^N}  
\overline{\mathscr{G}}^{(w), z} (t-s, \xi, y, \theta)
\overline{\mathscr{G}}^{(w), z} (s, x, \xi, \theta)  
d \xi  
\Bigr| \le C t^{-\tfrac{m(w)}{2}} 
\exp \Bigl( - \lambda 
\bigl| \Gamma_{t, (w)}^{-1}  \bar{v}^z (t, s, y, x, \theta)
\bigr|^2 \Bigr), 
\end{align}
for some $C, \lambda > 0$ independent of $t, s, x, y, z, \theta$, where  
\begin{align*}
\bar{v}^z (t, s, y, x, \theta) 
\equiv e^{s A_{z, \theta}} x + \int_0^s e^{(s-u)A_{z, \theta}}b_{z, \theta} du - e^{- (t-s) A_{z, \theta}} y 
+ \int_0^{t-s} e^{- u A_{z, \theta}} b_{z, \theta} du. 
\end{align*}
Noticing that 
\begin{align}
\bar{v}^z (t, s, y, x, \theta)  
% & = - e^{s A_{z, \theta}} 
% \Bigl( e^{-t A_{z, \theta}}y - x - \int_0^s e^{-u A_{z, \theta}} b_{z, \theta} du - \int_0^{t-s} e^{- (s+u) A_{z, \theta}} b_{z, \theta} du
% \Bigr) \nonumber \\ 
= - e^{s A_{z, \theta}} 
\Bigl( e^{-t A_{z, \theta}}y - x - \int_0^t e^{-u A_{z, \theta}} b_{z, \theta} du
\Bigr), 
\end{align}
we obtain the desired bound (\ref{eq:conv_2}). 
\qed 
\begin{lemma} \label{lemma:poly_bd}
Let $w \in \{\ref{eq:ellip}, \ref{eq:hypo} \}$. Let  $x \in \mathbb{R}^N$ be the initial state of SDE and Assumption \ref{ass:coeff} hold. For any $\alpha \in \{1, \ldots, N\}^k, \, k \in \mathbb{N}$ and $\lambda_1 > 0$, there exist constants $C, \lambda_2 > $ such that for all $(t, \xi, \theta) \in (0, 1) \times \mathbb{R}^N \times \Theta$, 
\begin{align}
\Bigl|
\prod_{1 \le i \le k} (\xi_{\alpha_i}-x_{\alpha_i})
\Bigr| 
\exp \Bigl(  - \lambda_1 \bigl| \Gamma_{t, (w)}^{-1} \bigl( \xi - x - V_0 (x, \theta) t  \bigr) \bigr|^2  \Bigr)  
\le C t^{\tfrac{k}{2}} 
\exp 
\Bigl(  - \lambda_2 \bigl| \Gamma_{t, (w)}^{-1} \bigl( \xi - x - V_0 (x, \theta) t  \bigr) \bigr|^2  \Bigr). 
\end{align}
\end{lemma}

\noindent 
{\textit{Proof of Lemma \ref{lemma:poly_bd}.}}
We have that: for all $1 \le j \le k$, 
\begin{align*}
& \bigl| \xi_{\alpha_j} - x_{\alpha_j}  \bigr| 
\exp 
\Bigl(  - \lambda_1 \bigl| \Gamma_{t, (w)}^{-1} \bigl( \xi - x - V_0 (x, \theta) t  \bigr) \bigr|^2  \Bigr) 
\\ 
& 
\le \bigl| \xi_{\alpha_j} - x_{\alpha_j} - V_0^{\alpha_j} (x, \theta) t \bigr| \, 
\exp 
\Bigl(  - \lambda_1 \bigl| \Gamma_{t, (w)}^{-1} \bigl( \xi - x - V_0 (x, \theta) t  \bigr) \bigr|^2  \Bigr)  
+ C_1 t  \exp 
\Bigl(  - \lambda_1 \bigl| \Gamma_{t, (w)}^{-1} \bigl( \xi - x - V_0 (x, \theta) t  \bigr) \bigr|^2  \Bigr)  
\\ 
& 
\le C_2  t^{\tfrac{1}{2} + \mathbf{1}_{1 \le \alpha_j \le N_S}} \, 
\exp 
\Bigl(  - \lambda_2 \bigl| \Gamma_{t, (w)}^{-1} \bigl( \xi - x - V_0 (x, \theta) t  \bigr) \bigr|^2  \Bigr)  
+ C_1 t  
\exp 
\Bigl(  - \lambda_1 \bigl| \Gamma_{t, (w)}^{-1} \bigl( \xi - x - V_0 (x, \theta) t  \bigr) \bigr|^2  \Bigr) \\ 
& \le C_3 t^{\tfrac{1}{2}} \, 
\exp 
\Bigl(  - \lambda_2 \bigl| \Gamma_{t, (w)}^{-1} \bigl( \xi - x - V_0 (x, \theta) t  \bigr) \bigr|^2  \Bigr)   
\end{align*}
for some constants $C_1, C_2, C_3, \lambda_2 >0$ independent of $t, \xi, \theta$, where we used \eqref{eq:exp_bd} in the second inequality. 
% 
% \begin{align}
% \tfrac{\bigl| \xi_{\alpha_j} - x_{\alpha_j} - V_0^{\alpha_j} (x, \theta) t \bigr|}{t^{1/2 + \mathbf{1}_{1 \le \alpha_j \le N_S}}}
% \exp 
% \Bigl(  - \lambda_1 \bigl| \Gamma_{t, (w)}^{-1} \bigl( \xi - x - V_0 (x, \theta) t  \bigr) \bigr|^2  \Bigr)   
% \le 
% C \exp 
% \Bigl(  - \lambda \bigl| \Gamma_{t, (w)}^{-1} \bigl( \xi - x - V_0 (x, \theta) t  \bigr) \bigr|^2  \Bigr), 
% \end{align}
% % 
% for some $C, \lambda$ independent of $t, \xi, x, \theta$.
\qed  \\ 

We also introduce the following two results:  
\begin{lemma} \label{lemma:grad_D}
Let Assumptions \ref{ass:diff}, \ref{ass:param}, \ref{ass:coeff} hold. Let $\psi^{(w)} \in \mathscr{F}^{(w)}, \, w \in \{\ref{eq:ellip}, \ref{eq:hypo} \}$ and $J$ be an non-negative integer.  
Then, for any $\alpha \in \mathbb{Z}_{\ge 0}^N$, there exist constants $C > 0$ and $q \ge 1$ such that for all $(t, \xi, z, \theta) \in (0,1) \times \mathbb{R}^N \times \mathscr{Z}_\kappa \times \Theta$,  
\begin{align} \label{eq:grad_J}
\bigl| 
\partial^{\alpha}  \mathscr{D}_{(J)}^{z, \theta}  
\bigl\{ \psi^{(w)}  (t, \cdot) \bigr\} (\xi) 
\bigr| \le C t^{- \| \alpha \|_w - Q_w (J)} 
\times \Bigl(1 + \bigl|\tfrac{\xi - z}{\sqrt{t}} \bigr|^q \Bigr) \times \bigl| \widetilde{\psi}^{(w)} (t, \xi) \bigr|, \qquad w \in \{\ref{eq:ellip}, \ref{eq:hypo} \}, 
\end{align}
where 
\begin{align} \label{eq:D_J}
\mathscr{D}_{(J)}^{z, \theta} \equiv \mathrm{ad}^J_{\mathscr{L}_\theta^{0, z}} \bigl( \widetilde{\mathscr{L}}^{z}_\theta \bigr)
\end{align} 
and $Q_w(J)$ is defined as: 
\begin{align} \label{eq:Q}
Q_{\ref{eq:ellip}}  (J) 
\equiv 
\begin{cases}
    \tfrac{1}{2}, & J = 0; \\[0.2cm] 
    \tfrac{1}{2} (J + 2),  & J \ge 1,   
\end{cases} 
\qquad  
Q_{\ref{eq:hypo}} (J) \equiv  
\begin{cases}
 \tfrac{1}{2}, & J = 0; \\[0.2cm]
 \tfrac{1}{2} \floor{\tfrac{3}{2} J} + 1, & J= 1, 2; \\[0.2cm]  
 \tfrac{1}{2} \floor{\tfrac{3}{2} J} + \tfrac{3}{2}, &  
 J \ge 3.  
\end{cases}
\end{align} 
Furthermore, (\ref{eq:grad_J}) holds with the differential operator $\mathscr{D}_{(J)}^{z, \theta}$ being replaced with $(\mathscr{D}_{(J)}^{z, \theta})^\ast$. 
\end{lemma} 
\begin{lemma} \label{lemma:grad_D_a}
Let Assumptions \ref{ass:diff}, \ref{ass:param} and \ref{ass:coeff} hold. Let $x \in \mathbb{R}^N$ be the initial state of SDE  
% $\kappa > 0$ be the constant defined in Assumption \ref{ass:coeff} 
and $\psi^{(w)} \in \mathscr{F}^{(w)}, \, w \in \{\ref{eq:ellip}, \ref{eq:hypo} \}$. Then for any $\alpha \in \mathbb{Z}_{\ge 0}^j,\, j \ge 2$, there exists a constant $C > 0$ such that for all $t \in (0,1)$, 
\begin{align} \label{eq:F}
\bigl| 
\mathscr{D}_\alpha^{z, \theta} \psi^{(w)} (t, x)
|_{z=x} 
\bigr| 
\le 
C t^{- |\alpha| - \tfrac{j}{2}}  
\times  
\bigl| \widetilde{\psi}^{(w)} (t, x) \bigr|. 
% \begin{cases}    
% C t^{- \tfrac{|\alpha|}{2} - j + \tfrac{\| \alpha \|_0}{2}} \times \bigl| \widetilde{\psi}^{(\ref{eq:ellip})} (t, x) \bigr|, 
% & w = \ref{eq:ellip}; \\[0.2cm] 
% C t^{- \tfrac{1}{2} \sum_{1 \le i \le j} \floor{\tfrac{3 \alpha_i}{2}} - \tfrac{3}{2}j + \tfrac{1}{2} \| \alpha \|_{0, 1, 2}} \times \bigl| \widetilde{\psi}^{(\ref{eq:hypo})} (t, x) \bigr|, 
% & w = \ref{eq:hypo}. 
% \end{cases} 
\end{align}
% 
%where we have set:
%\begin{align}
%F^{(w)} (\lambda) \equiv 
%\begin{cases}
%\tfrac{1}{2} | \lambda | +  k - \tfrac{1}{2} \| \lambda \|_0,  & w = \ref{eq:ellip},  \\ 
%\tfrac{1}{2} \sum_{i = 1}^{k} \floor{  \tfrac{3 \lambda_i }{2} }
%+  \tfrac{3}{2}k  -  \tfrac{1}{2} \| \lambda \|_{0,1,2},  & w = \ref{eq:hypo}, 
%\end{cases} 
%\qquad \lambda \in \mathbb{Z}_{\ge 0}^k, \  k \in \mathbb{N}. 
%\end{align} 
% 
\end{lemma}  

The proof of above two lemmas are postponed to Section \ref{app:pf_grad_D} after studying the expression of  differential operators $\mathscr{D}_{(J)}^{z, \theta}$ and $\mathscr{D}_\alpha^{z, \theta}$ in Section \ref{app:diff_step_1} and \ref{app:diff_step_2}, respectively. 

\subsection{Proof of Theorem \ref{thm:bd_r1}}
\label{app:pf_bd_r1} 
We introduce the following key result for proving Theorem \ref{thm:bd_r1}. 
\begin{proposition} \label{prop:bd_E_1}
Let Assumptions \ref{ass:diff}--\ref{ass:dim} hold. Let $x \in \mathbb{R}^N$ and $w \in 
\{ \ref{eq:ellip}, \ref{eq:hypo} \}$. For any $\alpha \in \mathbb{Z}_{\ge 0}^{N}$ and $j \in \mathbb{N}$, there exists a constant $C > 0$ such that for any $(\Delta, y, \theta) \in (0, 1) \times \mathbb{R}^N \times \Theta$ and for any $0 \le s_{j} \le s_{j-1} \le \cdots \le s_1 \le \Delta, \, j \in \mathbb{N}$, 

\begin{align}
\begin{aligned}
& 
\left| \partial^{\alpha}_y P_{s_j}^{\theta, (w)} \widetilde{\mathscr{L}}_\theta^z \bar{P}_{s_{j-1} - s_j}^{\theta, z, (w)}
\widetilde{\mathscr{L}}_\theta^z \bar{P}_{s_{j-2} - s_{j-1}}^{\theta, z, (w)} \cdots 
\widetilde{\mathscr{L}}_\theta^z p_{\Delta - s_1}^{\bar{X}^z, (w)} (\cdot, y ;\theta)  (x) |_{z=x}  
\right|    \\ 
& \quad 
\le C {\Delta^{ - \|\alpha\|_w  -  {1}/{2}}} 
\cdot \Bigl\{ \prod _{1 \le i \le j-1}  s_i^{-{1}/{2}} \Bigr\}
\cdot 
\Delta^{- {m(w)}/{2}}
\exp \Bigl( - \lambda 
\bigl| \Gamma_{\Delta,(w)}^{-1} 
\bigl( y - x -V_0 (x, \theta) \Delta \bigr) \bigr|^2 \Bigr),   
\end{aligned}  \label{eq:true_kernel} 
\end{align} 
where we interpret 
$\textstyle \prod _{1 \le i \le j-1}  s_i^{-{1}/{2}} \equiv 1$ for $j = 1$. 
\end{proposition}
% 
%\subsubsection{Proof of Lemma \ref{prop:bd_E_1}}

\noindent 
\textit{Proof of Proposition \ref{prop:bd_E_1}.}
We exploit the mathematical induction on $j \in \mathbb{N}$. \\ 

\noindent
\textbf{Step I.} We consider the case $j=1$. We have 
\begin{align}
\partial^{\alpha}_y P_{s_1}^{\theta, (w)} \widetilde{\mathscr{L}}_\theta^z p_{\Delta - s_1}^{\bar{X}^z, (w)} (\cdot, y ;\theta) (x) |_{z=x}  
= \int_{\mathbb{R}^N} \partial^{\alpha}_y  \widetilde{\mathscr{L}}_\theta^z p_{\Delta - s_1}^{\bar{X}^z, (w)} (\cdot, y ;\theta)  (\xi)
p_{s_1}^{X, (w)} (x, \xi; \theta) 
d \xi |_{z=x}. \label{eq:term_1}
\end{align}
To give a bound for (\ref{eq:term_1}) uniformly in $s_1 \in [0, \Delta]$, we separately consider the cases $s_1 \in [0, \Delta/2]$ and $s_1 \in [\Delta/2, \Delta]$. 
\\ 

\noindent 
\textbf{Step I-1.} Consider the case $s_1 \in [0, \Delta/2]$. We have that:
\begin{align}
& \Bigl| \partial_y^\alpha \widetilde{\mathscr{L}}_\theta^z p_{\Delta - s_1}^{\bar{X}^z, (w)} (\cdot, y ;\theta)  (\xi)
\bigr|_{z=x} \Bigr| 
\le C_1 \Bigl| \partial_\xi^\alpha  \widetilde{\mathscr{L}}_\theta^z  p_{\Delta - s_1}^{\bar{X}^z, (w)} (\cdot, y ;\theta)  (\xi) \bigr|_{z=x} \Bigr|  
\qquad \bigl( \because \eqref{eq:basic} \bigr)
\nonumber \\ 
& \le 
C_2 (\Delta - s_1)^{-\|\alpha \|_{w} - \tfrac{1}{2}} \times 
\Bigl(1 +  \bigl| \tfrac{\xi- x}{\sqrt{\Delta - s_1}} \bigr|^q \Bigr) 
\times \overline{\mathscr{G}}^{(w), x} (\Delta-s_1, \xi, y,  \theta) 
 \qquad \bigl( \because \mathrm{Lemma \, \ref{lemma:grad_D}} \, \mathrm{with} \, J = 0  \bigr) \nonumber \\
& 
\le C_3 \Delta^{-\|\alpha \|_{w} - \tfrac{1}{2}} \times
\Bigl( 1 +  \bigl| \tfrac{\xi- x}{\sqrt{\Delta}} \bigr|^q \Bigr) \times 
\overline{\mathscr{G}}^{(w), x} (\Delta-s_1, \xi, y, \theta), \qquad \bigl( \because 
s_1 \in [0, \Delta/2] \bigr)  
\label{eq:bd_0}
\end{align}
for some constants $C_1, C_2, C_3 > 0$ and $q \ge 1$, which are independent of $\Delta, s_1, \xi, y, \theta$. Making use of the bound (\ref{eq:bd_0}) and Proposition \ref{prop:density_bd} for $p_{s_1}^X (x, \xi; \theta)$, we have that: 
\begin{align*}
\bigl|  \eqref{eq:term_1} \bigr|  
& \le C_1 \Delta^{-\|\alpha \|_{w} - \tfrac{1}{2}} 
\times \int_{\mathbb{R}^N} 
\overline{\mathscr{G}}^{(w),x}(\Delta-s_1, \xi, y, \theta) 
\, {\mathscr{G}}^{(w)}(s_1, x, \xi, \theta) d\xi 
\qquad \bigl( \because \mathrm{Lemma \, \ref{lemma:poly_bd}} \bigr) \\ 
& \le C_2 \Delta^{-\|\alpha \|_{w} - \tfrac{1}{2}} \times
\mathscr{G}^{(w)} (\Delta, x, y, \theta). 
\qquad \big(\because \mathrm{Lemma \, \ref{lemma:conv}}, \, s_1 \in [0, \Delta/2] \bigr) 
\end{align*}
% 
% where $C_1, C_2, \lambda > 0$ are some constants independent of $y, \theta, \Delta, s_1$. 
% In particular, in the derivation of the second inequality, we have also exploited Proposition \ref{prop:density_bd} and Lemma \ref{lemma:poly_bd}.  

\noindent
\textbf{Step I-2.} Consider the case $s_1 \in [\Delta/2, \Delta]$.  We emphasise that the argument in \textbf{Step I-1} can not be employed here since the bound is obtained under 
$s_1 \in [0, \Delta/2]$, see \eqref{eq:bd_0}. To provide an appropriate bound, we apply the integration by parts (IBP) to remove the partial derivatives on the Gaussian kernel $p_{\Delta -s_1}^{\bar{X}^z, (w)}$ and then aim to derive an upper bound involving $(s_1)^{-K}$ for some $K > 0$ (not $(\Delta - s_1)^{-K}$), which is ultimately bounded by $\Delta^{-K}$ under $s_1 \in [\Delta/2, \Delta]$. We have that:  
\begin{align*}
\bigl|  
\eqref{eq:term_1}
\bigr| 
& = \Bigl| 
\int_{\mathbb{R}^N} 
\partial^{\alpha}_y  p_{\Delta - s_1}^{\bar{X}^z, (w)} (\xi, y; \theta) 
\bigl( \widetilde{\mathscr{L}}_\theta^{z} \bigr)^\ast 
\{ p_{s_1}^{X, (w)} (x, \cdot; \theta) \} (\xi)  d \xi
|_{z=x}  \Bigr|  \\ 
& 
\le C_1 
\int_{\mathbb{R}^N} 
p_{\Delta - s_1}^{\bar{X}^x, (w)} (\xi, y; \theta) 
\bigl| 
\partial^{\alpha}_{\xi} \bigl( 
\widetilde{\mathscr{L}}_\theta^{z} \bigr)^\ast 
\{ p_{s_1}^{X, (w)} (x, \cdot; \theta) \} (\xi) 
|_{z=x} \bigr| d \xi 
\qquad \bigl( \because \eqref{eq:basic}, \, \mathrm{IBP} \bigr) \\ 
& \le C_2 (s_1)^{- \| \alpha \|_w - \tfrac{1}{2}}
\int_{\mathbb{R}^N}  \overline{\mathscr{G}}^{(w), x} (\Delta-s_1, \xi, y, \theta) {\mathscr{G}}^{(w)} (s_1, x, \xi,  \theta)  d \xi \qquad \big(\because \mathrm{Lemmas \, \ref{lemma:poly_bd}, \, \ref{lemma:grad_D}} \bigr) 
\\ 
& \le C_3 \Delta^{- \| \alpha \|_w - \tfrac{1}{2}} \times  \mathscr{G}^{(w)} (\Delta, x, y, \theta) 
\qquad \big(\because \mathrm{Lemma \, \ref{lemma:conv}}, \, s_1 \in [\Delta/2, \Delta] \bigr),   
\end{align*}
for some constants $C_1, C_2, C_3 > 0$ independent of $y, \theta, \Delta, s_1$. Step I is now complete. \\ 

\noindent
\textbf{Step II.} We assume that the assertion holds for $j = 1, \ldots, k-1, \, k \ge 2$ and study the case $j = k$. To get an appropriate bound, we again consider the cases $s_1 \in [0, \Delta/2]$ and $s_1 \in [\Delta/2, \Delta]$ separately. \\ 

\noindent
\textbf{Step II-1.} Consider the case $s_1 \in [0, \Delta/2]$. We have 
\begin{align}
&  
\partial^{\alpha}_y \bigl( P_{s_k}^{\theta, (w)} \widetilde{\mathscr{L}}_\theta^z \bar{P}_{s_{k-1} - s_k}^{\theta, z, (w)}
\widetilde{\mathscr{L}}_\theta^z 
\cdots 
\bar{P}_{s_1 - s_2}^{\theta, z, (w)}
\widetilde{\mathscr{L}}_\theta^z 
\bigr) \bigl\{ p_{\Delta - s_1}^{\bar{X}^z, (w)} (\cdot, y ;\theta) \bigr\} (x) |_{z=x}  
\nonumber \\ 
& = 
\int_{\mathbb{R}^N}
\partial^{\alpha}_y \widetilde{\mathscr{L}}_\theta^z \bigl\{  p_{\Delta - s_1}^{\bar{X}^z, (w)} (\cdot, y; \theta) \bigr\} (\xi) 
\, 
\bigl( P_{s_k}^{\theta, (w)} \widetilde{\mathscr{L}}_\theta^z \bar{P}_{s_{k-1} - s_k}^{\theta, z, (w)}
\widetilde{\mathscr{L}}_\theta^z 
\cdots \widetilde{\mathscr{L}}_\theta^z \bigr) 
\{ p_{s_1 - s_2}^{\bar{X}^z} (\cdot, \xi; \theta)  \bigr\} (x)  d \xi
|_{z=x}.  \label{eq:term_2}
\end{align}
Noticing that the assumption of induction yields 
\begin{align}
\bigl| \bigl( P_{s_k}^{\theta, (w)} \widetilde{\mathscr{L}}_\theta^z \bar{P}_{s_{k-1} - s_k}^{\theta, z, (w)}
\widetilde{\mathscr{L}}_\theta^z 
\cdots \widetilde{\mathscr{L}}_\theta^z \bigr) 
\{ p_{s_1 - s_2}^{\bar{X}^z} (\cdot, \xi; \theta) \} (x) \bigr|_{z=x} 
\bigr| 
\le C s_1^{-1/2} \cdot  \Bigl\{ \prod_{2 \le i \le k-1} s_i^{-1/2} \Bigr\} 
\cdot \mathscr{G}^{(w)} (s_1, x, \xi, \theta), 
\label{eq:ass}
\end{align}
for some constant $C >0$ independent of $\{s_{i}\}_{i = 1 }^{k-1}$, $\xi$ and $\theta$, we have that:
\begin{align*}
& \bigl| \eqref{eq:term_2} \bigr| \\ 
& \le C_1 \Delta^{-\| \alpha \|_w - {1}/{2}} \cdot 
\Bigl\{ \prod_{1 \le i \le k-1} s_i^{-1/2} \Bigr\} \times \int_{\mathbb{R}^N} 
\Bigl(1 + \bigl| \tfrac{\xi - x}{\sqrt{\Delta}} \bigr|^q
\Bigr) 
\, \overline{\mathscr{G}}^{(w), x} (\Delta-s_1, \xi, y, \theta) 
\, \mathscr{G}^{(w)} (s_1, x, \xi, \theta)  d \xi
\qquad \bigl( \because \eqref{eq:bd_0}, \, \eqref{eq:ass} \bigr) \\
& \le C_2 \Delta^{-\| \alpha \|_w - 1/2} 
\cdot \Bigl\{ \prod_{1 \le i \le k-1} s_i^{-1/2} \Bigr\} 
\times \int_{\mathbb{R}^N}  
\overline{\mathscr{G}}^{(w), x} (\Delta-s_1, \xi, y, \theta) 
\, 
\mathscr{G}^{(w)} (s_1, x, \xi, \theta)  d \xi
\quad \bigl( \because \mathrm{Lemma \, \ref{lemma:poly_bd}} \bigr)
\\ 
& \le C_3 \Delta^{-\| \alpha \|_w - 1/2} 
\cdot 
\Bigl\{ \prod_{1 \le i \le k-1} s_i^{-1/2} \Bigr\} 
\times \mathscr{G}^{(w)} (\Delta, x, y, \theta) 
\qquad \bigl( \because \mathrm{Lemma \, \ref{lemma:conv}} \bigr)  
\end{align*}
for some constants $q \ge 1$ and $C_1, C_2, C_3 > 0$ independent of $\{s_i\}_{i = 1}^{k-1}, \Delta, y, \theta$. 
\\ 

\noindent 
\textbf{Step II-2.} 
We consider the case $s_1 \in [\Delta/2, \Delta]$. IBP together with (\ref{eq:basic}) yields that:
\begin{align}
\bigl| \eqref{eq:term_2} \bigr| 
\le C \int_{\mathbb{R}^N}
p_{\Delta - s_1}^{\bar{X}^x, (w)} (\xi, y; \theta) 
\times 
\bigl| 
\underbrace{ \partial^{\alpha}_\xi  
(\widetilde{\mathscr{L}}_\theta^z)^\ast 
\bigl( P_{s_k}^{\theta, (w)} \widetilde{\mathscr{L}}_\theta^z \bar{P}_{s_{k-1} - s_k}^{\theta, z, (w)}
\widetilde{\mathscr{L}}_\theta^z 
\cdots \widetilde{\mathscr{L}}_\theta^z p_{s_1 - s_2}^{\bar{X}^z, (w)} (x, \cdot; \theta) \bigr) (\xi)
|_{z=x}}_{=: (\flat)} \bigr| d \xi.  
\end{align}
Making use of the assumption of induction and Lemma \ref{lemma:grad_D} with 
$$
\psi^{(w)}(s_1, \xi) \equiv  
\bigl( P_{s_k}^{\theta, (w)} \widetilde{\mathscr{L}}_\theta^z \bar{P}_{s_{k-1} - s_k}^{\theta, z, (w)}
\widetilde{\mathscr{L}}_\theta^z 
\cdots \widetilde{\mathscr{L}}_\theta^z p_{s_1 - s_2}^{\bar{X}^z, (w)} (\cdot, \xi; \theta) \bigr) (x),
$$ we obtain: 
\begin{align}
% & \bigl| 
% \partial^{\alpha}_\xi  
% (\widetilde{\mathscr{L}}_\theta^z)^\ast 
% \bigl( P_{s_k}^{\theta, (w)} \widetilde{\mathscr{L}}_\theta^z \bar{P}_{s_{k-1} - s_k}^{\theta, z, (w)}
% \widetilde{\mathscr{L}}_\theta^z 
% \cdots \widetilde{\mathscr{L}}_\theta^z p_{s_1 - s_2}^{\bar{X}^z, (w)} (x, \cdot; \theta) \bigr) (\xi)
% |_{z=x} \bigr| \nonumber  \\
| (\flat) | & \le C_1 (s_1)^{- \| \alpha \|_w - \tfrac{1}{2}} 
\Bigl(1 + \bigl| \tfrac{\xi - x}{\sqrt{s_1}} \bigr|^q 
\Bigr) \times 
\Bigl\{ \prod_{1 \le i \le k-1} s_i^{-1/2} \Bigr\}
\times \mathscr{G}^{(w)} (s_1, x, \xi, \theta)  \qquad \bigl( \because \mathrm{Lemma \, \ref{lemma:grad_D}} \, \mathrm{with} \, J = 0 \bigr) 
\nonumber \\ 
& \le C_2 \Delta^{- \| \alpha \|_w - \tfrac{1}{2}} 
\times \Bigl\{ \prod_{1 \le i \le k-1} s_i^{-1/2} \Bigr\} \times \mathscr{G}^{(w)} (s_1, x, \xi, \theta)  
\qquad \bigl( \because \mathrm{Lemma \, \ref{lemma:poly_bd}}, \ s_1 \in [\Delta/2, \Delta] \bigr) \label{eq:term_3}
\end{align}
where $q \ge 1$ and $C_1, C_2 > 0$ are constants independent of $\xi, \theta, t$ and $\{ s_i \}_{i = 1}^{k-1}$. We thus obtain from Proposition \ref{prop:density_bd} and the upper bound (\ref{eq:term_3}) that: 
\begin{align*}
\bigl| \eqref{eq:term_2} \bigr|
& \le C \Delta^{- \| \alpha \|_w - \tfrac{1}{2}} 
\cdot \Bigl\{ \prod_{1 \le i \le k-1} s_i^{-1/2} \Bigr\} \times \int_{\mathbb{R}^N}
\overline{\mathscr{G}}^{(w), x} (\Delta - s_1, \xi, y, \theta) \mathscr{G}^{(w)} (s_1, x, \xi, \theta) 
d \xi \\
& \le C \Delta^{- \| \alpha \|_w - \tfrac{1}{2}} 
\cdot 
\Bigl\{ \prod_{1 \le i \le k-1} s_i^{-1/2} \Bigr\}
\times \mathscr{G}^{(w)} (\Delta, x, y, \theta).  \qquad \bigl( \because \mathrm{Lemma \, \ref{lemma:conv}} \bigr)   
\end{align*}
Thus, the assertion holds for the case $j = k$, and the proof of Proposition \ref{prop:bd_E_1} is now complete. 
\qed  \\  

\noindent 
\textit{Proof of Theorem \ref{thm:bd_r1}.}
Recall the definition of $\mathscr{R}_1^{M, (w)} (\Delta, x, y; \theta) $ in (\ref{eq:first_err}) in the main text. Making use of Proposition \ref{prop:bd_E_1} with $j = M$, 
% 
% $$
% P_{s_M}^{\theta, (w)} \widetilde{\mathscr{L}}_\theta^z \bar{P}_{s_{M-1} - s_M}^{\theta, z, (w)}
% \widetilde{\mathscr{L}}_\theta^z \bar{P}_{s_{M-2} - s_{M-1}}^{\theta, z, (w)} \cdots 
% \widetilde{\mathscr{L}}_\theta^z p_{\Delta - s_1}^{\bar{X}^z, (w)} (x, y ;\theta), $$
% 
we obtain:  
\begin{align*}
& \bigl| 
\mathscr{R}_1^{M, (w)} (\Delta, x, y; \theta) 
\bigr| 
= \left| 
% \int_0^{\Delta} ds_1 \int_0^{s_1} ds_2 \cdots \int_0^{s_{M-1}} ds_M  
\int_{0 \le s_M \le \cdots \le s_1 \le \Delta}
P_{s_M}^{\theta, (w)} \widetilde{\mathscr{L}}_\theta^z \bar{P}_{s_{M-1} - s_M}^{\theta, z, (w)}
\widetilde{\mathscr{L}}_\theta^z \bar{P}_{s_{M-2} - s_{M-1}}^{\theta, z, (w)} \cdots 
\widetilde{\mathscr{L}}_\theta^z p_{\Delta - s_1}^{\bar{X}^z, (w)} (x, y ;\theta) |_{z=x}
\, ds_M \cdots ds_1  
\right|
\\
& \quad \le C_1 \Delta^{-\tfrac{1}{2}}
% \int_0^{\Delta} ds_1 \int_0^{s_1} ds_2 \cdots \int_0^{s_{M-1}} ds_M \times 
\int_{0 \le s_M \le \cdots \le s_1 \le \Delta} 
\Bigl\{\prod_{1 \le i \le M-1} 
s_i^{-1/2} \Bigr\} 
\cdot {\Delta^{-m(w)/2}} 
\exp \Bigl( - \lambda \bigl| \Gamma_{\Delta, (w)}^{-1} \bigl( y - x - V_0 (x, \theta) \Delta \bigr) \bigr|^2 \Bigr) ds_M \cdots ds_1 
\\
& \quad \le {C_2}{\Delta^{M/2}} 
\cdot \Delta^{-m(w)/2}
\exp \Bigl( - \lambda \bigl| \Gamma_{\Delta, (w)}^{-1} \bigl( y - x - V_0 (x, \theta) \Delta \bigr) \bigr|^2 \Bigr), \qquad w \in \{\ref{eq:ellip}, \ref{eq:hypo} \}, 
\end{align*}
for some constants $C_1, C_2, \lambda > 0$ independent of $(\Delta, y, \theta) 
\in (0,1) \times  \mathbb{R}^N \times \Theta$. 
\qed 
\subsection{Proof of Theorem \ref{thm:bd_r2}} 
\label{app:pf_bd_r2}

For notational simplicity, we write: for $k \in \mathbb{N}$, $0 \le s_{j} \le s_{j-1} \le \cdots \le s_1 \le \Delta, \, j \in \mathbb{N}$ with $\Delta \in (0,1)$ and $(\xi, y, z, \theta) \in \mathbb{R}^N \times \mathbb{R}^N \times \mathscr{Z}_\kappa  \times \Theta$, 
\begin{align*}
{F}^{[k],z, (w)}_{s_{1:j}} (\xi, y; \theta) 
\equiv 
\bar{P}_{s_j}^{\theta, z, (w)} \mathscr{D}_{(k)}^{z, \theta}  
\bar{P}_{s_{j-1} - s_j}^{\theta, z, (w)}
\widetilde{\mathscr{L}}_\theta^z   
\bar{P}_{s_{j-2} - s_{j-1}}^{\theta, z, (w)} \cdots 
\bar{P}_{s_{1} - s_{2}}^{\theta, z, (w)} 
\widetilde{\mathscr{L}}_\theta^z 
\{ p_{\Delta - s_1}^{\bar{X}^z, (w)} (\cdot, y ;\theta)  \} (\xi), 
\quad w \in \{\ref{eq:ellip}, \ref{eq:hypo} \}, 
\end{align*}
where $s_{1:j} = \{s_i\}_{i =1}^j$ and $\mathscr{D}^{z, \theta}_{(k)}$ is defined in \eqref{eq:D_J}. We note that $\mathscr{D}_{(0)}^{z, \theta} \equiv \widetilde{\mathscr{L}}_\theta^z = \mathscr{L}_\theta - \mathscr{L}_\theta^{0, z}$ and in the above definition, the true intractable transition density $p^X$ is not involved.  
\begin{proposition} \label{prop:bd_E_2}
Let Assumptions \ref{ass:diff}--\ref{ass:dim} hold. Let $w \in \{\ref{eq:ellip}, \ref{eq:hypo} \}$  
and $\kappa > 0$ be the constant defined in Assumption \ref{ass:coeff}. 
Then, for any $J \in \mathbb{N}$, 
$j \in \mathbb{N}$ and multi-indices $\beta, \gamma \in \mathbb{Z}_{\ge 0}^{N}$, there exist constants $C > 0$ and $q \ge 1$ such that for any  $(\Delta, \xi, y, z, \theta) \in (0, 1) \times \mathbb{R}^N \times \mathbb{R}^N \times \mathscr{Z}_\kappa \times \Theta$ and for any $0 \le s_{j} \le s_{j-1} \le \cdots \le s_1 \le \Delta, \, j \in \mathbb{N}$, 
\begin{align}
% \sup_{
% \substack{z \in \mathscr{Z}_\kappa \\ 
% 0\le s_j \le \cdots \le s_1 \le \Delta }}
\bigl| 
\partial^\beta_{\xi} \partial^{\gamma}_y {F}^{[J],z, (w)}_{s_{1:j}} (\xi, y; \theta)   \bigr| 
\le
C \Delta^{- \| \beta  + \gamma \|_{w} 
- \tfrac{j-1}{2}  - Q_{w} (J)}
\cdot 
\Bigl( 
1 +  \bigl| \tfrac{\xi - z}{\sqrt{\Delta}} \bigr|^q 
\Bigr)
\cdot \overline{\mathscr{G}}^{(w), z} (\Delta, \xi, y, \theta), \label{eq:LDL_kernel} 
\end{align}  
where $Q_w (J)$ is defined in (\ref{eq:Q}). 
\end{proposition} 

\noindent 
\textit{Proof of Proposition \ref{prop:bd_E_2}}.
We exploit the mathematical induction on the integer $j$.  In what follows, we make use of notation $C$ as the positive constant satisfying the property in the statement. Note that `$C$' can be different from line to line. \\ 

\noindent 
\textbf{Step I.} 
We assume $j=1$. We have 
\begin{align}
\partial^\beta_\xi \partial^{\gamma}_y  
{F}^{[J], z, (w)}_{s_{1}} (\xi, y; \theta)  
& = \int_{\mathbb{R}^N} 
\mathscr{D}_{(J)}^{z, \theta}  
\bigl( \partial_y^\gamma p_{\Delta - s_1}^{\bar{X}^z, (w)} (\cdot, y ;\theta) 
\bigr) (\eta) \, 
\partial_{\xi}^\beta p_{s_1}^{\bar{X}^z, (w)} (\xi, \eta; \theta) d \eta. 
% \nonumber \\[0.2cm]
% & 
% = \sum_{\lambda \in \mathscr{I}_w (k)} 
% \int_{\mathbb{R}^N} \mathscr{W}_\lambda^{[k], z} (\xi, \theta) \, \partial_y^\gamma \partial_\xi^\lambda  p_{\Delta-s_1}^{\bar{X}^z, (w)} (\xi, y; \theta) 
% \partial_{\xi}^\beta p_{s_1}^{\bar{X}^z, (w)} (x, \xi; \theta) d \xi 
\label{eq:int_1} 
\end{align}  
% 
% where the set of multi-indices $\mathscr{I}_w(k)$ is defined in (\ref{eq:set_I}), and $\mathscr{W}_\lambda^{[k], z} : \mathbb{R}^N \times \Theta \to \mathbb{R}$ is defined in Lemma \ref{lemma:diff_step_1} and characterised as follows: under Assumptions \ref{ass:diff}, \ref{ass:param}, \ref{ass:coeff}, for any $k \in \mathbb{Z}_{\ge 0}$ and $\lambda \in \mathscr{I}_w (k)$, there exist constants $C >  0$ and $q \ge 0$ such that for all $(z, \xi, \theta) \in \mathscr{Z}_\kappa \times \mathbb{R}^N \times \Theta$
% % 
% \begin{align}
% \bigl| \mathscr{W}_\lambda^{[k], z} (\xi, \theta) \bigr|  \le C |\xi - z|^q. 
% \end{align} 
% 
To obtain an appropriate upper bound for (\ref{eq:int_1}), we separately treat the cases $s_1 \in [0, \Delta/2]$ and $s_1 \in [\Delta/2, \Delta]$ as we did in the proof of Proposition \ref{prop:bd_E_1}.  \\ 

\noindent 
\textbf{Step I-1.} 
We consider the case $s_1 \in [0, \Delta/2]$. (\ref{eq:basic}) and IBP yield that: 
\begin{align}
|\eqref{eq:int_1}|  
\le C \int_{\mathbb{R}^N} 
\bigl| \partial_\eta^{\beta}
\mathscr{D}_{(J)}^{z, \theta}  
\bigl( \partial_y^\gamma p_{\Delta - s_1}^{\bar{X}^z, (w)} (\cdot, y ;\theta) 
\bigr) (\eta) \bigr| 
\times 
\bigl|  p_{s_1}^{\bar{X}^z, (w)} (\xi, \eta; \theta) \bigr| d \eta.  
\end{align}
Due to (\ref{eq:LDL_density_bd}) in Proposition \ref{prop:density_bd} and (\ref{eq:basic}), for any $\alpha \in \mathbb{Z}_{\ge 0}^N$, 
\begin{align*} 
\bigl| 
\partial_\eta^{\alpha} 
\partial_y^\gamma p_{\Delta - s_1}^{\bar{X}^z, (w)} (\eta, y ;\theta) 
\bigr| \le C (\Delta-s_1)^{-\|\alpha\|_w -\|\gamma \|_w } \times 
% (\Delta-s_1)^{-\|\gamma \|_w} \times 
\overline{\mathscr{G}}^{(w), z} (\Delta - s_1, \eta, y, \theta),  
\end{align*}
and thus, we can apply Lemma \ref{lemma:grad_D} to $\psi^{(w)} (\Delta, \eta) = \partial_y^\gamma p_{\Delta - s_1}^{\bar{X}^z, (w)} (\eta, y ;\theta)$ and obtain: 
\begin{align*}
\bigl| \partial_\eta^{\beta}
\mathscr{D}_{(J)}^{z, \theta}  
\bigl( \partial_y^\gamma p_{\Delta - s_1}^{\bar{X}^z, (w)} (\cdot, y ;\theta) 
\bigr) (\eta)
\bigr| \le C (\Delta - s_1)^{-\|\beta \|_w  -  \|  \gamma\|_w - Q_w (J)}
\times 
\Bigl(1 + \bigl| \tfrac{\eta - z}{\sqrt{\Delta - s_1}} \bigr|^q \Bigr)
\times 
\overline{\mathscr{G}}^{(w), z} (\Delta - s_1, \eta, y, \theta).
\end{align*}
We then have that: 
\begin{align*}
& | \eqref{eq:int_1} |  \\ 
& \le C  (\Delta -s_1)^{- \| \beta + \gamma \|_w - Q_w (J)}   
\times 
\int_{\mathbb{R}^N} 
\Bigl(1 + \bigl| \tfrac{\eta - z}{\sqrt{\Delta - s_1}} \bigr|^q \Bigr) \times 
\overline{\mathscr{G}}^{(w), z} (\Delta -s_1, \eta, y, \theta) \, 
\overline{\mathscr{G}}^{(w), z} (s_1, \xi, \eta, \theta) d \eta 
\quad \bigl( \because \mathrm{Proposition \, \ref{prop:density_bd}} \bigr) \nonumber \\  
& \le C 
\Delta^{- \| \beta + \gamma \|_w - Q_w (J)}   
\times 
\int_{\mathbb{R}^N} 
\Bigl(1 + \bigl| \tfrac{\eta - z}{\sqrt{\Delta}} \bigr|^q \Bigr) \times 
\overline{\mathscr{G}}^{(w), z} (\Delta -s_1, \eta, y, \theta) \, 
\overline{\mathscr{G}}^{(w), z} (s_1, \xi, \eta, \theta) d \eta 
\quad \bigl( \because s_1  \in [0, \Delta/2] \bigr)  \nonumber \\  
& \le C \Delta^{- \| \beta + \gamma \|_w - Q_w (J)}   
\times \Bigl(1 + \bigl| \tfrac{\xi - z}{\sqrt{\Delta}} \bigr|^q \Bigr)  
\, \overline{\mathscr{G}}^{(w), z} (\Delta, \xi, y, \theta) \quad \bigl(  \because \mathrm{Lemmas \, \ref{lemma:conv}, \ref{lemma:poly_bd}} \bigr) 
\end{align*}
where $q \ge 1$ and $C_1, C_2 > 0$ are constants independent of $\Delta, \xi, y, \theta$. \\ 
% Making use of the notation $Q_w (J)$ defined in (\ref{eq:Q}), we have that: for $w \in \{\ref{eq:ellip}, \ref{eq:hypo}\}$, 
% % 
% \begin{align*}
% & \bigl| \partial_\eta^{\beta}
% \mathscr{D}_{(J)}^{z, \theta}  
% \bigl( \partial_y^\gamma p_{\Delta - s_1}^{\bar{X}^z, (w)} (\cdot, y ;\theta) 
% \bigr) (\eta)
% \bigr| \times 
% \bigl|  p_{s_1}^{\bar{X}^z, (w)} (\xi, \eta; \theta) \bigr| \nonumber \\ 
% & \le 
% C_1 (\Delta -s_1)^{- \| \beta + \gamma \|_w - Q_w (J)}   
% \times \Bigl(1 + \bigl| \tfrac{\eta - z}{\sqrt{\Delta - s_1}} \bigr|^q \Bigr) \times 
% \overline{\mathscr{G}}^{(w), z} (\Delta -s_1, \eta, y, \theta) \, 
% \overline{\mathscr{G}}^{(w), z} (s_1, \xi, \eta, \theta)
% \quad \bigl( \because \mathrm{Proposition} \, \ref{prop:density_bd} \bigr) \nonumber \\ 
% & \le C_2 \Delta^{- \| \beta + \gamma \|_w - Q_w (J)}   
% \times \Bigl(1 + \bigl| \tfrac{\xi - z}{\sqrt{\Delta}} \bigr|^q \Bigr) \times 
% \overline{\mathscr{G}}^{(w), z} (\Delta -s_1, \eta, y, \theta) \, 
% \overline{\mathscr{G}}^{(w), z} (s_1, \xi, \eta, \theta). 
% \quad \bigl(\because \mathrm{Lemma} \, \ref{lemma:poly_bd}, \ s_1 \in [0, \Delta/2] \bigr) 
% \end{align*}
% % 
% Therefore, Lemmas \ref{lemma:conv}-\ref{lemma:poly_bd} yield:  
% % 
% \begin{align*}
% \bigl| (\ref{eq:int_1})\bigr| 
% \le C \Delta^{- \| \beta + \gamma \|_w - Q_w (J)}   
% \times \Bigl(1 + \bigl| \tfrac{\eta - z}{\sqrt{\Delta}} \bigr|^q \Bigr)  
% \, \overline{\mathscr{G}}^{(w), z} (\Delta, \xi, y, \theta). 
% \end{align*} 

\noindent 
\textbf{Step I-2.} For the case $s_1 \in [\Delta/2, \Delta]$, we apply IBP to obtain: 
\begin{align*}
|\eqref{eq:int_1}|  
& \le
C \int_{\mathbb{R}^N} 
p_{\Delta -s_1}^{\bar{X}^z, (w)} (\eta, y; \theta)
\times 
\bigl| \partial_{\eta}^{\gamma} \bigl( \mathscr{D}_{(J)}^{z, \theta} \bigr)^\ast 
\{ \partial_{\xi}^\beta p_{s_1}^{\bar{X}^z, (w)} (\xi, \cdot ; \theta) \} (\eta) \bigr| d \eta. 
\end{align*} 
Then the same argument in \textbf{Step I-1} provides the desired bound, and we omit the detailed proof to avoid repetition. 
\\ 

\noindent 
\textbf{Step II.} Consider general $j \in \mathbb{N}$. We now assume that the assertion (\ref{eq:LDL_kernel}) holds for $j = 1, \ldots, i$ for $i \in \mathbb{N}$, and then consider the case $j = i+1$. We again consider the cases $s_1 \in [0, \Delta/2]$ and $s_1 \in [\Delta/2, \Delta]$, separately.  
\\

\noindent 
\textbf{Step II-1. $s_1 \in [0, \Delta/ 2]$.} We have that: 
\begin{align}
& \Bigl| \partial^\beta_\xi  \partial^{\gamma}_y  
{F}^{[J], z, (w)}_{s_{1:i+1}} (\xi, y; \theta)   \Bigr| \nonumber \\ 
& \le  C 
\int_{\mathbb{R}^N} 
\bigl| 
\underbrace{ \partial_{\eta}^\beta 
\mathscr{D}_{(J)}^{z, \theta}
\bigl( 
\partial^{\gamma}_y 
\bar{P}_{s_i - s_{i+1}}^{\theta, z, (w)} 
\widetilde{\mathscr{L}}_\theta^z 
\bar{P}_{s_{i-1} - s_i}^{\theta, z, (w)}  
\cdots 
\widetilde{\mathscr{L}}_\theta^z \bigr)   
\{ p_{\Delta -s_1}^{\bar{X}^z, (w)} (\cdot, y;  \theta) \} (\eta) 
}_{=: (\#1)}
\bigr| 
\times 
\bigl| p_{s_{i+1}}^{\bar{X}^z, (w)} (\xi, \eta ; \theta)
\bigr|  d \eta.  
\label{eq:step2-1-1}
\end{align}
Under the assumption of induction with $J = 0$, we have that: for any $\alpha \in \mathbb{Z}_{\ge 0}^N$, 
\begin{align*}
& \bigl| \partial^{\alpha}_\eta
\partial^{\gamma}_y 
\bar{P}_{s_i - s_{i+1}}^{\theta, z, (w)} 
\widetilde{\mathscr{L}}_\theta^z 
\bar{P}_{s_{i-1} - s_i}^{\theta, z, (w)}  
\cdots 
\widetilde{\mathscr{L}}_\theta^z  
\{ p_{\Delta -s_1}^{\bar{X}^z, (w)} (\cdot, y; \theta) \} (\eta) \bigr|  \\
& \le 
C 
(\Delta - s_{i+1})^{- \| \alpha  + \gamma \|_w 
- \tfrac{i-1}{2} - \tfrac{1}{2}} 
\times \Bigl(1 + 
\bigl| \tfrac{\eta - z}{\sqrt{\Delta - s_{i+1}}} \bigr|^q 
\Bigr) \times 
\overline{\mathscr{G}}^{(w), z} (\Delta - s_{i+1}, \eta, y, \theta). 
\end{align*}
Thus, making use of Lemma \ref{lemma:grad_D} with 
\begin{align*}
\psi^{(w)} (\Delta - s_{i+1}, \eta) 
& \equiv
\partial^{\gamma}_y 
\bar{P}_{s_i - s_{i+1}}^{\theta, z, (w)} 
\widetilde{\mathscr{L}}_\theta^z 
\bar{P}_{s_{i-1} - s_i}^{\theta, z, (w)}  
\cdots 
\widetilde{\mathscr{L}}_\theta^z  
\{ p_{\Delta -s_1}^{\bar{X}^z, (w)} (\cdot, y; \theta) \} (\eta) \\ 
& = 
\partial^{\gamma}_y 
\bar{P}_{s_i - s_{i+1}}^{\theta, z, (w)} 
\widetilde{\mathscr{L}}_\theta^z 
\bar{P}_{s_{i-1} - s_{i+1} - (s_i - s_{i+1})}^{\theta, z, (w)}  
\cdots 
\widetilde{\mathscr{L}}_\theta^z  
\{ p_{\Delta - s_{i+1} - (s_1 - s_{i+1})}^{\bar{X}^z, (w)} (\cdot, y; \theta) \} (\eta),  
\end{align*}
we have that  
\begin{align*}
%\bigl| 
%\partial_{\eta}^\beta 
%\mathscr{D}_{(J)}^{z, \theta}
%\bigl( 
%\partial^{\gamma}_y 
%\bar{P}_{s_i - s_{i+1}}^{\theta, z, (w)} 
%\widetilde{\mathscr{L}}_\theta^z 
%\bar{P}_{s_{i-1} - s_i}^{\theta, z, (w)}  
%\cdots 
%\widetilde{\mathscr{L}}_\theta^z \bigr)   
%\{ p_{\Delta -s_1}^{\bar{X}^z, (w)} (\cdot, y;  \theta) \} (\eta) 
%\bigr| 
%\\ 
| (\#1) | 
& \le C (\Delta - s_{i+1})^{-\| \beta + \gamma \|_w - \tfrac{i}{2} - Q_w (J)} \times 
\Bigl( 1 + \bigl| \tfrac{\eta - z}{\sqrt{\Delta - s_{i+1}}} \bigr|^q 
\Bigr)
\times 
\overline{\mathscr{G}}^{(w), z} (\Delta - s_{i+1}, \eta, y, \theta) \\ 
& \le 
C \Delta^{-\| \beta + \gamma \|_w - \tfrac{i}{2} - Q_w (J)} \times 
\Bigl( 1 + \bigl| \tfrac{\eta - z}{\sqrt{\Delta}} \bigr|^q 
\Bigr)
\times 
\overline{\mathscr{G}}^{(w), z} (\Delta - s_{i+1}, \eta, y, \theta). \qquad \bigl( \because 
0 \le s_{i+1} \le s_1 \le \tfrac{\Delta}{2} \bigr) 
\end{align*}
Therefore, application of Lemmas \ref{lemma:conv}, \ref{lemma:poly_bd} and Proposition \ref{prop:density_bd} to (\ref{eq:step2-1-1}) yield that: for $s_1 \in [0, \Delta/2]$
\begin{align*}
\Bigl| \partial^\beta_\xi  
\partial^{\gamma}_y  
{F}^{[J], z, (w)}_{s_{1:i+1}} (\xi, y; \theta)   \Bigr| 
\le C \Delta^{-\| \beta +  \gamma\|_w 
- Q_w (J) - \tfrac{i}{2}}
\times 
\Bigl( 
1 + \bigl| \tfrac{\xi- z}{\sqrt{\Delta}}\bigr|^q 
\Bigr) 
\times 
\overline{\mathscr{G}}^{(w), z} (\Delta, \xi,  y, \theta). 
\end{align*} 
\textbf{Step II-2}. $s_1 \in [\Delta/2, \Delta]$. IBP together with (\ref{eq:basic}) gives
\begin{align}
& 
\Bigl| \partial^\beta_\xi  \partial^{\gamma}_y  
{F}^{[J], z, (w)}_{s_{1:i+1}} (\xi, y; \theta)   \Bigr| \nonumber \\  
& \le 
C  \int_{\mathbb{R}^N} 
\bigl| p_{\Delta - s_1}^{\bar{X}^z, (w)} (\cdot, y ; \theta) (\xi) 
\bigr| 
\times 
\bigl|
\underbrace{\partial_\eta^\gamma 
(\widetilde{\mathscr{L}}_\theta^{z})^\ast_\eta
\bigl( \partial^\beta_\xi 
\bar{P}_{s_{i+1}}^{\theta, z, (w)} \mathscr{D}_{(J)}^{z, \theta} 
\bar{P}_{s_i - s_{i+1}}^{\theta, z, (w)} 
\widetilde{\mathscr{L}}_\theta^z  
\cdots 
\widetilde{\mathscr{L}}_\theta^z 
p_{s_1 - s_{2}}^{\bar{X}^z, (w)} (\cdot, \eta ; \theta) (\xi) \bigr) }_{=:(\# 2)}
\bigr| d \eta.   
\label{eq:step2-2-1}
\end{align} 
Under the assumption of induction, 
Lemma \ref{lemma:grad_D} with 
$$
\psi^{(w)} (s_1, \eta) = 
\partial^\beta_\xi 
\bar{P}_{s_{i+1}}^{\theta, z, (w)} \mathscr{D}_{(J)}^{z, \theta} 
\bar{P}_{s_i - s_{i+1}}^{\theta, z, (w)} 
\widetilde{\mathscr{L}}_\theta^z  
\cdots 
\widetilde{\mathscr{L}}_\theta^z 
p_{s_1 - s_{2}}^{\bar{X}^z, (w)} (\cdot, \eta ; \theta) (\xi)
$$ gives 
\begin{align}
%& \Bigl| 
%\partial_\eta^\gamma 
%(\widetilde{\mathscr{L}}_\theta^{z})^\ast_\eta
%\bigl( \partial^\beta_\xi 
%\bar{P}_{s_{i+1}}^{\theta, z, (w)} \mathscr{D}_{(J)}^{z, \theta} 
%\bar{P}_{s_i - s_{i+1}}^{\theta, z, (w)} 
%\widetilde{\mathscr{L}}_\theta^z  
%\cdots 
%\widetilde{\mathscr{L}}_\theta^z 
%p_{s_1 - s_{2}}^{\bar{X}^z, (w)} (\cdot, \eta ; \theta) (\xi) \bigr) 
%\Bigr|
%\nonumber  \\ 
| (\# 2) | &  
\le C 
(s_1)^{ - \| \beta +  \gamma \|_w - \tfrac{i-1}{2} - Q_w (J) - \tfrac{1}{2}}
\times \Bigl( 1 + \bigl| \tfrac{\xi -  z}{\sqrt{s_1}} \bigr|^q \Bigr)
\times 
\overline{\mathscr{G}}^{(w), z} (s_1, \xi, \eta, 
\theta)
\nonumber \\ 
& 
\le C 
\Delta^{ - \| \beta +  \gamma \|_w - \tfrac{i}{2} - Q_w (J)}
\times 
\Bigl( 1 + \bigl| \tfrac{\xi -  z}{\sqrt{\Delta}} \bigr|^q 
\Bigr)
\times 
\overline{\mathscr{G}}^{(w), z} (s_1, \xi, \eta, 
\theta). 
\qquad \bigl( \because s_1 \in [\Delta/2, \Delta] \bigr)
\label{eq:step2-2-2}
\end{align}
We thus obtain from (\ref{eq:step2-2-1}) and (\ref{eq:step2-2-2}) that: 
for $s_1 \in [\Delta/2, \Delta]$, 
\begin{align*}
& \Bigl| \partial^\beta_\xi 
\partial^{\gamma}_y  {F}^{[J], z, (w)}_{s_{1:i+1}} (\xi, y; \theta)  
\Bigr|   \\
& \le 
C \Delta^{ - \| \beta + \gamma \|_w - \tfrac{i}{2} - Q_w (J)}
\Bigl( 1 + \bigl| \tfrac{\xi -  z}{\sqrt{\Delta}} \bigr|^q 
\Bigr)
\int_{\mathbb{R}^N}
\overline{\mathscr{G}}^{(w), z} (\Delta - s_1, \eta, y, \theta) 
\, \overline{\mathscr{G}}^{(w), z} (s_1, \xi, \eta, \theta) 
d \eta 
\qquad \bigl( \because \mathrm{Proposition \, \ref{prop:density_bd}} \bigr) 
\\ 
& \le 
C \Delta^{ - \| \beta + \gamma \|_w - \tfrac{i}{2} - Q_w (J)} \times 
\Bigl( 1 + \bigl| \tfrac{\xi -  z}{\sqrt{\Delta}} \bigr|^q 
\Bigr) \times 
\overline{\mathscr{G}}^{(w), z} (\Delta, \xi, y, \theta).  
\qquad \bigl( \because \mathrm{Lemma \, \ref{lemma:conv}} \bigr)  
\end{align*} 
We have shown that the assertion holds with $j = i+1$, and thus the proof of Proposition \ref{prop:bd_E_2} is now complete.  
\qed 
\\

\noindent
\textit{Proof of  Theorem \ref{thm:bd_r2}.}
We first derive an analytic expression for the error term $\mathscr{R}_2$. We recall the definition of $\hat{\mathscr{T}}^{(j)} (\Delta, x, y; \theta)$ as: for $x, y \in \mathbb{R}^N$, $\theta \in \Theta$ and $0 \le s_j \le s_{j-1} \le \cdots \le s_1 \le \Delta$, $1 \le j \le M-1$, 
\begin{align}
\hat{\mathscr{T}}^{\, j, (w)} (\Delta, x, y; \theta) \equiv \bar{P}_{s_j}^{\theta,  z,  (w)} 
\widetilde{\mathscr{L}}_\theta^z 
\bar{P}_{s_{j-1} - s_j}^{\theta, z, (w)} 
\widetilde{\mathscr{L}}_\theta^z 
\cdots 
\widetilde{\mathscr{L}}_\theta^z \bar{P}_{s_1 - s_2}^{\theta, z, (w)} 
\widetilde{\mathscr{L}}_\theta^z 
p_{\Delta-s_1}^{\bar{X}^z, (w)} (\cdot, y; \theta) (x) \bigr|_{z=x}, \qquad w \in \{\ref{eq:ellip}, \ref{eq:hypo} \}. 
 \label{eq:target}
\end{align}
Recursive application of Lemma \ref{lemma:aux2} to (\ref{eq:target}) yields: for any $\beta^{[j]} \in \mathbb{Z}_{\ge 0}^j$, 
\begin{align}
\hat{\mathscr{T}}^{\, j, (w)} (\Delta, x, y; \theta)  
= \sum_{ \alpha \le \beta^{[j]} } 
\tfrac{\prod_{i = 1}^j (s_{j+1-i})^{\alpha_i}}{\alpha!}
\mathscr{D}_\alpha^{z, \theta} 
p_\Delta^{\bar{X}^z, (w)} 
( \cdot, y; \theta) (x) \bigr|_{z=x}
+ \mathscr{E}^{\beta^{[j]}, (w)}_{s_{1:j}} (\Delta, x, y; \theta), 
\end{align} 
where we have defined:  
\begin{align}
\mathscr{E}^{\gamma, (w)}_{s_{1:j}} (\Delta, x, y; \theta) \equiv 
\sum_{1 \le i \le j} \widetilde{\mathscr{E}}^{\gamma, (w)}_i (\Delta, x, y; \theta), \qquad \gamma \in \mathbb{Z}_{\ge  0}^j 
\label{eq:E_def}
\end{align}
with 
\begin{align}
& \widetilde{\mathscr{E}}^{\gamma, (w)}_i (\Delta, x, y; \theta) 
\equiv
\sum_{\substack{ \alpha \le (\gamma_1, \ldots, \gamma_{i-1}) }}
% \sum_{\substack{
%0 \le \alpha_1 \le \beta^{[j]}_1 \\ 
%\vdots \\ 
%0 \le \alpha_{i-1} \le \beta^{[j]}_{i-1} \\ 
%}} 
\biggl\{ 
\tfrac{(s_j)^{\alpha_1} \cdots (s_{j-(i-2)})^{\alpha_{i-1}} (s_{j-(i-1)})^{\gamma_i +1}}{\alpha_1! \cdots \alpha_{i-1}!  \gamma_i!} 
\label{eq:err_i} \\
&  \qquad \quad 
\times \int_0^1 (1-u)^{\gamma_i +1} du \, 
\underbrace{\mathscr{D}_{\alpha}^{z, \theta}  
\bar{P}_{us_{j-(i-1)}}^{\theta, z, (w)} 
\mathscr{D}_{(\gamma_i + 1)}^{z, \theta} 
\bar{P}_{s_{j-i} - us_{j-(i-1)}}^{\theta, z, (w)} 
\widetilde{\mathscr{L}}_\theta^z  
\cdots 
\widetilde{\mathscr{L}}_\theta^z 
p_{\Delta - s_1}^{\bar{X}^z, (w)} (\cdot, y; \theta)(x)
\Bigr|_{z=x}}_{=: \flat_u^{[i], (w)} (\alpha)}  \biggr\}, \quad 1 \le i \le j.  \nonumber 
\end{align} 
%  
%\begin{align}
%\begin{aligned}
%& \widetilde{\mathscr{E}}^{\gamma, (w)}_j (\Delta, x, y; \theta) 
%\\ 
%& \quad 
%\equiv 
%\sum_{\substack{\alpha \in \mathbb{Z}_{\ge 0}^{j-1} \\ \alpha \in \gamma^{(j-1)}} }
%\tfrac{(s_j)^{\alpha_1} \cdots (s_2)^{\alpha_{j-1}} (s_1)^{\gamma_j+1}}{
%\alpha! \, \gamma_j!}
%\int_0^1 (1-u)^{\gamma_j +1} du 
%\, \mathscr{D}_{\alpha}^{z, \theta}  
%\bar{P}_{us_{1}}^{\theta, z, (w)} 
%\mathscr{D}_{(\gamma_j +1)}^{z, \theta}
%p_{\Delta - us_1}^{\bar{X}^z, (w)} (\cdot, y; \theta) (x) |_{z=x}.
%\end{aligned}
%\label{eq:err_j}
%\end{align}
% 
Note that $
\flat_u^{[j], (w)} (\alpha) 
\equiv \mathscr{D}_{\alpha}^{z, \theta}  
\bar{P}_{us_{1}}^{\theta, z, (w)} 
\mathscr{D}_{(\gamma_j +1)}^{z, \theta}
p_{\Delta - us_1}^{\bar{X}^z, (w)} (\cdot, y; \theta) (x) |_{z=x} 
$ when $i = j$.  Then, the remainder term $\mathscr{R}_2$ introduced in (\ref{eq:err}) in the main text is expressed as follows: 
\begin{align}
\mathscr{R}_2^{\,j, \beta^{[j]}, (w)} (\Delta, x, y; \theta) \equiv 
\int_{0 \le s_j \le \cdots \le s_1 \le \Delta}
% \int_0^\Delta \int_0^{s_1} \cdots \int_0^{s_{j-1}} 
\mathscr{E}^{\beta^{[j]}, (w)}_{s_{1:j}}
(\Delta, x, y; \theta) ds_j \cdots ds_1.  
\end{align} 
To bound the residual $\mathscr{R}_2^{\, j, \beta^{[j]}, (w)} (\Delta, x, y; \theta)$, we study the upper bound of (\ref{eq:err_i}). 
%To simplify notation, we write 
%% 
%\begin{align*}
%\flat_u^{[i], (w)} (\alpha)  
%\equiv 
%\mathscr{D}_{\alpha}^{z, \theta}  
%\bar{P}_{us_{j-(i-1)}}^{\theta, z, (w)} 
%\mathscr{D}_{(\gamma_i + 1)}^{z, \theta} 
%\bar{P}_{s_{j-i} - us_{j-(i-1)}}^{\theta, z, (w)} 
%\widetilde{\mathscr{L}}_\theta^z  
%\cdots 
%\widetilde{\mathscr{L}}_\theta^z 
%p_{\Delta - s_1}^{\bar{X}^z, (w)} (\cdot, y; \theta)(x) \bigr|_{z=x}, 
%\qquad 1 \le i \le j, 
%\end{align*}
%% 
%with convention 
%% 
%$ 
%\flat_u^{[j], (w)} (\alpha) 
%\equiv \mathscr{D}_{\alpha}^{z, \theta}  
%\bar{P}_{us_{1}}^{\theta, z, (w)} 
%\mathscr{D}_{(\gamma_j +1)}^{z, \theta}
%p_{\Delta - us_1}^{\bar{X}^z, (w)} (\cdot, y; \theta) (x) |_{z=x} 
%$. 
% 
Given Proposition \ref{prop:bd_E_2},  we can use Lemma \ref{lemma:grad_D_a} with 
$$\psi^{(w)} (\Delta, x)
\equiv \bar{P}_{us_{j-(i-1)}}^{\theta, z, (w)} 
\mathscr{D}_{(\gamma_i + 1)}^{z, \theta} 
\bar{P}_{s_{j-i} - us_{j-(i-1)}}^{\theta, z, (w)} 
\widetilde{\mathscr{L}}_\theta^z  
\cdots 
\widetilde{\mathscr{L}}_\theta^z 
p_{\Delta - s_1}^{\bar{X}^z, (w)} (\cdot, y; \theta)(x)$$
and get 
\begin{align}
\sup_{u \in [0,1]} \bigl| \flat_u^{[i], (w)} (\alpha) \bigr|  
& \le C 
\Delta^{- | \alpha | - \tfrac{i-1}{2}  - \tfrac{j-i}{2} - Q_w(\gamma_i + 1)}
\cdot {\mathscr{G}}^{(w)} (\Delta, x, y, \theta) \nonumber \\
& =  C 
\Delta^{- | \alpha | - \tfrac{j-1}{2} - Q_w(\gamma_i + 1)}
\cdot {\mathscr{G}}^{(w)} (\Delta, x, y, \theta) 
\qquad 1 \le i \le j, 
\label{eq:bd_flat_i}
\end{align}
for some constant $C > 0$ independent of $\Delta, y, \theta, 
\{ s_k \}_{1 \le k \le j-(i-1)}$, where we recall the definition of $Q_w$ and $F_w$  in (\ref{eq:Q}) and (\ref{eq:F}), respectively. 
It follows from (\ref{eq:err_i}) (with $\gamma$ replaced by $\beta^{[j]} \in \mathbb{Z}_{\ge 0}^j$) and (\ref{eq:bd_flat_i}) that: 
\begin{align}
\begin{aligned}
& \bigl| \widetilde{\mathscr{E}}^{\beta^{[j]}, (w)}_i (\Delta, x, y; \theta)  \bigr| 
\nonumber \\ 
% \le C_1 
% \sum_{\substack{ \alpha \in \mathbb{Z}_{\ge 0}^{i-1} \\ \alpha \le \gamma^{(i-1)} }} 
% \tfrac{(s_j)^{\alpha_1} \cdots (s_{j-(i-2)})^{\alpha_{i-1}} (s_{j-(i-1)})^{\gamma_i +1}}{\alpha_1! \cdots \alpha_{i-1}!  \gamma_i!} 
% \sup_{u \in [0,1]} \bigl| \flat_u^{[i], (w)} (\alpha) \bigr| 
% \nonumber \\ 
& \le  C \sum_{\alpha \le 
(\beta^{[j]}_1, \ldots, \beta^{[j]}_{i-1})} 
\tfrac{(s_j)^{\alpha_1} \cdots (s_{j-(i-2)})^{\alpha_{i-1}} (s_{j-(i-1)})^{\beta^{[j]}_i +1}}{\alpha_1! \cdots \alpha_{i-1}!  \beta^{[j]}_i!} 
\cdot 
\Delta^{-  | \alpha | - \tfrac{j-1}{2} - Q_w(\beta^{[j]}_i + 1)}  
\cdot \mathscr{G}^{(w)} (\Delta, x, y, \theta),  
\end{aligned}
\end{align} 
and then 
\begin{align}
& \bigl| \mathscr{R}_2^{\, j, \beta^{[j]}, (w)} (\Delta, x, y; \theta) \bigr| 
\le 	\sum_{1 \le i \le j}  
 \int_{0 \le s_j \le \cdots \le s_1 \le \Delta} 
\bigl| 
\widetilde{\mathscr{E}}^{\beta^{[j]}, (w)}_i (\Delta, x, y; \theta) 
\bigr| 
ds_j \cdots ds_1 
\nonumber \\ 
% \le C_1  \sum_{i = 1}^j  \sum_{\substack{ \alpha \in \mathbb{Z}_{\ge 0}^{i-1} \\ \alpha \le \beta^{[j], (i-1)} }}  
% \Biggl\{ 
% \int_0^\Delta \int_0^{s_1} \cdots \int_0^{s_{j-1}}   
% \tfrac{(s_j)^{\alpha_1} \cdots (s_{j-(i-2)})^{\alpha_{i-1}} (s_{j-(i-1)})^{\beta^{[j]}_i +1}}{\alpha_1! \cdots \alpha_{i-1}!  \beta^{[j]}_i!} 
% \cdot 
% \Delta^{- F_w (
% \alpha) - \tfrac{j-i}{2} - Q_w(\beta^{[j]}_i + 1)} 
% ds_j \cdots ds_1  \Biggr\} 
% \times \mathscr{G}^{(w)} (\Delta, x, y, \theta)  
% \nonumber \\ 
%& \quad \le C_1 
%\sum_{1 \le i \le j} 
%\sum_{ \alpha \le (\beta^{[j]}_1, \ldots, \beta^{[j]}_{i-1}) }
%\Delta^{|\alpha| + (\beta^{[j]}_i + 1) + \tfrac{j + i}{2} - F_w (\alpha) - Q_w (\beta^{[j]}_{i}+1)
%} \times \mathscr{G}^{(w)} (\Delta, x, y, \theta)
%\nonumber \\ 
& \quad 
\le C  \sum_{1 \le i \le  j}  \Delta^{(\beta^{[j]}_i + 1) + \tfrac{j+1}{2} - Q_w (\beta^{[j]}_i  + 1)}  \times \mathscr{G}^{(w)} (\Delta, x, y, \theta) 
%\qquad \bigl(\because F_w (\alpha) \le | \alpha | + \tfrac{i-1}{2}; \, \mathrm{see} \, \eqref{eq:F_bd} \bigr)
\nonumber \\ 
& \quad 
\le C \times 
\begin{cases}
\sum_{1 \le i \le  j} \Delta^{\tfrac{\beta^{[j]}_i}{2} + \tfrac{j}{2}} \times \mathscr{G}^{(\ref{eq:ellip})}  (\Delta, x, y, \theta),  & w = \ref{eq:ellip};  \\ 
 \sum_{1 \le i \le  j} \Delta^{\beta_i^{[j]} - \tfrac{1}{2} 
 \floor{\tfrac{3}{2} (\beta_i^{[j]} + 1)} + \tfrac{j}{2} 
+ \tfrac{1}{2} \times \mathbf{1}_{\beta_i^{[j]} = 0,1}} 
\times \mathscr{G}^{(\ref{eq:hypo})} (\Delta, x, y, \theta), & w = \ref{eq:hypo},  
\end{cases} 
\qquad  \bigl( \because \eqref{eq:Q}  \bigr)  \label{eq:bd_R2}
%\\ 
%= 
%\begin{cases}
%C_3 \sum_{i = 1}^{j} \Delta^{\tfrac{\beta^{[j]}_i}{2} + \tfrac{j}{2}} \times \mathscr{G}^{(\ref{eq:ellip})}  (\Delta, x, y, \theta),  & w = \ref{eq:ellip};  \\ 
%C_3 \sum_{i = 1}^j \Delta^{ \tfrac{1}{2} 
%\bigl( 
%\bigl \lbrack \tfrac{\beta_i^{[j]}}{2}   \bigr \rbrack- 1
%\bigr)  + \tfrac{j}{2} 
%		+ \tfrac{1}{2} \times \mathbf{1}_{\beta_i^{[j]} = 0,1}} 
%	\times \mathscr{G}^{(\ref{eq:hypo})} (\Delta, x, y, \theta), & w = \ref{eq:hypo}. 
% 
\end{align}
where $C > 0$ is a  constant independent of $\Delta, y$ and $\theta$. 
For the case of $w = \ref{eq:hypo}$ in (\ref{eq:bd_R2}), we have further 
\begin{align*}
\beta_i^{[j]} - \tfrac{1}{2} 
\floor{\tfrac{3}{2} (\beta_i^{[j]} + 1) } 
+ \tfrac{j}{2} 
+ \tfrac{1}{2} \times \mathbf{1}_{\beta_i^{[j]} = 0,1} 
= \tfrac{1}{2} 
\bigl( 
\floor{\tfrac{\beta_i^{[j]}}{2}} 
- 1
\bigr)  + \tfrac{j}{2} 
+ \tfrac{1}{2} \times \mathbf{1}_{\beta_i^{[j]} = 0,1}  = 
\tfrac{1}{2} 
\bigl( 
\floor{\tfrac{\beta_i^{[j]}}{2}}
- \mathbf{1}_{\beta_i^{[j]} \ge 2}
\bigr)  + \tfrac{j}{2},  
\end{align*}
and the proof of Theorem \ref{thm:bd_r2} is now complete. \qed

\section{Expression of the differential operator} \label{app:DO} 
In this section, we study the expression of the operator 
\begin{align*}
\mathscr{D}_\alpha^{z, \theta}
\equiv \prod_{1 \le i \le  j} \bigl(\mathrm{ad}^{\alpha_i}_{\mathscr{L}_\theta^z} 
(\widetilde{\mathscr{L}}_\theta^z) \bigr), \qquad \alpha \in \mathbb{Z}_{\ge 0}^j, \ j \in \mathbb{N}, \  (z, \theta) \in \mathbb{R}^N \times \Theta 
\end{align*} 
introduced in (\ref{eq:D}) in the main text, where we recall: 
\begin{itemize}
	\item  $\widetilde{\mathscr{L}}_\theta^z = \mathscr{L}_\theta - \mathscr{L}_\theta^{0, z}$ with  $\mathscr{L}_\theta$ being the generator of SDE and $\mathscr{L}_\theta^{0, z}$ the generator of LDL scheme (\ref{eq:LDL_frozen}) defined as: 
	\begin{align} \label{eq:L_0_z}
		\mathscr{L}_\theta^{0, z} 
		&
		=  \sum_{1 \le i \le N} 
		\bigl[A_{z, \theta} x + b_{z, \theta} \bigr]_i 
		\partial_i  + \tfrac{1}{2} \sum_{1 \le i_1, i_2 \le N} \sum_{1 \le j \le d} V_j^{i_1} (z, \theta) V_j^{i_2} (z, \theta) 
		\partial_{i_1 i_2}. 
	\end{align}  
  \item For differential operators $D_1, D_2$,  
  $
  \mathrm{ad}_{D_1} (D_2) = [D_1, D_2] 
  \equiv D_1 D_2 - D_2 D_1. 
  $ 
  The $k$-times iteration of the commutator writes as
  $\mathrm{ad}_{D_1}^k (D_2) = [D_1, \mathrm{ad}_{D_1}^{k-1}, (D_2)]$, $k\ge 1$,  
  %$   \mathrm{ad}_{D_1}^k (D_2) = \underbrace{[D_1, [D_1, \ldots [D_1}_{k-\mathrm{times}}, D_2]]]$.    
  % \begin{gather*}
  	%     \mathrm{ad}_{D_1}^k (D_2) = \underbrace{[D_1, [D_1, \ldots [D_1}_{k-\mathrm{times}}, D_2]]]  
  	% \end{gather*}  
  % 
  with $\mathrm{ad}_{D_1}^0 (D_2)= D_2$.  
\end{itemize}

%which plays a key role in the construction of the CF expansion (\ref{eq:cf}) and obtaining the series expansion (\ref{eq:series}) as explained in Section \ref{sec:series}. 
To obtain the expression, we will proceed with this section as follows:
 \\

\noindent 
\textit{Step 1.} We derive the expression of $\mathscr{D}_{(J)}^{z, \theta} \equiv \mathrm{ad}^J_{\mathscr{L}_\theta^{0, z}} (\widetilde{\mathscr{L}}_\theta^z), \, J \in \mathbb{N}$, in Section \ref{app:diff_step_1}. \\ 
\noindent 
\textit{Step 2.} We derive the expression of $\mathscr{D}_\alpha^{z, \theta}$ in Section \ref{app:diff_step_2}. 
\\ 

\noindent
Throughout Section \ref{app:DO}, the following rule about the commutator of differential operators is critical: 

\begin{claim} \label{claim:comm}
Let $D_1$ and $D_2$ be linear differential operators defined as follows: for $\varphi \in C^\infty (\mathbb{R}^N, \mathbb{R})$, 
\begin{align}
    D_i \varphi (x) \equiv \sum_{ \substack{ \alpha \in \mathbb{Z}_{\ge 0}^N,  \, |\alpha| \le k_i }} h_\alpha^{[i]} (x) \partial^\alpha \varphi (x), 
    \qquad k_i \in \mathbb{N}, \ \  i = 1, 2,  
\end{align}    
with $h_\alpha^{[i]} \in C^\infty (\mathbb{R}^N, \mathbb{R})$. Then the commutator of $D_1$ and $D_2$, i.e. $[D_1, D_2]$, is a linear differential operator of at most order $k_1 + k_2 -1$. 
\end{claim}

The above claim is easily verified by noticing that the terms involving $\partial^\gamma, \, |\gamma| = k_1 + k_2$, are cancelled out due to the definition of commutator $[D_1, D_2] \equiv D_1 D_2 - D_2 D_1$. 
%Thus, the front differential operator in the products must act on the coefficients of the next operator at least once to generate non-zero coefficients $h_\alpha^{[i]}$, which leads to the above claim. 
Also, we notice that if all of the coefficients in $D_1$ and $D_2$ are constant, then $\lbrack D_1, D_2 \rbrack \varphi (x) = 0$. 
% An auxiliary result used in Step 1 is provided as Lemma \ref{lemma:ad1} in Section \ref{app:diff_aux}.  
% We use slight abuse of notation for the floor function as $\lbrack \alpha \rbrack = 0$ for $\alpha < 0$. 
% 
% Also, for the model class (\ref{eq:hypo}) and $\beta = (\beta_R, \beta_S) \in \mathbb{Z}_{\ge 0}^N \equiv \mathbb{Z}_{\ge 0}^d \times \mathbb{Z}_{\ge 0}^d$, 
% % 
% \begin{align}
% \| \beta \|_R \equiv | \beta_R |, \qquad   
% \| \beta \|_S \equiv | \beta_S |.  
% \end{align}
%
\subsection{Step 1. } 
\label{app:diff_step_1}

We introduce a set of multi-indices: for $J  \in \mathbb{Z}_{\ge 0}$, 
\begin{align}
\mathscr{I}_{w}(J)
= 
\begin{cases}
\bigl\{ \gamma \in \mathbb{Z}_{\ge 0}^N \ \bigl| \  1 \le | \gamma |  \le
J +2   \bigr\},  & w = \ref{eq:ellip};  \\[0.2cm] 
\bigl\{ \gamma  = (\gamma_S, \gamma_R)\in \mathbb{Z}_{\ge 0}^{N_S} \times \mathbb{Z}_{\ge 0}^{N_R} \ \bigl| \  1 \le | \gamma_R |  + 2 | \gamma_S | \le
J +2 \bigr\},  & w = \ref{eq:hypo}.    
\end{cases} 
\label{eq:set_I}
\end{align}
We also introduce the following set of functions: for the constant $\kappa > 0$ defined in Assumption \ref{ass:coeff}, 
\begin{align} \label{eq:S_lambda}
\mathscr{S}_\kappa \equiv 
\Bigl\{ f: \mathbb{R}^N \times \mathscr{Z}_\kappa \times \Theta \to \mathbb{R}  
\ \bigl|  \ \forall  \alpha \in \mathbb{Z}_{\ge 0}^N, \, \exists C > 0  \quad 
\mathrm{s.t.} \, \sup_{(\xi, z, \theta) \in \mathbb{R}^N \times \mathscr{Z}_\kappa \times \Theta } 
| \partial^\alpha_\xi f (\xi, z, \theta) | < C
\Bigr\},  
\end{align}
where we recall the definition of $\mathscr{Z}_\kappa$ in (\ref{eq:Z}). We make use of the notation 
\begin{align}
a_{i_1 i_2} (x, \theta) \equiv 
\bigl[ \sigma (x, \theta) \sigma (x, \theta)^\top\bigr]_{i_1 i_2}, \qquad (x, \theta)  \in \mathbb{R}^N \times \Theta, 
\end{align}
for $1 \le i_1, i_2 \le N$ with $\sigma \equiv [V_1, \ldots, V_d]$. Note that under model class (\ref{eq:hypo}), the matrix takes $0$ when either $i_1$ or $i_2$ is less than or equal to $N_S$.    
\begin{lemma} \label{lemma:diff_step_1}
Let Assumptions \ref{ass:diff}, \ref{ass:param} and \ref{ass:coeff} hold. Let $J \in \mathbb{N}$, $\varphi \in C^\infty (\mathbb{R}^N, \mathbb{R})$ and  $\kappa > 0$ be the constant defined in Assumption \ref{ass:coeff}. Also, recall the definition of $\mathscr{Z}_\kappa$ in (\ref{eq:Z}). 
\begin{enumerate}
\item For elliptic model \eqref{eq:ellip}, it holds that:
\begin{align} \label{eq:ad_exp}
\mathscr{D}_{(J)}^{z, \theta} \varphi (\xi) 
%\equiv \mathrm{ad}^{J}_{\mathscr{L}_{\theta}^{0, z}} \bigl( \widetilde{\mathscr{L}}_\theta^z \bigr) \varphi (\xi) 
= 
\sum_{\gamma \in \mathscr{I}_{\ref{eq:ellip}} (J)}
\mathscr{W}^{ z}_\gamma (\xi, \theta) \partial^\gamma \varphi (\xi),  
\qquad (\xi, z, \theta) \in \mathbb{R}^N \times \mathscr{Z}_\kappa \times \Theta, 
\end{align} 
where the coefficient $\mathscr{W}^{z}_\gamma : \mathbb{R}^N \times \Theta \to \mathbb{R}$ is explicitly given and takes the following form: 
\begin{align} \label{eq:W}
	\mathscr{W}^{z}_\gamma (\xi, \theta) 
	= h_{\gamma} (\xi,  z, \theta) + \sum_{i = 1}^N h_{\gamma}^i  (\xi, z, \theta) (\xi_i - z_i) 
	+ \sum_{i_1, i_2  = 1}^N h_{\gamma}^{i_1i_2} (\xi,  z, \theta) (\xi_{i_1} - z_{i_1}) (\xi_{i_2} - z_{i_2}),    
\end{align} 
for some $h_{\gamma}, h_{\gamma}^{i} ,h_{\gamma}^{i_1i_2} \in \mathscr{S}_\kappa$. 
%In particular, for any $\gamma \in \mathscr{I}_{\ref{eq:ellip}} (0)$, 
%% 
%\begin{align*}
%	h_{\gamma} (\xi, z , \theta) \equiv 0, \qquad \forall (\xi, z, \theta) \in \mathbb{R}^N \times \mathscr{Z}_\kappa \times \Theta. 
%\end{align*}
% 
% 
\item For hypo-elliptic model \eqref{eq:hypo}, it holds that:  
\begin{align} \label{eq:ad_exp_H}
\mathscr{D}_{(J)}^{z, \theta} \varphi (\xi) 
=  
\sum_{\gamma \in \mathscr{I}_{\ref{eq:hypo}} (J)}
\mathscr{W}^{z}_\gamma (\xi, \theta) \partial^\gamma \varphi (\xi),   
\qquad (\xi, z, \theta) \in \mathbb{R}^N \times \mathscr{Z}_\kappa \times \Theta, 
\end{align} 
where the coefficient $\mathscr{W}^{z}_\gamma : \mathbb{R}^N \times \Theta \to \mathbb{R}$ is explicitly given and takes the form of (\ref{eq:W}) with the function $h_\gamma$ satisfying:  
\begin{align*}
	h_{\gamma}  (\xi, z , \theta) \equiv 0, \qquad \forall (\xi, z, \theta) \in \mathbb{R}^N \times \mathscr{Z}_\kappa \times \Theta  
\end{align*} 
for any $\gamma \in \mathscr{I}_{\ref{eq:hypo}} (J), \, J = 1, 2$, satisfying $\| \gamma \|_{\ref{eq:hypo}} \, \bigl( \equiv \tfrac{3}{2} |\gamma_S| + \tfrac{1}{2} |\gamma_R| \bigr) =  \tfrac{1}{2} 
\floor{\tfrac{3J}{2}}
+ \tfrac{3}{2}$. 
%  
%\begin{align*}
%    \mathscr{W}_\gamma^{[0], z} (x, \theta) |_{z=x} = 0, 
%\end{align*}% 
%and for $J = 1, 2,$
%% 
%\begin{align}
%    \mathscr{W}^{[J], z}_\gamma (x, \theta) |_{z=x} = 0,   
%    \label{eq:zero}  
%\end{align} 
%if $\gamma \in \mathscr{I}_{\ref{eq:hypo}} (J)$ satisfies $\| \gamma \|_{\ref{eq:hypo}} = \tfrac{1}{2} \lbrack \tfrac{3J}{2} \rbrack
%+ \tfrac{3}{2}$. 
\end{enumerate}
\end{lemma} 

Before proceeding to proving Lemma \ref{lemma:diff_step_1}, we introduce the following auxiliary result: 

\label{app:diff_aux}
\begin{lemma} \label{lemma:ad1}
Let $\varphi \in C^3 (\mathbb{R}^N, \mathbb{R})$ and $(z, \xi, \theta) \in \mathbb{R}^N \times \mathbb{R}^N \times \Theta$. It holds that:   
\begin{align}
 & \bigl( \mathrm{ad}_{\mathscr{L}_{\theta}^{0, z}} (\widetilde{\mathscr{L}}_\theta^z) \bigr) \varphi (\xi)
% \equiv [\mathscr{L}_\theta^{0, z}, \widetilde{\mathscr{L}}^z_\theta] \varphi (\xi) \nonumber \\
= \tfrac{1}{2} \sum_{1 \le i_1, i_2, i_3 \le N} a_{i_1 i_2} (z, \theta) \partial_{i_1 i_2} V_0^{i_3} (\xi, \theta) \partial_{i_3} \varphi (\xi) \nonumber \\ 
& \qquad 
+ \sum_{1 \le i_1, i_2 \le N} 
\Bigl( \bigl\{V_0^{i_1} (z, \theta)  + \partial_z^\top V_0^{i_1} (z, \theta) \cdot (\xi - z) \bigr\}
\partial_{i_1} V_{0}^{i_2} (\xi, \theta)
- V_0^{i_1} (\xi, \theta) \partial_{i_1} V_0^{i_2} (z, \theta) 
\Bigr)  \partial_{i_2} \varphi (\xi)
\nonumber \\
& \qquad 
+ \tfrac{1}{2}
\sum_{1 \le i_1, i_2, i_3 \le N} 
\bigl\{ V_0^{i_1} (z, \theta) + \partial_z^\top V_0^{i_1} (z, \theta) \cdot (\xi - z) \bigr\}  \partial_{i_1} a_{i_2 i_3} (\xi, \theta) \partial_{i_2 i_3} \varphi (\xi)
\nonumber \\ 
& \qquad 
+ \sum_{1 \le i_1, i_2, i_3 \le N} 
\bigl\{ a_{i_1 i_2} (z, \theta) \partial_{i_1} V_0^{i_3} (\xi, \theta) 
- a_{i_1 i_2} (\xi, \theta)  \partial_{i_1} V_0^{i_3} (z, \theta)  \bigr\}   \partial_{i_2 i_3} \varphi (\xi) 
\label{eq:key_term}\\ 
& 
\qquad 
+ \tfrac{1}{4} \sum_{1 \le i_1, i_2, i_3, i_4 \le N} 
a_{i_1, i_2} (z, \theta) \partial_{i_1 i_2} a_{i_3 i_4} (\xi, \theta) \partial_{i_3 i_4} \varphi (\xi)  
\nonumber \\
& 
\qquad 
+ \tfrac{1}{2} \sum_{1 \le i_1, i_2, i_3, i_4 \le N} 
a_{i_1, i_2} (z, \theta) \partial_{i_1} a_{i_3 i_4} (\xi, \theta) \partial_{i_2, i_3 i_4} \varphi (\xi), \nonumber 
\end{align} 
and then 
\begin{align*}
& \bigl( \mathrm{ad}_{\mathscr{L}_{\theta}^{0, z}} (\widetilde{\mathscr{L}}_\theta^z) \bigr) \varphi (\xi) 
|_{z = x}  
= 
\tfrac{1}{2} \sum_{1 \le i_1, i_2, i_3 \le N} a_{i_1 i_2} (x, \theta) \partial_{i_1 i_2} V_0^{i_3} (x, \theta) \partial_{i_3} \varphi (x) 
\\ 
& \qquad 
+ \tfrac{1}{2}
\sum_{1 \le i_1, i_2, i_3 \le N} V_0^{i_1} (x, \theta)  \partial_{i_1} a_{i_2 i_3} (x, \theta) \partial_{i_2 i_3} \varphi (x) 
+ \tfrac{1}{4} \sum_{1 \le i_1, i_2, i_3, i_4 \le N}   
a_{i_1, i_2} (x, \theta) \partial_{i_1 i_2} a_{i_3 i_4} (x, \theta) \partial_{i_3 i_4} \varphi (x)  \\
&  \qquad 
+ \tfrac{1}{2} \sum_{1 \le i_1, i_2, i_3, i_4 \le N} 
a_{i_1, i_2} (x, \theta) \partial_{i_1} a_{i_3 i_4} (x, \theta) \partial_{i_2 i_3 i_4} \varphi (x).  
\end{align*}
\end{lemma}

\noindent 
\textit{Proof of Lemma \ref{lemma:ad1}.} This is obtained via straightforward algebraic computation of differential operators.  
%and we omit the detailed calculation. 
\qed  \\ 

\noindent
\textit{Proof of Lemma \ref{lemma:diff_step_1}.} 
We show only the hypo-elliptic model (\ref{eq:hypo}) case because the elliptic case follows a similar argument. 
% We show the assertion by separately considering the model classes \ref{eq:ellip} and \ref{eq:hypo}.  
\\ 

% \noindent  
% \textit{Proof for elliptic model (\ref{eq:ellip})}. Since both $\mathscr{L}_\theta^{0, z}$ and $\widetilde{\mathscr{L}}_\theta^z := \mathscr{L}_\theta - \mathscr{L}_\theta^{0, z}$ are second order differential operators, we have from the above claim that: 
% % 
% \begin{align}
% \bigl( \mathrm{ad}_{\mathscr{L}_{\theta}^{0, z}} (\widetilde{\mathscr{L}}_\theta^z) \bigr) \varphi (\xi) 
% \equiv [\mathscr{L}_\theta^{0, z}, \widetilde{\mathscr{L}}_\theta^z] \varphi (\xi) = \sum_{ \substack{ \gamma \in \mathbb{Z}_{\ge 0}^N \\
% |\gamma| \le 3 }} \mathscr{W}_{\gamma}^{z} (\xi, \theta) \partial^\gamma \varphi (\xi), \quad (z, \xi, \theta) \in \mathbb{R}^N \times \mathbb{R}^N \times \Theta,   
% \end{align} 
% % 
% where $\mathscr{W}_\gamma^{z} : \mathbb{R}^N \times \Theta \to \mathbb{R}$ is given in a closed-form. Iterative use of Claim \ref{claim:comm} leads to the assertion for $\bigl( \mathrm{ad}^J_{\mathscr{L}_{\theta}^{0, z}} (\widetilde{\mathscr{L}}_\theta^z) \bigr) \varphi (\xi)$ with a general $J \ge 2$. 
% \\ 

The model class \eqref{eq:hypo} requires a careful treatment because the differential operator contains the differentiation w.r.t. the smooth and components, which produces $\Delta^{- M}$ with $M$ varying for components when acting on the transition density of LDL scheme (\ref{eq:LDL})  (recall Lemma \ref{lemma:deriv_LDL} in the main text). We first notice that for both second order differential operators $\mathscr{L}_\theta^{0, z}$ and $\widetilde{\mathscr{L}}_\theta^z := \mathscr{L}_\theta - \mathscr{L}_\theta^{0, z}$, the differentiation of the smooth component is contained only in the first order differential part, and the second order derivatives are all taken w.r.t. the rough components. Precisely, for $\mathscr{L}_\theta^{0, z}$ and $\mathscr{L}_\theta$, the second order derivatives are given as: 
\begin{align}
\tfrac{1}{2} \sum_{N_S + 1 \le i_1, i_2 \le N} a_{i_1 i_2} (z, \theta) \partial_{\xi_{i_1}} \partial_{\xi_{i_2}}, 
\qquad  
\tfrac{1}{2} \sum_{N_S + 1 \le i_1, i_2 \le N} a_{i_1 i_2} (\xi, \theta) \partial_{\xi_{i_1}} \partial_{\xi_{i_2}}, 
\end{align} 
respectively, and do not involve the derivatives w.r.t. the smooth components.  We then study the following three cases separately: (i) $J = 1$, (ii) $J = 2$, (iii) $J \ge 3$.
\\ 

\noindent 
\textit{Proof for Case (i). $J=1$.} We have from Lemma \ref{lemma:ad1} that: 
\begin{align} \label{eq:J1}
\bigl( \mathrm{ad}_{\mathscr{L}_{\theta}^{0, z}} (\widetilde{\mathscr{L}}_\theta^z) \bigr) \varphi (\xi) 
& = \sum_{\substack{ \gamma  \in \mathbb{Z}_{\ge 0}^{N} \\ 
1 \le |\gamma_R| + 2 |\gamma_S| \le 3}} 
\mathscr{W}^{[1], z}_\gamma (\xi, \theta) \partial^\gamma \varphi (\xi),  
\end{align}
for some $\mathscr{W}_\gamma^{[1], z} : \mathbb{R}^N \times \Theta \to \mathbb{R}$ characterised as (\ref{eq:W}). In particular, for $\mathscr{W}_\gamma^{[1], z}$ with $\gamma \in \mathbb{Z}_{\ge 0}^{N}$ satisfying  
$|\gamma_S| = |\gamma_R| =1$, i.e., $\| \gamma \|_{\ref{eq:hypo}} = 2$, which corresponds to terms involved in term (\ref{eq:key_term}), we have 
\begin{align} \label{eq:W1_zero}
\mathscr{W}_\gamma^{[1], z} (\xi, \theta) 
= \sum_{i = 1}^N g_\gamma^i  (\xi, z, \theta) (\xi_i - z_i), 
\end{align} 
for some $g_\gamma^i \in \mathscr{S}_\kappa$, where we have performed Taylor expansion around $\xi = z$ for the following terms in (\ref{eq:key_term}): 
$ 
a_{i_1 i_2} (z, \theta) \partial_{i_1} V_0^{i_3} (\xi, \theta) - a_{i_1 i_2} (\xi, \theta) \partial_{i_1} V_0^{i_3} (z, \theta), \; N_S + 1 \le i_1, i_2 \le N, \; 1 \le i_3 \le N.    
$
Therefore, the assertion holds for $J = 1$. 
\\ 

\noindent 
\textit{Proof for Case (ii). $J=2$.}
Notice from $ \eqref{eq:J1}$ that: 
\begin{align}
\bigl( \mathrm{ad}_{\mathscr{L}_{\theta}^{0, z}} (\widetilde{\mathscr{L}}_\theta^z) \bigr) \varphi (\xi) 
& = \Biggl\{ 
\sum_{\substack{ \gamma  \in \mathbb{Z}_{\ge 0}^{N} \\ 
|\gamma_S|= 0, \, 1 \le |\gamma_R| \le 3}}
+ 
\sum_{\substack{ \gamma  \in \mathbb{Z}_{\ge 0}^{N} \\ 
|\gamma_S|= 1, \, 0 \le |\gamma_R| \le 1}} 
\Biggr\}
\mathscr{W}^{[1], z}_\gamma (\xi, \theta) \partial^\gamma \varphi (\xi).  
\end{align}
Then we have from Claim \ref{claim:comm} and (\ref{eq:L_0_z}) that:   
\begin{align}
& \bigl( \mathrm{ad}^2_{\mathscr{L}_{\theta}^{0, z}} (\widetilde{\mathscr{L}}_\theta^z) \bigr) \varphi (\xi)  
\equiv 
\mathscr{L}_\theta^{0, z} \bigl( \mathrm{ad}_{\mathscr{L}_{\theta}^{0, z}} (\widetilde{\mathscr{L}}_\theta^z) 
\bigr) \varphi (\xi) 
- \bigl( \mathrm{ad}_{\mathscr{L}_{\theta}^{0, z}} (\widetilde{\mathscr{L}}_\theta^z) 
\bigr) \mathscr{L}_\theta^{0, z} \varphi (\xi) 
\nonumber \\[0.2cm]  
& \qquad \quad 
\begin{aligned} 
& = 
\Biggl\{ 
\sum_{\substack{ \gamma  \in \mathbb{Z}_{\ge 0}^{N} \\ 
|\gamma_S|= 0, \, 1 \le |\gamma_R| \le 4}}
+ 
\sum_{\substack{ \gamma  \in \mathbb{Z}_{\ge 0}^{N} \\ 
|\gamma_S|= 1, \, 0 \le |\gamma_R| \le 2}} 
\Biggr\} 
\mathscr{W}^{[2], z}_\gamma (\xi, \theta) \partial^\gamma \varphi (\xi) 
\\[0.2cm] 
& \qquad 
- 
\Biggl\{ 
\sum_{ \substack{ \gamma \in \mathbb{Z}_{\ge 0}^{N} \\ 
|\gamma_S| = 0, \, 1 \le |\gamma_R| \le 3}}
+ 
\sum_{\substack{ \gamma \in   \mathbb{Z}_{\ge 0}^{N} \\ |\gamma_S| = 1, |\gamma_R | =  0, 1}}  
\Biggr\}  
\sum_{\substack{ 1 \le k_1 \le N \\ e_{k_1} \le \gamma }} 
\sum_{1 \le k_2 \le N} 
\mathscr{W}^{[1], z}_\gamma (\xi,  \theta)
\lbrack A_{z, \theta} \rbrack_{k_2 k_1}
% \partial_{z_{k_1}} V_0^{k_2} (z, \theta) 
\partial^{\gamma - e_{k_1} + e_{k_2}} \varphi (\xi), 
\end{aligned} \label{eq:J2}
\end{align} 
for some $\mathscr{W}^{[2], z}_\gamma: \mathbb{R}^N \times \Theta \to \mathbb{R}$ specified as (\ref{eq:W}) and $\mathscr{W}^{[1], z}_\gamma$ is defined as (\ref{eq:J1}). Due to  (\ref{eq:W1_zero}), for any $\gamma \in \mathbb{Z}_{\ge 0}^N$ satisfying $|\gamma_S| = |\gamma_R| = 1$, it holds that   
\begin{align}
\sum_{\substack{ 1 \le k_1 \le N \\ e_{k_1} \le \gamma}} 
\sum_{1 \le k_2 \le N} 
\mathscr{W}^{[1], z}_\gamma (\xi, \theta)
\, \lbrack A_{z, \theta} \rbrack_{k_2 k_1} \, 
= \sum_{\substack{ 1 \le k_1 \le N \\ e_{k_1} \le \gamma}} 
\sum_{1 \le k_2 \le N} 
\sum_{1 \le i \le N} g_\gamma^i (\xi, z, \theta) (\xi_i - z_i) \lbrack A_{z, \theta} \rbrack_{k_2 k_1} 
\end{align} 
with $g_\gamma^{i} \in \mathscr{S}_\kappa$ defined in \eqref{eq:W1_zero}. Then the term (\ref{eq:J2}) is expressed as:  
\begin{align}
\eqref{eq:J2} 
= 
\sum_{\substack{ \beta  \in \mathbb{Z}_{\ge 0}^{N}  
\\  1 \le |\beta_R|  + 2 | \beta_S | \le 4}} 
\mathscr{W}_\beta^{[2], z} (\xi, \theta) \partial^\beta \varphi (\xi),   
\end{align}
where $\mathscr{W}^{[2], z}_\beta : \mathbb{R}^N \times \Theta \to \mathbb{R}$ has the same property as $\mathscr{W}$ in Lemma \ref{lemma:diff_step_1}-(2). 
% 
% \begin{align}
% G_\beta^{[2], z} (x, \theta) |_{z=x} = 0 
% \end{align} 
% % 
% for $\beta \in \mathbb{Z}_{\ge 0}^N$ s.t. $|\beta_S| = 2$, $|\beta_R| = 0$, i.e., $\|\beta\|_{\ref{eq:hypo}} = 3$. Thus, $ \bigl( \mathrm{ad}^2_{\mathscr{L}_{\theta}^{0, z}} (\widetilde{\mathscr{L}}_\theta^z) \bigr) \varphi (x) $ is expressed as (\ref{eq:ad_exp_H}) with the coefficients specified as in Lemma \ref{lemma:diff_step_1}-(2).   
\\ 

\noindent
\textit{Proof for Case (iii). $J \ge 3$}.
We exploit the mathematical induction. 
% The case $J = 3$ is shown via the same argument used in the case (ii). 
We assume that the assertion holds for $J = 3, 4, \ldots, K$, $K \ge 3$ and consider $J = K+1$. We write  
$ 
    \mathscr{L}_\theta^{0, z}
     \equiv \mathscr{L}_\theta^{0, [1], z} + \mathscr{L}_\theta^{0, [2], z} 
$ 
with 
\begin{align*}
\mathscr{L}_\theta^{0, [1], z} 
= \sum_{1 \le i \le N} \bigl\{ V_0^i (z, \theta) + 
% \partial_z^\top V_0^i (z, \theta) \cdot (\xi -z) \bigr\} 
\sum_{1 \le j \le N} [A_{z, \theta}]_{ij} 
(\xi_j - z_j)  \bigr\} 
\partial_{\xi_i},  \quad 
\mathscr{L}_\theta^{0, [2], z}  
= \tfrac{1}{2} \sum_{N_S + 1 \le k_1, k_2 \le N} 
a_{k_1 k_2} (z, \theta)
\partial_{\xi_{k_1}} \partial_{\xi_{k_2}}. 
\end{align*} 
% 
% We introduce: 
% % 
% \begin{align}
% \mathcal{S}_K = 
% \Bigl\{ 
% \gamma \in \mathbb{Z}_{\ge 0}^N \, 
% | \,
% \| \gamma \|_\# \le K + 1  - \tfrac{1}{2} 
% \bigl\lbrack \tfrac{K-1}{2} \bigr\rbrack, \, |\gamma_R|  + 2 |\gamma_S| \le K+2 
% \Bigr\}, \qquad K \in \mathbb{N}.  
% \end{align} 
% 
From the assumption of the induction and Claim \ref{claim:comm},   
\begin{align*}
& \bigl( \mathrm{ad}^{K+1}_{\mathscr{L}_{\theta}^{0, z}} (\widetilde{\mathscr{L}}_\theta^z) \bigr) \varphi (\xi)  
%= \mathscr{L}_\theta^{0, z} \bigl( \mathrm{ad}^{K}_{\mathscr{L}_{\theta}^{0, z}} (\widetilde{\mathscr{L}}_\theta^z) \bigr) \varphi (\xi) 
%- \bigl( \mathrm{ad}^{K}_{\mathscr{L}_{\theta}^{0, z}} (\widetilde{\mathscr{L}}_\theta^z) \bigr) \mathscr{L}_\theta^{0, z} \varphi (\xi)
= \bigl[  \mathscr{L}_{\theta}^{0, z},  \mathrm{ad}^{K}_{\mathscr{L}_{\theta}^{0, z}} (\widetilde{\mathscr{L}}_\theta^z)  \bigr]
= T_1 + T_2 + T_3, 
\end{align*} 
% % 
where we have set: 
\begin{align}
T_1 & \equiv 
%\sum_{\substack{ \gamma \in \mathbb{Z}_{\ge 0}^N \\ 1 \le |\gamma_R| + 2 |\gamma_S| \le K + 2}} 
\sum_{\gamma \in \mathscr{I}_{\ref{eq:hypo}}(K)} 
\bigl( \mathscr{L}_\theta^{0, [1], z}
\{ \mathscr{W}_\gamma^{[K], z} (\cdot, \theta) \} (\xi)
\bigr)  \partial^\gamma \varphi (\xi); 
\label{eq:T1} \\ 
T_2 & \equiv 
\sum_{\gamma \in \mathscr{I}_{\ref{eq:hypo}}(K)} 
\Bigl( 
 \mathscr{L}_\theta^{0, [2], z} \bigl\{ \mathscr{W}_\gamma^{[K], z} (\cdot, \theta) \partial^\gamma \varphi (\cdot) \bigr\} (\xi) 
-  
 \mathscr{W}_\gamma^{[K], z} (\xi, \theta)  \mathscr{L}_\theta^{0, [2], z}  \partial^\gamma 
  \varphi (\xi)  
\Bigr);  \\ 
T_3 & \equiv   
-\sum_{\gamma \in \mathscr{I}_{\ref{eq:hypo}}(K)}  
\sum_{ \substack{ 1 \le k_1 \le N \\ 
e_{k_1} \le \gamma }} 
\sum_{1 \le k_2 \le N} \mathscr{W}_\gamma^{[K], z} (\xi, \theta) \times \lbrack A_{z, \theta} \rbrack_{k_2 k_1} 
\partial^{\gamma - e_{k_1} + e_{k_2}} \varphi (\xi). 
\end{align} 
We rewrite the expression of $T_i, \, i = 2, 3$. Due to Claim \ref{claim:comm}, the second term $T_2$ writes: 
\begin{align} \label{eq:T2}
T_2 = 
\sum_{\beta \in \mathscr{I}_{\ref{eq:hypo}} (K + 1)} 
\mathscr{W}_\beta^{z} (\xi, \theta) \partial^\beta \varphi (\xi), 
\end{align}
for some $\mathscr{W}_\beta^{z}: \mathbb{R}^N \times \Theta \to \mathbb{R}$ specified in the statement of Lemma \ref{lemma:diff_step_1}.  
%  satisfying that: $G_\beta^{[T_2], z} (x, \theta) = 0$ for any $(x, \theta) \in \mathbb{R}^N \times \Theta$ if 
% $$
% |\beta_R| + 2 | \beta_S | > K + 3, 
% \quad 
% \mathrm{or}
% \quad   
% \| \beta \|_\# > K + \tfrac{3}{2} - \tfrac{1}{2} \lbrack \tfrac{K-1}{2} \rbrack 
% = 
% \begin{cases} 
% K + 2 - \tfrac{1}{2} \lbrack \tfrac{K}{2} \rbrack - \tfrac{1}{2} & 
% \bigl( K = 2j + 1, \, j \ge 1 \bigr); \\ 
% K + 2 - \tfrac{1}{2} \lbrack \tfrac{K}{2} \rbrack 
% &\bigl( K = 2j, \, j \ge 1 \bigr). 
% \end{cases}  
% $$ 
Subsequently, we study the third term $T_3$. Set $\beta = (\beta_S, \beta_R) = \gamma - e_{k_1} + e_{k_2}$ for $1 \le k_1, k_2 \le N$ and $\gamma \in \mathscr{I}_{\ref{eq:hypo}} (K)$  (thus $1 \le |\gamma_R| + 2 | \gamma_S| \le K+2$). 
Then it holds that: 
\begin{align*}
	1 \le |\beta_R| + 2 | \beta_S | 
	\le  
	\begin{cases} 
		K + 1 & (1 \le k_1 \le N_S, \, N_S + 1 \le k_2 \le N); \\ 
		K + 2  & (1 \le k_1, k_2 \le N_S \ \mathrm{or} \ N_S + 1 \le k_1, k_2 \le N); \\   
		K + 3 & (N_S + 1 \le k_1 \le N, \, 
		1 \le k_2 \le N_S).  
	\end{cases} 
\end{align*}  
Thus, $T_3$ takes the same form as (\ref{eq:T2}).  
% 
% \begin{align} \label{eq:T3} 
% T_3 = 
% %\sum_{\substack{ \beta \in \mathbb{Z}_{\ge 0}^N \\ 
% %1 \le |\beta_R| + 2 | \beta_S |  \\le K + 3}}
% \sum_{\beta \in \mathscr{I}_{\ref{eq:hypo}}(K+1)} 
%  G_\beta^{[T_3], z} (x, \theta) \partial^\beta \varphi (x),  
% \end{align}
% % 
% for some $G_\beta^{[T_3], z} : \mathbb{R}^N \times \Theta \to \mathbb{R}$ given in a closed-form. 
% satisfying that: $G_\beta^{[T_3], z} (x, \theta) = 0$ for any $(x, \theta) \in \mathbb{R}^N \times \Theta$ if 
% $$ 
% | \beta_R | + 2 | \beta_S | > K + 3, 
% \quad 
% \mathrm{or} 
% \quad  
% \| \beta \|_\# > K + 2 - \tfrac{1}{2} \lbrack \tfrac{K}{2}
% \rbrack . 
% $$
% 
Therefore, the assertion holds for $J = K + 1$.  Proof of Lemma \ref{lemma:diff_step_1} is now complete. \qed 

\subsection{Step 2.}
\label{app:diff_step_2}
We introduce: for a multi-index $\alpha \in \mathbb{Z}_{\ge 0}^j, \, j \in \mathbb{N}$,
\begin{align}
\begin{aligned} \label{eq:set_mi}
\mathscr{J}_{\ref{eq:ellip}} (\alpha) 
& \equiv \Bigl\{  \gamma \in \mathbb{Z}_{\ge 0}^N \ \bigl| \  1 \le |\gamma| \le | \alpha |  + 2 j - \| \alpha \|_{0}  \Bigr\}; \\ 
\mathscr{J}_{\ref{eq:hypo}} (\alpha) 
& \equiv 
\Bigl\{  \gamma \in \mathbb{Z}_{\ge 0}^N \ \bigl| \  1 \le |\gamma_R | + 2 |\gamma_S|  \le | \alpha |  + 2 j  - \| \alpha \|_{0}, 
\  \| \gamma \|_{\ref{eq:hypo}} \le  F (\alpha) \Bigr\},  
\end{aligned} 
\end{align}
% \begin{align*}
% \mathscr{J}_w (\alpha) 
% = 
% \begin{cases}
% \bigl\{  \gamma \in \mathbb{Z}_{\ge 0}^N \ \Bigl| \  1 \le |\gamma| \le | \alpha |  + 2 j - \| \alpha \|_{0}  \bigr\}, & w = \ref{eq:ellip}; \\[0.2cm] 
% \bigl\{  \gamma \in \mathbb{Z}_{\ge 0}^N \ \Bigl| \  1 \le |\gamma_R | + 2 |\gamma_S|  \le | \alpha |  + 2 j  - \| \alpha \|_{0}, 
% \  \| \gamma \|_{\ref{eq:hypo}} \le \tfrac{1}{2} \sum_{1 \le i \le j} \floor{\tfrac{3 \alpha_i}{2}} + \tfrac{3}{2}j - \tfrac{1}{2} \| \alpha \|_{0, 1, 2} \bigr\}, & w = \ref{eq:hypo},  
% \end{cases}
% \end{align*} 
% 
with 
\begin{gather*}
F(\alpha) \equiv \tfrac{1}{2} \sum_{1 \le i \le j} \floor{\tfrac{3 \alpha_i}{2}} + \tfrac{3}{2}j - \tfrac{1}{2} \| \alpha \|_{0, 1, 2}, \quad 
\| \alpha \|_{0} \equiv \sum_{1 \le i \le j} \mathbf{1}_{\alpha_i = 0}, \quad 
\| \alpha \|_{0, 1, 2} \equiv \sum_{1 \le i \le j} \bigl\{ 2 \times \mathbf{1}_{\alpha_i = 0} +  \mathbf{1}_{\alpha_i \in \{ 1, 2\} } \bigr\}. 
\end{gather*}
We also introduce: for the initial state of SDE $x \in \mathbb{R}^N$, 
\begin{align*}
\mathscr{S}^x  
\equiv \Bigl\{  f : \mathbb{R}^N \times \Theta \to \mathbb{R} \,  \bigl|   \, 
\forall \alpha \in \mathbb{Z}_{\ge 0}^N, \, \exists C > 0 \ \mathrm{s.t.} \ \sup_{\theta \in \Theta} | \partial^\alpha_x f (x, \theta) |  \le C \Bigr\}.  
\end{align*}
% 
%Under consideration of model class (\ref{eq:hypo}), we make use of the following notation: for any multi-index $\lambda \in \mathbb{Z}_{\ge 0}^N$, $\lambda_S \equiv (\lambda_1, \ldots, \lambda_{N_S})$ and 
%$\lambda_R \equiv (\lambda_{N_S+1}, \ldots, \lambda_N )$. 
% 
\begin{lemma}  \label{lemma:diff_step_2}
Let Assumptions \ref{ass:diff}, \ref{ass:param} and \ref{ass:coeff} hold. 
Let $w \in \{\ref{eq:ellip}, \ref{eq:hypo} \}$, 
$\alpha \in \mathbb{Z}_{\ge 0}^j, \, j \in \mathbb{N}$, $(x, \theta) \in \mathbb{R}^N \times \Theta$ and $\varphi \in C^\infty (\mathbb{R}^N, \mathbb{R})$. 
It holds that: 
%$w \in \{ \ref{eq:ellip}, \ref{eq:hypo} \}$, 
% 
\begin{align} \label{eq:D_a}
\mathscr{D}_\alpha^{z, \theta} \varphi (x) |_{z=x} = 
\sum_{\gamma \in \mathscr{J}_w (\alpha)  }
\mathscr{V}_\gamma (x, \theta) \partial^\gamma \varphi (x),    
\end{align} 
where $\mathscr{V}_\gamma: \mathbb{R}^N \times \Theta \to \mathbb{R}^N$ is explicitly given by products of partial derivatives of SDE's coefficients ${V}_i (x, \theta), \, 0 \le i \le d$ w.r.t. $x$ 
and characterised as follows: 
$ \mathscr{V}_\gamma  (x, \theta) = 0$
for any $\gamma \in \mathscr{J}_{w} (\mathbf{0})$, and $\mathscr{V}_\gamma \in \mathscr{S}^x$ for any $\gamma \in \mathscr{J}_w (\alpha)$. 
\end{lemma} 
%In the above, we have set: for $\alpha \in \mathbb{Z}_{\ge 0}^{j}, \, j \in \mathbb{N}$, 
%% 
%\begin{align*}
%\| \alpha \|_0 \equiv \sum_{1 \le k \le j} \mathbf{1}_{\alpha_k = 0}, \qquad  
%\| \alpha \|_{0, 1, 2} 
%\equiv
%\sum_{1 \le k \le j}  \mathbf{1}_{\alpha_k \in  \{0, 1, 2 \}}.  
%\end{align*} 
% 
\begin{remark} \label{rem:D}
 We note that for the multi-index $\gamma \in \mathscr{J}_w (\alpha), \, \alpha \in \mathbb{Z}_{\ge 0}^j, \, j \in \mathbb{N}, \, w \in \{\ref{eq:ellip}, \ref{eq:hypo}\}$, it holds that: 
 \begin{align} \label{eq:bd_gamma}
   \| \gamma \|_{w} \le | \alpha | + \tfrac{j}{2}, \qquad  w \in \{ \ref{eq:ellip}, \ref{eq:hypo} \}. 
 \end{align} 
 Thus, we have from (\ref{eq:bd_gamma}) and Lemmas \ref{lemma:deriv_LDL} and \ref{lemma:diff_step_2} that: there exist constants $C, \lambda >0$ such that for all $(\Delta, y, \theta)  \in (0,1)  \times \mathbb{R}^N \times \Theta$, 
 \begin{align*}
\Bigl| \Delta^{|\alpha| + j} \cdot \mathscr{D}_{\alpha}^{z, \theta} 
\bigl\{ p_\Delta^{\bar{X}^z, (w)}  (\cdot, y; \theta) \bigr\} (x) |_{z=x} 
\Bigr| 
\le C \Delta^{\tfrac{j}{2}} 
\cdot \Delta^{-m(w)/2} \exp \Bigl(  - \lambda 
\bigl|  \Gamma_{\Delta,(w)}^{-1} \bigl( y - x - V_{0} (x, \theta) \Delta  \bigr) \bigr|^2  \Bigr),   
 \end{align*} 
where we recall the definition of $\Gamma_{\Delta, (w)}$ in (\ref{eq:Gamma}). We here check the inequality (\ref{eq:bd_gamma}) for the hypo-elliptic model \eqref{eq:hypo}, and the case of \eqref{eq:ellip} follows a similar argument. Since $\gamma \in  \mathscr{J}_{\mathrm{\ref{eq:hypo}}} (\alpha)$ satisfies $\| \gamma \|_{\ref{eq:hypo}} \le F (\alpha)$, it suffices to show that: for $\alpha \in \mathbb{Z}_{\ge 0}^j, \, j \in \mathbb{N}$, 
\begin{align} \label{eq:F_bd}
F(\alpha) \le | \alpha |  + \tfrac{j}{2}. 
\end{align}
When $j=1$, i.e., $\alpha \in \mathbb{Z}_{\ge 0}$, we have that: 
% for $\gamma \in \mathscr{J}_{\ref{eq:hypo}} (\alpha)$,  
	% 
\begin{itemize}
\item When $\alpha = 0$,  
\begin{align}
F(\alpha)
% \equiv \tfrac{1}{2} \floor{\tfrac{3 \alpha}{2}} 
% + \tfrac{3}{2} - \tfrac{1}{2} \| \alpha \|_{0, 1, 2}  
\le \tfrac{1}{2} \ \bigl( =  \alpha  + \tfrac{1}{2} \bigr); 
\end{align} 
\item When $\alpha = 1$ or $\alpha = 2$,
\begin{align}
% \| \gamma \|_{\ref{eq:hypo}}  \le \tfrac{1}{2} \floor{\tfrac{3 \alpha}{2}} 
% + \tfrac{3}{2} - \tfrac{1}{2} \| \alpha \|_{0, 1, 2}   
F (\alpha )= 
\begin{cases}
\tfrac{3}{2} \ \bigl(=  \alpha + \tfrac{1}{2} \bigr) & \alpha = 1;  \\[0.2cm] 
\tfrac{5}{2} \ \bigl(=  \alpha + \tfrac{1}{2} \bigr) & \alpha = 2, 
\end{cases}
\end{align} 
\item When $\alpha \ge 3$, 
\begin{align}
F (\alpha ) 
= 
\tfrac{1}{2} \floor{\tfrac{3 \alpha}{2}}
+ \tfrac{3}{2} 
= \tfrac{\alpha}{2}  + \tfrac{1}{2} 
\floor{\tfrac{\alpha}{2}} 
+ \tfrac{3}{2}
\le | \alpha | + \tfrac{1}{2}, 
\end{align}
since $\floor{ \tfrac{\alpha}{2} } \le  \alpha  - 2$ for $\alpha \ge 3$. 
\end{itemize}% 
Thus, (\ref{eq:F_bd}) holds for $\alpha \in \mathbb{Z}$. 
We now assume that (\ref{eq:F_bd}) holds for $\alpha \in \mathbb{Z}_{\ge 0}^j$ with $j \le k$ and consider the case $\alpha \in \mathbb{Z}_{\ge 0}^{k+1}$.  We write $\alpha^{(k)} = (\alpha_1, \ldots, \alpha_k)$. Then it holds that: 
\begin{align*}
F (\alpha) 
& = \tfrac{1}{2} \sum_{1 \le i \le k+1} \floor{\tfrac{3 \alpha_i}{2}}
+ \tfrac{3}{2}(k+1) - \tfrac{1}{2} \| \alpha \|_{0, 1, 2}   \\ 
&  \le  | \alpha^{(k)}| + \tfrac{k}{2}  + \tfrac{1}{2} \floor{ \tfrac{3 \alpha_{k+1}}{2} } 
+ \tfrac{3}{2} - \tfrac{1}{2} \| \alpha_{k+1} \|_{0, 1, 2} \qquad \bigl( \because \mathrm{assumption \, of \, induction} \bigr) 
\nonumber \\[0.2cm] 
& \le |\alpha^{(k)}| + \tfrac{k}{2}  + | \alpha_{k+1}| + \tfrac{1}{2} 
\quad \bigl( \because \mathrm{argument \, from \, the \, case} \, j =1 \bigr) \\[0.2cm]
&  = |\alpha| + \tfrac{k+1}{2}. 
\end{align*}
Thus, (\ref{eq:F_bd}) holds for $\alpha \in \mathbb{Z}_{\ge 0}^{k+1}$. We then conclude (\ref{eq:bd_gamma}) under model class \eqref{eq:hypo}. 
\end{remark}
\noindent
\textit{Proof of Lemma \ref{lemma:diff_step_2}.}  
We provide proofs separately for the model classes \ref{eq:ellip} and \ref{eq:hypo}.  
\\  

\noindent
\textit{Proof for elliptic model \eqref{eq:ellip}}. 
We consider the mathematical induction on $j \in \mathbb{N}$ for the multi-index $\alpha \in \mathbb{Z}_{\ge 0}^j$. The case $j = 1$ has been already proved in Lemma \ref{lemma:diff_step_1}. We assume that the assertion holds for $j = 1, \ldots, k, \, k \ge 2$ and consider the case $j = k+1$. For $\alpha \in \mathbb{Z}_{\ge 0}^{k+1}$, we write $\alpha^{(k)} 
= (\alpha_1, \ldots, \alpha_k)$. Then the assumption of induction and Lemma \ref{lemma:diff_step_1} yield: 
\begin{align}
\mathscr{D}_\alpha^{z, \theta} \varphi (x) |_{z=x} 
& =  
\mathscr{D}_{\alpha^{(k)}}^{z, \theta} 
\bigl( \mathrm{ad}_{\mathscr{L}_\theta^{0,z}}^{\alpha_{k+1}} (\widetilde{\mathscr{L}}_\theta^z) \bigr) \varphi (x) |_{z=x} 
\nonumber \\ 
& 
=
\sum_{\gamma \in \mathscr{J}_{\ref{eq:ellip}}(\alpha^{(k)})}
\sum_{\beta \in \mathscr{I}_{\ref{eq:ellip}}(\alpha_{k+1})} 
\sum_{\nu \le \gamma} 
\binom{\gamma}{\nu} 
\mathscr{V}_\gamma  (x, \theta) \,  
\partial^\nu \mathscr{W}_\beta^{z} (x, \theta) 
\partial^{\beta + \gamma - \nu} \varphi (x) |_{z=x},   \label{eq:D_ell}
\end{align} 
where 
$\mathscr{W}_\beta^{z} : \mathbb{R}^N \times \Theta \to \mathbb{R}$ and $\mathscr{V}_\gamma : \mathbb{R}^N \times \Theta \to \mathbb{R}$ are specified in Lemmas \ref{lemma:diff_step_1}-(1) and \ref{lemma:diff_step_2}, respectively. We note that when $\alpha_{k+1} = 0$, it holds that for all $\gamma \in \mathscr{J}_{\ref{eq:ellip}} (\alpha^{(k)}) $ and $\beta \in \mathscr{I}_{\ref{eq:ellip}} (\alpha_{k+1})$,   
\begin{align}
\mathscr{V}_\gamma  (x, \theta) \,  
\mathscr{W}_\beta^{z} (x, \theta)  
\partial^{\beta+ \gamma} \varphi (x)|_{z=x}  
= 0, 
\end{align}
since the term $\mathscr{W}_\beta^z (x, \theta)$ involves polynomials of $(x-z)$ that become $0$ when $z=x$. 
\begin{comment}
The above considers the case when $\mathscr{D}_{\alpha^{(k)}}^{z, \theta} \bigl( \mathscr{L}_\theta - \mathscr{L}_\theta^{0, z} \bigr)
\varphi (x) |_{z=x}$. In particular, if the derivatives in $\mathscr{D}_{\alpha^{(k)}}^{z, \theta}$ do not act on $\mathscr{L}_\theta - \mathscr{L}_\theta^{0, z}$, the term is always $0$.  
\end{comment}
% 
Thus, (\ref{eq:D_ell}) takes the form (\ref{eq:D_a}) with the coefficients satisfying the properties specified in Lemma \ref{lemma:diff_step_2}, and then the assertion holds for $j = k+1$. The proof for elliptic model \eqref{eq:ellip} is now complete. 
\\ 

\noindent 
\textit{Proof for hypo-elliptic model \eqref{eq:hypo}}. We again rely on the mathematical induction on $j = |\alpha|, \, \alpha \in \mathbb{Z}_{\ge 0}^j$. We first consider the case $j = 1$, i.e., $\alpha \equiv J \in \mathbb{Z}_{\ge 0}$. We have from Lemma \ref{lemma:diff_step_1}-(2) that:
\begin{align}
	\mathscr{D}_{\alpha}^{z, \theta} \varphi (x) |_{z=x} 
	= \sum_{\substack{ \gamma \in \mathbb{Z}_{\ge 0}^N \\ 
	1 \le |\gamma_R| + 2 |\gamma_S| \le J + 2}} \mathscr{W}_\gamma^{x} (x, \theta)  \partial^\gamma \varphi (x),  
\label{eq:j=1}
\end{align}
where the coefficient $\mathscr{W}_\gamma^{x} (x, \theta)$ is specified in Lemma \ref{lemma:diff_step_1}-(2). We then show that the summation of the right-hand side of (\ref{eq:j=1}) is taken for mutli-indices $\gamma \in \mathscr{J}_{\ref{eq:hypo}} (J)$. We notice that: 
% the multi-indices $\gamma \equiv (\gamma_S, \gamma_R) \in \mathbb{Z}_{\ge 0}^N$ with $1 \le |\gamma_R| + 2 | \gamma_S| \le J + 2$ satisfy: 
% 
\begin{align} 
\max_{\substack{ \gamma \in \mathbb{Z}_{\ge 0}^N \\ 
1 \le |\gamma_R| + 2 |\gamma_S| \le J + 2}} 
\| \gamma \|_{\ref{eq:hypo}} 
& = \tfrac{1}{2} 
\floor{ \tfrac{J}{2} }  
+ \tfrac{J}{2} + \tfrac{3}{2} 
= \tfrac{1}{2} \floor{ \tfrac{3J}{2} }  + \tfrac{3}{2}. 
\label{eq:max_mi}
\end{align} 
We also have from Lemma \ref{lemma:diff_step_1}-(2) that $W_\gamma^{x} (x, \theta)=0$ for any $\gamma \in \mathscr{I}_{\ref{eq:hypo}}(J)$,  $J = 1,2$ satisfying $\|\gamma\|_{\ref{eq:hypo}} = \tfrac{1}{2} \floor{\tfrac{3J}{2}} + \tfrac{3}{2}$. 
% This indicates that $\mathscr{W}_\gamma^x (x, \theta)$ can take non-zero values only for $\gamma \in \mathscr{I}_{\ref{eq:hypo}} (J)$ s.t. $\| \gamma \|_{\ref{eq:hypo}} \le $. 
Thus, (\ref{eq:j=1}) takes the form (\ref{eq:D_a}) with $\alpha \equiv J \in \mathbb{Z}_{\ge 0}$. 

Subsequently, we assume that (\ref{eq:D_a}) holds for multi-indices $\alpha \in \mathbb{Z}_{\ge 0}^j$ with $j = 1, \ldots, k$. We then consider $j =k+1$  
%  
%% 
%\begin{align*}
%J (\alpha) & \equiv |\alpha| + 2j - \| \alpha \|_{0, 1, 2}, \qquad \alpha \in \mathbb{Z}_{\ge 0}^j, \ j 
%\in \mathbb{N}; \\ 
%K (\beta) & \equiv |\beta_R| + 2 |\beta_S|, \qquad \beta \in \mathbb{Z}_{\ge 0}^{N}. 
%\end{align*}
%  
and write $\alpha^{(k)} \equiv (\alpha_1, \ldots, \alpha_k)$ for $\alpha \in \mathbb{Z}_{\ge 0}^{k+1}$. 
% For notational simplicity, we define: 
% % 
% \begin{align}
% F (\alpha) = \tfrac{1}{2} \sum_{1 \le i \le j} \floor{\tfrac{3 \alpha_j}{2}} 
% + \tfrac{3}{2}j - \tfrac{1}{2} \| \alpha \|_{0, 1, 2} ,  
% \qquad  \alpha \in \mathbb{Z}_{\ge 0}^j, \, j \in \mathbb{N}. 
% \end{align} 
% 
 We consider the following two cases: (i) $\alpha_{k+1} = 0$, (ii) $\alpha_{k+1} \ge 1$. For the case (i), first we have: 
% The differential operator $\widetilde{\mathscr{L}}_\theta^z \equiv \mathscr{L}_\theta - \mathscr{L}_\theta^z $ is expressed as follows: for a test function $\varphi \in C^2 (\mathbb{R}^N)$ and $\xi \in \mathbb{R}^N$, 
 % 
 \begin{align} \label{eq:tildeL}
 	 & \widetilde{\mathscr{L}}_\theta^z  
 	  = \sum_{1 \le i \le N} \bigl\{V_0^i (\xi, \theta) -  V_0^i (z, \theta) - \partial_z^\top V_0^i (z, \theta) \cdot (\xi - z)   \bigr\} \partial_{\xi_i} 
 	+  \tfrac{1}{2} \sum_{N_S + 1 \le i_1, i_2 \le N}  \bigl\{ a_{i_1 i_2} (\xi, \theta)  - a_{i_1 i_2} (z, \theta)   \bigr\} \partial_{\xi_{i_1} \xi_{i_2}} 
 	\nonumber  \\ 
 	& = \sum_{1 \le i  \le N}  \!  
 	\sum_{1 \le j_1, j_2 \le N} h_i^{j_1 j_2} (\xi, z, \theta) (\xi_{j_1} - z_{j_1}) (\xi_{j_2} - z_{j_2}) \partial_{\xi_i} 
 +   \! \! \! \! \! \!  \sum_{N_S + 1 \le i_1, i_2  \le N} \sum_{1 \le j \le N} 
 	h_{i_1 i_2}^j (\xi, z, \theta ) (\xi_j - z_j)
 	\partial_{\xi_{i_1} \xi_{i_2}}, 
 \end{align} 
 for some functions $h_i^{j_1 j_2}, \, h_{i_1 i_2}^{j} \in \mathscr{S}_\kappa$  under Assumptions \ref{ass:diff}, \ref{ass:param} and \ref{ass:coeff}. Note that we have performed Taylor expansion around $\xi = z$ to obtain the last line.  
Then, (\ref{eq:tildeL}) and the assumption of mathematical induction yield that: 
\begin{align*}
\mathscr{D}_\alpha^{z, \theta} \varphi (x) |_{z=x} 
= \mathscr{D}_{\alpha^{(k)}}^{z, \theta} \widetilde{\mathscr{L}}_\theta^{z}  
\varphi (x) |_{z=x}  = U_1^{(\mathrm{i})}   + U_2^{(\mathrm{i})},  
\end{align*}
where we have set: 
\begin{align}
	U_1^{(\mathrm{i})}  
	& \equiv  \sum_{\gamma \in \mathscr{J}_{\ref{eq:hypo} (\alpha^{(k)})}} 
	\sum_{\nu \le \gamma}  \sum_{1 \le i \le N} 
	 \sum_{1 \le j_1, j_2 \le N} \mathscr{V}_{\gamma} (x, \theta)  
%	h_i^{j_1 j_2} (\xi_1, \theta)    
	 \partial^\nu_x \bigl( h_i^{j_1 j_2} (x, z, \theta) (x_{j_1} - z_{j_1}) (x_{j_2} - z_{j_2}) \bigr) 
	\partial^{\gamma + e_i - \nu}  \varphi (x) |_{z=x};  \label{eq:U1_i}\\ 
	U_2^{(\mathrm{i})}  
	& \equiv 
	 \sum_{\gamma \in \mathscr{J}_{\ref{eq:hypo} (\alpha^{(k)})}}  
	 \sum_{\nu \le \gamma}
     \sum_{N_S + 1 \le i_1, i_2 \le N} \sum_{1 \le j  \le N} 
     \mathscr{V}_{\gamma}  (x, \theta)  
%     h_{i_1 i_2}^j (\xi_2, \theta ) 
     \partial^\nu_x \bigl( h_{i_1 i_2}^{j} (x, z, \theta) (x_j - z_j) \bigr) \, 
     \partial^{\gamma + e_{i_1} + e_{i_2} - \nu } \varphi (x) |_{z=x}, 
\end{align} 
% 
%for some $\xi_1, \xi_2 \in \mathbb{R}^N$, 
where $\mathscr{V}_{\gamma} \in \mathscr{S}^x$. 
% and $h_{i}^{j_1j_2}, \, h_{i_1 i_2}^j: \mathbb{R}^N \times \Theta \to \mathbb{R}$ satisfying: $\textstyle \sup_{\theta \in \Theta} | h_{i}^{j_1, j_2} (\xi_1, \theta) | \le C_1$, $\textstyle \sup_{\theta \in \Theta} | h_{i_1 i_2}^j (\xi_1, \theta) | \le C_2$ for some constants $C_1, C_2 > 0$ independent of $\xi_1$ and $\xi_2$, respectively. 
  Note that $U_1^{(\mathrm{i})}  = 0$ when $|\nu| \le 1$ and $U_2^{(\mathrm{i})}  = 0$ when $|\nu|=0$. 
For the multi-index $\lambda \equiv \gamma + e_i - \nu$ with $|\nu| \ge 2$, appearing in the term $U_1^{(\mathrm{i})} $, we have that: 
\begin{align*}
	|\lambda_R| + 2 | \lambda_S| 
	& \le |\gamma_R |  + 2 | \gamma_S|  \qquad (\because |e_{i,R}| + 2 |e_{i,S}| \le 2,  \, | \nu | \ge 2)   \\ 
	& \le | \alpha^{(k)} |	+  2k  - \| \alpha^{(k)} \|_{0} 
	\qquad (\because \mathrm{assumption \, of \, induction}) \\ 
	& \le | \alpha | + 2 k + 1 - \| \alpha \|_0 
	\qquad (\because  |\alpha_{k+1}| =0, \, \| \alpha_{k+1} \|_{0} = 1) 
\end{align*}
and 
\begin{align*}
	\| \lambda \|_{\ref{eq:hypo}}
	= \| \gamma  \|_{\ref{eq:hypo}} + \| e_i  \|_{\ref{eq:hypo}} - \| \nu \|_{\ref{eq:hypo}} 
	& \le \| \gamma \|_{\ref{eq:hypo}} + \tfrac{1}{2}  \qquad \bigl(  \because  \| e_i  \|_{\ref{eq:hypo}} \le \tfrac{3}{2}, \,  \| \nu \|_{\ref{eq:hypo}} \ge 1 \bigr) \\ 
	& \le F (\alpha^{(k)}) + \tfrac{1}{2} = F (\alpha).  
\end{align*} 
On the other hand, the multi-index $\lambda = \gamma + e_{i_1} + e_{i_2} - \nu, \, N_S + 1 \le i_1, i_2 \le N$, in the term $U_2^{(\mathrm{ii})}$ with $|\nu| \ge 1$ satisfies: 
\begin{align*}
|\lambda_R| + 2 | \lambda_S| 
& \le | \alpha^{(k)} | + 2k  -  \|\alpha^{(k)} \|_0 + 1 
= |\alpha| + 2 (k+1) - \| \alpha \|_0 \qquad 
\bigl(  \because |\alpha_{k+1}| = 0, \, \| \alpha_{k+1} \|_0 = 1 \bigr)
\end{align*}  
and 
\begin{align*}
	\| \lambda \|_{\ref{eq:hypo}}
	= \| \gamma  \|_{\ref{eq:hypo}} + \| e_{i_1}  \|_{\ref{eq:hypo}}  + \| e_{i_2}  \|_{\ref{eq:hypo}} - \| \nu \|_{\ref{eq:hypo}} 
	& \le \| \gamma \|_{\ref{eq:hypo}} + \tfrac{1}{2}  
	\qquad \bigl(  \because \| e_{i_1}  \|_{\ref{eq:hypo}}  + \| e_{i_2}  \|_{\ref{eq:hypo}} \le 1,  \  \| \nu \|_{\ref{eq:hypo}} \ge \tfrac{1}{2}  \bigr) \\ 
	& \le F (\alpha^{(k)}) + \tfrac{1}{2} = F (\alpha).  
\end{align*}  
Thus, we conclude that the summation of the terms $U_1^{(\mathrm{i})}$ and $U_2^{(\mathrm{i})}$ are expressed as (\ref{eq:D_a}) with $\alpha \in \mathbb{Z}_{\ge 0}^{k+1}$ when $\alpha_{k+1} = 0$.  

We move on to the case (ii), i.e., $\alpha_{k+1} \ge 1$. Lemma \ref{lemma:diff_step_1}-(2) and the assuption of mathematical induction yield that: 
\begin{align*}
\mathscr{D}_\alpha^{z, \theta} \varphi (x) |_{z=x} 
= \mathscr{D}_{\alpha^{(k)}}^{z, \theta} 
\bigl( \mathrm{ad}_{\mathscr{L}_\theta^{0, z}}^{\alpha_{k+1}}  (\widetilde{\mathscr{L}}_\theta^z) \bigr)  \varphi (x) |_{z=x}  
= U_1^{(\mathrm{ii})} + U_2^{(\mathrm{ii})}, 
\end{align*} 
where 
\begin{align}
U_1^{(\mathrm{ii})} & \equiv 
\sum_{\gamma \in \mathscr{J}_{\ref{eq:hypo}}(\alpha^{(k)})}
\sum_{ \beta \in \mathscr{I}_{\ref{eq:hypo}} (\alpha_{k+1})}
\mathscr{V}_\gamma (x, \theta) \,  
\mathscr{W}_\beta^{z} (x, \theta) |_{z=x}
\partial^{\beta + \gamma} \varphi (x) ; \\ 
U_2^{(\mathrm{ii})} & \equiv 
\sum_{\gamma \in \mathscr{J}_{\ref{eq:hypo}}(\alpha^{(k)})}
\sum_{ \beta \in \mathscr{I}_{\ref{eq:hypo}} (\alpha_{k+1})} 
\sum_{\substack{\nu \le \gamma \\ \nu \neq \mathbf{0}}} 
\binom{\gamma}{\nu} 
\mathscr{V}_\gamma (x, \theta) \,  
\partial^\nu \mathscr{W}_\beta^{z} (x, \theta) |_{z=x} 
\, \partial^{\beta + \gamma - \nu} \varphi (x) ,
\end{align} 
with 
$\mathscr{W}_\beta^{z} : \mathbb{R}^N \times \Theta \to \mathbb{R}$ and $\mathscr{V}_\gamma : \mathbb{R}^N \times \Theta \to \mathbb{R}$ being specified in Lemmas \ref{lemma:diff_step_1}-(2) and \ref{lemma:diff_step_2}, respectively. 
% 
%We show that these two terms $U_1^{[ii]}, U_2^{[ii]}$ are expressed in the form (\ref{eq:D_a_H}) with $j = k+1$. 
For the term $U_1^{(\mathrm{ii})}$,  we have that: 
\begin{align*}
\mathscr{V}_\gamma (x, \theta) \,  
\mathscr{W}_\beta^{z} (x, \theta) |_{z=x}
\partial^{\beta + \gamma} \varphi (x) = 0
\end{align*} 
for any $\beta \in \mathscr{I}_{\ref{eq:hypo}} (\alpha_{k+1})$, $\alpha_{k+1} = 1, 2$, satisfying $\| \beta \|_{\ref{eq:hypo}} = \tfrac{1}{2} \floor{\tfrac{3 \alpha_{k+1}}{2} }  + \tfrac{3}{2}$ 
%since $\mathscr{W}_\beta^{[\alpha_{k+1}], z} (x, \theta) |_{z=x} = 0$
 due to Lemma \ref{lemma:diff_step_1}-(2), since every term in $\mathscr{W}_\beta^z (x, \theta)$ involves $(x-z)$ that becomes $0$ when $z=x$. Thus, the multiindex $\beta$ in $U_1^{(\mathrm{ii})}$ satisfies $\| \beta \|_{\ref{eq:hypo}} \le \tfrac{1}{2} \floor{ \tfrac{3 \alpha_{k+1}}{2} } + \tfrac{3}{2} - \tfrac{1}{2} \| \alpha_{k+1} \|_{0, 1, 2}$ and then the multi-index $\lambda \equiv  \beta + \gamma$ involved in $U_1^{(\mathrm{ii})}$ satisfies: 
\begin{align*}
|\lambda_R| + 2 | \lambda_S| 
= |\beta_R| + 2 | \beta_S| + |\gamma_R| + 2 | \gamma_S|  
& \le |\alpha_{k+1}| + 2 + | \alpha^{(k)}| + 2 k - \| \alpha^{(k)} \|_0 \\
& = | \alpha | + 2 (k+1) - \| \alpha \|_0 
\qquad \bigl(  \because  \| \alpha_{k + 1} \|_0  =  0 \bigr) 
\end{align*}
and 
\begin{align*}
\| \lambda \|_{\ref{eq:hypo}} = \| \beta \|_{\ref{eq:hypo}} + \| \gamma \|_{\ref{eq:hypo}} \le \tfrac{1}{2} 
\floor{ \tfrac{3 \alpha_{k+1}}{2} }+ \tfrac{3}{2} 
- \tfrac{1}{2} \| \alpha_{k+1} \|_{0, 1, 2}
+ \tfrac{1}{2} \sum_{i = 1}^{k} \floor{  \tfrac{3 \alpha_i}{2} }  + \tfrac{3}{2}k - \tfrac{1}{2} \| \alpha^{(k)} \|_{0, 1, 2} 
= F(\alpha).   
\end{align*}
The multi-index $\lambda \equiv \beta + \gamma - \nu$ defined in the term $U_2^{(\mathrm{ii})}$ satisfies:
\begin{align*}
|\lambda_R| + 2 |\lambda_S| 
& = |\beta_R| + 2 |\beta_S|
+ |\gamma_R| + 2 |\gamma_S|
- |\nu_R| - 2 |\nu_S|  \\ 
& \le  | \alpha_{k+1} | + 2  + | \alpha^{(k)} | + 2 k 
- \| \alpha^{(k)} \|_0 - 1 
\qquad \bigl( \because  | \nu_R | + 2 | \nu_S |  \ge 1, \, \nu \neq \mathbf{0} \bigr) \\ 
& = |\alpha | + 2 (k + 1) - \| \alpha \|_0 - 1  
\qquad \bigl( \because  \alpha_{k+1} \geq 1 )
\end{align*}  
and 
\begin{align*}
\| \lambda \|_{\ref{eq:hypo}} 
= \| \beta \|_{\ref{eq:hypo}}  + \| \gamma \|_{\ref{eq:hypo}}  - \| \nu \|_{\ref{eq:hypo}} 
& \le \tfrac{1}{2} \floor{ \tfrac{3 \alpha_{k+1}}{2} }
+ \tfrac{3}{2}  + F (\alpha^{(k)})  - \| \nu \|_{\ref{eq:hypo}} 
\qquad \bigl( \because \eqref{eq:max_mi} \, \& \, \mathrm{assumption \, of \, induction} \bigr) \\
& \le \tfrac{1}{2} \floor{ \tfrac{3 \alpha_{k+1}}{2} } 
+ \tfrac{3}{2} 
+ F (\alpha^{(k)})  - \tfrac{1}{2}  \qquad \bigl(  \because \| \nu \|_{\ref{eq:hypo}} \ge \tfrac{1}{2}, \, \nu \neq \mathbf{0}  \bigr) \\ 
& = \tfrac{1}{2} \sum_{1 \le i \le k+1} \floor{\tfrac{3 \alpha_i}{2}}
+ \tfrac{3}{2} (k+1) - \tfrac{1}{2} \| \alpha^{(k)} \|_{0, 1, 2} - \tfrac{1}{2} 
\\ 
& \le F (\alpha).  \qquad \big( \because 0 \le \| \alpha_{k+1} \|_{0, 1, 2} \le 1 \ \mathrm{for} \ \alpha_{k+1} \ge 1 \bigr)
\end{align*}
Thus, the summation of the terms $U_1^{(\mathrm{ii})}$  and $U_2^{(\mathrm{ii})}$ is written in the form (\ref{eq:D_a}) with $\alpha \in \mathbb{Z}_{\ge 0}^{k+1}$ satisfying $\alpha_{k+1} \ge 1$.  We conclude that the assertion holds with $j = k+1$ for hypo-elliptic model \eqref{eq:hypo}, and then  proof of Lemma \ref{lemma:diff_step_2} is now complete. \qed 
\subsection{Proof of Lemmas \ref{lemma:grad_D} and \ref{lemma:grad_D_a}}
\label{app:pf_grad_D}
\subsubsection{Proof of Lemma \ref{lemma:grad_D}}
We will show (\ref{eq:grad_J}) only for the hypo-elliptic model class (\ref{eq:hypo}), and a similar argument is applied for the elliptic class (\ref{eq:ellip}). After showing (\ref{eq:grad_J}) under model class (\ref{eq:hypo}), we briefly show that (\ref{eq:grad_J}) holds with the adjoint operators $(\widetilde{\mathscr{L}}_\theta^z)^\ast$ and $(\mathscr{D}_{(J)}^{z, \theta})^\ast$.  

We divide the proof into the following two cases: (i) $J  = 0$ and (ii) $J \ge 1$. 
\\ 

\noindent 
\textbf{Case (i).} $J = 0$.  \; We have from (\ref{eq:tildeL}) that:  
\begin{align}
\bigl| \partial^{\alpha}  \widetilde{\mathscr{L}}_\theta^{ z} \bigl\{ \psi^{(\ref{eq:hypo})}  (t, \cdot) \bigr\} (\xi) \bigr|  
& \le C_1 \sum_{\nu \le \alpha} \sum_{1 \le j_1, j_2 \le N} \sum_{1 \le i \le N}   
\Bigl| 
\partial_\xi^\nu 
\bigl( (\xi_{j_1} - z_{j_1}) (\xi_{j_2} - z_{j_2})  \bigr) 
\partial^{\alpha - \nu + e_{i}} \psi^{(\ref{eq:hypo})} (t, \xi) 
\Bigr| \nonumber \\  
& \quad + C_2 
\sum_{\nu \le \alpha} \sum_{1 \le j \le N} 
\sum_{N_S  +1 \le i_1, i_2 \le N}   
\Bigl| 
\partial_\xi^\nu
\bigl( (\xi_{j} - z_{j})  \bigr) 
\partial^{\alpha - \nu + e_{i_1}+ e_{i_2}} \psi^{(\ref{eq:hypo})} (t, \xi) 
\Bigr|.  
\end{align}
for some constants $C_1, C_2 > 0$ independent of . Noting that 
$$ \|e_i\|_{\ref{eq:hypo}} = 
\begin{cases}
    \tfrac{3}{2}, & 1 \le i \le N_S;  \\ 
    \tfrac{1}{2}, & N_S + 1 \le i \le N, 
\end{cases}$$
we obtain (\ref{eq:grad_J}) for $J = 0$. 
\\ 

\noindent 
\textbf{Case (ii)}.  $J \ge 1$. \;  We have from Lemma \ref{lemma:diff_step_1}-2 that
\begin{align}
\bigl| 
\partial^{\alpha}  \mathscr{D}_{(J)}^{z, \theta}
\bigl\{ \psi^{(\ref{eq:hypo})}  (t, \cdot) \bigr\} (\xi) \bigr| 
\le \sum_{\gamma \in \mathscr{I}_{\ref{eq:hypo}}(J)} 
\sum_{\nu \le \alpha} \binom{\alpha}{\nu} 
\bigl| \partial_{\xi}^\nu 
\mathscr{W}_\gamma^z (\xi, \theta) \bigr| \times 
\bigl| \partial_\xi^{\alpha - \nu + \gamma} \psi^{(\ref{eq:hypo}} (t, \xi) \bigr|.  
\end{align} 
We also have from Lemma \ref{lemma:diff_step_1} that: for any $\gamma \in \mathscr{I}_{\ref{eq:hypo}} (J)$ with $J = 1, 2$ (recall the definition of $\mathscr{I}_{\ref{eq:hypo}} (J)$ in (\ref{eq:set_I})), 
\begin{align}
\mathscr{W}_\gamma^z (\xi, \theta) 
= 
h_{\gamma} (\xi,  z, \theta) \times \mathbf{1}_{\| \gamma\|_{\ref{eq:hypo}} \le \tfrac{1}{2} \floor{\tfrac{3J}{2}} + 1} + \sum_{1 \le i \le N} h_{\gamma}^i  (\xi, z, \theta) (\xi_i - z_i) 
+ \sum_{1 \le i_1, i_2 \le N} h_{\gamma}^{i_1i_2} (\xi,  z, \theta) (\xi_{i_1} - z_{i_1}) (\xi_{i_2} - z_{i_2}),
\nonumber 
\end{align}
for some $h_{\gamma}, h_{\gamma}^i, h_{\gamma}^{i_1i_2} \in \mathscr{S}_\kappa$ under Assumptions \ref{ass:diff}, \ref{ass:param} and \ref{ass:coeff}. Since we have from (\ref{eq:max_mi}) that 
$$ \max_{\gamma \in \mathscr{I}_{\ref{eq:hypo}}(J)} \| \gamma \|_{\ref{eq:hypo}} = \tfrac{1}{2} \floor{\tfrac{3J}{2}} + \tfrac{3}{2}, $$
we obtain (\ref{eq:grad_J}) for $J = 1, 2$. The case of $J \ge 3$ is also obtained by noting that 
\begin{align}
\mathscr{W}_\gamma^z (\xi, \theta) 
= 
h_{\gamma} (\xi,  z, \theta) 
+ \sum_{1 \le i \le N} h_{\gamma}^i  (\xi, z, \theta) (\xi_i - z_i) 
+ \sum_{1 \le i_1, i_2  \le N} h_{\gamma}^{i_1i_2} (\xi,  z, \theta) (\xi_{i_1} - z_{i_1}) (\xi_{i_2} - z_{i_2}),
\nonumber 
\end{align} 
for any $\gamma \in \mathscr{I}_{\ref{eq:hypo}} (J)$ with $J \ge 3$. 
\\ 

\noindent
\textit{Case of adjoint operators.} We will briefly explain the case of adjoint operators. We first consider $J = 0$. The adjoint of $\widetilde{\mathscr{L}}_\theta^z$ is given in the following form: for $\varphi \in C^2 (\mathbb{R}^N)$ and $\xi \in \mathbb{R}^N$, 
\begin{align} \label{eq:adj}
\begin{aligned}
& \bigl( \widetilde{\mathscr{L}}_\theta^{z} \bigr)^\ast \varphi (\xi)
= \sum_{1 \le j \le N} \widetilde{h}^j (\xi, z, \theta) (\xi_j - z_j) \varphi (\xi)
+ \sum_{1 \le i \le N} \sum_{1 \le j_1, j_2 \le N} \widetilde{h}_{i}^{j_1 j_2} (\xi, z, \theta) 
(\xi_{j_1} - z_{j_1}) (\xi_{j_2} - z_{j_2}) \partial_{\xi_i} \varphi (\xi)  \\ 
& \quad + \sum_{N_S + 1 \le i \le N} \sum_{1 \le j \le N} \widetilde{h}_i^j (\xi, z, \theta) \partial_{\xi_i} \varphi (\xi)
+ \sum_{N_S + 1 \le i_1, i_2 \le N} 
\sum_{1 \le j \le N} \widetilde{h}_{i_1 i_2}^j (\xi, z, \theta) 
(\xi_j - z_j)
\partial_{\xi_{i_1} \xi_{i_2}} \varphi (\xi), 
\end{aligned}
\end{align}
for some $\widetilde{h}^j, \, \widetilde{h}^{j_1 j_2}_i \, \widetilde{h}^{j}_i, \, \widetilde{h}^j_{i_1 i_2}  \in \mathscr{S}_\kappa$. Thus, following the argument in \textbf{Case (i)}, we have (\ref{eq:grad_J}) under the adjoint operator $(\mathscr{D}_{(J)}^{z, \theta})^\ast$ with $J = 0$. For $J \ge 1$, we can formulate the adjoint of $\mathscr{D}_{(J)}^{z, \theta}$ from (\ref{eq:ad_exp_H}) and then obtain (\ref{eq:grad_J}) from a similar argument from the proof of \textbf{Case (ii)}. The proof of Lemma \ref{lemma:grad_D} is now complete.  
\qed 

\subsubsection{Proof of Lemma \ref{lemma:grad_D_a}} 
Making use of Lemma \ref{lemma:diff_step_2}, we immediately have that: for any $\alpha \in \mathbb{Z}_{\ge 0}^{j}, \, j \ge 2$, 
\begin{align*}
\bigl| \mathscr{D}_\alpha^{z, \theta} \psi^{(w)} (t, \cdot) (x) |_{z=x} \bigr| 
& = \bigl| \sum_{\gamma \in \mathscr{J}_w (\alpha)} 
\mathscr{V}_\gamma (x, \theta) \partial^\gamma \psi^{(w)} (t, \cdot) (x) \bigr| \\ 
& \le C \sum_{\gamma \in \mathscr{J}_w (\alpha)}  
t^{- \| \gamma \|_w }  \times \bigl| \widetilde{\psi}^{(w)} (t, x) \bigr| 
%\begin{cases}    
%C t^{- \tfrac{|\alpha|}{2} - j + \tfrac{\| \alpha \|_0}{2}} \times \bigl| \widetilde{\psi}^{(\ref{eq:ellip})} (t, x) \bigr|, 
%& w = \ref{eq:ellip}; \\[0.2cm] 
%C t^{- \tfrac{1}{2} \sum_{1 \le i \le j} \floor{\tfrac{3 \alpha_i}{2}} - \tfrac{3}{2}j + \tfrac{1}{2} \| \alpha \|_{0, 1, 2}} \times \bigl| \widetilde{\psi}^{(\ref{eq:hypo})} (t, x) \bigr|, 
%& w = \ref{eq:hypo};  
 \le C t^{- | \alpha | - \tfrac{j}{2}} \times \bigl| \widetilde{\psi}^{(w)} (t, x) \bigr|,   
%\end{cases}
\end{align*} 
for some constant $C > 0$ independent of $t \in (0,1)$ and $\theta$, where we made use of (\ref{eq:bd_gamma}) in the last inequality.  
\qed

\section{Supporting information for numerical experiments}  \label{app:supp_num}

\subsection{CF-expansion for FHN model} \label{app:cf_fhn}
We provide the closed-form expression of transition density expansion for the FHN model (\ref{eq:FHN}) used in the numerical experiment in Section \ref{sec:num}. Throughout this section, we set $\theta = (\varepsilon, \gamma, \beta, \sigma)$ and write the drift function and diffusion matrix of the FHN model (\ref{eq:FHN}) as follows: 
\begin{align}
b (x; \theta) 
% \equiv
% \begin{bmatrix}
%     b^1 (x; \theta) \\[0.1cm]
%     b^2 (x; \theta) 
% \end{bmatrix}
= 
\begin{bmatrix}
\tfrac{1}{\varepsilon} \bigl(x_1 - (x_1)^3 - x_2 + s \bigr) \\[0.2cm]
\gamma x_1 - x_2 + \beta 
\end{bmatrix}, 
\quad 
\Sigma = 
\begin{bmatrix}
0 & 0 \\ 
0 & \sigma
\end{bmatrix}, \qquad x \in \mathbb{R}^2.   
\end{align} 
\subsubsection{CF-expansion based on full drift linearlisation} \label{app:fhn_cf_ldl}
First, we recall the (full) LDL scheme is defined in (\ref{eq:ldl_fhn}) with $\iota = \mathrm{I}$ in the main text. Then, the LDL scheme starting from the point $x \in \mathbb{R}^2$ with coefficients being frozen at $z \in \mathbb{R}^2$ is given as: 
\begin{align}
\bar{X}_t^{(\mathrm{I}), x, z} 
= x + \int_0^t \bigl( A_{z, \theta}^{(\mathrm{I})} \bar{X}_s^{(\mathrm{I}), x, z}  
+ b_{z, \theta}^{(\mathrm{I})} \bigr) ds
+ \Sigma (\theta) B_t, \qquad t \ge 0, 
\label{eq:ldl_frozen}
\end{align} 
with $B_t = (B_{1, t}, B_{2, t}), \, t \ge 0$ and 
\begin{align*}
A^{(\mathrm{I})}_{z, \theta} 
% \begin{bmatrix} 
% \partial_{z_1} b^1 (z; \theta) & \partial_{z_2} b^1 (z; \theta) \\[0.1cm]
% \partial_{z_1} b^2 (z; \theta) & \partial_{z_2} b^2 (z; \theta) 
% \end{bmatrix} 
=  
\begin{bmatrix} 
\tfrac{1}{\varepsilon} \bigl( 1 - 3 (z_1)^2 \bigr) 
& - \tfrac{1}{\varepsilon} \\[0.2cm] 
\gamma & -1 
\end{bmatrix}, \qquad 
b_{z, \theta}^{(\mathrm{I})} 
= b (z; \theta) - A^{(\mathrm{I})}_{z, \theta} z. 
\end{align*} 
% 
% % 
The transition density of the scheme (\ref{eq:ldl_frozen})  writes: 
\begin{align*}
\bar{p}_\Delta^{(\mathrm{I}), z} (x, y; \theta) 
\equiv 
\tfrac{1}{\sqrt{(2\pi)^2 \det a^{(\mathrm{I})} (\Delta, z)}} 
\exp \Bigl( - \tfrac{1}{2} 
\bigl(y  - 
\mu^{(\mathrm{I}), z} (\Delta, x, \theta) \bigr)^\top 
a^{(\mathrm{I})} (\Delta, z)^{-1} 
\bigl(y - \mu^{(\mathrm{I}), z} (\Delta, x, \theta)  \bigr) \Bigr)  
\end{align*}
where 
\begin{align*} 
\mu^{(\mathrm{I}), z} (\Delta, x, \theta)
= z + e^{\Delta A^{(\mathrm{I})}_{z, \theta}} 
(x - z) + 
\Bigl( \int_0^\Delta e^{(\Delta - s)  A^{(\mathrm{I})}_{z, \theta}}  ds \Bigr) 
b (z; \theta), \quad 
a^{(\mathrm{I})} (\Delta, z) 
= \int_0^t e^{(t-u) A^{(\mathrm{I})}_{z, \theta}} 
\, \Sigma \Sigma^\top \, e^{(t-u) A^{(\mathrm{I}), \top}_{z, \theta} } du.  
\end{align*}
The generator of the scheme $\bar{X}_t^{(\mathrm{I}), x, z, \theta}$ is given by 
\begin{align}
\mathscr{L}^{0, z}_\theta \varphi (x) = \sum_{i= 1, 2} 
\bigl[ A_{z, \theta}^{(\mathrm{I})} x + b^{(\mathrm{I})}_{z, \theta}  \bigr]_i \partial_{x_i} \varphi (x)
+ \tfrac{1}{2}  \sigma^2 \partial_{x_2}^2 \varphi (x), \qquad (z, x) \in \mathbb{R}^2 \times \mathbb{R}^2, 
\end{align}
for a sufficiently smooth function $\varphi : \mathbb{R}^2 \to \mathbb{R}$. We also write $\widetilde{\mathscr{L}}_\theta^{z} \equiv \mathscr{L}_\theta - \mathscr{L}_\theta^{0, z}$. Following the argument in Section \ref{sec:CF_expansion} in the main text, we obtain a closed-form expansion of transition density in the form of (\ref{eq:err}): 
\begin{align}
 p_\Delta^{X} (x, y; \theta) 
& = 
\Bigl\{ 
\bar{p}_\Delta^{(\mathrm{I}), z} (x, y; \theta)  
+ \tfrac{\Delta^2}{2} [\mathscr{L}^{0, z}_\theta, \widetilde{\mathscr{L}}^z_\theta ] \, \bar{p}_\Delta^{(\mathrm{I}), z} (\cdot, y; \theta) (x) 
+ \tfrac{\Delta^3}{6} \bigl[\mathscr{L}^{0, z}_\theta, [\mathscr{L}^{0, z}_\theta, \widetilde{\mathscr{L}}^z_\theta ]  \bigr] \, \bar{p}_\Delta^{(\mathrm{I}), z} (\cdot, y; \theta) (x)   \nonumber \\[0.1cm]
& \qquad + \tfrac{\Delta^4}{24}
\bigl\{ \mathrm{ad}_{\mathscr{L}^{0, z}_\theta}^3 (\widetilde{\mathscr{L}}^z_\theta) \bigr\} \,  \bar{p}_\Delta^{(\mathrm{I}), z} (\cdot, y; \theta)  (x)  
+ \tfrac{\Delta^5}{120} \bigl\{ \mathrm{ad}_{\mathscr{L}^{0, z}_\theta}^4 (\widetilde{\mathscr{L}}^z_\theta) \bigr\}  
\, \bar{p}_\Delta^{(\mathrm{I}), z} (\cdot, y; \theta) (x) \nonumber  \\[0.1cm]
& \qquad + \tfrac{\Delta^6}{720} 
\bigl\{ \mathrm{ad}_{\mathscr{L}^{0, z}_\theta}^5 (\widetilde{\mathscr{L}}^z_\theta) \bigr\}  
\bar{p}_\Delta^{(\mathrm{I}), z} (\cdot, y; \theta)  (x) 
+ \tfrac{\Delta^7}{5040} 
\bigl\{ \mathrm{ad}_{\mathscr{L}^{0, z}_\theta}^6 (\widetilde{\mathscr{L}}^z_\theta) \bigr\}
\bar{p}_\Delta^{(\mathrm{I}), z} (\cdot, y; \theta)  (x)   \nonumber \\[0.1cm] 
& \qquad + \tfrac{\Delta^3}{6} 
[\mathscr{L}^{0, z}_\theta, \widetilde{\mathscr{L}}^z_\theta] \, \widetilde{\mathscr{L}}^z_\theta \, 
\bar{p}_\Delta^{(\mathrm{II}), z} (\cdot, y; \theta)  (x) 
\Bigr\} \Bigr|_{z=x}  + 
\mathscr{R} (\Delta, x, y; \theta) \nonumber  \\ 
& \equiv q_\Delta^{(\mathrm{I})} ( x, y; \theta) +  	\mathscr{R} (\Delta, x, y; \theta),  \nonumber 
\end{align}   
where the residual $\mathscr{R}$ satisfies: 
\begin{align} \label{eq:R}
\int_{\mathbb{R}^2} \bigl| 	\mathscr{R} (\Delta, x, y; \theta)  \bigr| dy 
\le C \Delta^{3}, 
\end{align} 
for some positive constant $C = C (x, \theta)$ independent of $\Delta > 0$. We have that: 
\begin{align}
q_\Delta^{(\mathrm{I})} ( x, y; \theta) 
= \bar{p}_\Delta^{(\mathrm{I}), z} (x, y; \theta) |_{z=x}  
\times \Bigl\{ 1 +  \sum_{k = 1}^5 \Delta^{k/2} \times e_k^{(\mathrm{I})} (\Delta, x , y; \theta) 
\Bigr\} 
+ \widetilde{\mathscr{R}} (\Delta, x, y; \theta), 
\end{align}
where $\widetilde{\mathscr{R}}$ has the same property as the above $\mathscr{R}$ and 
\begin{align*}
\begin{aligned} %\label{eq:weight_I}
& e_k^{(\mathrm{I})} (\Delta, x , y; \theta)  = 0, 
\quad k = 1,2; \\[0.2cm] 
& e_3^{(\mathrm{I})} (\Delta, x , y; \theta) = - \tfrac{\Delta^{3/2}}{6} \cdot \tfrac{6 x_1 (s - x_1 + x_1^3 + x_2)^2}{\varepsilon^3} \cdot H_{(1)}^{(\mathrm{I})} (\Delta, x, y; \theta) ; 
\\[0.2cm]  
& e_4^{(\mathrm{I})} (\Delta, x , y; \theta) = - \tfrac{\Delta^2}{24} \cdot \tfrac{18 \sigma ^2 x_1 \left(s + x_1^3 - x_1 + x_2 \right)}{\epsilon ^3} \cdot H_{(1,2)}^{(\mathrm{I})} (\Delta, x, y; \theta)  
- \tfrac{\Delta^3}{120} \cdot \tfrac{24 \sigma ^2 x_1 
(s + x_1^3 - x_1 + x_2)}{\epsilon ^4} \cdot H_{(1,1)}^{(\mathrm{I})} (\Delta, x, y ; \theta);  
\\[0.2cm]
& e_5^{(\mathrm{I})} (\Delta, x , y; \theta) = \tfrac{\Delta^{3/2}}{24} \cdot 6 \cdot \tfrac{-3 x_1 \epsilon  (s + x_1^3 - x_1 + x_2) (\alpha +\gamma  x_1 - x_2 ) 
+ ( s + x_1^3 - x_1 + x_2 )^3
-\sigma ^2 x_1 \epsilon}{\epsilon^4} 
\cdot H_{(1)}^{(\mathrm{I})} (\Delta, x, y; \theta)  
\\[0.1cm] 
& \qquad \qquad \qquad \qquad   
- \tfrac{\Delta^{5/2}}{120} \cdot \tfrac{18 x_1 \sigma^4}{\varepsilon^3} \cdot 
H_{(1, 2, 2)}^{(\mathrm{I})} (\Delta, x, y; \theta)   
- \tfrac{\Delta^{7/2}}{720} \cdot  
\tfrac{60 x_1 \sigma^4}{\varepsilon^4} \cdot H_{(1, 1, 2)}^{(\mathrm{I})} (\Delta, x, y; \theta) 
\\[0.1cm]
& \qquad \qquad \qquad \qquad   
- \tfrac{\Delta^{9/2}}{5040} \cdot \tfrac{60 x_1 \sigma^4}{\varepsilon^5} \cdot H_{(1,1,1)}^{(\mathrm{I})} (\Delta, x, y; \theta),   
\end{aligned}
\end{align*}
with $H$ defined as: for a multi-index 
$\alpha = \{1, 2\}^{\ell}, \, \ell = 1, 2, 3$, 
\begin{align}
H_\alpha^{(\mathrm{I})} (\Delta, x, y; \theta) 
= \frac{\partial_\alpha 
\bar{p}_\Delta^{(\mathrm{I}), z} (\cdot, y; \theta) (x) }{\bar{p}_\Delta^{(\mathrm{I}), z} (x, y; \theta)} 
\Bigr|_{z=x}.   
\end{align}
%
% Thus, the CF-expansion around the LDL scheme for the FHN model (\ref{eq:FHN}) takes the following form: 
% % 
% \begin{align*}
% p_\Delta^X (x, y; \theta) 
% = \bar{p}_\Delta^{(\mathrm{I}), z} (x, y; \theta) |_{z=x}  
% \times \Bigl\{ 1 +  \sum_{k = 1}^5 \Delta^{k/2} \times w_k^{(\mathrm{I})} (\Delta, x , y; \theta) 
% \Bigr\} 
% + \overline{\mathscr{R}} (\Delta, x, y; \theta), 
% \end{align*} 
% % 
% where the residual $\overline{\mathscr{R}}$ has the same property as $\mathscr{R}$.  
% 
% 
\subsubsection{CF-expansion based on modified LDL scheme}
To obtain the CF-expansion based on the modified LDL scheme (defined in (\ref{eq:ldl_fhn}) with $\iota = \mathrm{II}$), we make use of the same argument in the case of full drift linearlisation in Section \ref{app:fhn_cf_ldl}. 
The modified LDL scheme starting from the point $x \in \mathbb{R}^2$ with coefficients being frozen at $z \in \mathbb{R}^2$ is given as: 
\begin{align}
\bar{X}_t^{(\mathrm{II}), x, z} 
= x + \int_0^t \bigl( A_{z, \theta}^{(\mathrm{II})} \bar{X}_s^{(\mathrm{II}), x, z}  
+ b_{z, \theta}^{(\mathrm{II})} \bigr) ds
+ \Sigma B_t, \qquad t \ge 0, 
\label{eq:mldl_frozen_2}
\end{align} 
with $B_t = (B_{1, t}, B_{2, t}), \, t \ge 0$ and 
\begin{align*}
A^{(\mathrm{II})}_{z, \theta} 
% \equiv  
% \begin{bmatrix} 
% \partial_{z_1} b^1 (z; \theta) & \partial_{z_2} 
% b^1 (z; \theta) \\[0.1cm]
% 0  & \partial_{z_2} b^2 (z; \theta) 
% \end{bmatrix} 
=  
\begin{bmatrix} 
\tfrac{1}{\varepsilon} \bigl( 1 - 3 (z_1)^2 \bigr) 
& - \tfrac{1}{\varepsilon} \\[0.2cm] 
0 & -1 
\end{bmatrix}, \qquad 
b_{z, \theta}^{(\mathrm{II})} 
= b (z; \theta) - A^{(\mathrm{II})}_{z, \theta}  z. 
\end{align*}
% 
% We have that: 
% % 
% \begin{align}
% \bar{X}_\Delta^{(\mathrm{II}), x, z, \theta}  
% = \mu^{(\mathrm{II}), z} (\Delta, x, \theta) 
% + \int_0^\Delta e^{(\Delta - s) A^{(\mathrm{II})}_{z, \theta} } \Sigma (\theta) \, d B_t, 
% \end{align} 
% % 
% where 
% % 
% \begin{align}
% \mu^{(\mathrm{II}), z} (\Delta, x, \theta)
% = z + e^{\Delta A^{(\mathrm{II})}_{z, \theta}} 
% (x - z) + 
% \left( \int_0^\Delta e^{(\Delta - s)  A^{(\mathrm{II})}_{z, \theta}}  ds \right) V_0 (z; \theta). 
% \end{align}
% 
The transition density of the scheme (\ref{eq:mldl_frozen_2}) frozen at $z \in \mathbb{R}^2$ writes: 
\begin{align*}
\bar{p}_\Delta^{(\mathrm{II}), z} (x, y; \theta) 
\equiv 
\tfrac{1}{\sqrt{(2\pi)^2 \det a^{(\mathrm{II})} (\Delta, z)}} 
\exp \Bigl( - \tfrac{1}{2} 
\bigl(y  - 
\mu^{(\mathrm{II}), z} (\Delta, x, \theta) \bigr)^\top 
a^{(\mathrm{II})} (\Delta, z)^{-1} 
\bigl(y - \mu^{(\mathrm{II}), z} (\Delta, x, \theta)  \bigr) \Bigr)  
\end{align*}
where 
\begin{align*} 
\mu^{(\mathrm{II}), z} (\Delta, x, \theta)
= z + e^{\Delta A^{(\mathrm{II})}_{z, \theta}} 
(x - z) + 
\Bigl( \int_0^\Delta e^{(\Delta - s)  A^{(\mathrm{II})}_{z, \theta}}  ds \Bigr) b (z; \theta), \quad 
a^{(\mathrm{II})} (\Delta, z)
= \int_0^t e^{(t-u) A^{(\mathrm{II})}_{z, \theta}} 
\, \Sigma \Sigma^\top \, e^{(t-u) A^{(\mathrm{II}), \top}_{z, \theta} } du.    
\end{align*}
% 
% The generator associated with the scheme $\bar{X}_t^{(\mathrm{II}), x, z, \theta}$ is given by 
% % 
% \begin{align}
% \mathscr{L}^{0, z}_\theta \varphi (x) = \sum_{i= 1, 2} 
% \bigl[ A_{z, \theta}^{(\mathrm{II})} x + b^{(\mathrm{II})}_{z, \theta}  \bigr]_i \partial_{x_i} \varphi (x)
% + \tfrac{1}{2}  \sigma^2 \partial_{x_2}^2 \varphi (x), \qquad (z, x) \in \mathbb{R}^2 \times \mathbb{R}^2, 
% \end{align}
% % 
% for sufficiently smooth function $\varphi : \mathbb{R}^2 \to \mathbb{R}$. We then have that:
From a similar argument in Section \ref{app:fhn_cf_ldl}, we obtain the following CF-expansion around the modified LDL scheme: 
\begin{align}
p_\Delta^X (x, y; \theta) 
= 
\bar{p}_\Delta^{(\mathrm{II}), z} (x, y; \theta) |_{z=x}  
\times \Bigl\{ 1 +  \sum_{k = 1}^5 \Delta^{k/2} \times w_k^{(\mathrm{II})} (\Delta, x , y; \theta) 
\Bigr\} + \mathscr{R} (\Delta, x, y; \theta) 
\end{align}
where the residual $\mathscr{R}$ is characterised as (\ref{eq:R}) and the correction terms are given as: 
\begin{align*}
& e_k^{(\mathrm{II})} (\Delta, x , y; \theta)  = 0, 
\quad k = 1,2; \nonumber \\[0.2cm] 
& e_3^{(\mathrm{II})} (\Delta, x , y; \theta)  
= - \tfrac{\Delta^{1/2}}{2} \cdot \tfrac{(s - x_1 + x_1^3 + x_2) \gamma}{\varepsilon} \cdot  H_{(2)}^{(\mathrm{II})} (\Delta, x, y; \theta) \\[0.1cm] 
& \qquad \qquad \qquad \qquad 
- \tfrac{\Delta^{3/2}}{6} \cdot \tfrac{6 x_1 (s - x_1 + x_1^3 + x_2)^2 + 2 (s - x_1 + x_1^3 + x_2) \gamma \epsilon}{\varepsilon^3} \cdot 
H_{(1)}^{(\mathrm{II})} (\Delta, x, y; \theta); 
\\[0.2cm]  
& e_4^{(\mathrm{II})} (\Delta, x , y; \theta) 
= - \tfrac{\Delta}{6} \cdot \tfrac{\gamma \sigma^2}{\varepsilon} \cdot H_{(2,2)}^{(\mathrm{II})} (\Delta, x, y; \theta)  
- \tfrac{\Delta^2}{24} \cdot 
\tfrac{2 \sigma^2  (9 x_1 (s + x_1^3 - x_1 + x_2 ) + 2 \gamma \varepsilon ) }{\epsilon ^3} \cdot H_{(1,2)}^{(\mathrm{II})} (\Delta, x, y; \theta) \nonumber \\[0.1cm]  
& \qquad \qquad \qquad \qquad 
- \tfrac{\Delta^3}{120} \cdot \tfrac{4 \sigma ^2 
(6 x_1 (s + x_1^3-x_1+x_2) + \gamma \varepsilon)}{\epsilon^4} \cdot H_{(1,1)}^{(\mathrm{II})} (\Delta, x, y; \theta); 
\\[0.2cm]
& e_5^{(\mathrm{II})} (\Delta, x , y; \theta) 
= \tfrac{\Delta^{1/2}}{6}
\cdot \tfrac{\gamma  (3 x_1^2-1 ) 
(s + x_1^3 - x_1 + x_2) - \gamma  \epsilon  
( \alpha + 2s +2 x_1^3 + (\gamma -2) x_1+ x_2 )}{\epsilon^2} \cdot H_{(2)}^{(\mathrm{II})} (\Delta, x, y ; \theta)  
\\[0.1cm]  
& \qquad \qquad \qquad \qquad  
+  \tfrac{\Delta^{3/2}}{24} \cdot  
g (x) 
\cdot H_{(1)}^{(\mathrm{II})} (\Delta, x, y; \theta) 
- \tfrac{\Delta^{5/2}}{120} \cdot \tfrac{18 x_1 \sigma^4}{\varepsilon^3} H_{(1, 2, 2)}^{(\mathrm{II})} (\Delta, x, y; \theta)   \\[0.1cm]  
& \qquad \qquad \qquad \qquad  
- \tfrac{\Delta^{7/2}}{720} \cdot 
\tfrac{60 x_1 \sigma^4}{\varepsilon^4} \cdot H_{(1, 1, 2)}^{(\mathrm{II})} (\Delta, x, y; \theta)
- \tfrac{\Delta^{9/2}}{5040} \cdot \tfrac{60 x_1 \sigma^4}{\varepsilon^5} \cdot H_{(1, 1, 1)}^{(\mathrm{II})} (\Delta, x, y; \theta) 
\end{align*} 
with
\begin{align*}
g (x) & \equiv  
- \tfrac{3}{\epsilon ^4}   
\Bigl\{ \alpha  \gamma  \epsilon ^2-2 s^3-6 s^2 (x_1^3-x_1 + x_2) + s \bigl(\gamma  \epsilon ^2
-6 (x_1^3 - x_1 + x_2 )^2 +6 x \epsilon (\alpha +\gamma  x_1 - x_2) \bigr) \\ 
& \qquad -2 x_1^9+6 x_1^7-6 x_1^6 x_2 +6 x_1^5 (\gamma  \epsilon -1)+ 6 x_1^4 (\alpha  \epsilon - x_2 (\epsilon -2))+x_1^3 \bigl(\gamma  (\epsilon
-6) \epsilon -6 x_2^2+2 \bigr) \\
& \qquad 
+6 x_1^2 (x_2 (\gamma  \epsilon +\epsilon -1) - 
\alpha  \epsilon ) 
+ x_1 \bigl( \epsilon \,  
\bigl\{ (\gamma -1) \gamma  \epsilon +2
\sigma ^2 \bigr\} -6 x_2^2 (\epsilon -1) + 6 \alpha  x_2 \epsilon \bigr)-2 x_2^3 \Bigr\} 
\end{align*}
and $H^{(\mathrm{II})}$ being defined as: for a multi-index 
$\alpha = \{1, 2\}^{\ell}, \, \ell = 1, 2, 3$, 
\begin{align}
H_\alpha^{(\mathrm{II})} (\Delta, x, y; \theta) 
= \frac{\partial_\alpha 
\bar{p}_\Delta^{(\mathrm{II}), z} (\cdot, y; \theta) (x) }{\bar{p}_\Delta^{(\mathrm{II}), z} (x, y; \theta)} 
\Bigr|_{z=x}.   
\end{align} 
% 
% \subsection{Bayesian inference of FHN model}
% 
% \subsubsection{Definition of posteriors P0 and P1}
% The local Gaussian scheme \citep{glot:20, glot:21} used in our experiment of FHN model (\ref{eq:FHN}) is defined as follows: for an initial point $x \in (v, u) \in \mathbb{R}^2$ and step-size $t \ge 0$,  
% % 
% \begin{align}
% \bar{X}_t^{\mathrm{LG}, x} 
% = 
% \begin{bmatrix}
% b^1 (x) t + \tfrac{1}{2} \mathscr{L}_\theta b^1 (x) t^2 \\[0.1cm]
% b^2 (x) t 
% \end{bmatrix}
% + 
% \begin{bmatrix}
% \mathscr{L}_1
% \end{bmatrix}
% \end{align} 
% 
\subsection{MCMC results for Section \ref{sec:results}}
\label{sec:results_suppl}
\begin{table}[h]
 \caption{Computational cost and summary statistics of MCMC chains for the FHN model.}
 \label{table:hmc_suppl}
 \centering 
  \begin{tabular}{cccccc}
\toprule 
scheme 
& parameter  & ess-bulk & ess-tail & r-hat 
& time(s)/iter \\
\midrule
\multirow{4}{*}{P0 (benchmark)}
 & $\varepsilon$ & 3491 & 4069 & 1 & \multirow{4}{*}{4.745}  \\
 & $\gamma$ & 3166 & 3571 & 1 & \\
 & $\beta$ & 6529 & 4996 & 1  & \\
 & $\sigma$ & 3616 & 3985 & 1 & \\
 \midrule  
\multirow{4}{*}{P1 (modified LDL)} 
 & $\varepsilon$ & 3781 & 3748 & 1 & \multirow{4}{*}{0.237} \\
 & $\gamma$ & 3449 & 3397 & 1 &  \\
 & $\beta$ & 6468 & 4226 & 1 &   \\
 & $\sigma$ & 3807 & 3697 & 1 &  \\
 \midrule  
\multirow{4}{*}{P2 (CF-expansion, $J =3$)} 
 & $\varepsilon$ & 3111 & 3696 & 1 & \multirow{4}{*}{0.460} \\
 & $\gamma$ & 3128 & 4162 & 1 & \\
 & $\beta$ & 5577 & 5665 & 1 &  \\
 & $\sigma$ & 3144 & 3695 & 1 &  \\
 \midrule  
 \multirow{4}{*}{P1' (Local Gaussian)}  
 & $\varepsilon$ & 2425 & 3147 & 1 & \multirow{4}{*}{0.118}  \\
 & $\gamma$ & 2300 & 3032 & 1 &  \\
 & $\beta$ & 4335 & 5493 & 1  & \\
 & $\sigma$ & 2482 & 3018 & 1  &  \\
 \midrule  
\multirow{4}{*}{P2' (CF-expansion, $J =4$)}
 & $\varepsilon$ & 2715 & 2876 & 1 & \multirow{4}{*}{0.548}  \\
 & $\gamma$ & 2926 & 2809 & 1 &  \\
 & $\beta$ & 5216 & 4928 & 1 &  \\
 & $\sigma$ & 2799 & 2908 & 1 &   \\
\midrule
\end{tabular}
\end{table} 

\begin{figure}[h]
\centering
\subfigure[\footnotesize P1' (Local Gaussian scheme)]{\includegraphics[width=6.5cm]{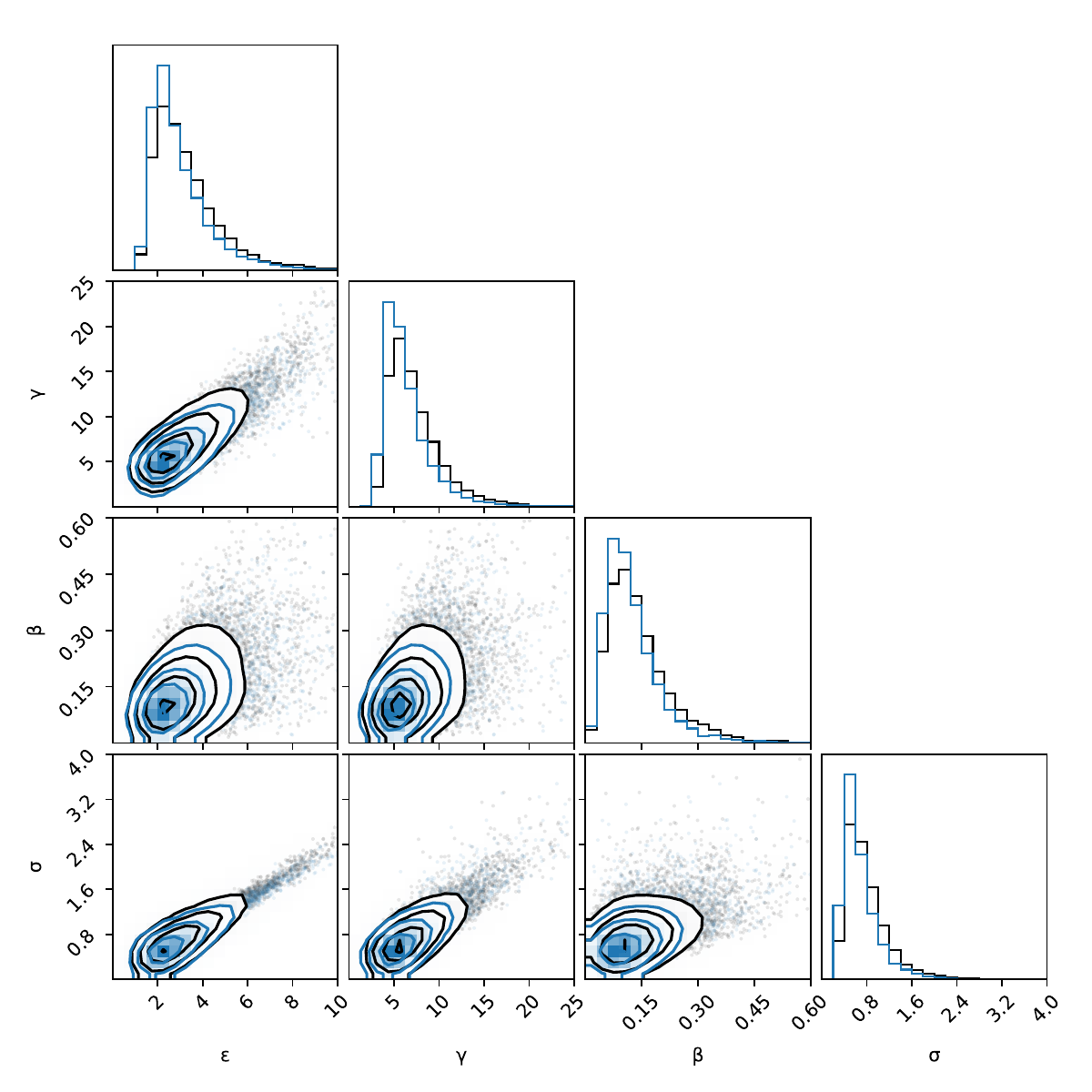} \label{fig:P3}}
%\hspace{1cm}
%\subfigure[\footnotesize P1 (Partial LDL scheme)]{\includegraphics[width=7cm]{FHN-Posterior-modified-LDLvsbenchmark_overlaid.pdf}  \label{fig:P1}}
%\subfigure[\footnotesize P2 (CF-expansion with $J = 3$)]{\includegraphics[width=7cm]{FHN-Posterior-DE-J-3vsbenchmark_overlaid.pdf} \label{fig:P2}}
\hspace{1cm} 
\subfigure[\footnotesize P2' (CF-expansion with $J = 4$)]{\includegraphics[width=6.5cm]{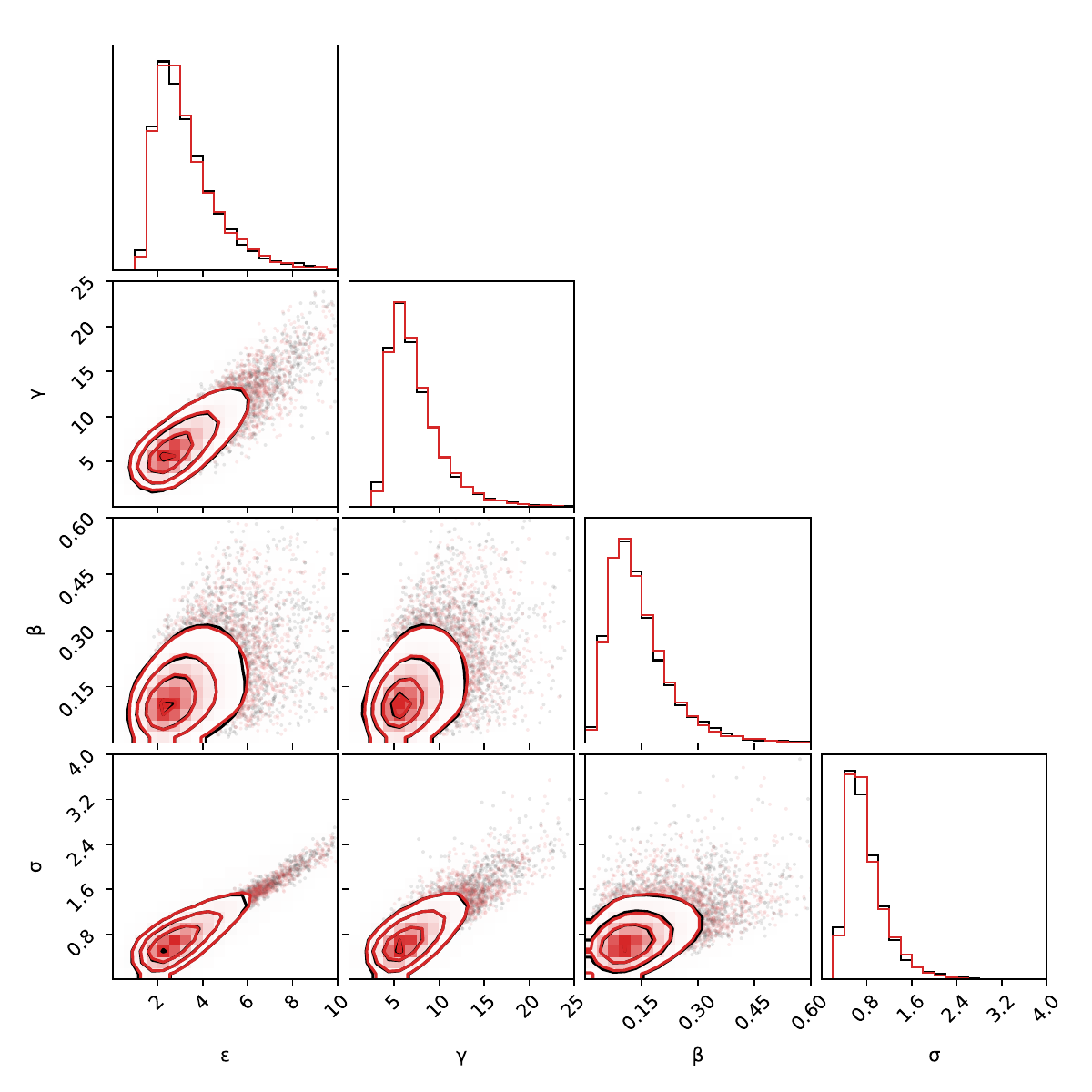}  \label{fig:P4}}
\caption{Posterior estimates. P0 (benchmark posterior) is overlaid in each figure in black.}
\label{fig:posterior2} 
\end{figure}

In addition to posteriors P0 (benchmark), P1 (partial LDL) and P2 (CF-expansion, $J=3$) from the main text, we also include here two more posteriors. Namely, P1' corresponds to the Gaussian approximation obtained via the local Gaussian scheme (\cite{glot:21}) and P2' to the CF-expansion with $J=4$.
Table \ref{table:hmc_suppl} provides average running time per iteration in secs from two chains and summary statistics that characterise the convergence of the chains, specifically, bulk effective sampling size ($ESS$), tail $ESS$ and $\hat{R}$ 
%with rank-normalisation and folding 
(analytic definitions are found in \cite{veh:21}). We note that for all schemes,
$\hat{R}<1.01$ and $ESS>400$, so the criteria recommended  in \cite{veh:21} are satisfied, thus we can conclude that the posteriors shown in Fig.~\ref{fig:posterior} in the main text and in Fig.~\ref{fig:posterior2} here are reliable.

%%%%%%%%%%%%%%%%%%%%%%%%%%%%%%%%%%%%%%%%%%%%%%
%%%% Main text entry area:

\bibliographystyle{plainnat} 
\bibliography{cf}     

%% or include bibliography directly:
% \begin{thebibliography}{}
% \bibitem[\protect\citeauthoryear{???}{???}]{b1}
% \end{thebibliography}

\end{document}